\documentclass{amsart}
\usepackage{amssymb}
\usepackage{amsfonts}
\usepackage[pagebackref,hypertexnames=false, colorlinks, citecolor=red, linkcolor=red]{hyperref}

\newtheorem{theorem}{Theorem}
\theoremstyle{plain}

\newtheorem{corollary}{Corollary}

\newtheorem{definition}{Definition}

\newtheorem{lemma}{Lemma}
\newtheorem{notation}{Notation}

\newtheorem{proposition}{Proposition}
\newtheorem{remark}{Remark}

\numberwithin{equation}{section}
\input{tcilatex}

\begin{document}
\title[The Corona Theorem in $\mathbb{C}^{n}$]{The Corona Theorem for the
Drury-Arveson Hardy space and other holomorphic Besov-Sobolev spaces on the
unit ball in $\mathbb{C}^{n}$}
\author[\c{S}. Costea]{\c{S}erban Costea}
\address{\c{S}. Costea\\
McMaster University\\
Department of Mathematics and Statistics\\
1280 Main Street West\\
Hamilton, Ontario L8S 4K1 Canada}
\email{secostea@math.mcmaster.ca}
\author[E. T. Sawyer]{Eric T. Sawyer$^{\dagger}$}
\address{E. T. Sawyer\\
McMaster University\\
Department of Mathematics and Statistics\\
1280 Main Street West\\
Hamilton, Ontario L8S 4K1 Canada}
\thanks{$\dagger$. Research supported in part by a grant from the National
Science and Engineering Research Council of Canada.}
\email{sawyer@mcmaster.ca}
\author[B. D. Wick]{Brett D. Wick$^{\ddagger}$}
\address{B. D. Wick\\
University of South Carolina\\
Department of Mathematics\\
LeConte College\\
1523 Greene Street\\
Columbia, SC 29208 USA}
\email{wick@math.sc.edu}
\thanks{$\ddagger$. Research supported in part by National Science
Foundation DMS Grant \# 0752703.}

\begin{abstract}
We prove that the multiplier algebra of the Drury-Arveson Hardy space $%
H_{n}^{2}$ on the unit ball in $\mathbb{C}^{n}$ has no corona in its maximal
ideal space, thus generalizing the Corona Theorem of L. Carleson to higher
dimensions. This result is obtained as a corollary of the Toeplitz corona
theorem and a new Banach space result: the Besov-Sobolev space $%
B_{p}^{\sigma }$ has the "baby corona property" for all $\sigma \geq 0$ and $%
1<p<\infty $. In addition we obtain infinite generator and semi-infinite
matrix versions of these theorems.
\end{abstract}

\maketitle
\tableofcontents

\section{Introduction}

In 1962 Lennart Carleson demonstrated in \cite{Car2} the absence of a corona
in the maximal ideal space of $H^{\infty }\left( \mathbb{D}\right) $ by
showing that if $\left\{ g_{j}\right\} _{j=1}^{N}$ is a finite set of
functions in $H^{\infty }\left( \mathbb{D}\right) $ satisfying 
\begin{equation}
\sum_{j=1}^{N}\left\vert g_{j}\left( z\right) \right\vert \geq
c>0,\;\;\;\;\;z\in \mathbb{D},  \label{coronadata}
\end{equation}%
then there are functions $\left\{ f_{j}\right\} _{j=1}^{N}$ in $H^{\infty
}\left( \mathbb{D}\right) $ with 
\begin{equation}
\sum_{j=1}^{N}f_{j}\left( z\right) g_{j}\left( z\right) =1,\;\;\;\;\;z\in 
\mathbb{D},  \label{coronasolutions}
\end{equation}%
In 1968 Fuhrmann \cite{Fuh} extended Carleson's corona theorem to the finite
matrix case. In 1980 Rosenblum \cite{Ros} and Tolokonnikov \cite{Tol} proved
the corona theorem for infinitely many generators $N=\infty $. This was
further generalized to the one-sided infinite matrix setting by Vasyunin in
1981 (see \cite{Tol3}). Finally Treil \cite{Trei} showed in 1988 that the
generalizations stop there by producing a counterexample to the two-sided
infinite matrix case.

H\"{o}rmander noted a connection between the corona problem and the Koszul
complex, and in the late 1970's Tom Wolff gave a simplified proof using the
theory of the $\overline{\partial }$ equation and Green's theorem (see \cite%
{Gar}). This proof has since served as a model for proving corona type
theorems for other Banach algebras.

While there is a large literature on corona theorems in one complex
dimension (see e.g. \cite{Nik}), progress in higher dimensions has been
limited. Indeed, apart from the simple cases in which the maximal ideal
space of the algebra can be identified with a compact subset of $\mathbb{C}%
^{n}$, no corona theorem has been proved until now in higher dimensions.
Instead, partial results have been obtained, such as the beautiful Toeplitz
corona theorem for Hilbert function spaces with a complete Nevanlinna-Pick
kernel, the $H^{p}$ corona theorem on the ball and polydisk, and results
restricting $N$ to $2$ generators in (\ref{coronadata}) (the case $N=1$ is
trivial). In particular, Varopoulos \cite{Var} published a lengthy classic
paper in an unsuccessful attempt to prove the corona theorem for the
multiplier algebra $H^{\infty }\left( \mathbb{B}_{n}\right) $ of the
classical Hardy space $H^{2}\left( \mathbb{B}_{n}\right) $ of holomorphic
functions on the ball with square integrable boundary values. His $BMO$
estimates for solutions with $N=2$ generators remain unimproved to this day.
We will discuss these partial results in more detail below.

Our main result is that the corona theorem, namely the absence of a corona
in the maximal ideal space, holds for the multiplier algebra $M_{H_{n}^{2}}$
of the Hilbert space $H_{n}^{2}$, the celebrated Drury-Arveson Hardy space
on the ball in $n$ dimensions.

\begin{theorem}
\label{DAcorona}If $\left\{ g_{j}\right\} _{j=1}^{N}$ is a finite set of
functions in $M_{H_{n}^{2}}$ satisfying (\ref{coronadata}), then there are
functions $\left\{ f_{j}\right\} _{j=1}^{N}$ in $M_{H_{n}^{2}}$ satisfying (%
\ref{coronasolutions}).
\end{theorem}

In many ways $H_{n}^{2}$, and not the more familiar space $H^{2}\left( 
\mathbb{B}_{n}\right) $, is the natural generalization to higher dimensions
of the classical Hardy space on the disk. For example, $H_{n}^{2}$ is
universal among Hilbert function spaces with the complete Pick property, and
its multiplier algebra $M_{H_{n}^{2}}$ is the correct home for the
multivariate von Neumann inequality (see e.g. \cite{AgMc}). See Arveson \cite%
{Arv} for more on the space $H_{n}^{2}$, including the model theory of $n$%
-contractions.

More generally, the corona theorem holds for the multiplier algebras $%
M_{B_{2}^{\sigma }\left( \mathbb{B}_{n}\right) }$ of the Besov-Sobolev
spaces $B_{2}^{\sigma }\left( \mathbb{B}_{n}\right) $, $0\leq \sigma \leq 
\frac{1}{2}$, on the unit ball $\mathbb{B}_{n}$ in $\mathbb{C}^{n}$. The
space $B_{2}^{\sigma }\left( \mathbb{B}_{n}\right) $ consists roughly of
those holomorphic functions $f$ whose derivatives of order $\frac{n}{2}%
-\sigma $ lie in the classical Hardy space $H^{2}\left( \mathbb{B}%
_{n}\right) =B_{2}^{\frac{n}{2}}\left( \mathbb{B}_{n}\right) $, and is
normed by%
\begin{equation*}
\left\Vert f\right\Vert _{B_{2}^{\sigma }\left( \mathbb{B}_{n}\right)
}=\left\{ \sum_{k=0}^{m-1}\left\vert f^{\left( k\right) }\left( 0\right)
\right\vert ^{2}+\int_{\mathbb{B}_{n}}\left\vert \left( 1-\left\vert
z\right\vert ^{2}\right) ^{m+\sigma }R^{m}f\left( z\right) \right\vert
^{2}d\lambda _{n}\left( z\right) \right\} ^{\frac{1}{2}},
\end{equation*}%
for some $m>\frac{n}{2}-\sigma $ where $R=\sum_{j=1}^{n}z_{j}\frac{\partial 
}{\partial z_{j}}$ is the radial derivative. In particular $H_{n}^{2}=B_{2}^{%
\frac{1}{2}}\left( \mathbb{B}_{n}\right) $. Finally, we also obtain
semi-infinite matrix versions of these results.

\begin{description}
\item[Note] Our techniques also yield BMO estimates for the $H^{\infty
}\left( \mathbb{B}_{n}\right) $ corona problem, which will appear elsewhere.
\end{description}

\section{The corona problem in $\mathbb{C}^{n}$\label{coronalift}}

Let $X$ be a Hilbert space of holomorphic functions in an open set $\Omega $
in $\mathbb{C}^{n}$ that is a reproducing kernel Hilbert space with a \emph{%
complete irreducible Nevanlinna-Pick kernel} (see \cite{AgMc} for the
definition). The following \emph{Toeplitz corona theorem} is due to Ball,
Trent and Vinnikov \cite{BaTrVi} (see also Ambrozie and Timotin \cite{AmTi}
and Theorem 8.57 in \cite{AgMc}).

For $f=\left( f_{\alpha }\right) _{\alpha =1}^{N}\in \oplus ^{N}X$ and $h\in
X$, define $\mathbb{M}_{f}h=\left( f_{\alpha }h\right) _{\alpha =1}^{N}$ and%
\begin{equation*}
\left\Vert f\right\Vert _{Mult\left( X,\oplus ^{N}X\right) }=\left\Vert 
\mathbb{M}_{f}\right\Vert _{X\rightarrow \oplus ^{N}X}=\sup_{\left\Vert
h\right\Vert _{X}\leq 1}\left\Vert \mathbb{M}_{f}h\right\Vert _{\oplus
^{N}X}.
\end{equation*}%
Note that $\max_{1\leq \alpha \leq N}\left\Vert \mathcal{M}_{f_{\alpha
}}\right\Vert _{M_{X}}\leq \left\Vert f\right\Vert _{Mult\left( X,\oplus
^{N}X\right) }\leq \sqrt{\sum_{\alpha =1}^{N}\left\Vert \mathcal{M}%
_{f_{\alpha }}\right\Vert _{M_{X}}^{2}}$.

\begin{description}
\item[Toeplitz corona theorem] \label{ToeCorThm}Let $X$ be a Hilbert
function space in an open set $\Omega $ in $\mathbb{C}^{n}$ with an
irreducible complete Nevanlinna-Pick kernel. Let $\delta >0$ and $N\in 
\mathbb{N}$. Then $g_{1},...,g_{N}\in M_{X}$ satisfy the following "baby
corona property"; for every $h\in X$, there are $f_{1},...,f_{N}\in X$ such
that 
\begin{eqnarray}
\left\Vert f_{1}\right\Vert _{X}^{2}+...+\left\Vert f_{N}\right\Vert
_{X}^{2} &\leq &\frac{1}{\delta }\left\Vert h\right\Vert _{X}^{2},
\label{corcon} \\
g_{1}\left( z\right) f_{1}\left( z\right) +...+g_{N}\left( z\right)
f_{N}\left( z\right) &=&h\left( z\right) ,\ \ \ \ \ z\in \Omega ,  \notag
\end{eqnarray}%
\emph{if and only if} $g_{1},...,g_{N}\in M_{X}$ satisfy the following
"multiplier corona property"; there are $\varphi _{1},...,\varphi _{N}\in
M_{X}$ such that%
\begin{eqnarray}
\left\Vert \varphi \right\Vert _{Mult\left( X,\oplus ^{N}X\right) } &\leq &1,
\label{multcorcon} \\
g_{1}\left( z\right) \varphi _{1}\left( z\right) +...+g_{N}\left( z\right)
\varphi _{N}\left( z\right) &=&\sqrt{\delta },\ \ \ \ \ z\in \Omega .  \notag
\end{eqnarray}
\end{description}

The \emph{baby corona theorem} is said to hold for $X$ if whenever $%
g_{1},...,g_{N}\in M_{X}$ satisfy 
\begin{equation}
\left\vert g_{1}\left( z\right) \right\vert ^{2}+...+\left\vert g_{N}\left(
z\right) \right\vert ^{2}\geq c>0,\ \ \ \ \ z\in \Omega ,  \label{lowerphi}
\end{equation}%
then $g_{1},...,g_{N}$ satisfy the baby corona property (\ref{corcon}). The
Toeplitz corona theorem thus provides a useful tool for reducing the
multiplier corona property (\ref{multcorcon}) to the more tractable, but
still very difficult, baby corona property (\ref{corcon}) for multiplier
algebras $M_{B_{p}^{\sigma }\left( \mathbb{B}_{n}\right) }$ of certain of
the Besov-Sobolev spaces $B_{p}^{\sigma }\left( \mathbb{B}_{n}\right) $ when 
$p=2$ - see below. The case of $M_{B_{p}^{\sigma }\left( \mathbb{B}%
_{n}\right) }$ when $p\neq 2$ must be handled by more classical methods and
remains largely unsolved.

\begin{remark}
\label{maxidealspace}A standard abstract argument applies to show that the
absence of a corona for the multiplier algebra $M_{X}$, i.e. the density of
the linear span of point evaluations in the maximal ideal space of $M_{X}$,
is equivalent to the following assertion: for each finite set $\left\{
g_{j}\right\} _{j=1}^{N}\subset M_{X}$ such that (\ref{lowerphi}) holds for
some $c>0$, there are $\left\{ \varphi _{j}\right\} _{j=1}^{N}\subset M_{X}$
and $\delta >0$ such that condition (\ref{multcorcon}) holds. See for
example Lemma 9.2.6 in \cite{Nik} or the proof of Criterion 3.5 on page 39
of \cite{Saw}.
\end{remark}

\subsection{The Baby Corona Theorem\label{The Baby Corona}}

To treat $N>2$ generators in (\ref{corcon}), it is just as easy to treat the
case $N=\infty $, and this has the advantage of not requiring bookkeeping of
constants depending on $N$. We will

\begin{enumerate}
\item use the Koszul complex for infinitely many generators, and

\item invert higher order forms in the $\overline{\partial }$ equation, and

\item devise new estimates for the Charpentier solution operators for these
equations including,

\begin{enumerate}
\item the use of sharp estimates on Euclidean expressions $\left\vert \left( 
\overline{w-z}\right) \frac{\partial }{\partial \overline{w}}f\right\vert $
in terms of the invariant derivative $\left\vert \widetilde{\nabla }%
f\right\vert $ (see Proposition \ref{threecrucial}),

\item the use of the exterior calculus together with the explicit form of
Charpentier's solution kernels in Theorems \ref{explicit} and \ref%
{explicitamel} to handle \emph{rogue} Euclidean factors $\overline{%
w_{j}-z_{j}}$ (see Section \ref{opest}), and

\item the application of generalized operator estimates of Schur type\ in
Lemma \ref{Zlemma} to obtain appropriate boundedness of solution operators.
\end{enumerate}
\end{enumerate}

In addition to these novel elements in the proof, we make crucial use of the
beautiful integration by parts formula of Ortega and Fabrega \cite{OrFa},
and in order to obtain $\ell ^{2}$-valued results, we use the clever
factorization of the Koszul complex in Andersson and Carlsson \cite{AnCa}
but adapted to $\ell ^{2}$.

\begin{notation}
For sequences $f\left( z\right) =\left( f_{i}\left( z\right) \right)
_{i=1}^{\infty }\in \ell ^{2}$ we will write%
\begin{equation*}
\left\vert f\left( z\right) \right\vert =\sqrt{\sum_{i=1}^{\infty
}\left\vert f_{i}\left( z\right) \right\vert ^{2}}.
\end{equation*}%
When considering sequences of vectors such as $\nabla ^{m}f\left( z\right)
=\left( \nabla ^{m}f_{i}\left( z\right) \right) _{i=1}^{\infty }$, the same
notation $\left\vert \nabla ^{m}f\left( z\right) \right\vert =\sqrt{%
\sum_{i=1}^{\infty }\left\vert \nabla ^{m}f_{i}\left( z\right) \right\vert
^{2}}$ will be used with $\left\vert \nabla ^{m}f_{i}\left( z\right)
\right\vert $ denoting the Euclidean length of the vector $\nabla
^{m}f_{i}\left( z\right) $. Thus the symbol $\left\vert \cdot \right\vert $
is used in at least three different ways; to denote the absolute value of a
complex number, the length of a finite vector in $\mathbb{C}^{N}$ and the
norm of a sequence in $\ell ^{2}$. Later it will also be used to denote the
Hilbert-Schmidt norm of a tensor, namely the square root of the sum of the
squares of the coefficients in the standard basis. In all cases the meaning
should be clear from the context.
\end{notation}

Recall that $B_{p}^{\sigma }\left( \mathbb{B}_{n};\ell ^{2}\right) $
consists of all $f=\left( f_{i}\right) _{i=1}^{\infty }\in H\left( \mathbb{B}%
_{n};\ell ^{2}\right) $ such that%
\begin{equation}
\left\Vert f\right\Vert _{B_{p}^{\sigma }\left( \mathbb{B}_{n};\ell
^{2}\right) }\equiv \sum_{k=0}^{m-1}\left\vert \nabla ^{k}f\left( 0\right)
\right\vert +\left( \int_{\mathbb{B}_{n}}\left\vert \left( 1-\left\vert
z\right\vert ^{2}\right) ^{m+\sigma }\nabla ^{m}f\left( z\right) \right\vert
^{p}d\lambda _{n}\left( z\right) \right) ^{\frac{1}{p}}<\infty ,
\label{defB}
\end{equation}%
for some $m>\frac{n}{p}-\sigma $. By Proposition \ref{Bequiv} below (see
also \cite{Bea}), the right side is finite for some $m>\frac{n}{p}-\sigma $
if and only if it is finite for all $m>\frac{n}{p}-\sigma $. As usual we
will write $B_{p}^{\sigma }\left( \mathbb{B}_{n}\right) $ for the
scalar-valued space.

We now state our baby corona theorem for the $\ell ^{2}$-valued Banach
spaces $B_{p}^{\sigma }\left( \mathbb{B}_{n};\mathbb{\ell }^{2}\right) $, $%
\sigma \geq 0$, $1<p<\infty $. Observe that for $\sigma <0$, $%
M_{B_{p}^{\sigma }\left( \mathbb{B}_{n}\right) }=B_{p}^{\sigma }\left( 
\mathbb{B}_{n}\right) $ is a subalgebra of $C\left( \overline{\mathbb{B}_{n}}%
\right) $ and so has no corona. The $N=2$ generator case of Theorem \ref%
{baby} when $\sigma \in \left[ 0,\frac{1}{p}\right) \cup \left( \frac{n}{p}%
,\infty \right) $ and $1<p<\infty $ is due to Ortega and Fabrega \cite{OrFa}%
, who also obtain the $N=2$ generator case when $\sigma =\frac{n}{p}$ and $%
1<p\leq 2$. See Theorem A in \cite{OrFa}. In \cite{OrFa2} Ortega and Fabrega
prove analogous results with scalar-valued Hardy-Sobolev spaces in place of
the Besov-Sobolev spaces.

Let $\left\Vert \mathbb{M}_{g}\right\Vert _{B_{p}^{\sigma }\left( \mathbb{B}%
_{n}\right) \rightarrow B_{p}^{\sigma }\left( \mathbb{B}_{n};\mathbb{\ell }%
^{2}\right) }$ denote the norm of the multiplication operator $\mathbb{M}%
_{g} $ from $B_{p}^{\sigma }\left( \mathbb{B}_{n}\right) $ to the $\mathbb{%
\ell }^{2}$-valued Besov-Sobolev space $B_{p}^{\sigma }\left( \mathbb{B}_{n};%
\mathbb{\ell }^{2}\right) $.

\begin{theorem}
\label{baby}Let $\delta >0$, $\sigma \geq 0$ and $1<p<\infty $. Then there
is a constant $C_{n,\sigma ,p,\delta }$ such that given $g=\left(
g_{i}\right) _{i=1}^{\infty }\in M_{B_{p}^{\sigma }\left( \mathbb{B}%
_{n}\right) \rightarrow B_{p}^{\sigma }\left( \mathbb{B}_{n};\mathbb{\ell }%
^{2}\right) }$ satisfying%
\begin{eqnarray*}
\left\Vert \mathbb{M}_{g}\right\Vert _{B_{p}^{\sigma }\left( \mathbb{B}%
_{n}\right) \rightarrow B_{p}^{\sigma }\left( \mathbb{B}_{n};\mathbb{\ell }%
^{2}\right) } &\leq &1, \\
\sum_{j=1}^{\infty }\left\vert g_{j}\left( z\right) \right\vert ^{2} &\geq
&\delta ^{2}>0,\ \ \ \ \ z\in \mathbb{B}_{n},
\end{eqnarray*}%
there is for each $h\in B_{p}^{\sigma }\left( \mathbb{B}_{n}\right) $ a
vector-valued function $f\in B_{p}^{\sigma }\left( \mathbb{B}_{n};\ell
^{2}\right) $ satisfying%
\begin{eqnarray}
\left\Vert f\right\Vert _{B_{p}^{\sigma }\left( \mathbb{B}_{n};\ell
^{2}\right) } &\leq &C_{n,\sigma ,p,\delta }\left\Vert h\right\Vert
_{B_{p}^{\sigma }\left( \mathbb{B}_{n}\right) },  \label{Ngen} \\
\sum_{j=1}^{\infty }g_{j}\left( z\right) f_{j}\left( z\right) &=&h\left(
z\right) ,\ \ \ \ \ z\in \mathbb{B}_{n}.  \notag
\end{eqnarray}
\end{theorem}

\begin{corollary}
\label{adult}Let $0\leq \sigma \leq \frac{1}{2}$. Then the Banach algebra $%
M_{B_{2}^{\sigma }\left( \mathbb{B}_{n}\right) }$ has no corona, i.e. (\ref%
{corcon}) implies (\ref{multcorcon}). In particular this includes Theorem %
\ref{DAcorona} that the multiplier algebra of the Drury-Arveson space $%
H_{n}^{2}=B_{2}^{\frac{1}{2}}\left( \mathbb{B}_{n}\right) $ has no corona
(the one-dimensional case is Carleson's corona theorem), and also includes
that the multiplier algebra of the $n$-dimensional Dirichlet space $\mathcal{%
D}\left( \mathbb{B}_{n}\right) =B_{2}^{0}\left( \mathbb{B}_{n}\right) $ has
no corona (the one-dimensional case here is due to Tolokonnikov \cite{Tol2}).
\end{corollary}

The corollary follows immediately from the finite generator case $p=2$ of
Theorem \ref{baby} and the Toeplitz corona theorem (and Remark \ref%
{maxidealspace}) since the spaces $B_{2}^{\sigma }\left( \mathbb{B}%
_{n}\right) $ have an irreducible complete Nevanlinna-Pick kernel when $%
0\leq \sigma \leq \frac{1}{2}$ (\cite{ArRoSa2}).

We also have a semi-infinite matricial corona theorem.

\begin{corollary}
Let $0\leq \sigma \leq \frac{1}{2}$. Let $\mathcal{H}_{1}$ be a finite $m$%
-dimensional Hilbert space and let $\mathcal{H}_{2}$ be an infinite
dimensional separable Hilbert space. Suppose that $F\in \mathcal{M}%
_{B_{2}^{\sigma }\left( \mathbb{B}_{n}\right) \left( \mathcal{H}%
_{1}\rightarrow \mathcal{H}_{2}\right) }$ satisfies $\delta ^{2}I_{m}\leq
F^{\ast }(z)F(z)\leq I_{m}$. Then there is $G\in \mathcal{M}_{B_{2}^{\sigma
}\left( \mathbb{B}_{n}\right) \left( \mathcal{H}_{2}\rightarrow \mathcal{H}%
_{1}\right) }$ such that 
\begin{eqnarray*}
G(z)F(z) &=&I_{m}, \\
\Vert G\Vert _{\mathcal{M}_{B_{2}^{\sigma }\left( \mathbb{B}_{n}\right)
\left( \mathcal{H}_{2}\rightarrow \mathcal{H}_{1}\right) }} &\leq &C_{\sigma
,n,\delta ,m}.
\end{eqnarray*}
\end{corollary}

This corollary follows immediately from the case $p=2$ of Theorem \ref{baby}
and the Toeplitz corona theorem together with Theorem (MCT) in Trent and
Zhang \cite{TrZh}. See \cite{TrZh} for the notation used here. We already
commented above on the special case of this corollary for the Hardy space $%
B_{2}^{\frac{1}{2}}\left( \mathbb{B}_{1}\right) =H^{2}\left( \mathbb{D}%
\right) $ on the disk. The case $m=1$ of this corollary for the classical
Dirichlet space $B_{2}^{0}\left( \mathbb{B}_{1}\right) =\mathcal{D}\left( 
\mathbb{D}\right) $ on the disk is due to Trent \cite{Tre2}. It would be of
interest to determine the dependence of the constants on $p$ and $\delta $
in Theorem \ref{baby}.

\subsubsection{Prior results\label{Prior results}}

In \cite{AnCa} Andersson and Carlsson solve the baby corona problem for $%
H^{2}\left( \mathbb{B}_{n}\right) $ and obtain the analogous (baby) $H^{p}$
corona theorem on the ball $\mathbb{B}_{n}$ for $1<p<\infty $ and with
constants independent of the number of generators (see also Amar \cite{Ama},
Andersson and Carlsson \cite{AnCa2},\cite{AnCa3} and Krantz and Li \cite%
{KrLi}). Partial results on the corona problem restricted to $N=2$
generators and $BMO$ in place of $L^{\infty }$ estimates have been obtained
for $H^{\infty }\left( \mathbb{B}_{n}\right) $ (the multiplier algebra of $%
H^{2}\left( \mathbb{B}_{n}\right) =B_{2}^{\frac{n}{2}}\left( \mathbb{B}%
_{n}\right) $) by N. Varopoulos\ \cite{Var} in 1977. This classical corona
problem remains open (Problem 19.3.7 in \cite{Rud}), along with the corona
problems for the multiplier algebras of $B_{2}^{\sigma }\left( \mathbb{B}%
_{n}\right) $, $\frac{1}{2}<\sigma <\frac{n}{2}$.

More recently\ in 2000 J. M. Ortega and J. Fabrega \cite{OrFa} obtain
partial results with $N=2$ generators in (\ref{corcon}) for the algebras $%
M_{B_{2}^{\sigma }\left( \mathbb{B}_{n}\right) }$ when $0\leq \sigma <\frac{1%
}{2}$, i.e. from the Dirichlet space $B_{2}^{0}\left( \mathbb{B}_{n}\right) $
up to but not including the Drury-Arveson Hardy space $H_{n}^{2}=B_{2}^{%
\frac{1}{2}}\left( \mathbb{B}_{n}\right) $. To handle $N=2$ generators they
exploit the fact that a $2\times 2$ antisymmetric matrix consists of just
one entry up to sign, so that as a consequence the form $\Omega _{1}^{2}$ in
the Koszul complex below is $\overline{\partial }$-closed. The paper \cite%
{OrFa} by Ortega and Fabrega has proved to be of enormous influence in our
work, as the basic groundwork and approach we use are set out there.

In \cite{TrWi} S. Treil and the third author obtain the $H^{p}$ corona
theorem on the polydisk $\mathbb{D}^{n}$ (see also Lin \cite{Lin} and Trent 
\cite{Tre}). The Hardy space $H^{2}\left( \mathbb{D}^{n}\right) $ on the
polydisk fails to have the complete Nevanlinna-Pick property, and
consequently the Toeplitz corona theorem only holds in a more complicated
sense that a family of kernels must be checked for positivity instead of
just one. As a result the corona theorem for the algebra $H^{\infty }\left( 
\mathbb{D}^{n}\right) $ on the polydisk remains open for $n\geq 2$. Finally,
even the baby corona problems, apart from that for $H^{p}$, remain open on
the polydisk.

\subsection{Plan of the paper\label{Plan of the paper}}

We will prove Theorem \ref{baby} using the Koszul complex and a
factorization of Andersson and Carlsson, an explicit calculation of
Charpentier's solution operators, and generalizations of the integration by
parts formulas of Ortega and Fabrega, together with new estimates for
boundedness of operators on certain real-variable analogues of the
holomorphic Besov-Sobolev spaces. Here is a brief plan of the proof.

We are given an infinite vector of multipliers $g=\left( g_{i}\right)
_{i=1}^{\infty }\in M_{B_{p}^{\sigma }\left( \mathbb{B}_{n}\right)
\rightarrow B_{p}^{\sigma }\left( \mathbb{B}_{n};\mathbb{\ell }^{2}\right) }$
that satisfy $\left\Vert \mathbb{M}_{g}\right\Vert _{B_{p}^{\sigma }\left( 
\mathbb{B}_{n}\right) \rightarrow B_{p}^{\sigma }\left( \mathbb{B}_{n};%
\mathbb{\ell }^{2}\right) }\leq 1$ and $\inf_{\mathbb{B}_{n}}\left\vert
g\right\vert \geq \delta >0$, and an element $h\in B_{p}^{\sigma }\left( 
\mathbb{B}_{n}\right) $. We wish to find $f=\left( f_{i}\right)
_{i=1}^{\infty }\in B_{p}^{\sigma }\left( \mathbb{B}_{n};\mathbb{\ell }%
^{2}\right) $ such that

\begin{enumerate}
\item $\mathcal{M}_{g}f=g\cdot f=h,$

\item $\overline{\partial }f=0,$

\item $\left\Vert f\right\Vert _{B_{p}^{\sigma }\left( \mathbb{B}_{n};%
\mathbb{\ell }^{2}\right) }\leq C_{n,\sigma ,p,\delta }\left\Vert
h\right\Vert _{B_{p}^{\sigma }\left( \mathbb{B}_{n}\right) }.$
\end{enumerate}

An obvious first attempt at a solution is%
\begin{equation*}
f=\frac{\overline{g}}{\left\vert g\right\vert ^{2}}h,
\end{equation*}%
since $f$ obviously satisfies (1), can be shown to satisfy (3), but fails to
satisfy (2) in general.

To rectify this we use the Koszul complex in Section \ref{The Koszul complex}%
, which employs \emph{any} solution to the $\overline{\partial }$ problem on
forms of bidegree $\left( 0,q\right) $, $1\leq q\leq n$, to produce a
correction term $\Lambda _{g}\Gamma _{0}^{2}$ so that%
\begin{equation*}
f=\frac{\overline{g}}{\left\vert g\right\vert ^{2}}h-\Lambda _{g}\Gamma
_{0}^{2}
\end{equation*}%
now satisfies (1) and (2), but (3) is now in doubt without specifying the
exact nature of the correction term $\Lambda _{g}\Gamma _{0}^{2}$.

In Section \ref{Charpentier's solution kernels for} we explicitly calculate
Charpentier's solution operators to the $\overline{\partial }$ equation for
use in solving the $\overline{\partial }$ problems arising in the Koszul
complex. These solution operators are remarkably simple in form and moreover
are superbly adapted for obtaining estimates in real-variable analogues of
the Besov-Sobolev spaces in the ball. In particular, the kernels $K\left(
w,z\right) $ of these solution operators involve expressions like 
\begin{equation}
\frac{\left( 1-w\overline{z}\right) ^{n-1-q}\left( 1-\left\vert w\right\vert
^{2}\right) ^{q}\overline{\left( w-z\right) }}{\bigtriangleup \left(
w,z\right) ^{n}},  \label{express}
\end{equation}%
where%
\begin{equation*}
\sqrt{\bigtriangleup \left( w,z\right) }=\left\vert P_{z}\left( w-z\right) +%
\sqrt{1-\left\vert z\right\vert ^{2}}Q_{z}\left( w-z\right) \right\vert
\end{equation*}%
is the length of the vector $w-z$ shortened by multiplying by $\sqrt{%
1-\left\vert z\right\vert ^{2}}$ its projection $Q_{z}\left( w-z\right) $
onto the orthogonal complement of the complex line through $z$. Also useful
is the identity $\sqrt{\bigtriangleup \left( w,z\right) }=\left\vert 1-w%
\overline{z}\right\vert \left\vert \varphi _{z}\left( w\right) \right\vert $
where $\varphi _{z}$ is the involutive automorphism of the ball that
interchanges $z$ and $0$; in particular this shows that $d\left( w,z\right) =%
\sqrt{\bigtriangleup \left( w,z\right) }$ is a quasimetric on the ball.

In Section \ref{Real variable analogues} we introduce real-variable
analogues $\Lambda _{p,m}^{\sigma }\left( \mathbb{B}_{n}\right) $ of the
Besov-Sobolev spaces $B_{p}^{\sigma }\left( \mathbb{B}_{n}\right) $ along
with $\ell ^{2}$-valued variants, that are based on the geometry inherent in
the complex structure of the ball and reflected in the solution kernels in (%
\ref{express}). In particular these norms involve modifications $D$\ of the
invariant derivative $\widetilde{\nabla }$ in the ball:%
\begin{equation*}
Df\left( w\right) =\left( 1-\left\vert w\right\vert ^{2}\right) P_{w}\nabla +%
\sqrt{1-\left\vert w\right\vert ^{2}}Q_{w}\nabla .
\end{equation*}%
Three crucial inequalities are then developed to facilitate the boundedness
of the Charpentier solution operators, most notably%
\begin{equation}
\left\vert \left( \overline{z-w}\right) ^{\alpha }\frac{\partial ^{m}}{%
\partial \overline{w}^{\alpha }}F\left( w\right) \right\vert \leq
C\bigtriangleup \left( w,z\right) ^{\frac{m}{2}}\left\vert \left(
1-\left\vert w\right\vert ^{2}\right) ^{-m}\overline{D}^{m}F\left( w\right)
\right\vert ,  \label{notable}
\end{equation}%
for $F\in H^{\infty }\left( \mathbb{B}_{n};\ell ^{2}\right) $, which
controls the product of Euclidean lengths with Euclidean derivatives on the
left, in terms of the product of the smaller length $\sqrt{\bigtriangleup
\left( w,z\right) }$ and the larger derivative $\left( 1-\left\vert
w\right\vert ^{2}\right) ^{-1}\overline{D}$ on the right. We caution the
reader that our definition of $\overline{D}^{m}$ is \emph{not} simply the
composition of $m$ copies of $\overline{D}$ - see Definition \ref{Dpowers}
below.

In Section \ref{Integration by parts} we recall the clever integration by
parts formulas of Ortega and Fabrega involving the left side of (\ref%
{notable}), and extend them to the Charpentier solution operators for higher
degree forms. If we differentiate (\ref{express}), the power of $%
\bigtriangleup \left( w,z\right) $ in the denominator can increase and the
integration by parts in Lemma \ref{IBP1'} below will temper this singularity
on the diagonal. On the other hand the radial integration by parts in
Corollary \ref{IBP2iter} below will temper singularities on the boundary of
the ball.

In Section \ref{Schur's Test section} we use Schur's Test to establish the
boundedness of positive operators with kernels of the form%
\begin{equation*}
\frac{\left( 1-\left\vert z\right\vert ^{2}\right) ^{a}\left( 1-\left\vert
w\right\vert ^{2}\right) ^{b}\sqrt{\bigtriangleup \left( w,z\right) }^{c}}{%
\left\vert 1-\overline{w}z\right\vert ^{a+b+c+n+1}}.
\end{equation*}%
The case $c=0$ is standard (see e.g. \cite{Zhu}) and the extension to the
general case follows from an automorphic change of variables. These results
are surprisingly effective in dealing with the ameliorated solution
operators of Charpentier.

Finally in Section \ref{opest} we put these pieces together to prove Theorem %
\ref{baby}.

The appendix collects technical proofs of formulas and modifications of
existing proofs in the literature that would otherwise interrupt the main
flow of the paper.

\section{Charpentier's solution kernels for $(0,q)$-forms on the ball\label%
{Charpentier's solution kernels for}}

In Theorem I.1 on page 127 of \cite{Cha}, Charpentier proves the following
formula for $(0,q)$-forms:

\begin{theorem}
\label{Ch}For $q\geq 0$ and all forms $f\left( \xi \right) \in C^{1}\left( 
\overline{\mathbb{B}_{n}}\right) $ of degree $\left( 0,q+1\right) $, we have
for $z\in \mathbb{B}_{n}$:%
\begin{equation*}
f\left( z\right) =C_{q}\int_{\mathbb{B}_{n}}\overline{\partial }f\left( \xi
\right) \wedge \mathcal{C}_{n}^{0,q+1}\left( \xi ,z\right) +c_{q}\overline{%
\partial }_{z}\left\{ \int_{\mathbb{B}_{n}}f\left( \xi \right) \wedge 
\mathcal{C}_{n}^{0,q}\left( \xi ,z\right) \right\} .
\end{equation*}
\end{theorem}

Here $\mathcal{C}_{n}^{0,q}\left( \xi ,z\right) $ is a $\left(
n,n-q-1\right) $-form in $\xi $ on the ball and a $\left( 0,q\right) $-form
in $z$ on the ball that is defined in Definition \ref{defpqkernel} below.
Using Theorem \ref{Ch}, we can solve $\overline{\partial }_{z}u=f$ for a $%
\overline{\partial }$-closed $(0,q+1)$-form $f$ as follows. Set 
\begin{equation*}
u(z)\equiv c_{q}\int_{\mathbb{B}_{n}}f(\xi )\wedge \mathcal{C}_{n}^{0,q}(\xi
,z)
\end{equation*}%
Taking $\overline{\partial }_{z}$ of this we see from Theorem \ref{Ch} and $%
\overline{\partial }f=0$ that 
\begin{equation*}
\overline{\partial }_{z}u=c_{q}\overline{\partial }_{z}\left( \int_{\mathbb{B%
}_{n}}f(\xi )\wedge \mathcal{C}_{n}^{0,q}(\xi ,z)\right) =f(z).
\end{equation*}

It is essential for our proof to explicitly compute the kernels $\mathcal{C}%
_{n}^{0,q}$ when $0\leq q\leq n-1$. The case $q=0$ is given in \cite{Cha}
and we briefly recall the setup. Denote by $\bigtriangleup :\mathbb{C}%
^{n}\times \mathbb{C}^{n}\rightarrow \left[ 0,\infty \right) $ the map%
\begin{equation*}
\bigtriangleup (w,z)\equiv \left\vert 1-w\overline{z}\right\vert ^{2}-\left(
1-\left\vert w\right\vert ^{2}\right) \left( 1-\left\vert z\right\vert
^{2}\right) .
\end{equation*}%
We compute that 
\begin{eqnarray}
\bigtriangleup (w,z) &=&1-2\func{Re}w\overline{z}+\left\vert w\overline{z}%
\right\vert ^{2}-\left\{ 1-\left\vert w\right\vert ^{2}-\left\vert
z\right\vert ^{2}+\left\vert w\right\vert ^{2}\left\vert z\right\vert
^{2}\right\}  \label{Dis} \\
&=&\left\vert w-z\right\vert ^{2}+\left\vert w\overline{z}\right\vert
^{2}-\left\vert w\right\vert ^{2}\left\vert z\right\vert ^{2}  \notag \\
&=&\left( 1-\left\vert z\right\vert ^{2}\right) \left\vert w-z\right\vert
^{2}+\left\vert z\right\vert ^{2}\left( \left\vert w-z\right\vert
^{2}-\left\vert w\right\vert ^{2}\right) +\left\vert w\overline{z}%
\right\vert ^{2}  \notag \\
&=&\left( 1-\left\vert z\right\vert ^{2}\right) \left\vert w-z\right\vert
^{2}+\left\vert z\right\vert ^{4}-2\func{Re}\left\vert z\right\vert ^{2}w%
\overline{z}+\left\vert w\overline{z}\right\vert ^{2}  \notag \\
&=&\left( 1-\left\vert z\right\vert ^{2}\right) \left\vert w-z\right\vert
^{2}+\left\vert \overline{z}(w-z)\right\vert ^{2},  \notag
\end{eqnarray}%
and by symmetry%
\begin{equation*}
\bigtriangleup (w,z)=\left( 1-\left\vert w\right\vert ^{2}\right) \left\vert
w-z\right\vert ^{2}+\left\vert \overline{w}(w-z)\right\vert ^{2}.
\end{equation*}%
We also have the standard identity%
\begin{equation}
\bigtriangleup \left( w,z\right) =\left\vert 1-z\overline{w}\right\vert
^{2}\left\vert \varphi _{w}\left( z\right) \right\vert ^{2},
\label{standardid}
\end{equation}%
where%
\begin{equation*}
\varphi _{w}\left( z\right) =\frac{P_{w}\left( w-z\right) +\sqrt{%
1-\left\vert w\right\vert ^{2}}Q_{w}\left( w-z\right) }{1-\overline{w}z}.
\end{equation*}%
Thus we also have%
\begin{eqnarray}
\bigtriangleup (w,z) &=&\left\vert P_{w}\left( z-w\right) +\sqrt{%
1-\left\vert w\right\vert ^{2}}Q_{w}\left( z-w\right) \right\vert ^{2}
\label{dimcomp} \\
&=&\left\vert P_{z}\left( z-w\right) +\sqrt{1-\left\vert z\right\vert ^{2}}%
Q_{z}\left( z-w\right) \right\vert ^{2}.  \notag
\end{eqnarray}%
It is convenient to combine the many faces of $\bigtriangleup \left(
w,z\right) $ in (\ref{Dis}), (\ref{standardid}) and (\ref{dimcomp}) in:%
\begin{eqnarray}
\bigtriangleup \left( w,z\right) &=&\left\vert 1-w\overline{z}\right\vert
^{2}-\left( 1-\left\vert w\right\vert ^{2}\right) \left( 1-\left\vert
z\right\vert ^{2}\right)  \label{manyfaces} \\
&=&\left( 1-\left\vert z\right\vert ^{2}\right) \left\vert w-z\right\vert
^{2}+\left\vert \overline{z}(w-z)\right\vert ^{2}  \notag \\
&=&\left( 1-\left\vert w\right\vert ^{2}\right) \left\vert w-z\right\vert
^{2}+\left\vert \overline{w}(w-z)\right\vert ^{2}  \notag \\
&=&\left\vert 1-w\overline{z}\right\vert ^{2}\left\vert \varphi _{w}\left(
z\right) \right\vert ^{2}  \notag \\
&=&\left\vert 1-w\overline{z}\right\vert ^{2}\left\vert \varphi _{z}\left(
w\right) \right\vert ^{2}  \notag \\
&=&\left\vert P_{w}\left( z-w\right) +\sqrt{1-\left\vert w\right\vert ^{2}}%
Q_{w}\left( z-w\right) \right\vert ^{2}  \notag \\
&=&\left\vert P_{z}\left( z-w\right) +\sqrt{1-\left\vert z\right\vert ^{2}}%
Q_{z}\left( z-w\right) \right\vert ^{2}.  \notag
\end{eqnarray}

To compute the kernels $\mathcal{C}_{n}^{0,q}$ we start with the
Cauchy-Leray form 
\begin{equation*}
\mu (\xi ,w,z)\equiv \frac{1}{(\xi (w-z))^{n}}\sum_{i=1}^{n}(-1)^{i-1}\xi
_{i}\left[ \wedge _{j\neq i}d\xi _{j}\right] \wedge _{i=1}^{n}d(w_{i}-z_{i}),
\end{equation*}%
which is a closed form on $\mathbb{C}^{n}\times \mathbb{C}^{n}\times \mathbb{%
C}^{n}$ since with $\zeta =w-z$, $\mu $ is a pullback of the form%
\begin{equation*}
\nu (\xi ,\zeta )\equiv \frac{1}{(\xi \zeta )^{n}}\sum_{i=1}^{n}(-1)^{i-1}%
\xi _{i}\left[ \wedge _{j\neq i}d\xi _{j}\right] \wedge _{i=1}^{n}d\zeta
_{i},
\end{equation*}%
which is easily computed to be closed (see e.g. 16.4.5 in \cite{Rud}).

One then lifts the form $\mu $ via a section $s$ to give a closed form on $%
\mathbb{C}^{n}\times \mathbb{C}^{n}$. Namely, for $s:\mathbb{C}^{n}\times 
\mathbb{C}^{n}\rightarrow \mathbb{C}^{n}$ one defines, 
\begin{equation*}
s^{\ast }\mu \left( w,z\right) \equiv \frac{1}{\left( s\left( w,z\right)
\left( w-z\right) \right) ^{n}}\sum_{i=1}^{n}(-1)^{i-1}s_{i}\left(
w,z\right) \left[ \wedge _{j\neq i}ds_{j}\right] \wedge _{i=1}^{n}d\left(
w_{i}-z_{i}\right) .
\end{equation*}%
Now we fix $s$ to be the following section used by Charpentier: 
\begin{equation}
s(w,z)\equiv \overline{w}(1-w\overline{z})-\overline{z}(1-\left\vert
w\right\vert ^{2}).  \label{section}
\end{equation}%
Simple computations \cite{OrFa} demonstrate that 
\begin{eqnarray}
s(w,z)(w-z) &=&\left\{ \overline{w}\left( 1-w\overline{z}\right) -\overline{z%
}\left( 1-\left\vert w\right\vert ^{2}\right) \right\} \left( w-z\right)
\label{sdotD} \\
&=&\left\{ \left( \overline{w}-\overline{z}\right) -\left( w\overline{z}%
\right) \overline{w}+\left\vert w\right\vert ^{2}\overline{z}\right\} \left(
w-z\right)  \notag \\
&=&\left\vert w-z\right\vert ^{2}-\left( w\overline{z}\right) \left(
\left\vert w\right\vert ^{2}-\overline{w}z\right) +\left\vert w\right\vert
^{2}\left( \overline{z}w-\left\vert z\right\vert ^{2}\right)  \notag \\
&=&\left\vert w-z\right\vert ^{2}-\left( w\overline{z}\right) \left\vert
w\right\vert ^{2}+\left\vert \overline{w}z\right\vert ^{2}+\left\vert
w\right\vert ^{2}\overline{z}w-\left\vert w\right\vert ^{2}\left\vert
z\right\vert ^{2}  \notag \\
&=&\left\vert w-z\right\vert ^{2}+\left\vert \overline{w}z\right\vert
^{2}-\left\vert w\right\vert ^{2}\left\vert z\right\vert ^{2}=\bigtriangleup
\left( w,z\right) ,  \notag
\end{eqnarray}%
by the second line in (\ref{Dis}).

\begin{definition}
\label{defCauchy}We define the Cauchy Kernel on $\mathbb{B}_{n}\times 
\mathbb{B}_{n}$ to be 
\begin{equation}
\mathcal{C}_{n}\left( w,z\right) \equiv s^{\ast }\mu (w,z)  \label{defCn}
\end{equation}%
for the section $s$ given in (\ref{section}) above.
\end{definition}

\begin{definition}
\label{defpqkernel}For $0\leq p\leq n$ and $0\leq q\leq n-1$ we let $%
\mathcal{C}_{n}^{p,q}$ be the component of $\mathcal{C}_{n}\left( w,z\right) 
$ that has bidegree $(p,q)$ in $z$ and bidegree $(n-p,n-q-1)$ in $w$.
\end{definition}

Thus if $\eta $ is a $(p,q+1)$-form in $w$, then $\mathcal{C}%
_{n}^{p,q}\wedge \eta $ is a $(p,q)$-form in $z$ and a multiple of the
volume form in $w$. We now prepare to give explicit formulas for
Charpentier's solution kernels $\mathcal{C}_{n}^{0,q}(w,z)$. First we
introduce some notation.

\begin{notation}
\label{perm}Let $\omega _{n}\left( z\right) =\bigwedge_{j=1}^{n}dz_{j}$. For 
$n$ a positive integer and $0\leq q\leq n-1$ let $P_{n}^{q}$ denote the
collection of all permutations $\nu $ on $\{1,\ldots ,n\}$ that map to $%
\{i_{\nu },J_{\nu },L_{\nu }\}$ where $J_{\nu }$ is an increasing
multi-index with $\mathnormal{card}(J_{\nu })=n-q-1$ and $\mathnormal{card}%
(L_{\nu })=q$. Let $\epsilon _{\nu }\equiv sgn\left( \nu \right) \in \left\{
-1,1\right\} $ denote the signature of the permutation $\nu $.
\end{notation}

Note that the number of increasing multi-indices of length $n-q-1$ is $\frac{%
n!}{\left( q+1\right) !(n-q-1)!}$, while the number of increasing
multi-indices of length $q$ are $\frac{n!}{q!(n-q)!}$. Since we are only
allowed certain combinations of $J_{\nu }$ and $L_{\nu }$ (they must have
disjoint intersection and they must be increasing multi-indices), it is
straightforward to see that the total number of permutations in $P_{n}^{q}$
that we are considering is $\frac{n!}{(n-q-1)!q!}$.

From \O vrelid \cite{Ovr} we obtain that Charpentier's kernel takes the
(abstract) form 
\begin{equation*}
\mathcal{C}_{n}^{0,q}(w,z)=\frac{1}{\bigtriangleup (w,z)^{n}}\sum_{\nu \in
P_{n}^{q}}sgn\left( \nu \right) s_{i_{\nu }}\bigwedge_{j\in J_{\nu }}%
\overline{\partial }_{w}s_{j}\bigwedge_{l\in L_{\nu }}\overline{\partial }%
_{z}s_{l}\wedge \omega _{n}(w).
\end{equation*}%
Fundamental for us will be the explicit formula for Charpentier's kernel
given in the next theorem. We are informed by Part 2 of Proposition I.1 in 
\cite{Cha} that $\mathcal{C}_{n}^{p,q}(w,z)=0$ for $w\in \partial \mathbb{B}%
_{n}$, and this serves as a guiding principle in the proof we give in the
appendix. It is convenient to isolate the following factor common to all
summands in the formula: 
\begin{equation}
\Phi _{n}^{q}\left( w,z\right) \equiv \frac{\left( 1-w\overline{z}\right)
^{n-1-q}\left( 1-\left\vert w\right\vert ^{2}\right) ^{q}}{\bigtriangleup
\left( w,z\right) ^{n}},\ \ \ \ \ 0\leq q\leq n-1.  \label{defPhi}
\end{equation}

\begin{theorem}
\label{explicit}Let $n$ be a positive integer and suppose that $0\leq q\leq
n-1$. Then%
\begin{equation}
\mathcal{C}_{n}^{0,q}\left( w,z\right) =\sum_{\nu \in P_{n}^{q}}\left(
-1\right) ^{q}\Phi _{n}^{q}\left( w,z\right) sgn\left( \nu \right) \left( 
\overline{w_{i_{\nu }}}-\overline{z_{i_{\nu }}}\right) \bigwedge_{j\in
J_{\nu }}d\overline{w_{j}}\bigwedge_{l\in L_{\nu }}d\overline{z_{l}}%
\bigwedge \omega _{n}\left( w\right) .  \label{gensolutionker}
\end{equation}
\end{theorem}

\begin{remark}
We can rewrite the formula for $\mathcal{C}_{n}^{0,q}\left( w,z\right) $ in (%
\ref{gensolutionker}) as%
\begin{equation}
\mathcal{C}_{n}^{0,q}\left( w,z\right) =\Phi _{n}^{q}\left( w,z\right)
\sum_{\left\vert J\right\vert =q}\sum_{k\notin J}\left( -1\right) ^{\mu
\left( k,J\right) }\left( \overline{z_{k}}-\overline{w_{k}}\right) d%
\overline{z}^{J}\wedge d\overline{w}^{\left( J\cup \left\{ k\right\} \right)
^{c}}\wedge \omega _{n}\left( w\right) ,  \label{gensolutionker'}
\end{equation}%
where $J\cup \left\{ k\right\} $ here denotes the increasing multi-index
obtained by rearranging the integers $\left\{ k,j_{1},...j_{q}\right\} $ as 
\begin{equation*}
J\cup \left\{ k\right\} =\left\{ j_{1},...j_{\mu \left( k,J\right)
-1},k,j_{\mu \left( k,J\right) },...j_{q}\right\} .
\end{equation*}%
Thus $k$ occupies the $\mu \left( k,J\right) ^{th}$ position in $J\cup
\left\{ k\right\} $. The notation $\left( J\cup \left\{ k\right\} \right)
^{c}$ refers to the increasing multi-index obtained by rearranging the
integers in $\left\{ 1,2,...n\right\} \setminus \left( J\cup \left\{
k\right\} \right) $. To see (\ref{gensolutionker'}), we note that in (\ref%
{gensolutionker}) the permutation $\nu $ takes the $n$-tuple $\left(
1,2,...n\right) $ to $\left( i_{\nu },J_{\nu },L_{\nu }\right) $. In (\ref%
{gensolutionker'}) the $n$-tuple $\left( k,\left( J\cup \left\{ k\right\}
\right) ^{c},J\right) $ corresponds to $\left( i_{\nu },J_{\nu },L_{\nu
}\right) $, and so $sgn\left( \nu \right) $ becomes in (\ref{gensolutionker'}%
) the signature of the permutation that takes $\left( 1,2,...n\right) $ to $%
\left( k,\left( J\cup \left\{ k\right\} \right) ^{c},J\right) $. This in
turn equals $\left( -1\right) ^{\mu \left( k,J\right) }$ with $\mu \left(
k,J\right) $ as above.
\end{remark}

We observe at this point that the functional coefficient in the summands in (%
\ref{gensolutionker}) looks like 
\begin{equation*}
(-1)^{q}\Phi _{n}^{q}\left( w,z\right) \left( \overline{w_{i_{\nu }}}-%
\overline{z_{i_{\nu }}}\right) =(-1)^{q}\frac{(1-w\overline{z}%
)^{n-q-1}(1-|w|^{2})^{q}}{\bigtriangleup (w,z)^{n}}\left( \overline{%
w_{i_{\nu }}}-\overline{z_{i_{\nu }}}\right) ,
\end{equation*}%
which behaves like a fractional integral operator of order $1$ in the
Bergman metric on the diagonal relative to invariant measure. See the
appendix for a proof of Theorem \ref{explicit}.

Finally, we will adopt the usual convention of writing 
\begin{equation*}
\mathcal{C}_{n}^{0,q}f\left( z\right) =\int_{\mathbb{B}_{n}}f\left( w\right)
\wedge \mathcal{C}_{n}^{0,q}\left( w,z\right) ,
\end{equation*}%
when we wish to view $\mathcal{C}_{n}^{0,q}$ as an operator taking $\left(
0,q+1\right) $-forms $f$ in $w$ to $\left( 0,q\right) $-forms $\mathcal{C}%
_{n}^{0,q}f$ in $z$.

\subsection{Ameliorated kernels\label{Ameliorated kernels}}

We now wish to define right inverses with improved behaviour at the
boundary. We consider the case when the right side $f$ of the $\overline{%
\partial }$ equation is a $\left( p,q+1\right) $-form in $\mathbb{B}_{n}$.

As usual for a positive integer $s>n$ we will "project" the formula $%
\overline{\partial }\mathcal{C}_{s}^{p,q}f=f$ in $\mathbb{B}_{s}$ for a $%
\overline{\partial }$-closed form $f$ in $\mathbb{B}_{s}$ to a formula $%
\overline{\partial }\mathcal{C}_{n,s}^{p,q}f=f$ in $\mathbb{B}_{n}$ for a $%
\overline{\partial }$-closed form $f$ in $\mathbb{B}_{n}$. To accomplish
this we define \emph{ameliorated} operators $\mathcal{C}_{n,s}^{p,q}$ by%
\begin{equation*}
\mathcal{C}_{n,s}^{p,q}=\mathsf{R}_{n}\mathcal{C}_{s}^{p,q}\mathsf{E}_{s},
\end{equation*}%
where for $n<s$, $\mathsf{E}_{s}$ ($\mathsf{R}_{n}$) is the extension
(restriction) operator that takes forms $\Omega =\sum \eta
_{I,J}dw^{I}\wedge d\overline{w}^{J}$ in $\mathbb{B}_{n}$ ($\mathbb{B}_{s}$)
and extends (restricts) them to $\mathbb{B}_{s}$ ($\mathbb{B}_{n}$) by%
\begin{eqnarray*}
\mathsf{E}_{s}\left( \sum \eta _{I,J}dw^{I}\wedge d\overline{w}^{J}\right)
&\equiv &\sum \left( \eta _{I,J}\circ R\right) dw^{I}\wedge d\overline{w}%
^{J}, \\
\mathsf{R}_{n}\left( \sum \eta _{I,J}dw^{I}\wedge d\overline{w}^{J}\right)
&\equiv &\sum_{I,J\subset \left\{ 1,2,...,n\right\} }\left( \eta _{I,J}\circ
E\right) dw^{I}\wedge d\overline{w}^{J}.
\end{eqnarray*}%
Here $R$ is the natural orthogonal projection from $\mathbb{C}^{s}$ to $%
\mathbb{C}^{n}$ and $E$ is the natural embedding of $\mathbb{C}^{n}$ into $%
\mathbb{C}^{s}$. In other words, we extend a form by taking the coefficients
to be constant in the extra variables, and we restrict a form by discarding
all wedge products of differentials involving the extra variables and
restricting the coefficients accordingly.

For $s>n$ we observe that the operator $\mathcal{C}_{n,s}^{p,q}$ has
integral kernel%
\begin{equation}
\mathcal{C}_{n,s}^{p,q}\left( w,z\right) \equiv \int_{\sqrt{1-\left\vert
w\right\vert ^{2}}\mathbb{B}_{s-n}}\mathcal{C}_{s}^{p,q}\left( \left(
w,w^{\prime }\right) ,\left( z,0\right) \right) dV\left( w^{\prime }\right)
,\ \ \ \ \ z,w\in \mathbb{B}_{n},  \label{amel}
\end{equation}%
where $\mathbb{B}_{s-n}$ denotes the unit ball in $\mathbb{C}^{s-n}$ with
respect to the orthogonal decomposition $\mathbb{C}^{s}=\mathbb{C}^{n}\oplus 
\mathbb{C}^{s-n}$, and $dV$ denotes Lebesgue measure. If $f\left( w\right) $
is a $\overline{\partial }$-closed form on $\mathbb{B}_{n}$ then $f\left(
w,w^{\prime }\right) =f\left( w\right) $ is a $\overline{\partial }$-closed
form on $\mathbb{B}_{s}$ and we have for $z\in \mathbb{B}_{n}$,%
\begin{eqnarray*}
f\left( z\right) &=&f\left( z,0\right) =\overline{\partial }\int_{\mathbb{B}%
_{s}}\mathcal{C}_{s}^{p,q}\left( \left( w,w^{\prime }\right) ,\left(
z,0\right) \right) f\left( w\right) dV\left( w\right) dV\left( w^{\prime
}\right) \\
&=&\overline{\partial }\int_{\mathbb{B}_{n}}\left\{ \int_{\sqrt{1-\left\vert
w\right\vert ^{2}}\mathbb{B}_{s-n}}\mathcal{C}_{s}^{p,q}\left( \left(
w,w^{\prime }\right) ,\left( z,0\right) \right) dV\left( w^{\prime }\right)
\right\} f\left( w\right) dV\left( w\right) \\
&=&\overline{\partial }\int_{\mathbb{B}_{n}}\mathcal{C}_{n,s}^{p,q}\left(
w,z\right) f\left( w\right) dV\left( w\right) .
\end{eqnarray*}%
We have proved that%
\begin{equation*}
\mathcal{C}_{n,s}^{p,q}f\left( z\right) \equiv \int_{\mathbb{B}_{n}}\mathcal{%
C}_{n,s}^{p,q}\left( w,z\right) f\left( w\right) dV\left( w\right)
\end{equation*}%
is a right inverse for $\overline{\partial }$ on $\overline{\partial }$%
-closed forms:

\begin{theorem}
For all $s>n$ and $\overline{\partial }$-closed forms $f$ in $\mathbb{B}_{n}$%
, we have%
\begin{equation*}
\overline{\partial }\mathcal{C}_{n,s}^{p,q}f=f\text{ in }\mathbb{B}_{n}.
\end{equation*}
\end{theorem}

We will use only the case $p=0$ of this theorem and from now on we restrict
our attention to this case. The operators $\mathcal{C}_{n,s}^{0,0}$ have
been computed in \cite{OrFa} and are given by%
\begin{equation}
\mathcal{C}_{n,s}^{0,0}f\left( z\right) =\int_{\mathbb{B}_{n}}%
\sum_{j=0}^{n-1}c_{n,j,s}\frac{\left( 1-\left\vert w\right\vert ^{2}\right)
^{s-n+j}\left( 1-\left\vert z\right\vert ^{2}\right) ^{j}}{\left( 1-%
\overline{w}z\right) ^{s-n+j}\left( 1-w\overline{z}\right) ^{j}}\mathcal{C}%
_{n}^{0,0}\left( w,z\right) \wedge f\left( w\right) ,  \label{OFformula}
\end{equation}%
where%
\begin{eqnarray*}
\mathcal{C}_{n}^{0,0}\left( w,z\right) &=&c_{0}\frac{\left( 1-w\overline{z}%
\right) ^{n-1}}{\left\{ \left\vert 1-w\overline{z}\right\vert ^{2}-\left(
1-\left\vert w\right\vert ^{2}\right) \left( 1-\left\vert z\right\vert
^{2}\right) \right\} ^{n}} \\
&&\times \sum_{j=1}^{n}\left( -1\right) ^{j-1}\left( \overline{w_{j}}-%
\overline{z_{j}}\right) \dbigwedge\limits_{k\neq j}d\overline{w_{k}}%
\dbigwedge\limits_{\ell =1}^{n}dw_{\ell }.
\end{eqnarray*}%
A similar result holds for the operators $\mathcal{C}_{n,s}^{0,q}$. Define%
\begin{eqnarray*}
\Phi _{n,s}^{q}\left( w,z\right) &=&\Phi _{n}^{q}\left( w,z\right) \left( 
\frac{1-\left\vert w\right\vert ^{2}}{1-\overline{w}z}\right)
^{s-n}\sum_{j=0}^{n-q-1}c_{j,n,s}\left( \frac{\left( 1-|w|^{2}\right) \left(
1-|z|^{2}\right) }{\left\vert 1-w\overline{z}\right\vert ^{2}}\right) ^{j} \\
&=&\frac{\left( 1-w\overline{z}\right) ^{n-1-q}\left( 1-\left\vert
w\right\vert ^{2}\right) ^{q}}{\bigtriangleup \left( w,z\right) ^{n}}\left( 
\frac{1-\left\vert w\right\vert ^{2}}{1-\overline{w}z}\right)
^{s-n}\sum_{j=0}^{n-q-1}c_{j,n,s}\left( \frac{\left( 1-|w|^{2}\right) \left(
1-|z|^{2}\right) }{\left\vert 1-w\overline{z}\right\vert ^{2}}\right) ^{j} \\
&=&\sum_{j=0}^{n-q-1}c_{j,n,s}\frac{\left( 1-w\overline{z}\right)
^{n-1-q-j}\left( 1-\left\vert w\right\vert ^{2}\right) ^{s-n+q+j}\left(
1-|z|^{2}\right) ^{j}}{\left( 1-\overline{w}z\right) ^{s-n+j}\bigtriangleup
\left( w,z\right) ^{n}}.
\end{eqnarray*}%
Note that the numerator and denominator are \emph{balanced} in the sense
that the sum of the exponents in the denominator minus the corresponding sum
in the numerator (counting $\bigtriangleup \left( w,z\right) $ double) is $%
s+n+j-\left( s+j-1\right) =n+1$, the exponent of the invariant measure of
the ball $\mathbb{B}_{n}$.

\begin{theorem}
\label{explicitamel}Suppose that $s>n$ and $0\leq q\leq n-1$. Then we have%
\begin{eqnarray*}
\mathcal{C}_{n,s}^{0,q}(w,z) &=&\mathcal{C}_{n}^{0,q}(w,z)\left( \frac{%
1-\left\vert w\right\vert ^{2}}{1-\overline{w}z}\right)
^{s-n}\sum_{j=0}^{n-q-1}c_{j,n,s}\left( \frac{\left( 1-|w|^{2}\right) \left(
1-|z|^{2}\right) }{\left\vert 1-w\overline{z}\right\vert ^{2}}\right) ^{j} \\
&=&\Phi _{n,s}^{q}\left( w,z\right) \sum_{\left\vert J\right\vert
=q}\sum_{k\notin J}\left( -1\right) ^{\mu \left( k,J\right) }\left( 
\overline{z_{k}}-\overline{w_{k}}\right) d\overline{z}^{J}\wedge d\overline{w%
}^{\left( J\cup \left\{ k\right\} \right) ^{c}}\wedge \omega _{n}\left(
w\right) .
\end{eqnarray*}
\end{theorem}

\textbf{Proof}: For $s>n$ recall that the kernels of the\emph{\ }ameliorated
operators $\mathcal{C}_{n,s}^{0,q}$ are given in (\ref{amel}). For ease of
notation, we will set $k=s-n$, so we have $\mathbb{C}^{s}=\mathbb{C}%
^{n}\oplus \mathbb{C}^{k}$. Suppose that $0\leq q\leq n-1$. Recall from (\ref%
{gensolutionker}) that%
\begin{eqnarray*}
\mathcal{C}_{s}^{0,q}\left( w,z\right) &=&\left( -1\right) ^{q}\frac{\left(
1-w\overline{z}\right) ^{s-q-1}\left( 1-\left\vert w\right\vert ^{2}\right)
^{q}}{\bigtriangleup \left( w,z\right) ^{s}} \\
&&\times \sum_{\nu \in P_{s}^{q}}sgn\left( \nu \right) \left( \overline{%
w_{i_{\nu }}}-\overline{z_{i_{\nu }}}\right) \bigwedge_{j\in J_{\nu }}d%
\overline{w_{j}}\bigwedge_{l\in L_{\nu }}d\overline{z_{l}}\bigwedge \omega
_{s}\left( w\right) \\
&=&\sum_{\nu \in P_{s}^{q}}\digamma _{s,i_{\nu }}^{q}\left( w,z\right)
\bigwedge_{j\in J_{\nu }}d\overline{w_{j}}\bigwedge_{l\in L_{\nu }}d%
\overline{z_{l}}\bigwedge \omega _{s}\left( w\right) .
\end{eqnarray*}%
where%
\begin{equation*}
\digamma _{s,i_{\nu }}^{q}\left( w,z\right) =\Phi _{s}^{q}\left( w,z\right)
\left( \overline{w_{i_{\nu }}}-\overline{z_{i_{\nu }}}\right) =\frac{\left(
1-w\overline{z}\right) ^{s-q-1}\left( 1-\left\vert w\right\vert ^{2}\right)
^{q}}{\bigtriangleup \left( w,z\right) ^{s}}\left( \overline{w_{i_{\nu }}}-%
\overline{z_{i_{\nu }}}\right) .
\end{equation*}

To compute the ameliorations of these kernels, we need only focus on the
functional coefficient $\digamma _{s,i_{\nu }}^{q}\left( w,z\right) $ of the
kernel. It is easy to see that the ameliorated kernel can only give a
contribution in the variables when $1\leq i_{\nu }\leq n$, since when $%
n+1\leq i_{\nu }\leq s$ the functional kernel becomes radial in certain
variables and thus reduces to zero upon integration.

Then for any $1\leq i\leq n$ the corresponding functional coefficient $%
\digamma _{s,i}^{q}\left( w,z\right) $ has amelioration $\digamma
_{n,s,i}^{q}\left( w,z\right) $ given by 
\begin{eqnarray*}
\digamma _{n,s,i}^{q}\left( w,z\right) &=&\int_{\sqrt{1-\left\vert
w\right\vert ^{2}}\mathbb{B}_{s-n}}\digamma _{s,i}^{q}\left( (w,w^{\prime
}),\left( z,0\right) \right) dV(w^{\prime }) \\
&=&\int_{\sqrt{1-\left\vert w\right\vert ^{2}}\mathbb{B}_{k}}\frac{\left( 1-w%
\overline{z}\right) ^{s-q-1}(1-|w|^{2}-\left\vert w^{\prime }\right\vert
^{2})^{q}\left( \overline{z}_{i}-\overline{w}_{i}\right) }{\bigtriangleup
(\left( w,w^{\prime }\right) ,\left( z,0\right) )^{s}}dV\left( w^{\prime
}\right) \\
&=&\left( \overline{z}_{i}-\overline{w}_{i}\right) \left( 1-w\overline{z}%
\right) ^{s-q-1}\int_{\sqrt{1-\left\vert w\right\vert ^{2}}\mathbb{B}_{k}}%
\frac{(1-|w|^{2}-\left\vert w^{\prime }\right\vert ^{2})^{q}}{\bigtriangleup
(\left( w,w^{\prime }\right) ,\left( z,0\right) )^{s}}dV\left( w^{\prime
}\right) .
\end{eqnarray*}%
Theorem \ref{explicitamel} is a thus a consequence of the following
elementary lemma, which will find application in Section \ref{Integration by
parts} below on integration by parts as well.

\begin{lemma}
\label{amelcoeff}We have 
\begin{eqnarray*}
&&\left( 1-w\overline{z}\right) ^{s-q-1}\int_{\sqrt{1-\left\vert
w\right\vert ^{2}}\mathbb{B}_{s-n}}\frac{(1-|w|^{2}-\left\vert w^{\prime
}\right\vert ^{2})^{q}}{\bigtriangleup (\left( w,w^{\prime }\right) ,\left(
z,0\right) )^{s}}dV(w^{\prime }) \\
&=&\frac{\pi ^{s-n}}{\left( s-n\right) !}\Phi _{n}^{q}\left( w,z\right)
\left( \frac{1-\left\vert w\right\vert ^{2}}{1-\overline{w}z}\right)
^{s-n}\sum_{j=0}^{n-q-1}c_{j,n,s}\left( \frac{\left( 1-|w|^{2}\right) \left(
1-|z|^{2}\right) }{\left\vert 1-w\overline{z}\right\vert ^{2}}\right) ^{j}.
\end{eqnarray*}
\end{lemma}

See the appendix for a proof of Lemma \ref{amelcoeff}.

\section{Integration by parts\label{Integration by parts}}

We begin with an integration by parts formula involving a covariant
derivative in \cite{OrFa} (Lemma 2.1 on page 57) that reduces the
singularity of the solution kernel on the diagonal at the expense of
differentiating the form. However, in order to prepare for a generalization
to higher order forms, we replace the covariant derivative with the notion
of $\overline{\mathcal{Z}_{z,w}}$-derivative defined in (\ref{defabs}) below.

Recall Charpentier's explicit solution $\mathcal{C}_{n}^{0,0}\eta $ to the $%
\overline{\partial }$ equation $\overline{\partial }\mathcal{C}%
_{n}^{0,0}\eta =\eta $ in the ball $\mathbb{B}_{n}$ when $\eta $ is a $%
\overline{\partial }$-closed $\left( 0,1\right) $-form with coefficients in $%
C\left( \overline{\mathbb{B}_{n}}\right) $: the kernel is given by%
\begin{equation*}
\mathcal{C}_{n}^{0,0}\left( w,z\right) =c_{0}\frac{\left( 1-w\overline{z}%
\right) ^{n-1}}{\bigtriangleup \left( w,z\right) ^{n}}\sum_{j=1}^{n}\left(
-1\right) ^{j-1}\left( \overline{w_{j}}-\overline{z_{j}}\right)
\dbigwedge\limits_{k\neq j}d\overline{w_{k}}\dbigwedge\limits_{\ell
=1}^{n}dw_{\ell },
\end{equation*}%
for $\left( w,z\right) \in \mathbb{B}_{n}\times \mathbb{B}_{n}$ where%
\begin{equation*}
\bigtriangleup \left( w,z\right) =\left\vert 1-w\overline{z}\right\vert
^{2}-\left( 1-\left\vert w\right\vert ^{2}\right) \left( 1-\left\vert
z\right\vert ^{2}\right) .
\end{equation*}%
Define the Cauchy operator $\mathcal{S}_{n}$ on $\partial \mathbb{B}%
_{n}\times \mathbb{B}_{n}$ with kernel 
\begin{equation*}
\mathcal{S}_{n}\left( \zeta ,z\right) =c_{1}\frac{1}{\left( 1-\overline{%
\zeta }z\right) ^{n}}d\sigma \left( \zeta \right) ,\ \ \ \ \ \left( \zeta
,z\right) \in \partial \mathbb{B}_{n}\times \mathbb{B}_{n}.
\end{equation*}

Let $\eta =\sum_{j=1}^{n}\eta _{j}d\overline{w_{j}}$ be a $\left( 0,1\right) 
$-form with smooth coefficients. Let $\overline{\mathcal{Z}}=\overline{%
\mathcal{Z}_{z,w}}$ be the vector field acting in the variable $w=\left(
w_{1},...,w_{n}\right) $ and parameterized by $z=\left(
z_{1},...,z_{n}\right) $ given by%
\begin{equation}
\overline{\mathcal{Z}}=\overline{\mathcal{Z}_{z,w}}=\sum_{j=1}^{n}\left( 
\overline{w_{j}}-\overline{z_{j}}\right) \frac{\partial }{\partial \overline{%
w_{j}}}.  \label{defZ}
\end{equation}%
It will usually be understood from the context what the acting variable $w$
and\ the parameter variable $z$ are in $\overline{\mathcal{Z}_{z,w}}$ and we
will then omit the subscripts and simply write $\overline{\mathcal{Z}}$ for $%
\overline{\mathcal{Z}_{z,w}}$.

\begin{definition}
For $m\geq 0$, define the $m^{th}$ order derivative $\overline{\mathcal{Z}}%
^{m}\eta $ of a $\left( 0,1\right) $-form $\eta =\sum_{k=1}^{n}\eta
_{k}\left( w\right) d\overline{w_{k}}$ to be the $\left( 0,1\right) $-form
obtained by componentwise differentiation holding monomials in $\overline{w}-%
\overline{z}$ fixed: 
\begin{equation}
\overline{\mathcal{Z}}^{m}\eta \left( w\right) =\sum_{k=1}^{n}\left( 
\overline{\mathcal{Z}}^{m}\eta _{k}\right) \left( w\right) d\overline{w_{k}}%
=\sum_{k=1}^{n}\left\{ \sum_{\left\vert \alpha \right\vert =m}^{n}\left( 
\overline{w-z}\right) ^{\alpha }\frac{\partial ^{m}\eta _{k}}{\partial 
\overline{w}^{\alpha }}\left( w\right) \right\} d\overline{w_{k}}.
\label{defabs}
\end{equation}
\end{definition}

\begin{lemma}
\label{IBP1}(cf. Lemma 2.1 of \cite{OrFa}) For all $m\geq 0$ and smooth $%
\left( 0,1\right) $-forms $\eta =\sum_{k=1}^{n}\eta _{k}\left( w\right) d%
\overline{w_{k}}$, we have the formula,%
\begin{eqnarray}
\mathcal{C}_{n}^{0,0}\eta \left( z\right) &\equiv &\int_{\mathbb{B}_{n}}%
\mathcal{C}_{n}^{0,0}\left( w,z\right) \wedge \eta \left( w\right)
\label{mformula} \\
&=&\sum_{j=0}^{m-1}c_{j}\int_{\partial \mathbb{B}_{n}}\mathcal{S}_{n}\left(
w,z\right) \left( \overline{\mathcal{Z}}^{j}\eta \right) \left[ \overline{%
\mathcal{Z}}\right] \left( w\right) d\sigma \left( w\right)  \notag \\
&&+c_{m}\int_{\mathbb{B}_{n}}\mathcal{C}_{n}^{0,0}\left( w,z\right) \wedge 
\overline{\mathcal{Z}}^{m}\eta \left( w\right) .  \notag
\end{eqnarray}
\end{lemma}

Here the $\left( 0,1\right) $-form $\overline{\mathcal{Z}}^{j}\eta $ acts on
the vector field $\overline{\mathcal{Z}}$ in the usual way:%
\begin{equation*}
\left( \overline{\mathcal{Z}}^{j}\eta \right) \left[ \overline{\mathcal{Z}}%
\right] =\left( \sum_{k=1}^{n}\overline{\mathcal{Z}}^{j}\eta _{k}\left(
w\right) d\overline{w_{k}}\right) \left( \sum_{i=1}^{n}\left( \overline{w_{i}%
}-\overline{z_{i}}\right) \frac{\partial }{\partial \overline{w_{i}}}\right)
=\sum_{k=1}^{n}\left( \overline{w_{k}}-\overline{z_{k}}\right) \overline{%
\mathcal{Z}}^{j}\eta _{k}\left( w\right) .
\end{equation*}%
We can also rewrite the final integral in (\ref{mformula}) as%
\begin{equation*}
\int_{\mathbb{B}_{n}}\mathcal{C}_{n}^{0,0}\left( w,z\right) \wedge \overline{%
\mathcal{Z}}^{m}\eta \left( w\right) =\int_{\mathbb{B}_{n}}\Phi
_{n}^{0}\left( w,z\right) \left( \overline{\mathcal{Z}}^{m}\eta \right) %
\left[ \overline{\mathcal{Z}}\right] \left( w\right) dV\left( w\right) .
\end{equation*}%
See the appendix for a proof of Lemma \ref{IBP1}.

\bigskip

We now extend Lemma \ref{IBP1} to $\left( 0,q+1\right) $-forms. Let%
\begin{equation*}
\eta =\sum_{\left\vert I\right\vert =q+1}\eta _{I}\left( w\right) d\overline{%
w}^{I}
\end{equation*}%
be a $\left( 0,q+1\right) $-form with smooth coefficients. Given a $\left(
0,q+1\right) $-form $\eta =\sum_{\left\vert I\right\vert =q+1}\eta _{I}d%
\overline{w}^{I}$ and an increasing sequence $J$ of length $\left\vert
J\right\vert =q$, we define the interior product $\eta \lrcorner d\overline{w%
}^{J}$ of $\eta $ and $d\overline{w}^{J}$ by 
\begin{equation}
\eta \lrcorner d\overline{w}^{J}=\sum_{\left\vert I\right\vert =q+1}\eta
_{I}d\overline{w}^{I}\lrcorner d\overline{w}^{J}=\sum_{k\notin J}\left(
-1\right) ^{\mu \left( k,J\right) }\eta _{J\cup \left\{ k\right\} }d%
\overline{w_{k}},  \label{defintJ}
\end{equation}%
since $d\overline{w}^{I}\lrcorner d\overline{w}^{J}=\left( -1\right) ^{\mu
\left( k,J\right) }d\overline{w_{k}}$ if $k\in I\setminus J$ is the $\mu
\left( k,J\right) ^{th}$ index in $I$, and $0$ otherwise. Recall the vector
field $\overline{\mathcal{Z}}$ defined in (\ref{defZ}). The key connection
between $\eta \lrcorner d\overline{w}^{J}$ and the vector field $\overline{%
\mathcal{Z}}$ is%
\begin{eqnarray}
\left( \eta \lrcorner d\overline{w}^{J}\right) \left( \overline{\mathcal{Z}}%
\right) &=&\left( \sum_{k=1}^{n}\left( -1\right) ^{\mu \left( k,J\right)
}\eta _{J\cup \left\{ k\right\} }d\overline{w_{k}}\right) \left(
\sum_{j=1}^{n}\left( \overline{w_{j}}-\overline{z_{j}}\right) \frac{\partial 
}{\partial \overline{w_{j}}}\right)  \label{key} \\
&=&\sum_{k=1}^{n}\left( \overline{w_{k}}-\overline{z_{k}}\right) \left(
-1\right) ^{\mu \left( k,J\right) }\eta _{J\cup \left\{ k\right\} }.  \notag
\end{eqnarray}

We now define an $m^{th}$ order derivative $\overline{\mathcal{D}}^{m}\eta $
of a $\left( 0,q+1\right) $-form $\eta $ using the interior product. In the
case $q=0$ we will have $\overline{\mathcal{D}}^{m}\eta =\left( \overline{%
\mathcal{Z}}^{m}\eta \right) \left[ \overline{\mathcal{Z}}\right] $ for a $%
\left( 0,1\right) $-form $\eta $.

\begin{remark}
We are motivated by the fact that the Charpentier kernel $\mathcal{C}%
_{n}^{0,q}\left( w,z\right) $ takes $\left( 0,q+1\right) $-forms in $w$ to $%
\left( 0,q\right) $-forms in $z$. Thus in order to express the solution
operator $\mathcal{C}_{n}^{0,q}$ in terms of a volume integral rather than
the integration of a form in $w$ and $z$, our definition of $\overline{%
\mathcal{D}}^{m}\eta $, even when $m=0$, must include an appropriate
exchange of $w$-differentials for $z$-differentials.
\end{remark}

\begin{definition}
\label{defDbar}Let $m\geq 0$. For a $\left( 0,q+1\right) $-form $\eta
=\sum_{\left\vert I\right\vert =q+1}\eta _{I}d\overline{w}^{I}$ in the
variable $w$, define the $\left( 0,q\right) $-form $\overline{\mathcal{D}}%
^{m}\eta $ in the variable $z$ by 
\begin{equation*}
\overline{\mathcal{D}}^{m}\eta \left( w\right) =\sum_{\left\vert
J\right\vert =q}\overline{\mathcal{Z}}^{m}\left( \eta \lrcorner d\overline{w}%
^{J}\right) \left[ \overline{\mathcal{Z}}\right] \left( w\right) d\overline{z%
}^{J}.
\end{equation*}
\end{definition}

Again it is usually understood what the acting and parameter variables are
in $\overline{\mathcal{D}}^{m}$ but we will write $\overline{\mathcal{D}%
_{z,w}}^{m}\eta \left( w\right) $ when this may not be the case. Note that
for a $\left( 0,q+1\right) $-form $\eta =\sum_{\left\vert I\right\vert
=q+1}\eta _{I}d\overline{w}^{I}$, we have%
\begin{equation*}
\eta =\sum_{\left\vert J\right\vert =q}\left( \eta \lrcorner d\overline{w}%
^{J}\right) \wedge d\overline{w}^{J},
\end{equation*}%
and using (\ref{defabs}) the above definition yields%
\begin{eqnarray}
&&\overline{\mathcal{D}}^{m}\eta \left( w\right)  \label{effectofDm} \\
&=&\sum_{\left\vert J\right\vert =q}\overline{\mathcal{Z}}^{m}\left( \eta
\lrcorner d\overline{w}^{J}\right) \left[ \overline{\mathcal{Z}}\right]
\left( w\right) d\overline{z}^{J}  \notag \\
&=&\sum_{\left\vert J\right\vert =q}\sum_{k=1}^{n}\left( \overline{w_{k}}-%
\overline{z_{k}}\right) \left( -1\right) ^{\mu \left( k,J\right) }\left( 
\overline{\mathcal{Z}}^{m}\eta _{J\cup \left\{ k\right\} }\right) \left(
w\right) d\overline{z}^{J}  \notag \\
&=&\sum_{\left\vert J\right\vert =q}\sum_{k=1}^{n}\left( \overline{w_{k}}-%
\overline{z_{k}}\right) \left( -1\right) ^{\mu \left( k,J\right) }\left\{
\sum_{\left\vert \alpha \right\vert =m}\left( \overline{w-z}\right) ^{\alpha
}\frac{\partial ^{m}\eta _{J\cup \left\{ k\right\} }}{\partial \overline{w}%
^{\alpha }}\left( w\right) \right\} d\overline{z}^{J}.  \notag
\end{eqnarray}%
Thus the effect of $\overline{\mathcal{D}}^{m}$ on a basis element $\eta
_{I}d\overline{w}^{I}$\ is to replace a differential $d\overline{w_{k}}$
from $d\overline{w}^{I}$ ($I=J\cup \left\{ k\right\} $) with the factor $%
\left( -1\right) ^{\mu \left( k,J\right) }\left( \overline{w_{k}}-\overline{%
z_{k}}\right) $ (and this is accomplished by acting a $\left( 0,1\right) $%
-form on $\overline{\mathcal{Z}}$), replace the remaining differential $d%
\overline{w}^{J}$ with $d\overline{z}^{J}$, and then to apply the
differential operator $\overline{\mathcal{Z}}^{m}$ to the coefficient $\eta
_{I}$. We will refer to the factor $\left( \overline{w_{k}}-\overline{z_{k}}%
\right) $ introduced above as a \emph{rogue} factor since it is not
associated with a derivative $\frac{\partial }{\partial \overline{w_{k}}}$
in the way that $\left( \overline{w-z}\right) ^{\alpha }$ is associated with 
$\frac{\partial ^{m}}{\partial \overline{w}^{\alpha }}$. The point of this
distinction will be explained in Section \ref{opest} on estimates for
solution operators.

The following lemma expresses $\mathcal{C}_{n}^{0,q}\eta \left( z\right) $
in terms of integrals involving $\overline{\mathcal{D}}^{j}\eta $ for $0\leq
j\leq m$. Note that the overall effect is to reduce the singularity of the
kernel on the diagonal by $m$ factors of $\sqrt{\bigtriangleup \left(
w,z\right) }$, at the cost of increasing by $m$ the number of derivatives
hitting the form $\eta $. Recall from (\ref{defPhi}) that%
\begin{equation*}
\Phi _{n}^{\ell }\left( w,z\right) \equiv \frac{\left( 1-w\overline{z}%
\right) ^{n-1-\ell }\left( 1-\left\vert w\right\vert ^{2}\right) ^{\ell }}{%
\bigtriangleup \left( w,z\right) ^{n}}.
\end{equation*}%
We define the operator $\Phi _{n}^{\ell }$ on forms $\eta $ by%
\begin{equation*}
\Phi _{n}^{\ell }\eta \left( z\right) =\int_{\mathbb{B}_{n}}\Phi _{n}^{\ell
}\left( w,z\right) \eta \left( w\right) dV\left( w\right) .
\end{equation*}

\begin{lemma}
\label{IBP1'}Let $q\geq 0$. For all $m\geq 0$ we have the formula,%
\begin{equation}
\mathcal{C}_{n}^{0,q}\eta \left( z\right) =\sum_{k=0}^{m-1}c_{k}\mathcal{S}%
_{n}\left( \overline{\mathcal{D}}^{j}\eta \right) \left( z\right)
+\sum_{\ell =0}^{q}c_{\ell }\Phi _{n}^{\ell }\left( \overline{\mathcal{D}}%
^{m}\eta \right) \left( z\right) .  \label{mformula'}
\end{equation}
\end{lemma}

The proof is simply a reprise of that of Lemma \ref{IBP1} complicated by the
algebra that reduces matters to $\left( 0,1\right) $-forms. See the appendix.

\subsection{The radial derivative\label{The radial derivative}}

Recall the radial derivative $R=\sum_{j=1}^{n}w_{j}\frac{\partial }{\partial
w_{j}}$ from (\ref{defradder}). Here is Lemma 2.2 on page 58 of \cite{OrFa}.
See the appendix for a proof.

\begin{lemma}
\label{IBP2}Let $b>-1$. For $\Psi \in C\left( \overline{\mathbb{B}_{n}}%
\right) \cap C^{\infty }\left( \mathbb{B}_{n}\right) $ we have%
\begin{eqnarray*}
&&\int_{\mathbb{B}_{n}}\left( 1-\left\vert w\right\vert ^{2}\right) ^{b}\Psi
\left( w\right) dV\left( w\right) \\
&=&\int_{\mathbb{B}_{n}}\left( 1-\left\vert w\right\vert ^{2}\right)
^{b+1}\left( \frac{n+b+1}{b+1}I+\frac{1}{b+1}R\right) \Psi \left( w\right)
dV\left( w\right) .
\end{eqnarray*}
\end{lemma}

\begin{remark}
Typically the above lemma is applied with%
\begin{equation*}
\Psi \left( w\right) =\frac{1}{\left( 1-\overline{w}z\right) ^{s}}\psi
\left( w,z\right)
\end{equation*}%
where $z$ is a parameter in the ball $\mathbb{B}_{n}$ and%
\begin{equation*}
R\Psi \left( w\right) =\frac{1}{\left( 1-\overline{w}z\right) ^{s}}R\psi
\left( w,z\right)
\end{equation*}%
since $\frac{1}{\left( 1-\overline{w}z\right) ^{s}}$ is antiholomorphic in $%
w $.
\end{remark}

We will also need to iterate Lemma \ref{IBP2}, and for this purpose it is
convenient to introduce for $m\geq 1$ the notation%
\begin{eqnarray*}
R_{b} &=&R_{b,n}=\frac{n+b+1}{b+1}I+\frac{1}{b+1}R, \\
R_{b}^{m} &=&R_{b+m-1}R_{b+m-2}...R_{b}=\prod_{k=1}^{m}R_{b+m-k}.
\end{eqnarray*}

\begin{corollary}
\label{IBP2iter}Let $b>-1$. For $\Psi \in C\left( \overline{\mathbb{B}_{n}}%
\right) \cap C^{\infty }\left( \mathbb{B}_{n}\right) $ we have%
\begin{eqnarray*}
&&\int_{\mathbb{B}_{n}}\left( 1-\left\vert w\right\vert ^{2}\right) ^{b}\Psi
\left( w\right) dV\left( w\right) \\
&=&\int_{\mathbb{B}_{n}}\left( 1-\left\vert w\right\vert ^{2}\right)
^{b+m}R_{b}^{m}\Psi \left( w\right) dV\left( w\right) .
\end{eqnarray*}
\end{corollary}

\begin{remark}
The important point in Corollary \ref{IBP2iter} is that combinations of
radial derivatives $R$ and the identity $I$ are played off against powers of 
$1-\left\vert w\right\vert ^{2}$. It will sometimes be convenient to write
this identity as%
\begin{equation*}
\int_{\mathbb{B}_{n}}F\left( w\right) dV\left( w\right) =\int_{\mathbb{B}%
_{n}}\mathcal{R}_{b}^{m}F\left( w\right) dV\left( w\right)
\end{equation*}%
where%
\begin{equation}
\mathcal{R}_{b}^{m}\equiv \left( 1-\left\vert w\right\vert ^{2}\right)
^{b+m}R_{b}^{m}\left( 1-\left\vert w\right\vert ^{2}\right) ^{-b},
\label{calR}
\end{equation}%
and provided that $\Psi \left( w\right) =\left( 1-\left\vert w\right\vert
^{2}\right) ^{-b}F\left( w\right) $ lies in $C\left( \overline{\mathbb{B}_{n}%
}\right) \cap C^{\infty }\left( \mathbb{B}_{n}\right) $.
\end{remark}

\subsection{Integration by parts in ameliorated kernels\label{Integration by
parts in ameliorated}}

We must now extend Lemma \ref{IBP1'} and Corollary \ref{IBP2iter} to the
ameliorated kernels $\mathcal{C}_{n,s}^{0,q}$ given by%
\begin{equation*}
\mathcal{C}_{n,s}^{0,q}=\mathsf{R}_{n}\mathcal{C}_{s}^{0,q}\mathsf{E}_{s}.
\end{equation*}%
Since Corollary \ref{IBP2iter} already applies to very general functions $%
\Psi \left( w\right) $, we need only consider an extension of Lemma \ref%
{IBP1'}. The procedure for doing this is to apply Lemma \ref{IBP1'} to $%
\mathcal{C}_{s}^{0,q}$ in $s$ dimensions, and then integrate out the
additional variables using Lemma \ref{amelcoeff}.

\begin{lemma}
\label{IBPamel}Suppose that $s>n$ and $0\leq q\leq n-1$. For all $m\geq 0$
and smooth $\left( 0,q+1\right) $-forms $\eta $\ in $\overline{\mathbb{B}_{n}%
}$\ we have the formula,%
\begin{equation*}
\mathcal{C}_{n,s}^{0,q}\eta \left( z\right)
=\sum_{k=0}^{m-1}c_{k,n,s}^{\prime }\mathcal{S}_{n,s}\left( \overline{%
\mathcal{D}}^{k}\eta \right) \left[ \overline{\mathcal{Z}}\right] \left(
z\right) +\sum_{\ell =0}^{q}c_{\ell ,n,s}\Phi _{n,s}^{\ell }\left( \overline{%
\mathcal{D}}^{m}\eta \right) \left( z\right) ,
\end{equation*}%
where the ameliorated operators $\mathcal{S}_{n,s}$ and $\Phi _{n,s}^{\ell }$
have kernels given by, 
\begin{eqnarray*}
\mathcal{S}_{n,s}\left( w,z\right) &=&c_{n,s}\frac{\left( 1-|w|^{2}\right)
^{s-n-1}}{\left( 1-\overline{w}z\right) ^{s}}=c_{n,s}\left( \frac{1-|w|^{2}}{%
1-\overline{w}z}\right) ^{s-n-1}\frac{1}{\left( 1-\overline{w}z\right) ^{n+1}%
}, \\
\Phi _{n,s}^{\ell }(w,z) &=&\Phi _{n}^{\ell }\left( w,z\right) \left( \frac{%
1-\left\vert w\right\vert ^{2}}{1-\overline{w}z}\right)
^{s-n}\sum_{j=0}^{n-\ell -1}c_{j,n,s}\left( \frac{\left( 1-|w|^{2}\right)
\left( 1-|z|^{2}\right) }{\left\vert 1-w\overline{z}\right\vert ^{2}}\right)
^{j}.
\end{eqnarray*}
\end{lemma}

\textbf{Proof}: Recall that for a smooth $\left( 0,q+1\right) $-form $\eta
\left( w\right) =\sum_{\left\vert I\right\vert =q+1}\eta _{I}d\overline{w}%
^{I}$ in $\overline{\mathbb{B}_{n}}$, the $\left( 0,q\right) $-form $%
\overline{\mathcal{D}^{m}}\mathsf{E}_{s}\eta $ is given by 
\begin{eqnarray*}
\overline{\mathcal{D}^{m}}\mathsf{E}_{s}\eta \left( w\right)
&=&\sum_{\left\vert J\right\vert =q}\overline{\mathcal{D}^{m}}\left( \eta
\lrcorner d\overline{w}^{J}\right) d\overline{z}^{J}=\sum_{\left\vert
J\right\vert =q}\overline{\mathcal{D}^{m}}\left( \sum_{k\notin J}\left(
-1\right) ^{\mu \left( k,J\right) }\eta _{J\cup \left\{ k\right\} }\left(
w\right) d\overline{w_{k}}\right) d\overline{z}^{J} \\
&=&\sum_{\left\vert J\right\vert =q}\overline{\mathcal{D}^{m}}\left(
\sum_{k\notin J}\left( -1\right) ^{\mu \left( k,J\right) }\eta _{J\cup
\left\{ k\right\} }\left( w\right) d\overline{w_{k}}\right) d\overline{z}^{J}
\\
&=&\sum_{\left\vert J\right\vert =q}\sum_{k\notin J}\left( -1\right) ^{\mu
\left( k,J\right) }\left( \sum_{\left\vert \alpha \right\vert =m}\overline{%
\left( w_{k}-z_{k}\right) }\overline{\left( w-z\right) ^{\alpha }}\frac{%
\partial ^{m}}{\partial \overline{w}^{\alpha }}\eta _{J\cup \left\{
k\right\} }\left( w\right) \right) ,
\end{eqnarray*}%
where $J\cup \left\{ k\right\} $ is a multi-index with entries in $\mathfrak{%
I}_{n}\equiv \left\{ 1,2,...,n\right\} $ since the coefficient $\eta _{I}$
vanishes if $I$ is not contained in $\mathfrak{I}_{n}$. Moreover, the
multi-index $\alpha $ lies in $\left( \mathfrak{I}_{n}\right) ^{m}$ since
the coefficients $\eta _{I}$ are constant in the variable $w^{\prime
}=\left( w_{n+1},...,w_{s}\right) $. Thus 
\begin{equation*}
\overline{\mathcal{D}_{\left( z,0\right) ,\left( w,w^{\prime }\right) }^{m}}%
\mathsf{E}_{s}\eta =\overline{\mathcal{D}_{z,w}^{m}}\eta =\overline{\mathcal{%
D}^{m}}\eta ,
\end{equation*}%
and we compute that%
\begin{eqnarray*}
&&\mathsf{R}_{n}\Phi _{s}^{\ell }\left( \overline{\mathcal{D}_{\left(
z,0\right) ,\left( w,w^{\prime }\right) }^{m}}\mathsf{E}_{s}\eta \right)
\left( z\right) \\
&=&\Phi _{s}^{\ell }\left( \overline{\mathcal{D}^{m}}\eta \right) \left(
\left( z,0\right) \right) \\
&=&\sum_{\left\vert J\right\vert =q}\sum_{k\in \mathfrak{I}_{n}\setminus
J}\left( -1\right) ^{\mu \left( k,J\right) }\sum_{\left\vert \alpha
\right\vert =m}\Phi _{s}^{\ell }\left( \overline{\left( w_{k}-z_{k}\right) }%
\overline{\left( w-z\right) ^{\alpha }}\frac{\partial ^{m}}{\partial 
\overline{w}^{\alpha }}\eta _{J\cup \left\{ k\right\} }\left( \left(
w,w^{\prime }\right) \right) \right) \left( \left( z,0\right) \right) ,
\end{eqnarray*}%
where $J\cup \left\{ k\right\} \subset \mathfrak{I}_{n}$ and $\alpha \in
\left( \mathfrak{I}_{n}\right) ^{m}$ and%
\begin{eqnarray*}
&&\Phi _{s}^{\ell }\left( \overline{\left( w_{k}-z_{k}\right) }\overline{%
\left( w-z\right) ^{\alpha }}\frac{\partial ^{m}}{\partial \overline{w}%
^{\alpha }}\eta _{J\cup \left\{ k\right\} }\left( w\right) \right) \left(
\left( z,0\right) \right) \\
&=&\int_{\mathbb{B}_{s}}\frac{\left( 1-w\overline{z}\right) ^{s-1-\ell
}\left( 1-\left\vert w\right\vert ^{2}-\left\vert w^{\prime }\right\vert
^{2}\right) ^{\ell }}{\bigtriangleup \left( \left( w,w^{\prime }\right)
,\left( z,0\right) \right) ^{s}}\overline{\left( w_{k}-z_{k}\right) }%
\overline{\left( w-z\right) ^{\alpha }}\frac{\partial ^{m}}{\partial 
\overline{w}^{\alpha }}\eta _{J\cup \left\{ k\right\} }\left( w\right)
dV\left( \left( w,w^{\prime }\right) \right) \\
&=&\int_{\mathbb{B}_{n}}\left\{ \left( 1-w\overline{z}\right) ^{s-\ell
-1}\int_{\mathbb{B}_{s-n}}\frac{\left( 1-\left\vert w\right\vert
^{2}-\left\vert w^{\prime }\right\vert ^{2}\right) ^{\ell }}{\bigtriangleup
\left( \left( w,w^{\prime }\right) ,\left( z,0\right) \right) ^{s}}dV\left(
w^{\prime }\right) \right\} \\
&&\times \overline{\left( w_{k}-z_{k}\right) }\overline{\left( w-z\right)
^{\alpha }}\frac{\partial ^{m}}{\partial \overline{w}^{\alpha }}\eta _{J\cup
\left\{ k\right\} }\left( w\right) dV\left( w\right) .
\end{eqnarray*}

By Lemma \ref{amelcoeff} the term in braces above equals%
\begin{equation*}
\frac{\pi ^{s-n}}{\left( s-n\right) !}\Phi _{n}^{\ell }\left( w,z\right)
\left( \frac{1-\left\vert w\right\vert ^{2}}{1-\overline{w}z}\right)
^{s-n}\sum_{j=0}^{n-\ell -1}c_{j,n,s}\left( \frac{\left( 1-|w|^{2}\right)
\left( 1-|z|^{2}\right) }{\left\vert 1-w\overline{z}\right\vert ^{2}}\right)
^{j},
\end{equation*}%
and now performing the sum $\sum_{\left\vert J\right\vert =q}\sum_{k\in 
\mathfrak{I}_{n}\setminus J}\left( -1\right) ^{\mu \left( k,J\right)
}\sum_{\left\vert \alpha \right\vert =m}$ yields%
\begin{equation}
\mathsf{R}_{n}\Phi _{s}^{\ell }\left( \overline{\mathcal{D}_{\left(
z,0\right) }^{m}}\mathsf{E}_{s}\eta \right) \left( z\right) =\Phi _{s}^{\ell
}\left( \overline{\mathcal{D}_{z}^{m}}\eta \right) \left( \left( z,0\right)
\right) =\Phi _{n,s}^{\ell }\left( \overline{\mathcal{D}_{z}^{m}}\eta
\right) \left( z\right) .  \label{secondpart}
\end{equation}%
An even easier calculation using formula (1) in 1.4.4 on page 14 of \cite%
{Rud} shows that%
\begin{equation}
\mathsf{R}_{n}\mathcal{S}_{s}\left( \mathsf{E}_{s}\overline{\mathcal{D}%
_{z}^{k}}\eta \right) \left( \left( z,0\right) \right) =\mathcal{S}%
_{s}\left( \overline{\mathcal{D}_{z}^{k}}\eta \right) \left( \left(
z,0\right) \right) =\mathcal{S}_{n,s}\left( \overline{\mathcal{D}_{z}^{k}}%
\eta \right) \left( z\right) ,  \label{firstpart}
\end{equation}%
and now the conclusion of Lemma \ref{IBPamel} follows from (\ref{secondpart}%
), (\ref{firstpart}), the definition $\mathcal{C}_{n,s}^{0,q}=\mathsf{R}_{n}%
\mathcal{C}_{s}^{0,q}\mathsf{E}_{s}$, and Lemma \ref{IBP1'}.

\section{The Koszul complex\label{The Koszul complex}}

Here we briefly review the algebra behind the Koszul complex as presented
for example in \cite{Lin} in the finite dimensional setting. A more detailed
treatment in that setting can be found in Section 5.5.3 of \cite{Saw}. Fix $%
h $ holomorphic as in (\ref{Ngen}). Now if $g=\left( g_{j}\right)
_{j=1}^{\infty }$ satisfies $\left\vert g\right\vert ^{2}=\sum_{j=1}^{\infty
}\left\vert g_{j}\right\vert ^{2}\geq \delta ^{2}>0$, let 
\begin{equation*}
\Omega _{0}^{1}=\frac{\overline{g}}{\left\vert g\right\vert ^{2}}=\left( 
\frac{\overline{g_{j}}}{\left\vert g\right\vert ^{2}}\right) _{j=1}^{\infty
}=\left( \Omega _{0}^{1}\left( j\right) \right) _{j=1}^{\infty },
\end{equation*}%
which we view as a $1$-tensor (in $\ell ^{2}=\mathbb{C}^{\infty }$) of $%
\left( 0,0\right) $-forms with components $\Omega _{0}^{1}\left( j\right) =%
\frac{\overline{g_{j}}}{\left\vert g\right\vert ^{2}}$. Then $f=\Omega
_{0}^{1}h$ satisfies $\mathcal{M}_{g}f=f\cdot g=h$, but in general fails to
be holomorphic. The Koszul complex provides a scheme which we now recall for
solving a sequence of $\overline{\partial }$ equations that result in a
correction term $\Lambda _{g}\Gamma _{0}^{2}$ that when subtracted from $f$
above yields a holomorphic solution to the second line in (\ref{Ngen}). See
below.

The $1$-tensor of $\left( 0,1\right) $-forms $\overline{\partial }\Omega
_{0}=\left( \overline{\partial }\frac{\overline{g_{j}}}{\left\vert
g\right\vert ^{2}}\right) _{j=1}^{\infty }=\left( \overline{\partial }\Omega
_{0}^{1}\left( j\right) \right) _{j=1}^{\infty }$ is given by%
\begin{equation*}
\overline{\partial }\Omega _{0}^{1}\left( j\right) =\overline{\partial }%
\frac{\overline{g_{j}}}{\left\vert g\right\vert ^{2}}=\frac{\left\vert
g\right\vert ^{2}\overline{\partial g_{j}}-\overline{g_{j}}\overline{%
\partial }\left\vert g\right\vert ^{2}}{\left\vert g\right\vert ^{4}}=\frac{1%
}{\left\vert g\right\vert ^{4}}\sum_{k=1}^{\infty }g_{k}\overline{\left\{
g_{k}\partial g_{j}-g_{j}\partial g_{k}\right\} }.
\end{equation*}%
and can be written as%
\begin{equation*}
\overline{\partial }\Omega _{0}^{1}=\Lambda _{g}\Omega _{1}^{2}\equiv \left[
\sum_{k=1}^{\infty }\Omega _{1}^{2}\left( j,k\right) g_{k}\right]
_{j=1}^{\infty },
\end{equation*}%
where the antisymmetric $2$-tensor $\Omega _{1}^{2}$ of $\left( 0,1\right) $%
-forms is given by%
\begin{equation*}
\Omega _{1}^{2}=\left[ \Omega _{1}^{2}\left( j,k\right) \right]
_{j,k=1}^{\infty }=\left[ \frac{\overline{\left\{ g_{k}\partial
g_{j}-g_{j}\partial g_{k}\right\} }}{\left\vert g\right\vert ^{4}}\right]
_{j,k=1}^{\infty }.
\end{equation*}%
and $\Lambda _{g}\Omega _{1}^{2}$ denotes its contraction by the vector $g$
in the final variable.

We can repeat this process and by induction we have%
\begin{equation}
\overline{\partial }\Omega _{q}^{q+1}=\Lambda _{g}\Omega _{q+1}^{q+2},\ \ \
\ \ 0\leq q\leq n,  \label{induction}
\end{equation}%
where $\Omega _{q}^{q+1}$ is an alternating $\left( q+1\right) $-tensor of $%
\left( 0,q\right) $-forms. Recall that $h$ is holomorphic. When $q=n$ we
have that $\Omega _{n}^{n+1}h$ is $\overline{\partial }$-closed and this
allows us to solve a chain of $\overline{\partial }$ equations%
\begin{equation*}
\overline{\partial }\Gamma _{q-2}^{q}=\Omega _{q-1}^{q}h-\Lambda _{g}\Gamma
_{q-1}^{q+1},
\end{equation*}%
for alternating $q$-tensors $\Gamma _{q-2}^{q}$ of $\left( 0,q-2\right) $%
-forms, using the ameliorated Charpentier solution operators $\mathcal{C}%
_{n,s}^{0,q}$ defined in (\ref{amel}) above (note that our notation
suppresses the dependence of $\Gamma $ on $h$). With the convention that $%
\Gamma _{n}^{n+2}\equiv 0$ we have 
\begin{eqnarray}
\overline{\partial }\left( \Omega _{q}^{q+1}h-\Lambda _{g}\Gamma
_{q}^{q+2}\right) &=&0,\ \ \ \ \ 0\leq q\leq n,  \label{antip} \\
\overline{\partial }\Gamma _{q-1}^{q+1} &=&\Omega _{q}^{q+1}h-\Lambda
_{g}\Gamma _{q}^{q+2},\ \ \ \ \ 1\leq q\leq n.  \notag
\end{eqnarray}

Now%
\begin{equation*}
f\equiv \Omega _{0}^{1}h-\Lambda _{g}\Gamma _{0}^{2}
\end{equation*}%
is holomorphic by the first line in (\ref{antip}) with $q=0$, and since $%
\Gamma _{0}^{2}$ is antisymmetric, we compute that $\Lambda _{g}\Gamma
_{0}^{2}\cdot g=\Gamma _{0}^{2}\left( g,g\right) =0$ and%
\begin{equation*}
\mathcal{M}_{g}f=f\cdot g=\Omega _{0}^{1}h\cdot g-\Lambda _{g}\Gamma
_{0}^{2}\cdot g=h-0=h.
\end{equation*}%
Thus $f=\left( f_{i}\right) _{i=1}^{\infty }$ is a vector of holomorphic
functions satisfying the second line in (\ref{Ngen}). The first line in (\ref%
{Ngen}) is the subject of the remaining sections of the paper.

\subsection{Wedge products and factorization of the Koszul complex\label%
{Wedge products}}

Here we record the remarkable factorization of the Koszul complex in
Andersson and Carlsson \cite{AnCa}. To describe the factorization we
introduce an exterior algebra structure on $\ell ^{2}=\mathbb{C}^{\infty }$.
Let $\left\{ e_{1},e_{2},...\right\} $ be the usual basis in $\mathbb{C}%
^{\infty }$, and for an increasing multiindex $I=\left( i_{1},...,i_{\ell
}\right) $ of integers in $\mathbb{N}$, define 
\begin{equation*}
e_{I}=e_{i_{1}}\wedge e_{i_{2}}\wedge ...\wedge e_{i_{\ell }},
\end{equation*}%
where we use $\wedge $ to denote the wedge product in the exterior algebra $%
\Lambda ^{\ast }\left( \mathbb{C}^{\infty }\right) $ of $\mathbb{C}^{\infty
} $, as well as for the wedge product on forms in $\mathbb{C}^{n}$. Note
that $\left\{ e_{I}:\left\vert I\right\vert =r\right\} $ is a basis for the
alternating $r$-tensors on $\mathbb{C}^{\infty }$.

If $f=\sum_{\left\vert I\right\vert =r}f_{I}e_{I}$ is an alternating $r$%
-tensor on $\mathbb{C}^{\infty }$ with values that are $\left( 0,k\right) $%
-forms in $\mathbb{C}^{n}$, which may be viewed as a member of the exterior
algebra of $\mathbb{C}^{\infty }\otimes \mathbb{C}^{n}$, and if $%
g=\sum_{\left\vert J\right\vert =s}g_{J}e_{J}$ is an alternating $s$-tensor
on $\mathbb{C}^{\infty }$ with values that are $\left( 0,\ell \right) $%
-forms in $\mathbb{C}^{n}$, then as in \cite{AnCa} we define the wedge
product $f\wedge g$ in the exterior algebra of $\mathbb{C}^{\infty }\otimes 
\mathbb{C}^{n}$ to be the alternating $\left( r+s\right) $-tensor on $%
\mathbb{C}^{\infty }$ with values that are $\left( 0,k+\ell \right) $-forms
in $\mathbb{C}^{n}$ given by%
\begin{eqnarray}
f\wedge g &=&\left( \sum_{\left\vert I\right\vert =r}f_{I}e_{I}\right)
\wedge \left( \sum_{\left\vert J\right\vert =s}g_{J}e_{J}\right)
\label{wedge} \\
&=&\sum_{\left\vert I\right\vert =r,\left\vert J\right\vert =s}\left(
f_{I}\wedge g_{J}\right) \left( e_{I}\wedge e_{J}\right)  \notag \\
&=&\sum_{\left\vert K\right\vert =r+s}\left( \pm \sum_{I+J=K}f_{I}\wedge
g_{J}\right) e_{K}.  \notag
\end{eqnarray}%
Note that we simply write the exterior product of an element from $\Lambda
^{\ast }\left( \mathbb{C}^{\infty }\right) $ with an element from $\Lambda
^{\ast }\left( \mathbb{C}^{n}\right) $ as juxtaposition, without writing an
explicit wedge symbol. This should cause no confusion since the basis we use
in $\Lambda ^{\ast }\left( \mathbb{C}^{\infty }\right) $ is $\left\{
e_{i}\right\} _{i=1}^{\infty }$, while the basis we use in $\Lambda ^{\ast
}\left( \mathbb{C}^{n}\right) $ is $\left\{ dz_{j},d\widehat{z_{j}}\right\}
_{j=1}^{n}$, quite different in both appearance and interpretation.

In terms of this notation we then have the following factorization in
Theorem 3.1 of Andersson and Carlsson \cite{AnCa}: 
\begin{equation}
\Omega _{0}^{1}\wedge \dbigwedge\limits_{i=1}^{\ell }\widetilde{\Omega
_{0}^{1}}=\left( \sum_{k_{0}=1}^{\infty }\frac{\overline{g_{k_{0}}}}{%
\left\vert g\right\vert ^{2}}e_{k_{0}}\right) \wedge
\dbigwedge\limits_{i=1}^{\ell }\left( \sum_{k_{i}=1}^{\infty }\frac{%
\overline{\partial g_{k_{i}}}}{\left\vert g\right\vert ^{2}}e_{k_{i}}\right)
=-\frac{1}{\ell +1}\Omega _{\ell }^{\ell +1},  \label{Omegaform}
\end{equation}%
where%
\begin{equation*}
\Omega _{0}^{1}=\left( \frac{\overline{g_{i}}}{\left\vert g\right\vert ^{2}}%
\right) _{i=1}^{\infty }\ \text{and }\widetilde{\Omega _{0}^{1}}=\left( 
\frac{\overline{\partial g_{i}}}{\left\vert g\right\vert ^{2}}\right)
_{i=1}^{\infty }.
\end{equation*}%
The factorization in \cite{AnCa} is proved in the finite dimensional case,
but this extends to the infinite dimensional case by continuity. Since the $%
\ell ^{2}$ norm is quasi-multiplicative on wedge products by Lemma 5.1 in 
\cite{AnCa} we have

\begin{equation}
\left\vert \Omega _{\ell }^{\ell +1}\right\vert ^{2}\leq C_{\ell }\left\vert
\Omega _{0}^{1}\right\vert ^{2}\left\vert \widetilde{\Omega _{0}^{1}}%
\right\vert ^{2\ell },\ \ \ \ \ 0\leq \ell \leq n,  \label{quasimult}
\end{equation}%
where the constant $C_{\ell }$ depends only on the number of factors $\ell $
in the wedge product, and \emph{not} on the underlying dimension of the
vector space (which is infinite for $\ell ^{2}=\mathbb{C}^{\infty }$).

It will be useful in the next section to consider also tensor products%
\begin{equation}
\widetilde{\Omega _{0}^{1}}\otimes \widetilde{\Omega _{0}^{1}}=\left(
\sum_{i=1}^{\infty }\frac{\overline{\partial g_{i}}}{\left\vert g\right\vert
^{2}}e_{i}\right) \otimes \left( \sum_{j=1}^{\infty }\frac{\overline{%
\partial g_{j}}}{\left\vert g\right\vert ^{2}}e_{j}\right)
=\sum_{i,j=1}^{\infty }\frac{\overline{\partial g_{i}}\otimes \overline{%
\partial g_{j}}}{\left\vert g\right\vert ^{4}}e_{i}\otimes e_{j},
\label{tensor}
\end{equation}%
and more generally $\mathcal{X}^{\alpha }\widetilde{\Omega _{0}^{1}}\otimes 
\mathcal{X}^{\beta }\widetilde{\Omega _{0}^{1}}$ where $\mathcal{X}^{m}$
denotes the vector derivative defined in Definition \ref{calX} below. We
will use the fact that the $\ell ^{2}$-norm is \emph{multiplicative} on
tensor products.

\section{An almost invariant holomorphic derivative\label{An almost
invariant}}

In this section we continue to consider $\ell ^{2}$-valued spaces. We refer
the reader to \cite{ArRoSa} for the definition of the Bergman tree $\mathcal{%
T}_{n}$ and the corresponding pairwise disjoint decomposition of the ball $%
\mathbb{B}_{n}$:%
\begin{equation*}
\mathbb{B}_{n}=\overset{\cdot }{\bigcup_{\alpha \in \mathcal{T}_{n}}}%
K_{\alpha },
\end{equation*}%
where the sets $K_{\alpha }$ are comparable to balls of radius one in the
Bergman metric $\beta $ on the ball $\mathbb{B}_{n}$: $\beta \left(
z,w\right) =\frac{1}{2}\ln \frac{1+\left\vert \varphi _{z}\left( w\right)
\right\vert }{1-\left\vert \varphi _{z}\left( w\right) \right\vert }$
(Proposition 1.21 in \cite{Zhu}). This decomposition gives an analogue in $%
\mathbb{B}_{n}$ of the standard decomposition of the upper half plane $%
\mathbb{C}_{+}$ into dyadic squares whose distance from the boundary $%
\partial \mathbb{C}_{+}$ equals their side length. We also recall from \cite%
{ArRoSa} the differential operator $D_{a}$ which on the Bergman kube $%
K_{\alpha }$, and provided $a\in K_{\alpha }$, is close to the invariant
gradient $\widetilde{\nabla }$, and which has the additional property that $%
D_{a}^{m}f\left( z\right) $ is holomorphic for $m\geq 1$ and $z\in K_{\alpha
}$ when $f$ is holomorphic. For our purposes the powers $D_{a}^{m}f$, $m\geq
1$, are easier to work with than the corresponding powers $\widetilde{\nabla 
}^{m}f$, which fail to be holomorphic. It is shown in \cite{ArRoSa} that $%
D_{a}^{m}$ can be used to define an equivalent norm on the Besov space $%
B_{p}\left( \mathbb{B}_{n}\right) =B_{p}^{0}\left( \mathbb{B}_{n}\right) $,
and it is a routine matter to extend this result to the Besov-Sobolev space $%
B_{p}^{\sigma }\left( \mathbb{B}_{n}\right) $ when $\sigma \geq 0$ and $%
m>2\left( \frac{n}{p}-\sigma \right) $. The further extension to $\ell ^{2}$
-valued functions is also routine.

We define%
\begin{equation*}
\nabla _{z}=\left( \frac{\partial }{\partial z_{1}},...,\frac{\partial }{%
\partial z_{n}}\right) \text{ and }\overline{\nabla _{z}}=\left( \frac{%
\partial }{\partial \overline{z_{1}}},...,\frac{\partial }{\partial 
\overline{z_{n}}}\right)
\end{equation*}%
so that the usual Euclidean gradient is given by the pair $\left( \nabla
_{z},\overline{\nabla _{z}}\right) $. Fix $\alpha \in \mathcal{T}_{n}$ and
let $a=c_{\alpha }$. Recall that the gradient with invariant length given by 
\begin{eqnarray*}
\widetilde{\nabla }f\left( a\right) &=&\left( f\circ \varphi _{a}\right)
^{\prime }\left( 0\right) =f^{\prime }\left( a\right) \varphi _{a}^{\prime
}\left( 0\right) \\
&=&-f^{\prime }\left( a\right) \left\{ \left( 1-\left\vert a\right\vert
^{2}\right) P_{a}+\left( 1-\left\vert a\right\vert ^{2}\right) ^{\frac{1}{2}%
}Q_{a}\right\}
\end{eqnarray*}%
fails to be holomorphic in $a$. To rectify this, we define as in \cite%
{ArRoSa}, 
\begin{eqnarray}
D_{a}f\left( z\right) &=&f^{\prime }\left( z\right) \varphi _{a}^{\prime
}\left( 0\right)  \label{defDa} \\
&=&-f^{\prime }\left( z\right) \left\{ \left( 1-\left\vert a\right\vert
^{2}\right) P_{a}+\left( 1-\left\vert a\right\vert ^{2}\right) ^{\frac{1}{2}%
}Q_{a}\right\} ,  \notag
\end{eqnarray}%
for $z\in \mathbb{B}_{n}$. Note that $\nabla _{z}\left( \overline{a}\cdot
z\right) =\overline{a}^{t}$ when we view $w\in \mathbb{B}_{n}$ as an $%
n\times 1$ complex matrix, and denote by $w^{t}$ the $1\times n$ transpose
of $w$. With this interpretation, we observe that $P_{a}z=\frac{\overline{a}z%
}{\left\vert a\right\vert ^{2}}a$ has derivative $P_{a}=P_{a}^{\prime }z=%
\frac{a\overline{a}^{t}}{\left\vert a\right\vert ^{2}}=\left\vert
a\right\vert ^{-2}\left[ a_{i}\overline{a_{j}}\right] _{1\leq i,j\leq n}$.

The next lemma from \cite{ArRoSa} shows that $D_{a}^{m}$ and $D_{b}^{m}$ are
comparable when $a$ and $b$ are close in the Bergman metric.

\begin{lemma}
\label{boxcomp}(Lemma 6.2 in \cite{ArRoSa}) Let $a,b\in \mathbb{B}_{n}$
satisfy $\beta \left( a,b\right) \leq C$. There is a positive constant $%
C_{m} $ depending only on $C$ and $m$ such that 
\begin{equation*}
C_{m}^{-1}\left\vert D_{b}^{m}f\left( z\right) \right\vert \leq \left\vert
D_{a}^{m}f\left( z\right) \right\vert \leq C_{m}\left\vert D_{b}^{m}f\left(
z\right) \right\vert ,
\end{equation*}%
for all $f\in H\left( \mathbb{B}_{n};\ell ^{2}\right) $.
\end{lemma}

We remind the reader that $\left\vert D_{a}^{m}f\left( z\right) \right\vert =%
\sqrt{\sum_{i=1}^{\infty }\left\vert D_{a}^{m}f_{i}\left( z\right)
\right\vert ^{2}}$ if $f=\left( f_{i}\right) _{i=1}^{\infty }$. The scalar
proof in \cite{ArRoSa} is easily extended to $\ell ^{2}$-valued $f$.

\begin{definition}
\label{treeseminorm}(see \cite{ArRoSa}) Suppose $\sigma \geq 0$, $1<p<\infty 
$ and $m\geq 1$. We define a \textquotedblleft tree
semi-norm\textquotedblright\ $\left\Vert \cdot \right\Vert _{B_{p,m}^{\sigma
}\left( \mathbb{B}_{n};\ell ^{2}\right) }^{\ast }$ by 
\begin{equation}
\left\Vert f\right\Vert _{B_{p,m}^{\sigma }\left( \mathbb{B}_{n};\ell
^{2}\right) }^{\ast }=\left( \sum_{\alpha \in \mathcal{T}_{n}}\int_{B_{d}%
\left( c_{\alpha },C_{2}\right) }\left\vert \left( 1-\left\vert z\right\vert
^{2}\right) ^{\sigma }D_{c_{\alpha }}^{m}f\left( z\right) \right\vert
^{p}d\lambda _{n}\left( z\right) \right) ^{\frac{1}{p}}.  \label{norms'}
\end{equation}
\end{definition}

We now recall the invertible \textquotedblleft radial\textquotedblright\
operators $R^{\gamma ,t}:H\left( \mathbb{B}_{n}\right) \rightarrow H\left( 
\mathbb{B}_{n}\right) $ given in \cite{Zhu} by 
\begin{equation*}
R^{\gamma ,t}f\left( z\right) =\sum_{k=0}^{\infty }\frac{\Gamma \left(
n+1+\gamma \right) \Gamma \left( n+1+k+\gamma +t\right) }{\Gamma \left(
n+1+\gamma +t\right) \Gamma \left( n+1+k+\gamma \right) }f_{k}\left(
z\right) ,
\end{equation*}%
provided neither $n+\gamma $ nor $n+\gamma +t$ is a negative integer, and
where $f\left( z\right) =\sum_{k=0}^{\infty }f_{k}\left( z\right) $ is the
homogeneous expansion of $f$. This definition is easily extended to $f\in
H\left( \mathbb{B}_{n};\ell ^{2}\right) $. If the inverse of $R^{\gamma ,t}$
is denoted $R_{\gamma ,t}$, then Proposition 1.14 of \cite{Zhu} yields 
\begin{eqnarray}
R^{\gamma ,t}\left( \frac{1}{\left( 1-\overline{w}z\right) ^{n+1+\gamma }}%
\right) &=&\frac{1}{\left( 1-\overline{w}z\right) ^{n+1+\gamma +t}},
\label{Zhuidentity} \\
R_{\gamma ,t}\left( \frac{1}{\left( 1-\overline{w}z\right) ^{n+1+\gamma +t}}%
\right) &=&\frac{1}{\left( 1-\overline{w}z\right) ^{n+1+\gamma }},  \notag
\end{eqnarray}%
for all $w\in \mathbb{B}_{n}$. Thus for any $\gamma $, $R^{\gamma ,t}$ is
approximately differentiation of order $t$. The next proposition shows that
the derivatives $R^{\gamma ,m}f\left( z\right) $ are \textquotedblleft $%
L^{p} $ norm equivalent\textquotedblright\ to $\left\{ f\left( 0\right)
,...,\nabla ^{m-1}f\left( 0\right) ,\nabla ^{m}f\left( z\right) \right\} $
for $m$ large enough. The scalar case $\sigma =0$ is Proposition 2.1 in \cite%
{ArRoSa} and follows from Theorems 6.1 and Theorem 6.4 of \cite{Zhu}. The
extension to $\sigma \geq 0$ and $\ell ^{2}$-valued $f$ is routine. See the
appendix and also \cite{Bea}.

\begin{proposition}
\label{Bequiv}Suppose that $\sigma \geq 0$, $0<p<\infty $, $n+\gamma $ is
not a negative integer, and $f\in H\left( \mathbb{B}_{n};\ell ^{2}\right) $.
Then the following four conditions are equivalent: 
\begin{eqnarray*}
&&\left( 1-\left\vert z\right\vert ^{2}\right) ^{m+\sigma }\nabla
^{m}f\left( z\right) \in L^{p}\left( d\lambda _{n};\ell ^{2}\right) \,\text{%
for \emph{some }}m>\frac{n}{p}-\sigma ,m\in \mathbb{N}, \\
&&\left( 1-\left\vert z\right\vert ^{2}\right) ^{m+\sigma }\nabla
^{m}f\left( z\right) \in L^{p}\left( d\lambda _{n};\ell ^{2}\right) \text{
for \emph{all }}m>\frac{n}{p}-\sigma ,m\in \mathbb{N}, \\
&&\left( 1-\left\vert z\right\vert ^{2}\right) ^{m+\sigma }R^{\gamma
,m}f\left( z\right) \in L^{p}\left( d\lambda _{n};\ell ^{2}\right) \text{
for \emph{some }}m>\frac{n}{p}-\sigma ,m+n+\gamma \notin -\mathbb{N}, \\
&&\left( 1-\left\vert z\right\vert ^{2}\right) ^{m+\sigma }R^{\gamma
,m}f\left( z\right) \in L^{p}\left( d\lambda _{n};\ell ^{2}\right) \,\text{%
for \emph{all }}m>\frac{n}{p}-\sigma ,m+n+\gamma \notin -\mathbb{N}.
\end{eqnarray*}%
Moreover, with $\psi \left( z\right) =1-\left\vert z\right\vert ^{2}$, we
have for $1<p<\infty $, 
\begin{eqnarray*}
&&C^{-1}\left\Vert \psi ^{m_{1}+\sigma }R^{\gamma ,m_{1}}f\right\Vert
_{L^{p}\left( d\lambda _{n};\ell ^{2}\right) } \\
\;\;\;\;\; &\leq &\sum_{k=0}^{m_{2}-1}\left\vert \nabla ^{k}f\left( 0\right)
\right\vert +\left( \int_{\mathbb{B}_{n}}\left\vert \left( 1-\left\vert
z\right\vert ^{2}\right) ^{m_{2}+\sigma }\nabla ^{m_{2}}f\left( z\right)
\right\vert ^{p}d\lambda _{n}\left( z\right) \right) ^{\frac{1}{p}} \\
\;\;\;\;\; &\leq &C\left\Vert \psi ^{m_{1}+\sigma }R^{\gamma
,m_{1}}f\right\Vert _{L^{p}\left( d\lambda _{n};\ell ^{2}\right) }
\end{eqnarray*}%
for all $m_{1},m_{2}>\frac{n}{p}-\sigma $, $m_{1}+n+\gamma \notin -\mathbb{N}
$, $m_{2}\in \mathbb{N}$, and where the constant $C$ depends only on $\sigma 
$, $m_{1}$, $m_{2}$, $n$, $\gamma $ and $p$.
\end{proposition}

There is one further equivalent norm involving the radial derivative 
\begin{equation}
Rf\left( z\right) =z\cdot \nabla f\left( z\right) =\sum_{j=1}^{n}z_{j}\frac{%
\partial f}{\partial z_{j}}\left( z\right) ,  \label{defradder}
\end{equation}%
and its iterates $R^{k}=R\circ R\circ ...\circ R$ ($k$ times).

\begin{proposition}
\label{Bequivrad}Suppose that $\sigma \geq 0$, $0<p<\infty $ and $f\in
H\left( \mathbb{B}_{n};\ell ^{2}\right) $. Then 
\begin{eqnarray*}
&&\sum_{k=0}^{m_{1}}\left( \int_{\mathbb{B}_{n}}\left\vert \left(
1-\left\vert z\right\vert ^{2}\right) ^{m_{1}+\sigma }R^{k}f\left( z\right)
\right\vert ^{p}d\lambda _{n}\left( z\right) \right) ^{\frac{1}{p}} \\
&\approx &\sum_{k=0}^{m_{2}-1}\left\vert \nabla ^{k}f\left( 0\right)
\right\vert +\left( \int_{\mathbb{B}_{n}}\left\vert \left( 1-\left\vert
z\right\vert ^{2}\right) ^{m_{2}+\sigma }\nabla ^{m_{2}}f\left( z\right)
\right\vert ^{p}d\lambda _{n}\left( z\right) \right) ^{\frac{1}{p}}
\end{eqnarray*}%
for all $m_{1},m_{2}>\frac{n}{p}-\sigma $, $m_{1}+n+\gamma \notin -\mathbb{N}
$, $m_{2}\in \mathbb{N}$, and where the constants in the equivalence depend
only on $\sigma $, $m_{1}$, $m_{2}$, $n$ and $p$.
\end{proposition}

The seminorms $\left\Vert \cdot \right\Vert _{B_{p,m}^{\sigma }\left( 
\mathbb{B}_{n};\ell ^{2}\right) }^{\ast }$ turn out to be independent of $%
m>2\left( \frac{n}{p}-\sigma \right) $. We will obtain this fact as a
corollary of the equivalence of the standard norm\ in (\ref{defB}) with the
corresponding norm in Proposition \ref{Bequiv} using the \textquotedblleft
radial\textquotedblright\ derivative $R^{0,m}$. Note that the restriction $%
m>2\left( \frac{n}{p}-\sigma \right) $ is dictated by the fact that $%
\left\vert D_{c_{\alpha }}^{m}f\left( z\right) \right\vert $ involves the
factor $\left( 1-\left\vert z\right\vert ^{2}\right) ^{\frac{m}{2}}$ times $%
m^{th}$ order tangential derivatives of $f$, and so we must have that $%
\left( 1-\left\vert z\right\vert ^{2}\right) ^{\left( \frac{m}{2}+\sigma
\right) p}d\lambda _{n}\left( z\right) $ is a finite measure, i.e. $\left( 
\frac{m}{2}+\sigma \right) p-n-1>-1$. The case scalar $\sigma =0$ of the
following lemma is Lemma 6.4 in \cite{ArRoSa}.

\begin{lemma}
\label{radinv}Let $1<p<\infty $, $\sigma \geq 0$ and $m>2\left( \frac{n}{p}%
-\sigma \right) $. Denote by $B_{\beta }\left( c,C\right) $ the ball center $%
c$ radius $C$ in the Bergman metric $\beta $. Then for $f\in H\left( \mathbb{%
B}_{n};\ell ^{2}\right) $, 
\begin{eqnarray}
&&\left\Vert f\right\Vert _{B_{p,m}^{\sigma }\left( \mathbb{B}_{n};\mathbb{%
\ell }^{2}\right) }^{\ast }+\sum_{j=0}^{m-1}\left\vert \nabla ^{j}f\left(
0\right) \right\vert  \label{norms} \\
\;\;\;\;\; &\equiv &\left( \sum_{\alpha \in \mathcal{T}_{n}}\int_{B_{\beta
}\left( c_{\alpha },C_{2}\right) }\left\vert \left( 1-\left\vert
z\right\vert ^{2}\right) ^{\sigma }D_{c_{\alpha }}^{m}f\left( z\right)
\right\vert ^{p}d\lambda _{n}\left( z\right) \right) ^{\frac{1}{p}%
}+\sum_{j=0}^{m-1}\left\vert \nabla ^{j}f\left( 0\right) \right\vert  \notag
\\
\;\;\;\;\; &\approx &\left( \int_{\mathbb{B}_{n}}\left\vert \left(
1-\left\vert z\right\vert ^{2}\right) ^{m+\sigma }R^{\sigma ,m}f\left(
z\right) \right\vert ^{p}d\lambda _{n}\left( z\right) \right) ^{\frac{1}{p}%
}+\sum_{j=0}^{m-1}\left\vert \nabla ^{j}f\left( 0\right) \right\vert
=\left\Vert f\right\Vert _{B_{p,m}^{\sigma }\left( \mathbb{B}_{n};\mathbb{%
\ell }^{2}\right) }.  \notag
\end{eqnarray}
\end{lemma}

See the appendix for an adaptation of the proof in \cite{ArRoSa} to the case 
$\sigma \geq 0$ and $\ell ^{2}$-valued $f$.

\bigskip

We will also need to know that the pointwise multipliers in $%
M_{B_{p}^{\sigma }\left( \mathbb{B}_{n}\right) \rightarrow B_{p}^{\sigma
}\left( \mathbb{B}_{n};\ell ^{2}\right) }$ are bounded. Indeed, standard
arguments show the following.

\begin{lemma}
\label{boundedmult}%
\begin{equation}
M_{B_{p}^{\sigma }\left( \mathbb{B}_{n}\right) \rightarrow B_{p}^{\sigma
}\left( \mathbb{B}_{n};\ell ^{2}\right) }\subset H^{\infty }\left( \mathbb{B}%
_{n};\ell ^{2}\right) \cap B_{p}^{\sigma }\left( \mathbb{B}_{n};\ell
^{2}\right) .  \label{multarebounded}
\end{equation}
\end{lemma}

\textbf{Proof}: If $\varphi \in M_{B_{p}^{\sigma }\left( \mathbb{B}%
_{n}\right) \rightarrow B_{p}^{\sigma }\left( \mathbb{B}_{n};\ell
^{2}\right) }$, then $\varphi \in B_{p}^{\sigma }\left( \mathbb{B}_{n};\ell
^{2}\right) $ since $1\in B_{p}^{\sigma }\left( \mathbb{B}_{n}\right) $, and%
\begin{equation*}
\mathbb{M}_{\varphi }:B_{p}^{\sigma }\left( \mathbb{B}_{n}\right)
\rightarrow B_{p}^{\sigma }\left( \mathbb{B}_{n};\ell ^{2}\right) \text{ and 
}\mathbb{M}_{\varphi }^{\ast }:B_{p}^{\sigma }\left( \mathbb{B}_{n};\ell
^{2}\right) ^{\ast }\rightarrow B_{p}^{\sigma }\left( \mathbb{B}_{n}\right)
^{\ast }.
\end{equation*}%
If $e_{z}$ denotes point evaluation at $z\in \mathbb{B}_{n}$, $x\in \ell
^{2} $ and $f\in B_{p}^{\sigma }\left( \mathbb{B}_{n}\right) $, then the
calculation%
\begin{eqnarray*}
\left\langle f,\mathbb{M}_{\varphi }^{\ast }\left( xe_{z}\right)
\right\rangle _{B_{p}^{\sigma }\left( \mathbb{B}_{n}\right) }
&=&\left\langle \mathbb{M}_{\varphi }f,xe_{z}\right\rangle _{B_{p}^{\sigma
}\left( \mathbb{B}_{n};\ell ^{2}\right) }=\sum_{i=1}^{\infty }\left\langle
\varphi _{i}f,x_{i}e_{z}\right\rangle _{B_{p}^{\sigma }\left( \mathbb{B}%
_{n}\right) } \\
&=&\sum_{i=1}^{\infty }\overline{x_{i}}\varphi _{i}\left( z\right) f\left(
z\right) =\sum_{i=1}^{\infty }\overline{x_{i}}\varphi _{i}\left( z\right)
\left\langle f,e_{z}\right\rangle _{B_{p}^{\sigma }\left( \mathbb{B}%
_{n}\right) } \\
&=&\sum_{i=1}^{\infty }\left\langle f,\overline{\varphi _{i}\left( z\right) }%
x_{i}e_{z}\right\rangle _{B_{p}^{\sigma }\left( \mathbb{B}_{n}\right)
}=\left\langle f,\left\langle x,\varphi \left( z\right) \right\rangle _{\ell
^{2}}e_{z}\right\rangle _{B_{p}^{\sigma }\left( \mathbb{B}_{n}\right) },
\end{eqnarray*}%
shows that 
\begin{equation*}
\mathbb{M}_{\varphi }^{\ast }\left( xe_{z}\right) =\left\langle x,\varphi
\left( z\right) \right\rangle _{\mathbb{\ell }^{2}}e_{z}.
\end{equation*}%
This yields%
\begin{eqnarray*}
\left\vert \left\langle x,\varphi \left( z\right) \right\rangle _{\mathbb{%
\ell }^{2}}\right\vert \left\Vert e_{z}\right\Vert _{B_{p}^{\sigma }\left( 
\mathbb{B}_{n}\right) ^{\ast }} &=&\left\Vert \mathbb{M}_{\varphi }^{\ast
}\left( xe_{z}\right) \right\Vert _{B_{p}^{\sigma }\left( \mathbb{B}%
_{n}\right) ^{\ast }} \\
&\leq &\left\Vert \mathbb{M}_{\varphi }^{\ast }\right\Vert _{B_{p}^{\sigma
}\left( \mathbb{B}_{n};\ell ^{2}\right) ^{\ast }\rightarrow B_{p}^{\sigma
}\left( \mathbb{B}_{n}\right) ^{\ast }}\left\Vert xe_{z}\right\Vert
_{B_{p}^{\sigma }\left( \mathbb{B}_{n};\ell ^{2}\right) ^{\ast }} \\
&=&\left\Vert \mathbb{M}_{\varphi }\right\Vert _{B_{p}^{\sigma }\left( 
\mathbb{B}_{n}\right) \rightarrow B_{p}^{\sigma }\left( \mathbb{B}_{n};\ell
^{2}\right) }\left\vert x\right\vert \left\Vert e_{z}\right\Vert
_{B_{p}^{\sigma }\left( \mathbb{B}_{n}\right) ^{\ast }},
\end{eqnarray*}%
which gives%
\begin{equation*}
\left\vert \varphi \left( z\right) \right\vert =\sup_{x\neq 0}\frac{%
\left\vert \left\langle x,\varphi \left( z\right) \right\rangle _{\mathbb{%
\ell }^{2}}\right\vert }{\left\vert x\right\vert }\leq \left\Vert \mathbb{M}%
_{\varphi }\right\Vert _{B_{p}^{\sigma }\left( \mathbb{B}_{n}\right)
\rightarrow B_{p}^{\sigma }\left( \mathbb{B}_{n};\ell ^{2}\right) },\ \ \ \
\ z\in \mathbb{B}_{n},
\end{equation*}%
and completes the proof of Lemma \ref{boundedmult}.

\bigskip

In order to deal with functions $f$ on $\mathbb{B}_{n}$\ that are not
necessarily holomorphic, we use a notion of higher order derivative $D^{m}$
introduced in \cite{ArRoSa} that is based on iterating $D_{a}$ rather than $%
\widetilde{\nabla }$.

\begin{definition}
\label{Dpowers}For $m\in \mathbb{N}$ and $f\in C^{\infty }\left( \mathbb{B}%
_{n};\ell ^{2}\right) $ smooth in $\mathbb{B}_{n}$\ we define $\Theta
^{m}f\left( a,z\right) =D_{a}^{m}f\left( z\right) $ for $a,z\in \mathbb{B}%
_{n}$, and then set%
\begin{equation*}
D^{m}f\left( z\right) =\Theta ^{m}f\left( z,z\right) =D_{z}^{m}f\left(
z\right) ,\ \ \ \ \ z\in \mathbb{B}_{n}.
\end{equation*}
\end{definition}

Note that in this definition, we iterate the operator $D_{z}$ holding $z$
fixed, and then evaluate the result at the same $z$. If we combine Lemmas %
\ref{boxcomp} and \ref{radinv} we obtain that for $f\in H\left( \mathbb{B}%
_{n};\mathbb{\ell }^{2}\right) $,%
\begin{equation*}
\left\Vert f\right\Vert _{B_{p,m}^{\sigma }\left( \mathbb{B}_{n};\ell
^{2}\right) }\approx \sum_{j=0}^{m-1}\left\vert \nabla ^{j}f\left( 0\right)
\right\vert +\left( \int_{\mathbb{B}_{n}}\left\vert \left( 1-\left\vert
z\right\vert ^{2}\right) ^{\sigma }D^{m}f\left( z\right) \right\vert
^{p}d\lambda _{n}\left( z\right) \right) ^{\frac{1}{p}}.
\end{equation*}

\subsection{Real variable analogues of Besov-Sobolev spaces\label{Real
variable analogues}}

In order to handle the operators arising from integration by parts formulas
below, we will need yet more general equivalent norms on $B_{p,m}^{\sigma
}\left( \mathbb{B}_{n};\ell ^{2}\right) $.

\begin{definition}
\label{calX}We denote by $\mathcal{X}^{m}$ the vector of all differential
operators of the form $X_{1}X_{2}...X_{m}$ where each $X_{i}$ is either $%
1-\left\vert z\right\vert ^{2}$ times the identity operator $I$, the
operator $\overline{D}$, or the operator $\left( 1-\left\vert z\right\vert
^{2}\right) R$. Just as in Definition \ref{Dpowers}, we calculate the
products $X_{1}X_{2}...X_{m}$ by composing $\overline{D_{a}}$ and $\left(
1-\left\vert a\right\vert ^{2}\right) R$ and then setting $a=z$ at the end.
Note that $\overline{D_{a}}$ and $\left( 1-\left\vert a\right\vert
^{2}\right) R$ commute since the first is an antiholomorphic derivative and
the coefficient $z$ in $R=z\cdot \nabla $ is holomorphic. Similarly we
denote by $\mathcal{Y}^{m}$ the corresponding products of $\left(
1-\left\vert z\right\vert ^{2}\right) I$, $D$ (instead of $\overline{D}$)
and $\left( 1-\left\vert z\right\vert ^{2}\right) R$.
\end{definition}

In the iterated derivative $\mathcal{X}^{m}$ we are differentiating only
with the \emph{antiholomorphic} derivative $\overline{D}$ or the \emph{%
holomorphic} derivative $R$. When $f$ is holomorphic, we thus have $\mathcal{%
X}^{m}f\sim \left\{ \left( 1-\left\vert z\right\vert ^{2}\right)
^{m}R^{k}f\right\} _{k=0}^{m}$. The reason we allow $1-\left\vert
z\right\vert ^{2}$ times the identity $I$ to occur in $\mathcal{X}^{m}$ is
that this produces a norm (as opposed to just a seminorm) without including
the term $\sum_{k=0}^{m-1}\left\vert \nabla ^{k}f\left( 0\right) \right\vert 
$. We define the norm $\left\Vert \cdot \right\Vert _{B_{p,m}^{\sigma
}\left( \mathbb{B}_{n};\mathbb{\ell }^{2}\right) }$ for \emph{smooth} $f$ on
the ball $\mathbb{B}_{n}$ by%
\begin{equation*}
\left\Vert f\right\Vert _{B_{p,m}^{\sigma }\left( \mathbb{B}_{n};\ell
^{2}\right) }\equiv \left( \sum_{k=0}^{m}\int_{\mathbb{B}_{n}}\left\vert
\left( 1-\left\vert z\right\vert ^{2}\right) ^{m+\sigma }R^{k}f\left(
z\right) \right\vert ^{p}d\lambda _{n}\left( z\right) \right) ^{\frac{1}{p}},
\end{equation*}%
and note that provided $m+\sigma >\frac{n}{p}$, this provides an equivalent
norm for the Besov-Sobolev space $B_{p}^{\sigma }\left( \mathbb{B}_{n};\ell
^{2}\right) $ of holomorphic functions on $\mathbb{B}_{n}$. These
considerations motivate the following two definitions of a \emph{%
real-variable} analogue of the norm $\left\Vert \cdot \right\Vert
_{B_{p,m}^{\sigma }\left( \mathbb{B}_{n};\ell ^{2}\right) }$.

\begin{definition}
We define the \emph{norms} $\left\Vert \cdot \right\Vert _{\Lambda
_{p,m}^{\sigma }\left( \mathbb{B}_{n};\ell ^{2}\right) }$ and $\left\Vert
\cdot \right\Vert _{\Phi _{p,m}^{\sigma }\left( \mathbb{B}_{n};\ell
^{2}\right) }$\ for $f=\left( f_{i}\right) _{i=1}^{\infty }$ smooth on the
ball $\mathbb{B}_{n}$ by 
\begin{eqnarray}
\left\Vert f\right\Vert _{\Lambda _{p,m}^{\sigma }\left( \mathbb{B}_{n};\ell
^{2}\right) } &\equiv &\left( \int_{\mathbb{B}_{n}}\left\vert \left(
1-\left\vert z\right\vert ^{2}\right) ^{\sigma }\mathcal{X}^{m}f\left(
z\right) \right\vert ^{p}d\lambda _{n}\left( z\right) \right) ^{\frac{1}{p}},
\label{defLambda} \\
\left\Vert f\right\Vert _{\Phi _{p,m}^{\sigma }\left( \mathbb{B}_{n};\ell
^{2}\right) } &\equiv &\left( \int_{\mathbb{B}_{n}}\left\vert \left(
1-\left\vert z\right\vert ^{2}\right) ^{\sigma }\mathcal{Y}^{m}f\left(
z\right) \right\vert ^{p}d\lambda _{n}\left( z\right) \right) ^{\frac{1}{p}}.
\notag
\end{eqnarray}
\end{definition}

It is \emph{not} true that either of the norms $\left\Vert \cdot \right\Vert
_{\Lambda _{p,m}^{\sigma }\left( \mathbb{B}_{n};\ell ^{2}\right) }$ or $%
\left\Vert \cdot \right\Vert _{\Phi _{p,m}^{\sigma }\left( \mathbb{B}%
_{n};\ell ^{2}\right) }$ are independent of $m$ for large $m$ when acting on
smooth functions. However, Lemmas \ref{boxcomp} and \ref{radinv} show the
equivalence of norms when restricted to holomorphic vector functions:

\begin{lemma}
Let $1<p<\infty $, $\sigma \geq 0$ and $m>2\left( \frac{n}{p}-\sigma \right) 
$. Then for $f$ a holomorphic vector function we have 
\begin{equation}
\left\Vert f\right\Vert _{B_{p,m}^{\sigma }\left( \mathbb{B}_{n};\ell
^{2}\right) }\approx \left\Vert f\right\Vert _{\Lambda _{p,m}^{\sigma
}\left( \mathbb{B}_{n};\ell ^{2}\right) }\approx \left\Vert f\right\Vert
_{\Phi _{p,m}^{\sigma }\left( \mathbb{B}_{n};\ell ^{2}\right) .}
\label{allnormsequal}
\end{equation}
\end{lemma}

The norms $\left\Vert \cdot \right\Vert _{\Lambda _{p,m}^{\sigma }\left( 
\mathbb{B}_{n};\ell ^{2}\right) }$ arise in the integration by parts in
iterated Charpentier kernels in Section \ref{opest}, while the norms $%
\left\Vert \cdot \right\Vert _{\Phi _{p,m}^{\sigma }\left( \mathbb{B}%
_{n};\ell ^{2}\right) }$ are useful for estimating the holomorphic function $%
g$ in the Koszul complex. For this latter purpose we will use the following
multilinear inequality whose scalar version is, after translating notation,
Theorem 3.5 in \cite{OrFa}. The extension to $\ell ^{2}$-valued functions is
routine but again, for the convenience of the reader, we give a detailed
proof in the appendix.

\begin{proposition}
\label{multilinear}Suppose that $1<p<\infty $, $0\leq \sigma <\infty $, $%
M\geq 1$, $m>2\left( \frac{n}{p}-\sigma \right) $ and $\alpha =\left( \alpha
_{0},...,\alpha _{M}\right) \in \mathbb{Z}_{+}^{M+1}$ with $\left\vert
\alpha \right\vert =m$. For $g\in M_{B_{p}^{\sigma }\left( \mathbb{B}%
_{n}\right) \rightarrow B_{p}^{\sigma }\left( \mathbb{B}_{n};\ell
^{2}\right) }$ and $h\in B_{p}^{\sigma }\left( \mathbb{B}_{n}\right) $ we
have,%
\begin{eqnarray*}
&&\int_{\mathbb{B}_{n}}\left( 1-\left\vert z\right\vert ^{2}\right)
^{p\sigma }\left\vert \left( \mathcal{Y}^{\alpha _{1}}g\right) \left(
z\right) \right\vert ^{p}...\left\vert \left( \mathcal{Y}^{\alpha
_{M}}g\right) \left( z\right) \right\vert ^{p}\left\vert \left( \mathcal{Y}%
^{\alpha _{0}}h\right) \left( z\right) \right\vert ^{p}d\lambda _{n}\left(
z\right) \\
&&\ \ \ \ \ \ \ \ \ \ \leq C_{n,M,\sigma ,p}\left( \left\Vert \mathbb{M}%
_{g}\right\Vert _{B_{p}^{\sigma }\left( \mathbb{B}_{n}\right) \rightarrow
B_{p}^{\sigma }\left( \mathbb{B}_{n};\ell ^{2}\right) }^{Mp}\right)
\left\Vert h\right\Vert _{B_{p}^{\sigma }\left( \mathbb{B}_{n}\right) }^{p}.
\end{eqnarray*}
\end{proposition}

\begin{remark}
The inequalities for $M=1$ in Proposition \ref{multilinear} actually
characterize multipliers $g$ in the sense that a function $g\in
B_{p}^{\sigma }\left( \mathbb{B}_{n};\ell ^{2}\right) \cap H^{\infty }\left( 
\mathbb{B}_{n};\ell ^{2}\right) $ is in $M_{B_{p}^{\sigma }\left( \mathbb{B}%
_{n}\right) \rightarrow B_{p}^{\sigma }\left( \mathbb{B}_{n};\ell
^{2}\right) }$ \emph{if and only if} the inequalities with $M=1$ in
Proposition \ref{multilinear} hold. This follows from noting that each term
in the Leibniz expansion of $\mathcal{Y}^{m}\left( gh\right) $ occurs on the
left side of the display above with $M=1$.
\end{remark}

\subsubsection{Three crucial inequalities\label{Three crucial inequalities}}

In order to establish appropriate inequalities for the Charpentier solution
operators, we will need to control terms of the form $\left( \overline{z-w}%
\right) ^{\alpha }\frac{\partial ^{m}}{\partial \overline{w}^{\alpha }}%
F\left( w\right) $, $\overline{D_{\left( z\right) }^{m}}\bigtriangleup
\left( w,z\right) $, $\overline{D_{\left( z\right) }^{m}}\left\{ \left( 1-w%
\overline{z}\right) ^{k}\right\} $and $R_{\left( z\right) }^{m}\left\{
\left( 1-\overline{w}z\right) ^{k}\right\} $ inside the integral for $T$ as
given in the integration by parts formula in Lemma \ref{IBP1'}\ above. Here
we are using the subscript $\left( z\right) $ in parentheses to indicate the
variable being differentiated. This is to avoid confusion with the notation $%
D_{a}$ introduced in (\ref{defDa}). We collect the necessary estimates in
the following proposition.

\begin{proposition}
\label{threecrucial}For $z,w\in \mathbb{B}_{n}$ and $m\in \mathbb{N}$, we
have the following three crucial estimates:%
\begin{equation}
\left\vert \left( \overline{z-w}\right) ^{\alpha }\frac{\partial ^{m}}{%
\partial \overline{w}^{\alpha }}F\left( w\right) \right\vert \leq C\left( 
\frac{\sqrt{\bigtriangleup \left( w,z\right) }}{1-\left\vert w\right\vert
^{2}}\right) ^{m}\left\vert \overline{D}^{m}F\left( w\right) \right\vert ,\
\ \ \ \ F\in H\left( \mathbb{B}_{n};\ell ^{2}\right) ,m=\left\vert \alpha
\right\vert .  \label{moddelta}
\end{equation}%
\begin{eqnarray}
\left\vert D_{\left( z\right) }\bigtriangleup \left( w,z\right) \right\vert
&\leq &C\left\{ \left( 1-\left\vert z\right\vert ^{2}\right) \bigtriangleup
\left( w,z\right) ^{\frac{1}{2}}+\bigtriangleup \left( w,z\right) \right\} ,
\label{rootD} \\
\left\vert \left( 1-\left\vert z\right\vert ^{2}\right) R_{\left( z\right)
}\bigtriangleup \left( w,z\right) \right\vert &\leq &C\left( 1-\left\vert
z\right\vert ^{2}\right) \sqrt{\bigtriangleup \left( w,z\right) },  \notag
\end{eqnarray}%
\begin{eqnarray}
\left\vert D_{\left( z\right) }^{m}\left\{ \left( 1-\overline{w}z\right)
^{k}\right\} \right\vert &\leq &C\left\vert 1-\overline{w}z\right\vert
^{k}\left( \frac{1-\left\vert z\right\vert ^{2}}{\left\vert 1-\overline{w}%
z\right\vert }\right) ^{\frac{m}{2}},  \label{Dbound} \\
\left\vert \left( 1-\left\vert z\right\vert ^{2}\right) ^{m}R_{\left(
z\right) }^{m}\left\{ \left( 1-\overline{w}z\right) ^{k}\right\} \right\vert
&\leq &C\left\vert 1-\overline{w}z\right\vert ^{k}\left( \frac{1-\left\vert
z\right\vert ^{2}}{\left\vert 1-\overline{w}z\right\vert }\right) ^{m}. 
\notag
\end{eqnarray}
\end{proposition}

\textbf{Proof}: To prove (\ref{moddelta}) we view $D_{a}$ as a
differentiation operator in the variable $w$ so that%
\begin{equation*}
D_{a}=-\nabla _{w}\left\{ \left( 1-\left\vert a\right\vert ^{2}\right) P_{a}+%
\sqrt{1-\left\vert a\right\vert ^{2}}Q_{a}\right\} .
\end{equation*}%
A basic calculation is then:%
\begin{eqnarray*}
&&\left( 1-\overline{a}z\right) \varphi _{a}\left( z\right) \cdot \left(
D_{a}\right) ^{t} \\
&&\ \ \ \ \ =\left\{ P_{a}\left( z-a\right) +\sqrt{1-\left\vert a\right\vert
^{2}}Q_{a}\left( z-a\right) \right\} \\
&&\ \ \ \ \ \ \ \ \ \ \ \ \ \ \ \cdot \left\{ \left( 1-\left\vert
a\right\vert ^{2}\right) P_{a}\nabla _{w}+\sqrt{1-\left\vert a\right\vert
^{2}}Q_{a}\nabla _{w}\right\} \\
&&\ \ \ \ \ =P_{a}\left( z-a\right) \left( 1-\left\vert a\right\vert
^{2}\right) P_{a}\nabla _{w} \\
&&\ \ \ \ \ \ \ \ \ \ \ \ \ \ \ +\sqrt{1-\left\vert a\right\vert ^{2}}%
Q_{a}\left( z-a\right) \sqrt{1-\left\vert a\right\vert ^{2}}Q_{a}\nabla _{w}
\\
&&\ \ \ \ \ =\left( 1-\left\vert a\right\vert ^{2}\right) \left( z-a\right)
\cdot \nabla _{w}.
\end{eqnarray*}%
From this we conclude the inequality%
\begin{eqnarray*}
\left\vert \left( z_{i}-a_{i}\right) \frac{\partial }{\partial w_{i}}F\left(
w\right) \right\vert &\leq &\left\vert \left( z-a\right) \cdot \nabla
F\left( w\right) \right\vert \\
&\leq &\left\vert \frac{1-\overline{a}z}{1-\left\vert a\right\vert ^{2}}%
\varphi _{a}\left( z\right) \right\vert \left\vert D_{a}F\left( w\right)
\right\vert \\
&=&\frac{\sqrt{\bigtriangleup \left( a,z\right) }}{1-\left\vert a\right\vert
^{2}}\left\vert D_{a}F\left( w\right) \right\vert ,
\end{eqnarray*}%
as well as its conjugate%
\begin{equation*}
\left\vert \overline{\left( z_{i}-a_{i}\right) }\frac{\partial }{\partial 
\overline{w_{i}}}F\left( w\right) \right\vert \leq C\frac{\sqrt{%
\bigtriangleup \left( a,z\right) }}{1-\left\vert a\right\vert ^{2}}%
\left\vert \overline{D_{a}}F\left( w\right) \right\vert .
\end{equation*}%
Moreover, we can iterate this inequality to obtain%
\begin{equation*}
\left\vert \left( \overline{z-a}\right) ^{\alpha }\frac{\partial ^{m}}{%
\partial \overline{w}^{\alpha }}F\left( w\right) \right\vert \leq C\left( 
\frac{\sqrt{\bigtriangleup \left( a,z\right) }}{1-\left\vert a\right\vert
^{2}}\right) ^{m}\left\vert \left( \overline{D_{a}}\right) ^{m}F\left(
w\right) \right\vert ,
\end{equation*}%
for a multi-index of length $m$. With $a=w$ this becomes the first estimate (%
\ref{moddelta}).

\bigskip

To see the second estimate (\ref{rootD}), recall from (\ref{defDa}) that%
\begin{equation*}
D_{a}f\left( z\right) =-\left\{ \left( 1-\left\vert a\right\vert ^{2}\right)
P_{a}\nabla f\left( z\right) +\left( 1-\left\vert a\right\vert ^{2}\right) ^{%
\frac{1}{2}}Q_{a}\nabla f\left( z\right) \right\} .
\end{equation*}%
We let $a=z$. By the unitary invariance of 
\begin{equation*}
\bigtriangleup \left( w,z\right) =\left\vert 1-\overline{w}z\right\vert
^{2}-\left( 1-\left\vert z\right\vert ^{2}\right) \left( 1-\left\vert
w\right\vert ^{2}\right) ,
\end{equation*}%
we may assume that $z=\left( \left\vert z\right\vert ,0,...,0\right) $. Then
we have%
\begin{eqnarray*}
\frac{\partial }{\partial z_{j}}\bigtriangleup \left( w,z\right) &=&\frac{%
\partial }{\partial z_{j}}\left\{ \left( 1-\overline{w}z\right) \left( 1-%
\overline{z}w\right) -\left( 1-\overline{z}z\right) \left( 1-\left\vert
w\right\vert ^{2}\right) \right\} \\
&=&-\overline{w_{j}}\left( 1-\overline{z}w\right) +\overline{z_{j}}\left(
1-\left\vert w\right\vert ^{2}\right) \\
&=&\left( \overline{z_{j}}-\overline{w_{j}}\right) +\overline{w_{j}}\left( 
\overline{z}w\right) -\overline{z_{j}}\left\vert w\right\vert ^{2} \\
&=&\left( \overline{z_{j}}-\overline{w_{j}}\right) \left( 1-\left\vert
z\right\vert ^{2}\right) +\overline{z_{j}}\left\vert z\right\vert ^{2}-%
\overline{w_{j}}\left\vert z\right\vert ^{2}+\overline{w_{j}}\left( 
\overline{z}w\right) -\overline{z_{j}}\left\vert w\right\vert ^{2} \\
&=&\left( \overline{z_{j}}-\overline{w_{j}}\right) \left( 1-\left\vert
z\right\vert ^{2}\right) +\overline{z_{j}}\left( \left\vert z\right\vert
^{2}-\left\vert w\right\vert ^{2}\right) +\overline{w_{j}}\left( \overline{z}%
\left( w-z\right) \right) .
\end{eqnarray*}

Now $Q_{z}\nabla f=\left( 0,\frac{\partial f}{\partial z_{2}},...,\frac{%
\partial f}{\partial z_{n}}\right) $ and thus a typical term in $Q_{z}\nabla
\bigtriangleup $ is $\frac{\partial }{\partial z_{j}}\bigtriangleup \left(
w,z\right) $ with $j\geq 2$. From $z=\left( \left\vert z\right\vert
,0,...,0\right) $ and $j\geq 2$ we have $z_{j}=0$ and so 
\begin{equation*}
\frac{\partial }{\partial z_{j}}\bigtriangleup \left( w,z\right) =\left( 
\overline{z_{j}}-\overline{w_{j}}\right) \left( 1-\left\vert z\right\vert
^{2}\right) -\left( \overline{z_{j}}-\overline{w_{j}}\right) \left( 
\overline{z}\left( w-z\right) \right) ,\ \ \ \ \ j\geq 2.
\end{equation*}%
Now (\ref{manyfaces}) implies%
\begin{equation}
\bigtriangleup \left( w,z\right) =\left( 1-\left\vert z\right\vert
^{2}\right) \left\vert w-z\right\vert ^{2}+\left\vert \overline{z}%
(w-z)\right\vert ^{2},  \label{Del}
\end{equation}%
which together with the above shows that%
\begin{eqnarray}
\sqrt{1-\left\vert z\right\vert ^{2}}\left\vert Q_{z}\nabla \bigtriangleup
\left( w,z\right) \right\vert &\leq &C\left\vert z-w\right\vert \left(
1-\left\vert z\right\vert ^{2}\right) ^{\frac{3}{2}}  \label{Qest} \\
&&+C\sqrt{1-\left\vert z\right\vert ^{2}}\left\vert z-w\right\vert
\left\vert \overline{z}\left( w-z\right) \right\vert  \notag \\
&\leq &C\left( 1-\left\vert z\right\vert ^{2}\right) \bigtriangleup \left(
w,z\right) ^{\frac{1}{2}}+C\bigtriangleup \left( w,z\right) .  \notag
\end{eqnarray}

As for $P_{z}\nabla D=\left( \frac{\partial f}{\partial z_{1}}%
,0,...,0\right) $ we use (\ref{Del}) to obtain 
\begin{eqnarray*}
\left\vert P_{z}\nabla \bigtriangleup \left( w,z\right) \right\vert
&=&\left\vert \left( \overline{z_{1}}-\overline{w_{1}}\right) \left(
1-\left\vert z\right\vert ^{2}\right) +\overline{z_{1}}\left( \left\vert
z\right\vert ^{2}-\left\vert w\right\vert ^{2}\right) +\overline{w_{1}}%
\overline{z}\left( w-z\right) \right\vert \\
&\leq &\left\vert z-w\right\vert \left( 1-\left\vert z\right\vert
^{2}\right) +\left\vert \left\vert z\right\vert ^{2}-\left\vert w\right\vert
^{2}\right\vert +\left\vert \overline{z}\left( w-z\right) \right\vert \\
&\leq &C\sqrt{\bigtriangleup \left( w,z\right) }+2\left\vert \left\vert
z\right\vert -\left\vert w\right\vert \right\vert .
\end{eqnarray*}%
However,%
\begin{eqnarray*}
\bigtriangleup \left( w,z\right) &\geq &\left( 1-\left\vert w\right\vert
\left\vert z\right\vert \right) ^{2}-\left( 1-\left\vert z\right\vert
^{2}\right) \left( 1-\left\vert w\right\vert ^{2}\right) \\
&=&1-2\left\vert w\right\vert \left\vert z\right\vert +\left\vert
w\right\vert ^{2}\left\vert z\right\vert ^{2}-\left\{ 1-\left\vert
z\right\vert ^{2}-\left\vert w\right\vert ^{2}+\left\vert z\right\vert
^{2}\left\vert w\right\vert ^{2}\right\} \\
&=&\left\vert z\right\vert ^{2}+\left\vert w\right\vert ^{2}-2\left\vert
w\right\vert \left\vert z\right\vert =\left( \left\vert z\right\vert
-\left\vert w\right\vert \right) ^{2}
\end{eqnarray*}%
and so altogether we have the estimate%
\begin{equation}
\left\vert P_{z}\nabla \bigtriangleup \left( w,z\right) \right\vert \leq C%
\sqrt{\bigtriangleup \left( w,z\right) }.  \label{Pest}
\end{equation}%
Combining (\ref{Qest}) and (\ref{Pest}) with the definition (\ref{defDa})
completes the proof of the first line in (\ref{rootD}). The second line in (%
\ref{rootD}) follows from (\ref{Pest}) since $R_{\left( z\right)
}=P_{z}\nabla $.

\bigskip

To prove the third estimate (\ref{Dbound}) we compute:%
\begin{eqnarray*}
D_{\left( z\right) }\left( 1-\overline{w}z\right) ^{k} &=&k\left( 1-%
\overline{w}z\right) ^{k-1}D_{\left( z\right) }\left( 1-\overline{w}z\right)
\\
&=&k\left( 1-\overline{w}z\right) ^{k-1}\left\{ \left( 1-\left\vert
z\right\vert ^{2}\right) P_{z}\nabla +\sqrt{1-\left\vert z\right\vert ^{2}}%
Q_{z}\nabla \right\} \left( 1-\overline{w}z\right) \\
&=&-k\left( 1-\overline{w}z\right) ^{k-1}\left\{ \left( 1-\left\vert
z\right\vert ^{2}\right) P_{z}\overline{w}+\sqrt{1-\left\vert z\right\vert
^{2}}Q_{z}\overline{w}\right\} ; \\
R_{\left( z\right) }\left( 1-\overline{w}z\right) ^{k} &=&k\left( 1-%
\overline{w}z\right) ^{k-1}\left( -\overline{w}z\right) .
\end{eqnarray*}%
Since $\left\vert w\right\vert ^{2}+\left\vert a\right\vert ^{2}\leq 2$ we
have%
\begin{eqnarray*}
\left\vert Q_{z}\overline{w}\right\vert ^{2} &=&\left\vert Q_{z}\left( 
\overline{w}-\overline{z}\right) \right\vert ^{2}\leq \left\vert \overline{w}%
-\overline{z}\right\vert ^{2}, \\
&=&\left\vert w\right\vert ^{2}+\left\vert z\right\vert ^{2}-2\func{Re}%
\left( w\overline{z}\right) \\
&\leq &2\func{Re}\left( 1-w\overline{z}\right) \leq 2\left\vert 1-w\overline{%
z}\right\vert ,
\end{eqnarray*}%
which yields%
\begin{eqnarray*}
\left\vert D_{\left( z\right) }\left\{ \left( 1-\overline{w}z\right)
^{k}\right\} \right\vert &\leq &C\left\vert 1-\overline{w}z\right\vert
^{k}\left\{ \frac{\left( 1-\left\vert z\right\vert ^{2}\right) +\sqrt{\left(
1-\left\vert z\right\vert ^{2}\right) \left\vert 1-\overline{w}z\right\vert }%
}{\left\vert 1-\overline{w}z\right\vert }\right\} \\
&\leq &C\left\vert 1-\overline{w}z\right\vert ^{k}\sqrt{\frac{1-\left\vert
z\right\vert ^{2}}{\left\vert 1-\overline{w}z\right\vert }}.
\end{eqnarray*}%
Iteration then yields (\ref{Dbound}).

\section{Schur's Test\label{Schur's Test section}}

Here we characterize boundedness of the positive operators that arise as
majorants of the solution operators below. The case $c=0$ of the following
lemma is Theorem 2.10 in \cite{Zhu}.

\begin{lemma}
\label{Zlemma}Let $a,b,c,t\in \mathbb{R}$. Then the operator%
\begin{equation*}
T_{a,b,c}f\left( z\right) =\int_{\mathbb{B}_{n}}\frac{\left( 1-\left\vert
z\right\vert ^{2}\right) ^{a}\left( 1-\left\vert w\right\vert ^{2}\right)
^{b}\left( \sqrt{\bigtriangleup \left( w,z\right) }\right) ^{c}}{\left\vert
1-w\overline{z}\right\vert ^{n+1+a+b+c}}f\left( w\right) dV\left( w\right)
\end{equation*}%
is bounded on $L^{p}\left( \mathbb{B}_{n};\left( 1-\left\vert w\right\vert
^{2}\right) ^{t}dV\left( w\right) \right) $ if and only if $c>-2n$ and 
\begin{equation}
-pa<t+1<p\left( b+1\right) .  \label{indexcondition}
\end{equation}
\end{lemma}

We sketch the proof for the case $c\neq 0$ when $p=2$ and $t=-n-1$. Let $%
\psi _{\varepsilon }\left( \zeta \right) =\left( 1-\left\vert \zeta
\right\vert ^{2}\right) ^{\varepsilon }$ and recall that $\sqrt{%
\bigtriangleup \left( w,z\right) }=\left\vert 1-w\overline{z}\right\vert
\left\vert \varphi _{z}\left( w\right) \right\vert $. We compute conditions
on $a$, $b$, $c$ and $\varepsilon $ such that we have%
\begin{equation*}
T_{a,b,c}\psi _{\varepsilon }\left( z\right) \leq C\psi _{\varepsilon
}\left( z\right) \text{ and }T_{a,b,c}^{\ast }\psi _{\varepsilon }\left(
w\right) \leq C\psi _{\varepsilon }\left( w\right) ,\ \ \ \ \ z,w\in \mathbb{%
B}_{n},
\end{equation*}%
where $T_{a,b,c}^{\ast }$ denotes the dual relative to $L^{2}\left( \lambda
_{n}\right) $. For this we take $\varepsilon \in \mathbb{R}$ and compute%
\begin{equation*}
T_{a,b,c}\psi _{\varepsilon }\left( z\right) =\int_{\mathbb{B}_{n}}\frac{%
\left( 1-\left\vert z\right\vert ^{2}\right) ^{a}\left( 1-\left\vert
w\right\vert ^{2}\right) ^{n+1+b+\varepsilon }\left\vert \varphi _{z}\left(
w\right) \right\vert ^{c}}{\left\vert 1-w\overline{z}\right\vert ^{n+1+a+b}}%
d\lambda _{n}\left( w\right) .
\end{equation*}

Note that the integral defining $T_{a,b,c}\psi _{\varepsilon }\left(
z\right) $ is finite if and only if $\varepsilon >-b-1$. Now in this
integral make the change of variable $w=\varphi _{z}\left( \zeta \right) $
and use that $\lambda _{n}$ is invariant to obtain%
\begin{equation*}
T_{a,b,c}\psi _{\varepsilon }\left( z\right) =\int_{\mathbb{B}_{n}}\frac{%
\left( 1-|z|^{2}\right) ^{a}\left( 1-\left\vert \varphi _{z}\left( \zeta
\right) \right\vert ^{2}\right) ^{n+1+b+\varepsilon }\left\vert \zeta
\right\vert ^{c}}{\left\vert 1-\overline{\varphi _{z}\left( \zeta \right) }%
z\right\vert ^{n+1+a+b}\left( 1-|\zeta |^{2}\right) ^{n+1}}dV\left( \zeta
\right) .
\end{equation*}%
Plugging the identities%
\begin{eqnarray*}
1-\varphi _{z}\left( \zeta \right) \overline{z} &=&1-\left\langle \varphi
_{z}\left( \zeta \right) ,\varphi _{z}\left( 0\right) \right\rangle =\frac{%
1-|z|^{2}}{1-\zeta \overline{z}}, \\
1-\left\vert \varphi _{z}\left( \zeta \right) \right\vert ^{2}
&=&1-\left\langle \varphi _{z}\left( \zeta \right) ,\varphi _{z}\left( \zeta
\right) \right\rangle =\frac{\left( 1-\left\vert z\right\vert ^{2}\right)
\left( 1-\left\vert \zeta \right\vert ^{2}\right) }{\left\vert 1-\zeta 
\overline{z}\right\vert ^{2}},
\end{eqnarray*}
into the formula for $T_{a,b,c}\psi _{\varepsilon }\left( z\right) $ we
obtain%
\begin{equation*}
T_{a,b,c}\psi _{\varepsilon }\left( z\right) =\psi _{\varepsilon }\left(
z\right) \int_{\mathbb{B}_{n}}\frac{\left( 1-\left\vert \zeta \right\vert
^{2}\right) ^{b+\varepsilon }\left\vert \zeta \right\vert ^{c}}{\left\vert
1-\zeta \overline{z}\right\vert ^{n+1+b-a+2\varepsilon }}dV\left( \zeta
\right) .
\end{equation*}

Then from Theorem 1.12 in \cite{Zhu} we obtain that%
\begin{equation*}
\sup_{z\in \mathbb{B}_{n}}\int_{\mathbb{B}_{n}}\frac{\left( 1-\left\vert
\zeta \right\vert ^{2}\right) ^{\alpha }}{\left\vert 1-\zeta \overline{z}%
\right\vert ^{\beta }}dV\left( \zeta \right) <\infty
\end{equation*}%
if and only if $\beta -\alpha <n+1$. Provided $c>-2n$ it is now easy to see
that we also have%
\begin{equation*}
\sup_{z\in \mathbb{B}_{n}}\int_{\mathbb{B}_{n}}\frac{\left( 1-\left\vert
\zeta \right\vert ^{2}\right) ^{\alpha }\left\vert \zeta \right\vert ^{c}}{%
\left\vert 1-\zeta \overline{z}\right\vert ^{\beta }}dV\left( \zeta \right)
<\infty
\end{equation*}%
if and only if $\beta -\alpha <n+1$. It now follows from the above that 
\begin{equation*}
T_{a,b,c}\psi _{\varepsilon }\left( z\right) \leq C\psi _{\varepsilon
}\left( z\right) ,\ \ \ \ \ z\in \mathbb{B}_{n},
\end{equation*}%
if and only if%
\begin{equation*}
-b-1<\varepsilon <a.
\end{equation*}

Arguing as above and provided $c>-2n$, we obtain%
\begin{equation*}
T_{a,b,c}^{\ast }\psi _{\varepsilon }\left( w\right) \leq C\psi
_{\varepsilon }\left( w\right) ,\ \ \ \ \ w\in \mathbb{B}_{n},
\end{equation*}%
if and only if%
\begin{equation*}
-a+n<\varepsilon <b+n+1.
\end{equation*}

Altogether then there is $\varepsilon \in \mathbb{R}$ such that $h=\sqrt{%
\psi _{\varepsilon }}$ is a Schur function for $T_{a,b,c}$ on $L^{2}\left(
\lambda _{n}\right) $ in Schur's Test (as given in Theorem 2.9 on page 51 of 
\cite{Zhu}) if and only if%
\begin{equation*}
\max \left\{ -a+n,-b-1\right\} <\min \left\{ a,b+n+1\right\} .
\end{equation*}%
This is equivalent to $-2a<-n<2\left( b+1\right) $, which is (\ref%
{indexcondition}) in the case $p=2,t=-n-1$. This completes the proof (in
this case) that (\ref{indexcondition}) implies the boundedness of $T_{a,b,c}$
on $L^{2}\left( \lambda _{n}\right) $. The converse is easy - see for
example the argument for the case $c=0$ on page 52 of \cite{Zhu}.

See the appendix for a more detailed proof of Lemma \ref{Zlemma}.

\begin{remark}
\label{Zlemma'}We will also use the trivial consequence of Lemma \ref{Zlemma}
that the operator%
\begin{equation*}
T_{a,b,c,d}f\left( z\right) =\int_{\mathbb{B}_{n}}\frac{\left( 1-\left\vert
z\right\vert ^{2}\right) ^{a}\left( 1-\left\vert w\right\vert ^{2}\right)
^{b}\left( \sqrt{\bigtriangleup \left( w,z\right) }\right) ^{c}}{\left\vert
1-w\overline{z}\right\vert ^{n+1+a+b+c+d}}f\left( w\right) dV\left( w\right)
\end{equation*}%
is bounded on $L^{p}\left( \mathbb{B}_{n};\left( 1-\left\vert w\right\vert
^{2}\right) ^{t}dV\left( w\right) \right) $ if $c>-2n$, $d\leq 0$ and (\ref%
{indexcondition}) holds. This is simply because $\left\vert 1-w\overline{z}%
\right\vert \leq 2$.
\end{remark}

\section{Operator estimates\label{opest}}

We must show that $f=\Omega _{0}^{1}h-\Lambda _{g}\Gamma _{0}^{2}\in
B_{p}^{\sigma }\left( \mathbb{B}_{n};\ell ^{2}\right) $ where $\Gamma
_{0}^{2}$ is an antisymmetric $2$-tensor of $\left( 0,0\right) $-forms that
solves%
\begin{equation*}
\overline{\partial }\Gamma _{0}^{2}=\Omega _{1}^{2}h-\Lambda _{g}\Gamma
_{1}^{3},
\end{equation*}%
and inductively where $\Gamma _{q}^{q+2}$ is an alternating $\left(
q+2\right) $-tensor of $\left( 0,q\right) $-forms that solves%
\begin{equation*}
\overline{\partial }\Gamma _{q}^{q+2}=\Omega _{q+1}^{q+2}h-\Lambda
_{g}\Gamma _{q+1}^{q+3},
\end{equation*}%
up to $q=n-1$ (since $\Gamma _{n}^{n+2}=0$ and the $\left( 0,n\right) $-form 
$\Omega _{n}^{n+1}$ is $\overline{\partial }$-closed). Using the Charpentier
solution operators $\mathcal{C}_{n,s}^{0,q}$ on $\left( 0,q+1\right) $-forms
we then get

\begin{eqnarray*}
f &=&\Omega _{0}^{1}h-\Lambda _{g}\Gamma _{0}^{2} \\
&=&\Omega _{0}^{1}h-\Lambda _{g}\mathcal{C}_{n,s_{1}}^{0,0}\left( \Omega
_{1}^{2}h-\Lambda _{g}\Gamma _{1}^{3}\right) \\
&=&\Omega _{0}^{1}h-\Lambda _{g}\mathcal{C}_{n,s_{1}}^{0,0}\left( \Omega
_{1}^{2}h-\Lambda _{g}\mathcal{C}_{n,s_{2}}^{0,1}\left( \Omega
_{2}^{3}h-\Lambda _{g}\Gamma _{2}^{4}\right) \right) \\
&&\vdots \\
&=&\Omega _{0}^{1}h-\Lambda _{g}\mathcal{C}_{n,s_{1}}^{0,0}\Omega
_{1}^{2}h+\Lambda _{g}\mathcal{C}_{n,s_{1}}^{0,0}\Lambda _{g}\mathcal{C}%
_{n,s_{2}}^{0,1}\Omega _{2}^{3}h-\Lambda _{g}\mathcal{C}_{n,s_{1}}^{0,0}%
\Lambda _{g}\mathcal{C}_{n,s_{2}}^{0,1}\Lambda _{g}\mathcal{C}%
_{n,s_{3}}^{0,2}\Omega _{3}^{4}h-... \\
&&+\left( -1\right) ^{n}\Lambda _{g}\mathcal{C}_{n,s_{1}}^{0,0}...\Lambda
_{g}\mathcal{C}_{n,s_{n}}^{0,n-1}\Omega _{n}^{n+1}h \\
&\equiv &\mathcal{F}^{0}+\mathcal{F}^{1}+...+\mathcal{F}^{n}.
\end{eqnarray*}

The goal is to establish%
\begin{equation*}
\left\Vert f\right\Vert _{B_{p}^{\sigma }\left( \mathbb{B}_{n};\mathbb{\ell }%
^{2}\right) }\leq C_{n,\sigma ,p,\delta }\left( g\right) \left\Vert
h\right\Vert _{B_{p}^{\sigma }\left( \mathbb{B}_{n}\right) },
\end{equation*}%
which we accomplish by showing that%
\begin{equation}
\left\Vert \mathcal{F}^{\mu }\right\Vert _{B_{p,m_{1}}^{\sigma }\left( 
\mathbb{B}_{n};\mathbb{\ell }^{2}\right) }\leq C_{n,\sigma ,p,\delta }\left(
g\right) \left\Vert h\right\Vert _{\Lambda _{p,m_{\mu }}^{\sigma }\left( 
\mathbb{B}_{n}\right) },\ \ \ \ \ 0\leq \mu \leq n,  \label{accomplishFmu}
\end{equation}%
for a choice of integers $m_{\mu }$ satisfying 
\begin{equation*}
\frac{n}{p}-\sigma <m_{1}<m_{2}<...<m_{\ell }<...<m_{n}.
\end{equation*}%
Recall that we defined both of the norms $\left\Vert F\right\Vert
_{B_{p,m_{\mu }}^{\sigma }\left( \mathbb{B}_{n};\mathbb{\ell }^{2}\right) }$
and $\left\Vert F\right\Vert _{\Lambda _{p,m_{\mu }}^{\sigma }\left( \mathbb{%
B}_{n};\mathbb{\ell }^{2}\right) }$ for smooth vector functions $F$ in the
ball $\mathbb{B}_{n}$.

\begin{description}
\item[Note on constants] We often indicate via subscripts, such as $n,\sigma
,p,\delta $, the important parameters on which a given constant $C$ depends,
especially when the constant appears in a basic inequality. However, at
times in mid-argument, we will often revert to suppressing some or all of
the subscripts in the interests of readability.
\end{description}

The norms $\left\Vert \cdot \right\Vert _{\Lambda _{p,m}^{\sigma }\left( 
\mathbb{B}_{n};\mathbb{\ell }^{2}\right) }$ in (\ref{defLambda}) above will
now be used to estimate the composition of Charpentier solution operators in
each function%
\begin{equation*}
\mathcal{F}^{\mu }=\Lambda _{g}\mathcal{C}_{n,s_{1}}^{0,0}...\Lambda _{g}%
\mathcal{C}_{n,s_{\mu }}^{0,\mu -1}\Omega _{\mu }^{\mu +1}h
\end{equation*}%
as follows. More precisely we will use the specialized variants of the
seminorms given by%
\begin{equation*}
\left\Vert F\right\Vert _{\Lambda _{p,m^{\prime },m^{\prime \prime
}}^{\sigma }\left( \mathbb{B}_{n};\mathbb{\ell }^{2}\right) }^{p}\equiv
\int_{\mathbb{B}_{n}}\left\vert \left( 1-\left\vert z\right\vert ^{2}\right)
^{\sigma }\left\{ \left( 1-\left\vert z\right\vert ^{2}\right) ^{m^{\prime
}}R^{m^{\prime }}\right\} \overline{D}^{m^{\prime \prime }}F\left( z\right)
\right\vert ^{p}d\lambda _{n}\left( z\right) ,
\end{equation*}%
where we take $m^{\prime \prime }$ derivatives in $\overline{D}$ followed by 
$m^{\prime }$ derivatives in the invariant radial operator $\left(
1-\left\vert z\right\vert ^{2}\right) R$. Recall from Definition \ref{calX}
that $\mathcal{X}^{m}$ denotes the vector of all differential operators of
the form $X_{1}X_{2}...X_{m}$ where each $X_{i}$ is either $I$, $\overline{D}
$, or $\left( 1-\left\vert z\right\vert ^{2}\right) R$, and where by
definition $1-\left\vert z\right\vert ^{2}$ is held constant in composing
operators. It will also be convenient at times to use the notation%
\begin{equation}
\mathcal{R}^{m}\equiv \left( 1-\left\vert z\right\vert ^{2}\right)
^{m}\left( R^{k}\right) _{k=0}^{m},  \label{shorthand}
\end{equation}%
which should cause no confusion with the related operators $\mathcal{R}%
_{b}^{m}$ in (\ref{calR}) introduced in the remark following Corollary \ref%
{IBP2iter}. Note that $\mathcal{R}^{m}$ is simply $\mathcal{X}^{m}$ when
none of the operators $\overline{D}$ appear. We will make extensive use the
multilinear estimate in Proposition \ref{multilinear}.

Let us fix our attention on the function $\mathcal{F}^{\mu }=\mathcal{F}%
_{0}^{\mu }$ and write%
\begin{eqnarray*}
\mathcal{F}_{0}^{\mu } &=&\Lambda _{g}\mathcal{C}_{n,s_{1}}^{0,0}\left\{
\Lambda _{g}\mathcal{C}_{n,s_{2}}^{0,1}...\Lambda _{g}\mathcal{C}_{n,s_{\mu
}}^{0,\mu -1}\Omega _{\mu }^{\mu +1}h\right\} =\Lambda _{g}\mathcal{C}%
_{n,s_{1}}^{0,0}\left\{ \mathcal{F}_{1}^{\mu }\right\} , \\
\mathcal{F}_{1}^{\mu } &=&\Lambda _{g}\mathcal{C}_{n,s_{2}}^{0,1}\left\{
\Lambda _{g}\mathcal{C}_{n,s_{3}}^{0,2}...\Lambda _{g}\mathcal{C}_{n,s_{\mu
}}^{0,\mu -1}\Omega _{\mu }^{\mu +1}h\right\} =\Lambda _{g}\mathcal{C}%
_{n,s_{2}}^{0,1}\left\{ \mathcal{F}_{2}^{\mu }\right\} , \\
\mathcal{F}_{q}^{\mu } &=&\Lambda _{g}\mathcal{C}_{n,s_{q+1}}^{0,q}\left\{ 
\mathcal{F}_{q+1}^{\mu }\right\} ,\ \ \ \ \ etc,
\end{eqnarray*}%
where $\mathcal{F}_{q}^{\mu }$ is a $\left( 0,q\right) $-form. We now
perform the integration by parts in Lemma \ref{IBPamel} in each iterated
Charpentier operator $\mathcal{F}_{q}^{\mu }=\Lambda _{g}\mathcal{C}%
_{n,s_{q+1}}^{0,q}\left\{ \mathcal{F}_{q+1}^{\mu }\right\} $ to obtain%
\begin{eqnarray}
\mathcal{F}_{q}^{\mu } &=&\Lambda _{g}\mathcal{C}_{n,s_{q+1}}^{0,q}\mathcal{F%
}_{q+1}^{\mu }  \label{Fterms} \\
&=&\sum_{j=0}^{m_{q+1}^{\prime }-1}c_{j,n,s_{q+1}}^{\prime }\Lambda _{g}%
\mathcal{S}_{n,s_{q+1}}\left( \overline{\mathcal{D}}^{j}\mathcal{F}%
_{q+1}^{\mu }\right) \left( z\right)  \notag \\
&&+\sum_{\ell =0}^{\mu }c_{\ell ,n,s_{q+1}}\Lambda _{g}\Phi
_{n,s_{q+1}}^{\ell }\left( \overline{\mathcal{D}}^{m_{q+1}^{\prime }}%
\mathcal{F}_{q+1}^{\mu }\right) \left( z\right) .  \notag
\end{eqnarray}

Now we compose these formulas for $\mathcal{F}_{k}^{\mu }$ to obtain an
expression for $\mathcal{F}^{\mu }$ that is a complicated sum of
compositions of the individual operators in (\ref{Fterms}) above. For now we
will concentrate on the main terms $\Lambda _{g}\Phi _{n,s_{k+1}}^{\mu
}\left( \overline{\mathcal{D}}^{m_{k+1}^{\prime }}\mathcal{F}_{k+1}^{\mu
}\right) $ that arise in the second sum above when $\ell =\mu $. We will see
that the same considerations apply to any of the other terms in (\ref{Fterms}%
). Recall from Lemma \ref{IBPamel} that the "boundary" operators $\mathcal{S}%
_{n,s_{q+1}}$ are projections of operators on $\partial \mathbb{B}_{s_{q}}$
to the ball $\mathbb{B}_{n}$ and have (balanced) kernels even simpler than
those of the operators $\Phi _{n,s_{q+1}}^{\ell }$. The composition of these
main terms is%
\begin{eqnarray}
&&\left( \Lambda _{g}\Phi _{n,s_{1}}^{\mu }\overline{\mathcal{D}}%
^{m_{1}^{\prime }}\right) \mathcal{F}_{1}^{\mu }  \label{fmucomp} \\
&=&\left( \Lambda _{g}\Phi _{n,s_{1}}^{\mu }\overline{\mathcal{D}}%
^{m_{1}^{\prime }}\right) \left( \Lambda _{g}\Phi _{n,s_{2}}^{\mu }\overline{%
\mathcal{D}}^{m_{2}^{\prime }}\right) \mathcal{F}_{2}^{\mu }  \notag \\
&=&\left( \Lambda _{g}\Phi _{n,s_{1}}^{\mu }\overline{\mathcal{D}}%
^{m_{1}^{\prime }}\right) \left( \Lambda _{g}\Phi _{n,s_{2}}^{\mu }\overline{%
\mathcal{D}}^{m_{2}^{\prime }}\right) ...\left( \Lambda _{g}\Phi _{n,s_{\mu
}}^{\mu }\overline{\mathcal{D}}^{m_{\mu }^{\prime }}\right) \Omega _{\mu
}^{\mu +1}h.  \notag
\end{eqnarray}

At this point we would like to take absolute values inside all of these
integrals and use the crucial inequalities in Proposition \ref{threecrucial}
to obtain a composition of positive operators of the type considered in
Lemma \ref{Zlemma}. However, there is a difficulty in using the crucial
inequality (\ref{moddelta}) to estimate the derivative $\overline{\mathcal{D}%
}^{m}$ on $\left( 0,q+1\right) $-forms $\eta $ given by (\ref{effectofDm}):%
\begin{equation*}
\overline{\mathcal{D}}^{m}\eta \left( z\right) =\sum_{\left\vert
J\right\vert =q}\sum_{k\notin J}\sum_{\left\vert \alpha \right\vert
=m}\left( -1\right) ^{\mu \left( k,J\right) }\overline{\left(
w_{k}-z_{k}\right) }\overline{\left( w-z\right) ^{\alpha }}\frac{\partial
^{m}}{\partial \overline{w}^{\alpha }}\eta _{J\cup \left\{ k\right\} }\left(
w\right) .
\end{equation*}%
The problem is that the factor $\overline{\left( w_{k}-z_{k}\right) }$ has
no derivative $\frac{\partial }{\partial \overline{w_{k}}}$ naturally
associated with it, as do the other factors in $\overline{\left( w-z\right)
^{\alpha }}$. We refer to the factor $\overline{\left( w_{k}-z_{k}\right) }$
as a \emph{rogue} factor, as it requires special treatment in order to apply
(\ref{moddelta}). Note that we cannot simply estimate $\overline{\left(
w_{k}-z_{k}\right) }$ by $\left\vert w-z\right\vert $ because this is much
larger in general than the estimate $\sqrt{\bigtriangleup \left( w,z\right) }
$ obtained in (\ref{moddelta}) (where the difference in size between $%
\left\vert w-z\right\vert $ and $\sqrt{\bigtriangleup \left( w,z\right) }$
is compensated by the difference in size between $\frac{\partial }{\partial 
\overline{w_{k}}}$ and $\overline{D}$).

\bigskip

We now describe how to circumvent this difficulty in the composition of
operators in (\ref{fmucomp}). Let us write each $\overline{\mathcal{D}}%
^{m_{q+1}^{\prime }}\mathcal{F}_{q+1}^{\mu }$ as%
\begin{equation*}
\sum_{\left\vert J\right\vert =q}\sum_{k\notin J}\sum_{\left\vert \alpha
\right\vert =m_{q+1}^{\prime }}\left( -1\right) ^{\mu \left( k,J\right) }%
\overline{\left( w_{k}-z_{k}\right) }\overline{\left( w-z\right) ^{\alpha }}%
\frac{\partial ^{m}}{\partial \overline{w}^{\alpha }}\left( \mathcal{F}%
_{q+1}^{\mu }\right) _{J\cup \left\{ k\right\} }\left( w\right) ,
\end{equation*}%
where $\left( \mathcal{F}_{q+1}^{\mu }\right) _{J\cup \left\{ k\right\} }$
is the coefficient of the form $\mathcal{F}_{q+1}^{\mu }$ with differential $%
d\overline{w}^{J\cup \left\{ k\right\} }$. We now replace each of these sums
with just one of the summands, say%
\begin{equation}
\overline{\left( w_{k}-z_{k}\right) }\overline{\left( w-z\right) ^{\alpha }}%
\frac{\partial ^{m}}{\partial \overline{w}^{\alpha }}\left( \mathcal{F}%
_{q+1}^{\mu }\right) _{J\cup \left\{ k\right\} }\left( w\right) .
\label{justone}
\end{equation}%
Here the factor $\overline{\left( w_{k}-z_{k}\right) }$ is a \emph{rogue}
factor, not associated with a corresponding derivative $\frac{\partial }{%
\partial \overline{w_{k}}}$. We will refer to $k$ as the \emph{rogue }index
associated with the \emph{rogue} factor when it is not convenient to
explicitly display the variables.

The key fact in treating the \emph{rogue} factor $\overline{\left(
w_{k}-z_{k}\right) }$ is that its presence in (\ref{justone}) means that the
coefficient $\left( \mathcal{F}_{q+1}^{\mu }\right) _{I}$ $\ $of the form $%
\mathcal{F}_{q+1}^{\mu }$ that multiplies it \emph{must} have $k$ in the
multi-index $I$. Since $\mathcal{F}_{q+1}^{\mu }=\Lambda _{g}\mathcal{C}%
_{n,s_{q+2}}^{0,q+1}\left\{ \mathcal{F}_{q+2}^{\mu }\right\} $, the form of
the ameliorated Charpentier kernel $\mathcal{C}_{n,s_{q+2}}^{0,q+1}$ in
Theorem \ref{explicitamel} shows that the coefficients of $\mathcal{C}%
_{n,s_{q+2}}^{0,q+1}\left( w,z\right) $ that multiply the \emph{rogue}
factor \emph{must} have the differential $d\overline{z_{k}}$ in them. In
turn, this means that the differential $d\overline{w_{k}}$ must be \emph{%
missing} in the coefficient of $\mathcal{C}_{n,s_{q+2}}^{0,q+1}\left(
w,z\right) $, and hence finally that the coefficients $\left( \mathcal{F}%
_{q+2}^{\mu }\right) _{H}$ with multi-index $H$ that survive the wedge
products in the integration \emph{must} have $k\in H$. This observation can
be repeated, and we now derive an important consequence.

Returning to (\ref{fmucomp}), each summand in $\overline{\mathcal{D}}%
^{m_{q+1}^{\prime }}\mathcal{F}_{q+1}^{\mu }$ has a \emph{rogue} factor with
associated \emph{rogue} index $k_{q+1}$. Thus the function in (\ref{fmucomp}%
) is a sum of terms of the form%
\begin{eqnarray*}
&&\left( \Lambda _{g}\Phi _{n,s_{1}}^{\mu }\overline{\left(
w_{k_{1}}-z_{k_{1}}\right) }\overline{\mathcal{Z}}^{m_{1}^{\prime }}\right)
\circ \left( \Lambda _{g}\Phi _{n,s_{2}}^{\mu }\overline{\left(
w_{k_{2}}-z_{k_{2}}\right) }\overline{\mathcal{Z}}^{m_{2}^{\prime }}\right)
_{I_{1}}\circ \\
&&\ \ \ \ \ \ \ \ \ \ ...\circ \left( \Lambda _{g}\Phi _{n,s_{\nu }}^{\nu }%
\overline{\left( w_{k_{\nu }}-z_{k_{\nu }}\right) }\overline{\mathcal{Z}}%
^{m_{\nu }^{\prime }}\right) _{I_{\nu -1}}\circ \\
&&\ \ \ \ \ \ \ \ \ \ ...\circ \left( \Lambda _{g}\Phi _{n,s_{\mu }}^{\mu -1}%
\overline{\left( w_{k_{\mu }}-z_{k_{\mu }}\right) }\overline{\mathcal{Z}}%
^{m_{\mu }^{\prime }}\right) _{I_{\mu -1}}\circ \left( \Omega _{\mu }^{\mu
+1}h\right) _{I_{\mu }},
\end{eqnarray*}%
where the subscript $I_{\nu }$ on the form $\Lambda _{g}\Phi _{n,s_{\nu
}}^{\nu }\overline{\left( w_{k_{\nu }}-z_{k_{\nu }}\right) }\overline{%
\mathcal{Z}}^{m_{\nu }^{\prime }}$ indicates that we are composing with the
component of $\Lambda _{g}\Phi _{n,s_{\nu }}^{\nu }\overline{\left(
w_{k_{\nu }}-z_{k_{\nu }}\right) }\overline{\mathcal{Z}}^{m_{\nu }^{\prime
}} $ corresponding to the multi-index $I_{\nu -1}$, i.e. the component with
the differential $d\overline{z}^{I_{\nu -1}}$. The notation will become
exceedingly unwieldy if we attempt to identify the different variables
associated with each of the iterated integrals, so we refrain from this in
general. The considerations of the previous paragraph now show that we must
have $\left\{ k_{1}\right\} =I_{1}$, $\left\{ k_{2}\right\} \cup I_{1}=I_{2}$
and more generally 
\begin{equation*}
\left\{ k_{\nu }\right\} \cup I_{\nu -1}=I_{\nu },\ \ \ \ \ 1<\nu \leq \mu .
\end{equation*}%
In particular we see that the associated \emph{rogue} indices $%
k_{1},k_{2},...k_{\mu }$ are all distinct and that as sets%
\begin{equation*}
\left\{ k_{1},k_{2},...,k_{\mu }\right\} =I_{\mu }.
\end{equation*}

If we denote by $\zeta $ the variable in the final form $\Omega _{\mu }^{\mu
+1}h$, we can thus write each \emph{rogue} factor $\overline{\left(
w_{k_{\nu }}-z_{k_{\nu }}\right) }$ as 
\begin{equation*}
\overline{\left( w_{k_{\nu }}-z_{k_{\nu }}\right) }=\overline{\left(
w_{k_{\nu }}-\zeta _{k_{\nu }}\right) }-\overline{\left( z_{k_{\nu }}-\zeta
_{k_{\nu }}\right) },
\end{equation*}%
and since $k_{\nu }\in I_{\mu }$, there is a factor of the form $\frac{%
\partial }{\partial \overline{\zeta _{k_{\nu }}}}\frac{\partial ^{\left\vert
\beta \right\vert }g_{i}}{\partial \overline{\zeta }^{\beta }}$ in each
summand of the component $\left( \Omega _{\mu }^{\mu +1}h\right) _{I_{\mu }}$
of $\Omega _{\mu }^{\mu +1}h$. So we are able to associate the \emph{rogue}
factor $\overline{\left( w_{k_{\nu }}-z_{k_{\nu }}\right) }$ with
derivatives of $g$ as follows:%
\begin{equation}
\left\{ \overline{\left( w_{k_{\nu }}-\zeta _{k_{\nu }}\right) }\frac{%
\partial }{\partial \overline{\zeta _{k_{\nu }}}}\right\} \frac{\partial
^{\left\vert \beta \right\vert }g_{i}}{\partial \overline{\zeta }^{\beta }}%
-\left\{ \overline{\left( z_{k_{\nu }}-\zeta _{k_{\nu }}\right) }\frac{%
\partial }{\partial \overline{\zeta _{k_{\nu }}}}\right\} \frac{\partial
^{\left\vert \gamma \right\vert }g_{j}}{\partial \overline{\zeta }^{\gamma }}%
.  \label{spreadrogue}
\end{equation}

Thus it is indeed possible to

\begin{enumerate}
\item apply the radial integration by parts in Corollary \ref{IBP2iter},

\item then take absolute values and $\ell ^{2}$-norms inside all the
integrals,

\item and then apply the crucial inequalities in Proposition \ref%
{threecrucial}.
\end{enumerate}

One of the difficulties remaining after this is that we are now left with
additional factors of the form%
\begin{equation*}
\frac{\sqrt{\bigtriangleup \left( w,\zeta \right) }}{1-\left\vert
w\right\vert ^{2}},\frac{\sqrt{\bigtriangleup \left( z,\zeta \right) }}{%
1-\left\vert z\right\vert ^{2}}
\end{equation*}%
that result from an application of (\ref{moddelta}) to the derivatives in (%
\ref{spreadrogue}). These factors are still \emph{rogue} in the sense that
the variable pairs occurring in them, namely $\left( w,\zeta \right) $ and $%
\left( z,\zeta \right) $, do not consist of consecutive variables in the
iterated integrals of (\ref{fmucomp}). This is rectified by using the fact
that $d\left( w,z\right) =\sqrt{\bigtriangleup \left( w,z\right) }$ is a
quasimetric, which in turn follows from the identity%
\begin{equation*}
\sqrt{\bigtriangleup \left( w,z\right) }=\left\vert 1-w\overline{z}%
\right\vert \left\vert \varphi _{z}\left( w\right) \right\vert =\delta
\left( w,z\right) ^{2}\rho \left( w,z\right) ,
\end{equation*}%
where $\rho \left( w,z\right) =\left\vert \varphi _{z}\left( w\right)
\right\vert $ is the invariant pseudohyperbolic metric on the ball
(Corollary 1.22 in \cite{Zhu}) and where $\delta \left( w,z\right)
=\left\vert 1-w\overline{z}\right\vert ^{\frac{1}{2}}$ satisfies the
triangle inequality on the ball (Proposition 5.1.2 in \cite{Rud}). Using the
quasi-subadditivity of $d\left( w,z\right) $ we can, with some care,
redistribute appropriate factors back to the iterated integrals where they
can be favourably estimated using Lemma \ref{Zlemma}. It is simplest to
illustrate this procedure in specific cases, so we defer further discussion
of this point until we treat in detail the cases $\mu =0,1,2$ below. We
again emphasize that all of the above observations regarding \emph{rogue}
factors in (\ref{fmucomp}) apply equally well to the \emph{rogue} factors in
the other terms $\Phi _{n,s_{q+1}}^{\ell }\left( \overline{\mathcal{D}}%
^{m_{q}^{\prime }}\mathcal{F}_{q+1}^{\mu }\right) \left( z\right) $ in (\ref%
{Fterms}), as well as to the boundary terms $\mathcal{S}_{n,s_{q+1}}\left( 
\overline{\mathcal{D}}^{j}\mathcal{F}_{q+1}^{\mu }\right) \left( z\right) $
in (\ref{Fterms}).

The other difficulty remaining is that in order to obtain a favourable
estimate using Lemma \ref{Zlemma} for the iterated integrals resulting from
the bullet items above, it is necessary to generate additional powers of $%
\left( 1-\left\vert z\right\vert ^{2}\right) $ (we are using $z$ as a
generic variable in the iterated integrals here). This is accomplished by
applying the radial integrations by parts in Corollary \ref{IBP2iter} to the 
\emph{previous} iterated integral. Of course such a possibility is
impossible for the first of the iterated integrals, but there we are only
applying the radial derivative $R$ thanks to the fact that our candidate $f$
from the Koszul complex is holomorphic. As a result, we see from (\ref{rootD}%
) that $\left( 1-\left\vert z\right\vert ^{2}\right) R$, unlike $D$,
generates positive powers of $1-\left\vert z\right\vert ^{2}$ even when
acting on $\bigtriangleup \left( w,z\right) $. This procedure is also best
illustrated in specific cases and will be treated in the next subsection.

So ignoring these technical issues for the moment, the integrals that result
from taking absolute values and $\ell ^{2}$-norms inside (\ref{fmucomp}) are
now estimated using Lemma \ref{Zlemma} and Remark \ref{Zlemma'}. Note that
we only use \emph{scalar-valued} Schur estimates since all the integrals to
which Lemma \ref{Zlemma} and Remark \ref{Zlemma'} are applied have positive
integrands. Here is the rough idea. Suppose that $\left\{
T_{1},T_{2},...,T_{\mu }\right\} $ is a collection of Charpentier solution
operators and that for a sequence of large integers $\left\{ m_{1}^{\prime
},m_{1}^{\prime \prime },m_{2}^{\prime },,m_{2}^{\prime \prime }...,m_{\mu
+1}^{\prime },m_{\mu +1}^{\prime \prime }\right\} $, we have the inequalities%
\begin{equation}
\left\Vert T_{j}F\right\Vert _{\Lambda _{p,m_{j}^{\prime },m_{j}^{\prime
\prime }}^{\sigma }\left( \mathbb{B}_{n};\mathbb{\ell }^{2}\right) }\leq
C_{j}\left\Vert F\right\Vert _{\Lambda _{p,m_{j+1}^{\prime },m_{j+1}^{\prime
\prime }}^{\sigma }\left( \mathbb{B}_{n};\mathbb{\ell }^{2}\right) },\ \ \ \
\ 1\leq j\leq \ell +1,  \label{seminormineq}
\end{equation}%
for the class of smooth functions $F$ that arise as $TG$ for some
Charpentier solution operator $T$ and some smooth $G$. Then we can estimate $%
\left\Vert T_{1}\circ T_{2}\circ ...\circ T_{\mu }\Omega \right\Vert
_{B_{p,m}^{\sigma }\left( \mathbb{B}_{n};\mathbb{\ell }^{2}\right) }$ by%
\begin{eqnarray*}
&&\left\Vert T_{1}\circ T_{2}\circ ...\circ T_{\ell }\Omega \right\Vert
_{\Lambda _{p,m_{1}^{\prime },m_{1}^{\prime \prime }}^{\sigma }\left( 
\mathbb{B}_{n};\mathbb{\ell }^{2}\right) } \\
&\leq &C_{1}\left\Vert T_{2}\circ ...\circ T_{\ell }\Omega \right\Vert
_{\Lambda _{p,m_{2}^{\prime },m_{2}^{\prime \prime }}^{\sigma }\left( 
\mathbb{B}_{n};\mathbb{\ell }^{2}\right) } \\
&\leq &C_{1}C_{2}\left\Vert T_{3}\circ ...\circ T_{\ell }\Omega \right\Vert
_{\Lambda _{p,m_{3}^{\prime },m_{3}^{\prime \prime }}^{\sigma }\left( 
\mathbb{B}_{n};\mathbb{\ell }^{2}\right) } \\
&\leq &C_{1}C_{2}...C_{\ell }\left\Vert \Omega \right\Vert _{\Lambda
_{p,m_{\ell +1}^{\prime },m_{\ell +1}^{\prime \prime }}^{\sigma }\left( 
\mathbb{B}_{n};\mathbb{\ell }^{2}\right) }.
\end{eqnarray*}%
Finally we will show that if $\Omega $ is one of the forms $\Omega
_{q}^{q+1} $ in the Koszul complex, then 
\begin{equation*}
\left\Vert \Omega \right\Vert _{\Lambda _{p,m_{\ell +1}^{\prime },m_{\ell
+1}^{\prime \prime }}^{\sigma }\left( \mathbb{B}_{n};\mathbb{\ell }%
^{2}\right) }\leq \left\Vert \Omega \right\Vert _{\Lambda _{p,m_{\ell
+1}^{\prime }+m_{\ell +1}^{\prime \prime }}^{\sigma }\left( \mathbb{B}_{n};%
\mathbb{\ell }^{2}\right) }\leq C_{n,\sigma ,p,\delta }\left( g\right)
\left\Vert h\right\Vert _{B_{p,m}^{\sigma }\left( \mathbb{B}_{n}\right) },
\end{equation*}%
and so altogether this proves that%
\begin{equation*}
\left\Vert f\right\Vert _{B_{p}^{\sigma }\left( \mathbb{B}_{n};\mathbb{\ell }%
^{2}\right) }\leq C_{n,\sigma ,p,\delta }\left( g\right) \left\Vert
h\right\Vert _{B_{p,m}^{\sigma }\left( \mathbb{B}_{n}\right) }.
\end{equation*}

\bigskip

We now make some brief comments on how to obtain the inequalities in (\ref%
{seminormineq}). Complete details will be given in the cases $\mu =0,1,2$
below, and the general case $0\leq \mu \leq n$ is no different than these
three cases. We note that from (\ref{gensolutionker}) the kernel of $%
\mathcal{C}_{n}^{0,q}$ typically looks like a sum of terms%
\begin{equation}
\frac{\left( 1-w\overline{z}\right) ^{n-1-q}\left( 1-\left\vert w\right\vert
^{2}\right) ^{q}}{\bigtriangleup \left( w,z\right) ^{n}}\left( \overline{%
z_{j}}-\overline{w_{j}}\right)  \label{typ}
\end{equation}%
times a wedge product of differentials in which the differential $d\overline{%
w_{j}}$ is missing. We again emphasize that the \emph{rogue} factor $\left( 
\overline{z_{j}}-\overline{w_{j}}\right) $ cannot simply be estimated by $%
\left\vert \overline{z_{j}}-\overline{w_{j}}\right\vert $ as the formula (%
\ref{manyfaces}) shows that%
\begin{equation*}
\sqrt{\bigtriangleup \left( w,z\right) }=\left\vert P_{z}\left( z-w\right) +%
\sqrt{1-\left\vert z\right\vert ^{2}}Q_{z}\left( z-w\right) \right\vert
\end{equation*}%
can be much smaller than $\left\vert z-w\right\vert $. As we mentioned
above, it is possible to exploit the fact that any surviving term in the
form $\Omega _{\mu }^{\mu +1}$ must then involve the derivative $\frac{%
\partial }{\partial \overline{w_{j}}}$ hitting a component of $g$. This
permits us to absorb part of the complex tangential component of $z-w$ into
the almost invariant derivative $D$ which is larger than the usual gradient
in the complex tangential directions. This results in a good estimate for
the \emph{rogue} factor $\left( \overline{z_{j}}-\overline{w_{j}}\right) $
in (\ref{typ}) based on the smaller quantity $\sqrt{\bigtriangleup \left(
w,z\right) }$. We have already integrated by parts to write (\ref{typ}) as
(recall that the factors $\overline{z_{j}}-\overline{w_{j}}$ are already
incorporated into $\overline{\mathcal{D}_{z}^{m}}\eta \left( w\right) $) 
\begin{equation*}
\int_{\mathbb{B}_{n}}\frac{\left( 1-w\overline{z}\right) ^{n-1-q}\left(
1-\left\vert w\right\vert ^{2}\right) ^{q}}{\bigtriangleup \left( w,z\right)
^{n}}\overline{\mathcal{D}^{m}}\eta \left( w\right) dV\left( w\right) ,
\end{equation*}%
plus boundary terms which we ignore for the moment. Then we use the three
crucial inequalities (\ref{moddelta}), (\ref{rootD}) and (\ref{Dbound});%
\begin{eqnarray*}
\left\vert \left( \overline{z_{j}}-\overline{w_{j}}\right) \overline{%
\mathcal{D}_{z,w}^{m}}\Omega _{\ell }^{\ell +1}\left( w\right) \right\vert
&\leq &\left( \frac{\sqrt{\bigtriangleup \left( w,z\right) }}{1-\left\vert
w\right\vert ^{2}}\right) ^{m+1}\left\vert \overline{D^{m}}\widehat{\Omega
_{\ell }^{\ell +1}}\left( w\right) \right\vert , \\
\left\vert D_{\left( z\right) }\bigtriangleup \left( w,z\right) \right\vert
&\leq &C\left( 1-\left\vert z\right\vert ^{2}\right) \bigtriangleup \left(
w,z\right) ^{\frac{1}{2}}+\bigtriangleup \left( w,z\right) , \\
\left\vert \left( 1-\left\vert z\right\vert ^{2}\right) R_{\left( z\right)
}\bigtriangleup \left( w,z\right) \right\vert &\leq &C\left( 1-\left\vert
z\right\vert ^{2}\right) \bigtriangleup \left( w,z\right) ^{\frac{1}{2}}, \\
\left\vert D_{\left( z\right) }^{m}\left\{ \left( 1-\overline{w}z\right)
^{k}\right\} \right\vert &\leq &C\left\vert 1-\overline{w}z\right\vert
^{k}\left( \frac{1-\left\vert z\right\vert ^{2}}{\left\vert 1-\overline{w}%
z\right\vert }\right) ^{\frac{m}{2}} \\
\left\vert \left( 1-\left\vert z\right\vert ^{2}\right) ^{m}R_{\left(
z\right) }^{m}\left\{ \left( 1-\overline{w}z\right) ^{k}\right\} \right\vert
&\leq &C\left\vert 1-\overline{w}z\right\vert ^{k}\left( \frac{1-\left\vert
z\right\vert ^{2}}{\left\vert 1-\overline{w}z\right\vert }\right) ^{m},
\end{eqnarray*}%
to help show that the resulting iterated kernels can be factored (after
accounting for all \emph{rogue} factors $\overline{z_{j}}-\overline{w_{j}}$)
into operators that satisfy the hypotheses of Lemma \ref{Zlemma} or Remark %
\ref{Zlemma'} above.

\begin{definition}
\label{defhat}The expression $\widehat{\Omega _{\ell }^{\ell +1}}$ denotes
the form $\Omega _{\ell }^{\ell +1}$ but with every occurrence of the
derivative $\frac{\partial }{\partial \overline{w_{j}}}$ replaced by the
derivative $\overline{D_{j}}$.
\end{definition}

Recall that each summand of $\Omega _{\ell }^{\ell +1}$ includes a product
of exactly $\ell $ distinct derivatives $\frac{\partial }{\partial \overline{%
w_{j}}}$ applied to components of $g$. Thus the entries of $\overline{D^{m}}%
\widehat{\Omega _{\ell }^{\ell +1}}\left( w\right) $ consist of $m+\ell $
derivatives distributed among components of $g$. Using the factorization of $%
\Omega _{\ell }^{\ell +1}$ in (\ref{Omegaform}), we obtain the corresponding
factorization for $\widehat{\Omega _{\ell }^{\ell +1}}$:%
\begin{equation}
\Omega _{0}^{1}\wedge \dbigwedge\limits_{i=1}^{\ell }\widehat{\Omega _{0}^{1}%
}=-\frac{1}{\ell +1}\widehat{\Omega _{\ell }^{\ell +1}},  \label{genformfact}
\end{equation}%
where $\Omega _{0}^{1}=\left( \frac{\overline{g_{i}}}{\left\vert
g\right\vert ^{2}}\right) _{i=1}^{\infty }\ $and $\widehat{\Omega _{0}^{1}}%
=\left( \frac{\overline{Dg_{i}}}{\left\vert g\right\vert ^{2}}\right)
_{i=1}^{\infty }.$

It is important for this purpose of using Lemma \ref{Zlemma} and Remark \ref%
{Zlemma'} to first apply the integration by parts Lemma \ref{IBP1'}\ to
temper the singularity due to negative powers of $\bigtriangleup \left(
w,z\right) $, and to use the integration by parts Corollary \ref{IBP2iter}
to infuse enough powers of $\left( 1-\left\vert w\right\vert ^{2}\right) $
for use in the subsequent iterated integral.

Finally it follows from Lemma \ref{radinv}, Proposition \ref{Bequivrad} and
Proposition \ref{multilinear} together with the factorization (\ref%
{Omegaform}) that%
\begin{equation}
\left\Vert \left( 1-\left\vert z\right\vert ^{2}\right) ^{\sigma }\mathcal{X}%
^{m}\widehat{\Omega _{\mu }^{\mu +1}}h\left( z\right) \right\Vert
_{L^{p}\left( \lambda _{n};\mathbb{\ell }^{2}\right) }\leq C\left\Vert 
\mathbb{M}_{g}\right\Vert _{B_{p}^{\sigma }\left( \mathbb{B}_{n}\right)
\rightarrow B_{p}^{\sigma }\left( \mathbb{B}_{n};\mathbb{\ell }^{2}\right)
}^{m+\mu }\left\Vert h\right\Vert _{B_{p}^{\sigma }\left( \mathbb{B}%
_{n}\right) }.  \label{Omegabound}
\end{equation}%
We defer the proof of (\ref{Omegabound}) until Subsubsection \ref{The
estimate for F0} when further calculations are available.

\begin{remark}
At this point we observe from (\ref{accomplishFmu}) that the exponent $m+\mu 
$ in (\ref{Omegabound}) is at most $m_{n}+n$, and thus we may take $\kappa
=m_{n}+n$. We leave it to the interested reader to estimate the size of $%
m_{n}$.
\end{remark}

Taking into account all of the above, the conclusion is that with $\kappa
=m_{n}+n$, 
\begin{equation*}
\left\Vert f\right\Vert _{B_{p}^{\sigma }\left( \mathbb{B}_{n};\mathbb{\ell }%
^{2}\right) }\leq C_{n,\sigma ,p,\delta }\left\Vert \mathbb{M}%
_{g}\right\Vert _{B_{p}^{\sigma }\left( \mathbb{B}_{n}\right) \rightarrow
B_{p}^{\sigma }\left( \mathbb{B}_{n};\mathbb{\ell }^{2}\right) }^{\kappa
}\left\Vert h\right\Vert _{B_{p}^{\sigma }\left( \mathbb{B}_{n}\right) }.
\end{equation*}

As the arguments described above are rather complicated we illustrate them
by considering the three cases $\mu =0,1,2$ in complete detail in the next
subsection before proceeding to the general case.

\subsection{Estimates in special cases\label{specialcases}}

Here we prove the estimates (\ref{accomplishFmu}) for $\mu =0,1,2$. Recall
that%
\begin{eqnarray*}
\mathcal{F}^{0} &=&\Omega _{0}^{1}h, \\
\mathcal{F}^{1} &=&\Lambda _{g}\mathcal{C}_{n,s_{1}}^{0,0}\Omega _{1}^{2}h,
\\
\mathcal{F}^{2} &=&\Lambda _{g}\mathcal{C}_{n,s_{1}}^{0,0}\Lambda _{g}%
\mathcal{C}_{n,s_{2}}^{0,1}\Omega _{2}^{3}h.
\end{eqnarray*}%
To obtain the estimate for $\mathcal{F}^{0}$ we use the multilinear
inequality in Proposition \ref{multilinear}.

In estimating $\mathcal{F}^{1}$ we confront for the first time a \emph{rogue}
factor $\overline{z_{k}-w_{k}}$ that we must associate with a derivative $%
\frac{\partial }{\partial \overline{w_{k}}}$ occurring in each surviving
summand of the $k^{th}$ component of the form $\Omega _{1}^{2}$. After
applying the integration by parts formula in \ref{IBPamel} as in \cite{OrFa}%
, we use the crucial inequalities in Proposition \ref{threecrucial} and the
Schur type operator estimates in Lemma \ref{Zlemma} with $c=0$ to obtain the
desired estimates. Finally we must also deal with the boundary terms in the
integration by parts formula for ameliorated Charpentier kernels in Lemma %
\ref{IBPamel}. This requires using the radial derivative integration by
parts formula in Corollary \ref{IBP2iter} as in \cite{OrFa}, and also
requires dealing with the corresponding \emph{rogue} factors.

The final trick in the proof arises in estimating $\mathcal{F}^{2}$. This
time there are two iterated integrals each with a \emph{rogue} factor. The
problematic \emph{rogue} factor $\overline{z_{k}-\zeta _{k}}$ occurs in the 
\emph{first} of the iterated integrals since there is \emph{no} derivative $%
\frac{\partial }{\partial \overline{\zeta _{k}}}$ hitting the second
iterated integral with which to associate the \emph{rogue} factor $\overline{%
z_{k}-\zeta _{k}}$. Instead we decompose the factor as $\overline{z_{k}-w_{k}%
}-\overline{\zeta _{k}-w_{k}}$ and associate each of these summands with a
derivative $\frac{\partial }{\partial \overline{w_{k}}}$ already occurring
in $\Omega _{2}^{3}$. Then we can apply the crucial inequality (\ref%
{moddelta}) and use the fact that $\sqrt{\bigtriangleup \left( w,z\right) }$
is a quasimetric to redistribute the estimates appropriately. As a result of
this redistribution we are forced to use Lemma \ref{Zlemma} with $c=\pm 1$
this time as well as $c=0$. In applying the Schur type estimates in Lemma %
\ref{Zlemma} to the \emph{second} iterated integral, we require a
sufficiently large power of $\left( 1-\left\vert w\right\vert ^{2}\right) $
to be carried over from the first iterated integral. To ensure this we again
use the radial derivative integration by parts formula in Corollary \ref%
{IBP2iter}.

The estimate (\ref{accomplishFmu}) for general $\mu $ involves no new ideas.
There are now $\mu $ rogue terms and we need to apply Lemma \ref{Zlemma}
with $c=0,\pm 1,...,\pm \left( \mu -1\right) $. With this noted the
arguments needed are those used above in the cases $\mu =0,1,2$.

\subsubsection{The estimate for $\mathcal{F}^{0}$\label{The estimate for F0}}

We begin with the estimate%
\begin{eqnarray*}
\left\Vert \mathcal{F}^{0}\right\Vert _{B_{p,m}^{\sigma }\left( \mathbb{B}%
_{n};\ell ^{2}\right) } &=&\left\Vert \Omega _{0}^{1}h\right\Vert
_{B_{p,m}^{\sigma }\left( \mathbb{B}_{n};\ell ^{2}\right) } \\
&\leq &C_{n,\sigma ,p,\delta }\left\Vert \mathbb{M}_{g}\right\Vert
_{B_{p}^{\sigma }\left( \mathbb{B}_{n}\right) \rightarrow B_{p}^{\sigma
}\left( \mathbb{B}_{n};\ell ^{2}\right) }^{m}\left\Vert h\right\Vert
_{B_{p,m}^{\sigma }\left( \mathbb{B}_{n}\right) },
\end{eqnarray*}%
for $m+\sigma >\frac{n}{p}$. However, for later use we prove instead the
more general estimate with $\mathcal{X}$ in place of $R$, except that $m$
must then be chosen twice as large:%
\begin{eqnarray}
&&\int_{\mathbb{B}_{n}}\left\vert \left( 1-\left\vert z\right\vert
^{2}\right) ^{\sigma }\mathcal{X}^{m}\left( \Omega _{0}^{1}h\right) \left(
z\right) \right\vert ^{p}d\lambda _{n}\left( z\right)  \label{twice} \\
&&\ \ \ \ \ \leq C_{n,\sigma ,p,\delta }\left\Vert \mathbb{M}_{g}\right\Vert
_{B_{p}^{\sigma }\left( \mathbb{B}_{n}\right) \rightarrow B_{p}^{\sigma
}\left( \mathbb{B}_{n};\ell ^{2}\right) }^{mp}\left\Vert h\right\Vert
_{B_{p}^{\sigma }\left( \mathbb{B}_{n}\right) }^{p},  \notag
\end{eqnarray}%
for $m>2\left( \frac{n}{p}-\sigma \right) $. Recall that $\mathcal{X}^{m}$
is the differential operator of order $m$ given in Definition \ref{calX}
that is adapted to the complex geometry of the unit ball $\mathbb{B}_{n}$.
It will be in estimating iterated Charpentier integrals below that the
derivatives $R^{m}$ and $\overline{\mathcal{D}^{m}}$ will arise from
integration by parts in the previous iterated integral, and this will
require estimates using $\mathcal{X}^{m}$.

By Leibniz' rule for $\mathcal{X}^{m}$ we have%
\begin{equation*}
\mathcal{X}^{m}\left( \Omega _{0}^{1}h\right) =\sum_{k=0}^{m}c_{k}\left( 
\mathcal{X}^{k}\Omega _{0}^{1}\right) \left( \mathcal{X}^{m-k}h\right) ,
\end{equation*}%
and%
\begin{equation}
\mathcal{X}^{k}\left( \Omega _{0}^{1}\right) =\mathcal{X}^{k}\left( \frac{%
\overline{g}}{\left\vert g\right\vert ^{2}}\right) =\sum_{\ell
=0}^{k}c_{\ell }\left( \mathcal{X}^{k-\ell }\overline{g}\right) \left( 
\mathcal{X}^{\ell }\left\vert g\right\vert ^{-2}\right) .  \label{diffOmega}
\end{equation}%
It suffices to prove%
\begin{eqnarray*}
&&\int_{\mathbb{B}_{n}}\left\vert \left( 1-\left\vert z\right\vert
^{2}\right) ^{\sigma }\left( \sum_{k=0}^{m}\sum_{\ell =0}^{k}c_{k}c_{\ell
}\left( \mathcal{X}^{k-\ell }\overline{g}\right) \left( \mathcal{X}^{\ell
}\left\vert g\right\vert ^{-2}\right) \left( \mathcal{X}^{m-k}h\right)
\right) \right\vert ^{p}d\lambda _{n} \\
&&\ \ \ \ \ \ \ \ \ \ \leq C_{n,\sigma ,p,\delta }\left\Vert \mathbb{M}%
_{g}\right\Vert _{B_{p}^{\sigma }\left( \mathbb{B}_{n}\right) \rightarrow
B_{p}^{\sigma }\left( \mathbb{B}_{n};\ell ^{2}\right) }^{mp}\left\Vert
h\right\Vert _{B_{p}^{\sigma }\left( \mathbb{B}_{n}\right) }^{p},
\end{eqnarray*}%
and hence%
\begin{eqnarray}
&&\int_{\mathbb{B}_{n}}\left( 1-\left\vert z\right\vert ^{2}\right)
^{p\sigma }\left\vert \mathcal{X}^{k-\ell }\overline{g}\right\vert
^{p}\left\vert \mathcal{X}^{\ell }\left\vert g\right\vert ^{-2}\right\vert
^{p}\left\vert \mathcal{X}^{m-k}h\right\vert ^{p}d\lambda _{n}
\label{suffprove} \\
&&\ \ \ \ \ \ \ \ \ \ \leq C_{n,\sigma ,p,\delta }\left\Vert \mathbb{M}%
_{g}\right\Vert _{B_{p}^{\sigma }\left( \mathbb{B}_{n}\right) \rightarrow
B_{p}^{\sigma }\left( \mathbb{B}_{n};\ell ^{2}\right) }^{mp}\left\Vert
h\right\Vert _{B_{p}^{\sigma }\left( \mathbb{B}_{n}\right) }^{p},  \notag
\end{eqnarray}%
for each fixed $0\leq \ell \leq k\leq m$.

Now we can profitably estimate both $\left\vert \mathcal{X}%
^{m-k}h\right\vert $ and $\left\vert \mathcal{X}^{k-\ell }\overline{g}%
\right\vert $ as they are, but we must be more careful with $\left\vert 
\mathcal{X}^{\ell }\left\vert g\right\vert ^{-2}\right\vert $. In the case $%
\ell =1$, we assume for convenience that $\mathcal{X}$ annihilates $g_{i}$\
(if not it will annihilate $\overline{g_{i}}$ unless $\mathcal{X}=I$) and
obtain, 
\begin{eqnarray*}
\left\vert \mathcal{X}\left\vert g\right\vert ^{-2}\right\vert ^{2}
&=&\left\vert -\left\vert g\right\vert ^{-4}\sum_{i=1}^{\infty }g_{i}%
\mathcal{X}\overline{g_{i}}\right\vert ^{2} \\
&\leq &\left\vert g\right\vert ^{-8}\left( \sum_{i=1}^{\infty }\left\vert
g_{i}\right\vert ^{2}\right) \left( \sum_{i=1}^{\infty }\left\vert \mathcal{X%
}\overline{g_{i}}\right\vert ^{2}\right) \leq \left\vert g\right\vert
^{-6}\sum_{i=1}^{\infty }\left\vert \mathcal{X}\overline{g_{i}}\right\vert
^{2}.
\end{eqnarray*}%
Similarly when $\ell =2$,%
\begin{eqnarray*}
\left\vert \mathcal{X}^{2}\left\vert g\right\vert ^{-2}\right\vert ^{2}
&=&\left\vert -\left\vert g\right\vert ^{-4}\sum_{i=1}^{\infty }g_{i}%
\mathcal{X}^{2}\overline{g_{i}}+2\left\vert g\right\vert ^{-6}\sum_{i\neq
j}\left( g_{i}\mathcal{X}\overline{g_{i}}\right) \left( g_{j}\mathcal{X}%
\overline{g_{j}}\right) \right\vert ^{2} \\
&\leq &2\left\vert g\right\vert ^{-6}\sum_{i=1}^{\infty }\left\vert \mathcal{%
X}^{2}\overline{g_{i}}\right\vert ^{2}+4\left\vert g\right\vert ^{-8}\left(
\sum_{i=1}^{\infty }\left\vert \mathcal{X}\overline{g_{i}}\right\vert
^{2}\right) ^{2},
\end{eqnarray*}%
and the general case is%
\begin{eqnarray}
&&\left\vert \mathcal{X}^{\ell }\left\vert g\right\vert ^{-2}\right\vert ^{2}
\label{diffneg} \\
&\leq &C_{\ell }\left\vert g\right\vert ^{-6}\sum_{i=1}^{\infty }\left\vert 
\mathcal{X}^{\ell }g_{i}\right\vert ^{2}+C_{\ell -1}\left\vert g\right\vert
^{-8}\left( \sum_{i=1}^{\infty }\left\vert \mathcal{X}^{\ell -1}\overline{%
g_{i}}\right\vert ^{2}\right) \left( \sum_{i=1}^{\infty }\left\vert \mathcal{%
X}\overline{g_{i}}\right\vert ^{2}\right)  \notag \\
&&+...+C_{0}\left\vert g\right\vert ^{-4-2\ell }\left( \sum_{i=1}^{\infty
}\left\vert \mathcal{X}\overline{g_{i}}\right\vert ^{2}\right) ^{\ell } 
\notag \\
&=&\sum_{1\leq \alpha _{1}\leq \alpha _{2}\leq ...\leq \alpha _{M}:\alpha
_{1}+\alpha _{2}+...+\alpha _{M}=\ell }c_{\alpha }\left\vert g\right\vert
^{-4-2\ell }\prod_{m=1}^{M}\left( \sum_{i=1}^{\infty }\left\vert \mathcal{X}%
^{\alpha _{m}}\overline{g_{i}}\right\vert ^{2}\right) .  \notag
\end{eqnarray}

We can ignore the powers of $\left\vert g\right\vert $ since $\left\vert
g\right\vert $ is bounded above and below by Lemma \ref{boundedmult} and the
hypotheses of Theorem \ref{baby}. Fixing $\alpha $ we see that the left side
of (\ref{suffprove}) is thus at most%
\begin{equation*}
C_{n,\sigma ,p,\delta }\int_{\mathbb{B}_{n}}\left( 1-\left\vert z\right\vert
^{2}\right) ^{p\sigma }\left\vert \mathcal{X}^{k-\ell }\overline{g}%
\right\vert ^{p}\left\vert \mathcal{Y}^{m-k}h\right\vert ^{p}\left(
\prod_{j=1}^{M}\left\vert \mathcal{X}^{\alpha _{j}}\overline{g}\right\vert
^{p}\right) d\lambda _{n}.
\end{equation*}%
Since $\left\vert \mathcal{X}^{k-\ell }\overline{g}\right\vert
^{2}=\sum_{i=1}^{\infty }\left\vert \mathcal{X}^{k-\ell }\overline{g_{i}}%
\right\vert ^{2}$ and $k-\ell $ could vanish (unlike the exponents $\alpha
_{\ell }$ which are positive), we see that altogether after renumbering, it
suffices to prove%
\begin{eqnarray}
&&\int_{\mathbb{B}_{n}}\left( 1-\left\vert z\right\vert ^{2}\right)
^{p\sigma }\left\vert \mathcal{Y}^{\alpha _{1}}h\right\vert ^{p}\left\vert 
\mathcal{Y}^{\alpha _{2}}g\right\vert ^{p}...\left\vert \mathcal{Y}^{\alpha
_{M}}g\right\vert ^{p}d\lambda _{n}  \label{suffprove'} \\
&&\ \ \ \ \ \ \ \ \ \ \leq C_{n,\sigma ,p,\delta }\left\Vert \mathbb{M}%
_{g}\right\Vert _{B_{p}^{\sigma }\left( \mathbb{B}_{n}\right) \rightarrow
B_{p}^{\sigma }\left( \mathbb{B}_{n};\ell ^{2}\right) }^{Mp}\left\Vert
h\right\Vert _{B_{p}^{\sigma }\left( \mathbb{B}_{n}\right) }^{p}  \notag
\end{eqnarray}%
for each fixed $\alpha =\left( \alpha _{1},\alpha _{2},...,\alpha
_{M}\right) $ with $M\geq 2$, $\left\vert \alpha \right\vert =m$ and at most
one of $\alpha _{2},...,\alpha _{M}$ is zero. We have used here that $%
\left\vert \overline{D}\overline{g}\right\vert =\left\vert Dg\right\vert $.
Now Proposition \ref{multilinear} yields (\ref{suffprove'}) for each $0\leq
k\leq m$ and $\left\vert \alpha \right\vert =m-k$. Summing these estimates
completes the proof of (\ref{twice}).

We can now prove the more general inequality (\ref{Omegabound}). Indeed,
using the factorization (\ref{Omegaform}) of $\widehat{\Omega _{\mu }^{\mu
+1}}$ together with the Leibniz formula gives%
\begin{eqnarray*}
\mathcal{X}^{m}\left( \widehat{\Omega _{\mu }^{\mu +1}}h\right) &=&\mathcal{X%
}^{m}\left( \Omega _{0}^{1}\wedge \left( \widehat{\Omega _{0}^{1}}\right)
^{\mu }h\right) \\
&=&\sum_{\alpha \in \mathbb{Z}_{+}^{\mu +2}:\left\vert \alpha \right\vert
=m}\left( \mathcal{X}^{\alpha _{0}}\Omega _{0}^{1}\right) \wedge
\dbigwedge\limits_{j=1}^{\mu }\left( \mathcal{X}^{\alpha _{j}}\widehat{%
\Omega _{0}^{1}}\right) \left( \mathcal{X}^{\alpha _{\mu +1}}h\right) \\
&=&\sum_{\alpha \in \mathbb{Z}_{+}^{\mu +2}:\left\vert \alpha \right\vert
=m}\left\{ \left( \mathcal{X}^{\alpha _{0}}\Omega _{0}^{1}\right) \wedge
\dbigwedge\limits_{j=1}^{\mu }\left( \mathcal{X}^{\alpha _{j}+1}\Omega
_{0}^{1}\right) \right\} \left( \mathcal{X}^{\alpha _{\mu +1}}h\right) ,
\end{eqnarray*}%
where we have used that $\widehat{\Omega _{0}^{1}}$ already has an $\mathcal{%
X}$ derivative in each summand, and so $\mathcal{X}^{\alpha _{j}}\widehat{%
\Omega _{0}^{1}}$ can be written as $\mathcal{X}^{\alpha _{j}+1}\Omega
_{0}^{1}$. Now use (\ref{diffOmega}) and (\ref{diffneg}) to see that $%
\left\vert \mathcal{X}^{m}\left( \widehat{\Omega _{\mu }^{\mu +1}}h\right)
\right\vert $ is controlled by a tensor product of at most $m+\mu $ factors,
and then apply Proposition \ref{multilinear} as above to complete the proof
of (\ref{Omegabound}).

\subsubsection{The estimate for $\mathcal{F}^{1}$\label{The estimate for F1}}

The estimate in (\ref{accomplishFmu}) with $\mu =1$ will follow from (\ref%
{Omegabound}) and the estimate 
\begin{eqnarray}
&&\left\Vert \left( 1-\left\vert z\right\vert ^{2}\right) ^{\sigma }\mathcal{%
Y}^{m_{1}}\left( \Lambda _{g}\mathcal{C}_{n,s}^{0,0}\Omega _{1}^{2}h\right)
\right\Vert _{L^{p}\left( \lambda _{n}\right) }^{p}  \label{want'''} \\
&&\ \ \ \ \ \ \ \ \ \ \leq C\int_{\mathbb{B}_{n}}\left\vert \left(
1-\left\vert z\right\vert ^{2}\right) ^{\sigma }\mathcal{X}^{m_{2}}\left( 
\widehat{\Omega _{1}^{2}}h\right) \left( z\right) \right\vert ^{p}d\lambda
_{n}\left( z\right) ,  \notag
\end{eqnarray}%
where as in Definition \ref{defhat}, we define $\widehat{\Omega _{1}^{2}}$
to be $\Omega _{1}^{2}$ with $\partial $ replaced by $D$ throughout: 
\begin{equation*}
\widehat{\Omega _{1}^{2}}=\sum_{j,k=1}^{N}\frac{\overline{\left\{
g_{k}Dg_{j}-g_{j}Dg_{k}\right\} }}{\left\vert g\right\vert ^{4}}e_{j}\wedge
e_{k},
\end{equation*}%
and where $Dh=\sum_{k=1}^{n}\left( D_{k}h\right) dz_{k}$ and $D_{k}$ is the $%
k^{th}$ component of $D$. We are using here the following observation
regarding the interior product $\Omega _{1}^{2}h\lrcorner d\overline{w_{k}}$:%
\begin{eqnarray}
&&\text{For each summand of }\Omega _{1}^{2}h\lrcorner d\overline{w_{k}}%
\text{, there is a unique }1\leq i\leq N\text{ so that }  \label{gi} \\
&&\ \ \ \ \ \ \ \ \ \ \frac{\partial g_{i}}{\partial \overline{w_{k}}}\text{
occurs as a factor in the summand.}  \notag
\end{eqnarray}

We rewrite (\ref{want'''}) as%
\begin{eqnarray}
&&\left\Vert \left( 1-\left\vert z\right\vert ^{2}\right) ^{\sigma }\mathcal{%
R}^{m_{1}^{\prime \prime }}D^{m_{1}^{\prime }}\left( \Lambda _{g}\mathcal{C}%
_{n,s}^{0,0}\Omega _{1}^{2}h\right) \right\Vert _{L^{p}\left( \lambda
_{n}\right) }^{p}  \label{want'} \\
&&\ \ \ \ \ \ \ \ \ \ \leq C\int_{\mathbb{B}_{n}}\left\vert \left(
1-\left\vert z\right\vert ^{2}\right) ^{\sigma }\mathcal{R}^{m_{2}^{\prime
\prime }}\overline{D^{m_{2}^{\prime }}}\left( \widehat{\Omega _{1}^{2}}%
h\right) \left( z\right) \right\vert ^{p}d\lambda _{n}\left( z\right) , 
\notag
\end{eqnarray}%
where $\mathcal{R}^{m}=\left( 1-\left\vert z\right\vert ^{2}\right)
^{m}\left( R^{k}\right) _{k=0}^{m}$ as in (\ref{shorthand}). As mentioned
above, we only need to prove the case $m_{1}^{\prime \prime }=0$ since (\ref%
{accomplishFmu}) only requires that we estimate $\left\Vert \mathcal{F}%
^{1}\right\Vert _{B_{p,m}^{\sigma }\left( \mathbb{B}_{n}\right) }$. However,
when considering the estimate for $\mathcal{F}^{2}$ in (\ref{accomplishFmu})
we will no longer have the luxury of using the norm $\left\Vert \cdot
\right\Vert _{B_{p,m}^{\sigma }\left( \mathbb{B}_{n}\right) }$ in the second
iterated integral occuring there, and so we will consider the more general
case now in preparation for what comes later. As we will see however, it is
necessary to choose $m_{1}^{\prime }$ sufficiently large in order to obtain (%
\ref{want'}). It is useful to recall that the operator $\left( 1-\left\vert
z\right\vert ^{2}\right) R$ is "smaller" than $\overline{D}$ in the sense
that%
\begin{eqnarray*}
\overline{D} &=&\left( 1-\left\vert z\right\vert ^{2}\right) P_{z}\overline{%
\nabla }+\sqrt{1-\left\vert z\right\vert ^{2}}Q_{z}\overline{\nabla }, \\
\left( 1-\left\vert z\right\vert ^{2}\right) R &=&\left( 1-\left\vert
z\right\vert ^{2}\right) P_{z}\nabla .
\end{eqnarray*}

To prove (\ref{want'}) we will ignore the contraction $\Lambda _{g}$ since
if derivatives hit $g$ in the contraction, the estimates are similar if not
easier. Note also that $\left\vert \Lambda _{g}F\right\vert \leq \left\vert
g\right\vert \left\vert F\right\vert $ for the contraction $\Lambda _{g}F$
of any tensor $F$.

We will also initially suppose that $m_{1}^{\prime \prime }=0$ and later
take $m_{1}^{\prime \prime }$ sufficiently large. Now we apply Lemma \ref%
{IBPamel} to $\mathcal{C}_{n,s}^{0,0}\Omega _{1}^{2}h$ and obtain%
\begin{eqnarray}
\mathcal{C}_{n,s}^{0,0}\Omega _{1}^{2}h\left( z\right)  &=&c_{0}\mathcal{C}%
_{n,s}^{0,0}\left( \overline{\mathcal{D}}^{m_{2}^{\prime }}\Omega
_{1}^{2}h\right) \left( z\right) +boundary\ terms  \label{GammaIBP} \\
&=&\int_{\mathbb{B}_{n}}\Phi _{n,s}^{0}\left( w,z\right) \overline{\mathcal{D%
}}^{m_{2}^{\prime }}\left( \Omega _{1}^{2}h\right) dV\left( w\right)   \notag
\\
&&+\ boundary\ terms.  \notag
\end{eqnarray}%
A typical term above looks like%
\begin{equation}
\int_{\mathbb{B}_{n}}\left( \frac{1-\left\vert w\right\vert ^{2}}{1-%
\overline{w}z}\right) ^{s-n}\frac{\left( 1-w\overline{z}\right) ^{n-1}}{%
\bigtriangleup \left( w,z\right) ^{n}}\overline{\mathcal{D}}^{m_{2}^{\prime
}}\left( \Omega _{1}^{2}h\right) dV\left( w\right)   \label{integralabove}
\end{equation}%
where we are discarding the sum of (balanced) factors $\left( \frac{\left(
1-|w|^{2}\right) \left( 1-|z|^{2}\right) }{\left\vert 1-w\overline{z}%
\right\vert ^{2}}\right) ^{j}$ for $1\leq j\leq n-1$ in Lemma \ref{IBPamel}\
that turn out to only help with the estimates. This can be seen from (\ref%
{Dbound}) and its trivial counterpart%
\begin{equation*}
\left\vert D_{\left( z\right) }^{m}\left\{ \left( 1-\left\vert z\right\vert
^{2}\right) ^{k}\right\} \right\vert +\left\vert \left( 1-\left\vert
z\right\vert ^{2}\right) ^{m}R_{\left( z\right) }^{m}\left\{ \left(
1-\left\vert z\right\vert ^{2}\right) ^{k}\right\} \right\vert \leq C\left(
1-\left\vert z\right\vert ^{2}\right) ^{k}.
\end{equation*}

Recall from the general discussion above that in the integral (\ref%
{integralabove}) there are \emph{rogue} factors $\overline{z_{k}-w_{k}}$ in $%
\overline{\mathcal{D}}^{m_{2}^{\prime }}\left( \Omega _{1}^{2}h\right)
\left( w\right) $ that must be associated with a $\frac{\partial }{\partial 
\overline{w_{k}}}$ derivative that hits some factor of each summand in the $%
k^{th}$ component $\Omega _{1}^{2}\lrcorner d\overline{w_{k}}$ of $\Omega
_{1}^{2}\approx \overline{\left\{ g_{i}\partial g_{j}-g_{j}\partial
g_{i}\right\} }$. Thus we can apply (\ref{moddelta}) to the components of $%
\Omega _{1}^{2}h\left( z\right) $ to obtain%
\begin{eqnarray}
&&\left\vert \overline{\mathcal{D}}^{m_{2}^{\prime }}\Omega _{1}^{2}h\left(
z\right) \right\vert  \label{moddest} \\
&\approx &\left\vert \sum_{k=1}^{n}\sum_{\left\vert \alpha \right\vert
=m_{2}^{\prime }}^{n}\left( \overline{w_{k}}-\overline{z_{k}}\right) \left( 
\overline{w-z}\right) ^{\alpha }\frac{\partial ^{m_{2}^{\prime }}}{\partial 
\overline{w}^{\alpha }}\left( \Omega _{1}^{2}h\lrcorner d\overline{w_{k}}%
\right) \right\vert  \notag \\
&\leq &C\left( \frac{\sqrt{\bigtriangleup \left( w,z\right) }}{1-\left\vert
w\right\vert ^{2}}\right) ^{m_{2}^{\prime }+1}\left\vert \overline{%
D^{m_{2}^{\prime }}}\left( \widehat{\Omega _{1}^{2}}h\right) \left( w\right)
\right\vert .  \notag
\end{eqnarray}%
Thus we get%
\begin{eqnarray}
&&\left( 1-\left\vert z\right\vert ^{2}\right) ^{\sigma }\left\vert
D^{m_{1}^{\prime }}\mathcal{C}_{n,s}^{0,0}\Omega _{1}^{2}h\left( z\right)
\right\vert  \label{defS} \\
&\leq &\int_{\mathbb{B}_{n}}\left( 1-\left\vert z\right\vert ^{2}\right)
^{\sigma }\left\vert D_{\left( z\right) }^{m_{1}^{\prime }}\left\{ \frac{%
\left( 1-\left\vert w\right\vert ^{2}\right) ^{s-n}\left( 1-w\overline{z}%
\right) ^{n-1}}{\left( 1-\overline{w}z\right) ^{s-n}\bigtriangleup \left(
w,z\right) ^{n}}\right\} \right\vert  \notag \\
&&\times \left( \frac{\sqrt{\bigtriangleup \left( w,z\right) }}{1-\left\vert
w\right\vert ^{2}}\right) ^{m_{2}^{\prime }+1}\left\vert \overline{%
D^{m_{2}^{\prime }}}\left( \widehat{\Omega _{1}^{2}}h\right) \left( w\right)
\right\vert dV\left( w\right)  \notag \\
&\equiv &S_{m_{1}^{\prime },m_{2}^{\prime }}^{s}f\left( z\right) ,  \notag
\end{eqnarray}%
where%
\begin{equation}
f\left( w\right) =\left( 1-\left\vert w\right\vert ^{2}\right) ^{\sigma
}\left\vert \overline{D^{m_{2}^{\prime }}}\left( \widehat{\Omega _{1}^{2}}%
h\right) \left( w\right) \right\vert .  \label{deff}
\end{equation}

Now we iterate the estimate (\ref{rootD}),%
\begin{equation*}
\left\vert D_{\left( z\right) }\bigtriangleup \left( w,z\right) \right\vert
\leq C\left( 1-\left\vert z\right\vert ^{2}\right) \bigtriangleup \left(
w,z\right) ^{\frac{1}{2}}+\bigtriangleup \left( w,z\right) ,
\end{equation*}%
to obtain%
\begin{eqnarray}
&&\left\vert D_{\left( z\right) }^{m_{1}^{\prime }}\left\{ \frac{\left(
1-\left\vert w\right\vert ^{2}\right) ^{s-n}\left( 1-w\overline{z}\right)
^{n-1}}{\left( 1-\overline{w}z\right) ^{s-n}\bigtriangleup \left( w,z\right)
^{n}}\right\} \right\vert  \label{OKterm} \\
&\leq &\frac{\left( 1-\left\vert z\right\vert ^{2}\right) ^{m_{1}^{\prime
}}\left( 1-\left\vert w\right\vert ^{2}\right) ^{s-n}\bigtriangleup \left(
w,z\right) ^{\frac{m_{1}^{\prime }}{2}}}{\left\vert 1-w\overline{z}%
\right\vert ^{s-2n+1}\bigtriangleup \left( w,z\right) ^{n+m_{1}^{\prime }}} 
\notag \\
&&+...+\frac{\left( 1-\left\vert w\right\vert ^{2}\right) ^{s-n}}{\left\vert
1-w\overline{z}\right\vert ^{s-2n+1}\bigtriangleup \left( w,z\right) ^{n}}%
+OK,  \notag
\end{eqnarray}%
where the terms in $OK$ are obtained when some of the derivatives $D$ hit
the factor $\frac{1}{\left( 1-\overline{w}z\right) ^{s-n}}$ or factors $%
D\bigtriangleup \left( w,z\right) $ already in the numerator. Leaving the $%
OK $ terms for later, we combine all the estimates above to get that if we
plug the first term on the right in (\ref{OKterm}) into the left side of (%
\ref{want'}), then the result is dominated by

\begin{eqnarray*}
&&\int_{\mathbb{B}_{n}}\frac{\left( 1-\left\vert z\right\vert ^{2}\right)
^{m_{1}^{\prime }+\sigma }\left( 1-\left\vert w\right\vert ^{2}\right)
^{s-n-m_{2}^{\prime }-1-\sigma }\bigtriangleup \left( w,z\right) ^{\frac{%
m_{1}^{\prime }+m_{2}^{\prime }+1}{2}}}{\left\vert 1-w\overline{z}%
\right\vert ^{s-2n+1}\bigtriangleup \left( w,z\right) ^{n+m_{1}^{\prime }}}%
f\left( w\right) dV\left( w\right) \\
&=&\int_{\mathbb{B}_{n}}\frac{\left( 1-\left\vert z\right\vert ^{2}\right)
^{m_{1}^{\prime }+\sigma }\left( 1-\left\vert w\right\vert ^{2}\right)
^{s-n-1-m_{2}^{\prime }-\sigma }}{\left\vert 1-w\overline{z}\right\vert
^{s-2n+1}}\sqrt{\bigtriangleup \left( w,z\right) }^{m_{2}^{\prime
}-m_{1}^{\prime }-2n+1}f\left( w\right) dV\left( w\right) .
\end{eqnarray*}%
Now for convenience choose $m_{2}^{\prime }=m_{1}^{\prime }+2n-1$ so that
the factor of $\sqrt{\bigtriangleup \left( w,z\right) }$ disappears. We then
get%
\begin{equation}
\left( 1-\left\vert z\right\vert ^{2}\right) ^{\sigma }\left\vert
D^{m_{1}^{\prime }}\mathcal{C}_{n,s}^{0,0}\Omega _{1}^{2}h\left( z\right)
\right\vert \leq \int_{\mathbb{B}_{n}}\frac{\left( 1-\left\vert z\right\vert
^{2}\right) ^{m_{1}^{\prime }+\sigma }\left( 1-\left\vert w\right\vert
^{2}\right) ^{s-3n-m_{1}^{\prime }-\sigma }}{\left\vert 1-w\overline{z}%
\right\vert ^{s-2n+1}}f\left( w\right) dV\left( w\right) .  \label{altL}
\end{equation}

Lemma \ref{Zlemma} shows that the operator%
\begin{equation*}
T_{a,b,0}f\left( z\right) =\int_{\mathbb{B}_{n}}\frac{\left( 1-\left\vert
z\right\vert ^{2}\right) ^{a}\left( 1-\left\vert w\right\vert ^{2}\right)
^{b}}{\left\vert 1-w\overline{z}\right\vert ^{n+1+a+b}}f\left( w\right)
dV\left( w\right)
\end{equation*}%
is bounded on $L^{p}\left( \mathbb{B}_{n};\left( 1-\left\vert w\right\vert
^{2}\right) ^{t}dV\left( w\right) \right) $ if and only if%
\begin{equation*}
-pa<t+1<p\left( b+1\right) .
\end{equation*}%
We apply this lemma with $t=-n-1$, $a=m_{1}^{\prime }+\sigma $ and $%
b=s-3n-m_{1}^{\prime }-\sigma $. Note that the sum of the exponents in the
numerator and denominator of (\ref{altL}) are equal if we write the integral
in terms of invariant measure $d\lambda _{n}\left( w\right) =\left(
1-\left\vert w\right\vert ^{2}\right) ^{-n-1}dV\left( w\right) $. We
conclude that $S_{m_{1}^{\prime },m_{2}^{\prime }}^{s}$ is bounded on $%
L^{p}\left( d\lambda _{n}\right) $ provided $T$ is, and that this latter
happens if and only if 
\begin{equation*}
-p\left( m_{1}^{\prime }+\sigma \right) <-n<p\left( s-3n+1-m_{1}^{\prime
}-\sigma \right) .
\end{equation*}%
This requires $m_{1}^{\prime }+\sigma >\frac{n}{p}$ and $s>3n-1+m_{1}^{%
\prime }+\sigma -\frac{n}{p}$.

\begin{remark}
Suppose instead that we choose $m_{2}^{\prime }$ above to be a positive
integer satisfying $c=m_{2}^{\prime }-m_{1}^{\prime }-2n+1>-2n$. Then we
would be dealing with the operator $T_{a,b,c}$ where $a=m_{1}^{\prime
}+\sigma $ and 
\begin{equation*}
b=s-n-1-m_{2}^{\prime }-\sigma =s-3n-c-m_{1}^{\prime }-\sigma .
\end{equation*}%
By Lemma \ref{Zlemma}, $T_{a,b,c}$ is bounded on $L^{p}\left( d\lambda
_{n}\right) $ if and only if%
\begin{equation*}
-p\left( m_{1}^{\prime }+\sigma \right) <-n<p\left( s-3n+1-c-m_{1}^{\prime
}-\sigma \right) ,
\end{equation*}%
i.e. $m_{1}^{\prime }+\sigma >\frac{n}{p}$ and $s>c+3n-1+m_{1}^{\prime
}+\sigma -\frac{n}{p}$. Thus we can use \emph{any} value of $c>-2n$ provided
we choose $m_{2}^{\prime }\geq m_{1}^{\prime }$ and $s$ large enough.
\end{remark}

Now we turn to the second displayed term on the right side of (\ref{OKterm})
which leads to the operator $T_{a,b,0}$ with $a=\sigma $, $b=s-3n-\sigma $.
This time we will \emph{not} in general have the required boundedness
condition $\sigma >\frac{n}{p}$. It is for this reason that we must return
to (\ref{want'}) and insist that $m_{1}^{\prime \prime }$ be chosen
sufficiently large that $m_{1}^{\prime \prime }+\sigma >\frac{n}{p}$. For
convenience we let $m_{1}^{\prime }=0$ for now. Indeed, it follows from the
second line in the crucial inequality (\ref{rootD}) that the second
displayed term on the right side of (\ref{OKterm}) is%
\begin{equation*}
\frac{\left( 1-\left\vert z\right\vert ^{2}\right) ^{m_{1}^{\prime \prime
}}\left( 1-\left\vert w\right\vert ^{2}\right) ^{s-n}\bigtriangleup \left(
w,z\right) ^{\frac{m_{1}^{\prime \prime }}{2}}}{\left\vert 1-w\overline{z}%
\right\vert ^{s-2n+1}\bigtriangleup \left( w,z\right) ^{n+m_{1}^{\prime
\prime }}}+better\ terms.
\end{equation*}%
Using this expression and choosing $m_{2}^{\prime }=m_{1}^{\prime \prime
}+2n-1$ so that the term $\sqrt{\bigtriangleup \left( w,z\right) }$
disappears from the ensuing integral, we obtain the following analogue of (%
\ref{altL}):%
\begin{eqnarray*}
&&\left( 1-\left\vert z\right\vert ^{2}\right) ^{\sigma }\left( 1-\left\vert
z\right\vert ^{2}\right) ^{m_{1}^{\prime \prime }}\left\vert \mathcal{R}%
^{m_{1}^{\prime \prime }}\mathcal{C}_{n,s}^{0,0}\Omega _{1}^{2}h\left(
z\right) \right\vert \\
&\leq &\int_{\mathbb{B}_{n}}\frac{\left( 1-\left\vert z\right\vert
^{2}\right) ^{m_{1}^{\prime \prime }+\sigma }\left( 1-\left\vert
w\right\vert ^{2}\right) ^{s-3n-m_{1}^{\prime \prime }-\sigma }}{\left\vert
1-w\overline{z}\right\vert ^{s-2n+1}}f\left( w\right) dV\left( w\right) .
\end{eqnarray*}%
The corresponding operator $T_{a,b,0}$ has $a=m_{1}^{\prime \prime }+\sigma $
and $b=s-3n-m_{1}^{\prime \prime }-\sigma $ and is bounded on $L^{p}\left(
\lambda _{n}\right) $ when $-p\left( m_{1}^{\prime \prime }+\sigma \right)
<-n<p\left( s-3n+1-m_{1}^{\prime \prime }-\sigma \right) $. Thus there is no
unnecessary restriction on $\sigma $ if $m_{1}^{\prime \prime }$ and $s$ are
chosen appropriately large. Note that the only difference between this
operator $T_{a,b,0}$ and the previous one is that $m_{1}^{\prime }$ has been
replaced by $m_{1}^{\prime \prime }$.

The above arguments are easily modified to handle the general case of (\ref%
{want'}) provided $m_{1}^{\prime \prime }+\sigma >\frac{n}{p}$ and $s$ is
chosen sufficiently large.

\bigskip

Now we return to consider the $OK$ terms in (\ref{OKterm}). For this we use
the inequality (\ref{Dbound}):%
\begin{equation*}
\left\vert D_{\left( z\right) }^{m}\left\{ \left( 1-\overline{w}z\right)
^{k}\right\} \right\vert \leq C\left\vert 1-\overline{w}z\right\vert
^{k}\left( \frac{1-\left\vert z\right\vert ^{2}}{\left\vert 1-\overline{w}%
z\right\vert }\right) ^{\frac{m}{2}}.
\end{equation*}%
We ignore the derivative $\left( 1-\left\vert z\right\vert ^{2}\right) R$ as
the second line in (\ref{Dbound}) shows that it satisfies a better estimate.
We also write $m_{1}$ and $m_{2}$ in place of $m_{1}^{\prime }$ and $%
m_{2}^{\prime }$ now. As a result, one of the extremal $OK$ terms in (\ref%
{OKterm}) is%
\begin{equation*}
\frac{\left( 1-\left\vert z\right\vert ^{2}\right) ^{\frac{m_{1}}{2}}\left(
1-\left\vert w\right\vert ^{2}\right) ^{s-n}}{\left\vert 1-w\overline{z}%
\right\vert ^{s-2n+1+\frac{m_{1}}{2}}\bigtriangleup \left( w,z\right) ^{n}},
\end{equation*}%
which when combined with the other estimates leads to the integral operator%
\begin{equation*}
\int_{\mathbb{B}_{n}}\frac{\left( 1-\left\vert z\right\vert ^{2}\right) ^{%
\frac{m_{1}}{2}+\sigma }\left( 1-\left\vert w\right\vert ^{2}\right)
^{s-n-1-m_{2}-\sigma }}{\left\vert 1-w\overline{z}\right\vert ^{s-2n+1+\frac{%
m_{1}}{2}}}\sqrt{\bigtriangleup \left( w,z\right) }^{m_{2}-2n-1}f\left(
w\right) dV\left( w\right) .
\end{equation*}%
This is $T_{a,b,c}$ with $a=\frac{m_{1}}{2}+\sigma $, $b=s-n-1-m_{2}-\sigma $
and $c=m_{2}-2n-1$. This is bounded on $L^{p}\left( \lambda _{n}\right) $
provided $m_{2}\geq 2$ and%
\begin{equation*}
-p\left( \frac{m_{1}}{2}+\sigma \right) <-n<p\left( s-n-m_{2}-\sigma \right)
,
\end{equation*}%
i.e. $\frac{m_{1}}{2}+\sigma >\frac{n}{p}$ and $s>n+m_{2}+\sigma -\frac{n}{p}
$. The intermediate $OK$ terms are handled similarly. Note that the crux of
the matter is that all of the positive operators have the form $T_{a,b,c}$,
and moreover, if $s$ and the $m^{\prime }s$ are chosen appropriately large,
then $T_{a,b,c}$ is bounded on $L^{p}\left( \lambda _{n}\right) $.

\subsubsection{Boundary terms for $\mathcal{F}^{1}$\label{Boundary terms for
F1}}

Now we turn to estimating the boundary terms in (\ref{GammaIBP}). A typical
term is%
\begin{equation}
\mathcal{S}_{n,s}\left( \overline{\mathcal{D}}^{k}\left( \Omega
_{1}^{2}h\right) \right) \left[ \overline{\mathcal{Z}}\right] \left(
z\right) =\int_{\mathbb{B}_{n}}\frac{\left( 1-\left\vert w\right\vert
^{2}\right) ^{s-n-1}}{\left( 1-\overline{w}z\right) ^{s}}\overline{\mathcal{D%
}}^{k}\left( \Omega _{1}^{2}h\right) \left[ \overline{\mathcal{Z}}\right]
\left( w\right) dV\left( w\right) ,  \label{typterm}
\end{equation}%
with $0\leq k\leq m-1$ upon appealing to Lemma \ref{IBPamel}.

We now apply the operator $\left( 1-\left\vert z\right\vert ^{2}\right)
^{m_{1}+\sigma }R^{m_{1}}$ to the integral in the right side of (\ref%
{typterm}) and using Proposition \ref{threecrucial} we obtain that the
absolute value of the result is dominated by%
\begin{eqnarray*}
&&\int_{\mathbb{B}_{n}}\frac{\left( 1-\left\vert z\right\vert ^{2}\right)
^{m_{1}+\sigma }\left( 1-\left\vert w\right\vert ^{2}\right) ^{s-n-1}}{%
\left\vert 1-\overline{w}z\right\vert ^{s+m_{1}}}\left( \frac{\sqrt{%
\bigtriangleup \left( w,z\right) }}{1-\left\vert w\right\vert ^{2}}\right)
^{k+1}\left\vert \overline{D}^{k}\left( \widehat{\Omega _{1}^{2}}h\right)
\right\vert dV\left( w\right) \\
&=&\int_{\mathbb{B}_{n}}\frac{\left( 1-\left\vert z\right\vert ^{2}\right)
^{m_{1}+\sigma }\left( 1-\left\vert w\right\vert ^{2}\right)
^{s-n-2-k-\sigma }\sqrt{\bigtriangleup \left( w,z\right) }^{k+1}}{\left\vert
1-\overline{w}z\right\vert ^{s+m_{1}}}\left\vert \left( 1-\left\vert
w\right\vert ^{2}\right) ^{\sigma }\overline{D}^{k}\left( \widehat{\Omega
_{1}^{2}}h\right) \left( w\right) \right\vert dV\left( w\right) .
\end{eqnarray*}%
The operator in question here is $T_{a,b,c}$ with $a=m_{1}+\sigma $, $%
b=s-n-2-k-\sigma $ and $c=k+1$ since%
\begin{equation*}
a+b+c+n+1=s+m_{1}.
\end{equation*}%
Lemma \ref{Zlemma} applies to prove the desired boundedness on $L^{p}\left(
\lambda _{n}\right) $ provided $m_{1}+\sigma >\frac{n}{p}$.

However, if $k$ fails to satisfy $k+1>2\left( \frac{n}{p}-\sigma \right) $,
then the derivative $D^{k+1}\Omega $ cannot be used to control the norm $%
\left\Vert \Omega \right\Vert _{B_{p}^{\sigma }\left( \mathbb{B}_{n}\right)
} $. To compensate for a small $k$, we must then apply Corollary \ref%
{IBP2iter} to the right side of (\ref{typterm}) (which for fixed $z$ is in $%
C\left( \overline{\mathbb{B}_{n}}\right) \cap C^{\infty }\left( \mathbb{B}%
_{n}\right) $) before differentiating and taking absolute values inside the
integral. This then leads to operators of the form%
\begin{eqnarray*}
&&\left( 1-\left\vert z\right\vert ^{2}\right) ^{m_{1}+\sigma
}R^{m_{1}}\left\{ \int_{\mathbb{B}_{n}}\frac{\left( 1-\left\vert
w\right\vert ^{2}\right) ^{s-n-1}}{\left( 1-\overline{w}z\right) ^{s}}\right.
\\
&&\ \ \ \ \ \times \left. \left( 1-\left\vert w\right\vert ^{2}\right)
^{m}R^{m}\left[ \overline{\mathcal{D}}^{k}\left( \Omega _{1}^{2}h\right)
\left( w\right) \right] dV\left( w\right) \right\} ,
\end{eqnarray*}%
which are dominated by 
\begin{eqnarray*}
&&\int_{\mathbb{B}_{n}}\frac{\left( 1-\left\vert z\right\vert ^{2}\right)
^{m_{1}+\sigma }\left( 1-\left\vert w\right\vert ^{2}\right) ^{s-n-1}}{%
\left\vert 1-\overline{w}z\right\vert ^{s+m_{1}}} \\
&&\times \left( \frac{\sqrt{\bigtriangleup \left( w,z\right) }}{1-\left\vert
w\right\vert ^{2}}\right) ^{k+1}\left\vert \mathcal{R}^{m}\overline{D}%
^{k}\left( \widehat{\Omega _{1}^{2}}h\right) \left( w\right) \right\vert
dV\left( w\right) ,
\end{eqnarray*}%
which is%
\begin{eqnarray*}
&&\int_{\mathbb{B}_{n}}\frac{\left( 1-\left\vert z\right\vert ^{2}\right)
^{m_{1}+\sigma }\left( 1-\left\vert w\right\vert ^{2}\right)
^{s-n-2-k-\sigma }\sqrt{\bigtriangleup \left( w,z\right) }^{k+1}}{\left\vert
1-\overline{w}z\right\vert ^{s+m_{1}}} \\
&&\ \ \ \ \ \ \ \ \ \ \ \ \ \ \ \ \ \ \ \ \times \left\vert \left(
1-\left\vert w\right\vert ^{2}\right) ^{\sigma }\mathcal{R}^{m}\overline{D}%
^{k}\left( \widehat{\Omega _{1}^{2}}h\right) \left( w\right) \right\vert
dV\left( w\right) .
\end{eqnarray*}%
This latter operator is $T_{a,b,c}H\left( z\right) $ with 
\begin{equation*}
a=m_{1}+\sigma ,b=s-n-2-k-\sigma ,c=k+1
\end{equation*}%
and $H\left( w\right) =\left\vert \left( 1-\left\vert w\right\vert
^{2}\right) ^{\sigma }R_{b^{\prime }}^{m}\overline{D}^{k}\left( \widehat{%
\Omega _{1}^{2}}h\right) \left( w\right) \right\vert $. Note that for $%
m>2\left( \frac{n}{p}-\sigma \right) $ we do indeed now have $\left\Vert
H\right\Vert _{L^{p}\left( \lambda _{n}\right) }\approx \left\Vert \widehat{%
\Omega _{1}^{2}}h\right\Vert _{B_{p}^{\sigma }\left( \mathbb{B}_{n}\right) }$%
. The operator here is the same as that above and so Lemma \ref{Zlemma}
applies to prove the desired boundedness on $L^{p}\left( \lambda _{n}\right) 
$.

\subsubsection{The estimate for $\mathcal{F}^{2}$\label{The estimate for F
two}}

Our next task is to obtain the estimate (\ref{accomplishFmu}) for $\mu =2$,
and for this we will show that%
\begin{eqnarray}
&&\int_{\mathbb{B}_{n}}\left\vert \left( 1-\left\vert z\right\vert
^{2}\right) ^{m_{1}+\sigma }R^{m_{1}}\Lambda _{g}\mathcal{C}%
_{n,s_{1}}^{0,0}\Lambda _{g}\mathcal{C}_{n,s_{2}}^{0,1}\Omega
_{2}^{3}\right\vert ^{p}d\lambda _{n}\left( z\right)  \label{want} \\
&\leq &C\int_{\mathbb{B}_{n}}\left\vert \left( 1-\left\vert z\right\vert
^{2}\right) ^{\sigma }\left( 1-\left\vert z\right\vert ^{2}\right)
^{m_{3}^{\prime \prime }}R^{m_{3}^{\prime \prime }}\overline{D}%
^{m_{3}^{\prime }}\left( \widehat{\Omega _{2}^{3}}h\right) \left( z\right)
\right\vert ^{p}d\lambda _{n}\left( z\right) .  \notag
\end{eqnarray}%
Unlike the previous argument we will have to deal with a \emph{rogue} term $%
\left( \overline{z_{2}}-\overline{\xi _{2}}\right) $ this time where there
is no derivative $\frac{\partial }{\partial \overline{\xi _{2}}}$ to
associate to the factor $\left( \overline{z_{2}}-\overline{\xi _{2}}\right) $%
. Again we ignore the contractions $\Lambda _{g}$. Then we use Lemma \ref%
{IBPamel} to perform integration by parts $m_{2}^{\prime }$ times in the
first iterated integral and $m_{3}^{\prime }$ times in the second iterated
integral. We also use Corollary \ref{IBP2iter} to perform integration by
parts in the \emph{radial} derivative $m_{2}^{\prime \prime }$ times in the 
\emph{first} iterated integral (for fixed $z$, $\mathcal{C}%
_{n,s_{2}}^{0,1}\Omega _{2}^{3}\in C\left( \overline{\mathbb{B}_{n}}\right)
\cap C^{\infty }\left( \mathbb{B}_{n}\right) $ by standard estimates \cite%
{Cha}), so that the additional factor $\left( 1-\left\vert \xi \right\vert
^{2}\right) ^{m_{2}^{\prime \prime }}$ can be used crucially in the second
iterated integral, and also $m_{3}^{\prime \prime }$ times in the \emph{%
second} iterated integral for use in acting on $\Omega _{2}^{3}$.

Recall from Lemma \ref{IBPamel} that%
\begin{eqnarray*}
&&\mathcal{C}_{n,s}^{0,q}\eta \left( z\right) =boundary\ terms\ \text{%
(depending on }m\text{)} \\
&&\ \ \ \ \ \ \ \ \ \ \ \ \ \ \ \ \ \ \ \ +\sum_{\ell =0}^{q}\int_{\mathbb{B}%
_{n}}\frac{\left( 1-w\overline{z}\right) ^{n-1-\ell }\left( 1-\left\vert
w\right\vert ^{2}\right) ^{\ell }}{\bigtriangleup \left( w,z\right) ^{n}}%
\left( \frac{1-\left\vert w\right\vert ^{2}}{1-\overline{w}z}\right) ^{s-n}
\\
&&\ \ \ \ \ \ \ \ \ \ \ \ \ \ \ \ \ \ \ \ \ \ \ \ \ \times \left(
\sum_{j=0}^{n-\ell -1}c_{j,\ell ,n,s}\left[ \frac{\left( 1-\left\vert
w\right\vert ^{2}\right) \left( 1-\left\vert z\right\vert ^{2}\right) }{%
\left\vert 1-w\overline{z}\right\vert ^{2}}\right] ^{j}\right) \overline{%
\mathcal{D}}^{m}\eta \left( z\right) .
\end{eqnarray*}%
Recall also that that $\overline{\mathcal{D}}^{m}$ already has the rogue
terms built in, as can be seen from (\ref{effectofDm}). Now we use the right
side above with $q=\ell =j=0$ to substitute for $\mathcal{C}_{n,s_{1}}^{0,0}$%
, and the right side above with $q=\ell =1$ and $j=0$ to substitute for $%
\mathcal{C}_{n,s_{2}}^{0,1}$. Then a typical part of the resulting kernel of
the operator $\mathcal{C}_{n,s_{1}}^{0,0}\mathcal{C}_{n,s_{2}}^{0,1}\Omega
_{2}^{3}\left( z\right) $ is 
\begin{eqnarray}
&&\int_{\mathbb{B}_{n}}\frac{\left( 1-\xi \overline{z}\right) ^{n-1}}{%
\bigtriangleup \left( \xi ,z\right) ^{n}}\left( \frac{1-\left\vert \xi
\right\vert ^{2}}{1-\overline{\xi }z}\right) ^{s_{1}-n}\left( \overline{z_{2}%
}-\overline{\xi _{2}}\right)  \label{prep} \\
&&\times \left( 1-\left\vert \xi \right\vert ^{2}\right) ^{m_{2}^{\prime
}}R^{m_{2}^{\prime }}\overline{\mathcal{D}}^{m_{2}^{\prime \prime }}\int_{%
\mathbb{B}_{n}}\frac{\left( 1-w\overline{\xi }\right) ^{n-2}\left(
1-|w|^{2}\right) }{\bigtriangleup \left( w,\xi \right) ^{n}}\left( \frac{%
1-\left\vert w\right\vert ^{2}}{1-\overline{w}\xi }\right) ^{s_{2}-n}  \notag
\\
&&\times \left( \overline{w_{1}}-\overline{\xi _{1}}\right) \left(
1-\left\vert w\right\vert ^{2}\right) ^{m_{3}^{\prime }}R^{m_{3}^{\prime }}%
\overline{\mathcal{D}}^{m_{3}^{\prime \prime }}\left( \Omega
_{2}^{3}h\right) \left( w\right) dV\left( w\right) dV\left( \xi \right) , 
\notag
\end{eqnarray}%
where we have arbitrarily chosen $\left( \overline{z_{2}}-\overline{\xi _{2}}%
\right) $ and $\left( \overline{w_{1}}-\overline{\xi _{1}}\right) $ as the 
\emph{rogue} factors.

\begin{remark}
\label{conjugatediff}It is important to note that the differential operators 
$\overline{\mathcal{D}_{\zeta }^{m_{2}}}$ are conjugate in the variable $z$
and hence \emph{vanish} on the kernels of the boundary terms $\mathcal{S}%
_{n,s}\left( \overline{\mathcal{D}}^{k}\Omega _{2}^{3}h\right) \left(
z\right) $ in the integration by parts formula (\ref{mformula'}) associated
to the Charpentier solution operator $\mathcal{C}_{n,s_{2}}^{0,1}$ since
these kernels are \emph{holomorphic}. As a result the operator $\overline{%
\mathcal{D}}^{m_{2}^{\prime }}$ hits only the factor $\overline{\mathcal{D}}%
^{k}\Omega _{2}^{3}h$ and a typical term is%
\begin{equation*}
\overline{\left( z_{i}-\zeta _{i}\right) }\frac{\partial }{\partial 
\overline{z_{i}}}\left\{ \overline{\left( w_{i}-z_{i}\right) }\Omega
_{2}^{3}h\right\} =-\overline{\left( z_{i}-\zeta _{i}\right) }\Omega
_{2}^{3}h,
\end{equation*}%
where the derivative $\frac{\partial }{\partial \overline{w_{i}}}$ must
occur in each surviving term in $\Omega _{2}^{3}h$, and this term which is
then handled like the \emph{rogue} terms.
\end{remark}

Now we recall the factorization (\ref{Omegaform}) with $\ell =2$,%
\begin{equation*}
\Omega _{2}^{3}=-4\Omega _{0}^{1}\wedge \widetilde{\Omega _{0}^{1}}\wedge 
\widetilde{\Omega _{0}^{1}},
\end{equation*}%
and that $\Omega _{2}^{3}\left( w\right) $ must have both derivatives $\frac{%
\partial g}{\partial \overline{w_{1}}}$ and $\frac{\partial g}{\partial 
\overline{w_{2}}}$ occurring in it, one surviving in each of the factors $%
\widetilde{\Omega _{0}^{1}}$, along with other harmless powers of $g$ that
we ignore. Thus we may replace $\widetilde{\Omega _{0}^{1}}\wedge \widetilde{%
\Omega _{0}^{1}}$ with $\frac{\partial }{\partial \overline{w_{2}}}\Omega
_{0}^{1}\wedge \frac{\partial }{\partial \overline{w_{1}}}\Omega _{0}^{1}$.
If we use 
\begin{equation*}
\overline{z_{2}}-\overline{\xi _{2}}=\left( \overline{z_{2}}-\overline{w_{2}}%
\right) -\left( \overline{\xi _{2}}-\overline{w_{2}}\right) ,
\end{equation*}%
we can write the above iterated integral as%
\begin{eqnarray*}
&&\int_{\mathbb{B}_{n}}\frac{\left( 1-\xi \overline{z}\right) ^{n-1}}{%
\bigtriangleup \left( \xi ,z\right) ^{n}}\left( \frac{1-\left\vert \xi
\right\vert ^{2}}{1-\overline{\xi }z}\right) ^{s_{1}-n} \\
&&\times \int_{\mathbb{B}_{n}}\left( 1-\left\vert \xi \right\vert
^{2}\right) ^{m_{2}^{\prime \prime }}R^{m_{2}^{\prime \prime }}\overline{%
\mathcal{D}}^{m_{2}^{\prime }}\left\{ \frac{\left( 1-w\overline{\xi }\right)
^{n-2}\left( 1-|w|^{2}\right) }{\bigtriangleup \left( w,\xi \right) ^{n}}%
\left( \frac{1-\left\vert w\right\vert ^{2}}{1-\overline{w}\xi }\right)
^{s_{2}-n}\right\} \\
&&\times \left[ \left( 1-\left\vert w\right\vert ^{2}\right) ^{m_{3}^{\prime
\prime }}R^{m_{3}^{\prime \prime }}\left( \overline{\xi _{2}}-\overline{w_{2}%
}\right) \frac{\partial }{\partial \overline{w_{2}}}\overline{\mathcal{D}}%
^{m_{3}^{\prime }-\ell }\Omega _{0}^{1}\right] \\
&&\wedge \left[ \left( 1-\left\vert w\right\vert ^{2}\right) ^{m_{3}^{\prime
\prime }}R^{m_{3}^{\prime \prime }}\left( \overline{\xi _{1}}-\overline{w_{1}%
}\right) \frac{\partial }{\partial \overline{w_{1}}}\overline{\mathcal{D}}%
^{\ell }\Omega _{0}^{1}\right] dV\left( w\right) dV\left( \xi \right)
\end{eqnarray*}%
minus%
\begin{eqnarray*}
&&\int_{\mathbb{B}_{n}}\frac{\left( 1-\xi \overline{z}\right) ^{n-1}}{%
\bigtriangleup \left( \xi ,z\right) ^{n}}\left( \frac{1-\left\vert \xi
\right\vert ^{2}}{1-\overline{\xi }z}\right) ^{s_{1}-n} \\
&&\times \int_{\mathbb{B}_{n}}\left( 1-\left\vert \xi \right\vert
^{2}\right) ^{m_{2}^{\prime \prime }}R^{m_{2}^{\prime \prime }}\overline{%
\mathcal{D}}^{m_{2}^{\prime }}\left\{ \frac{\left( 1-w\overline{\xi }\right)
^{n-2}\left( 1-|w|^{2}\right) }{\bigtriangleup \left( w,\xi \right) ^{n}}%
\left( \frac{1-\left\vert w\right\vert ^{2}}{1-\overline{w}\xi }\right)
^{s_{2}-n}\right\} \\
&&\times \left[ \left( 1-\left\vert w\right\vert ^{2}\right) ^{m_{3}^{\prime
\prime }}R^{m_{3}^{\prime \prime }}\left( \overline{z_{2}}-\overline{w_{2}}%
\right) \frac{\partial }{\partial \overline{w_{2}}}\overline{\mathcal{D}}%
^{m_{3}^{\prime }-\ell }\Omega _{0}^{1}\right] \\
&&\wedge \left[ \left( 1-\left\vert w\right\vert ^{2}\right) ^{m_{3}^{\prime
\prime }}R^{m_{3}^{\prime \prime }}\left( \overline{\xi _{1}}-\overline{w_{1}%
}\right) \frac{\partial }{\partial \overline{w_{1}}}\overline{\mathcal{D}}%
^{\ell }\Omega _{0}^{1}\right] dV\left( w\right) dV\left( \xi \right) ,
\end{eqnarray*}%
where we have temporarily ignored the wedge products with terms that do not
include derivatives of $g$, as these terms are bounded and so harmless.

Now we apply $\left( 1-\left\vert z\right\vert ^{2}\right) ^{\sigma }\left(
1-\left\vert z\right\vert ^{2}\right) ^{m_{1}^{\prime \prime
}}R^{m_{1}^{\prime \prime }}D^{m_{1}^{\prime }}$ to these operators. Using
the crucial inequalities in Proposition \ref{threecrucial} together with the
factorization (\ref{genformfact}) with $\ell =2$,%
\begin{equation*}
\widehat{\Omega _{2}^{3}}=-4\Omega _{0}^{1}\wedge \widehat{\Omega _{0}^{1}}%
\wedge \widehat{\Omega _{0}^{1}},
\end{equation*}%
the result of this application on the first integral is\ then dominated by%
\begin{eqnarray}
&&\int_{\mathbb{B}_{n}}\frac{\left( 1-\left\vert z\right\vert ^{2}\right)
^{\sigma }\left\vert 1-\xi \overline{z}\right\vert ^{n-1}}{\bigtriangleup
\left( \xi ,z\right) ^{m_{1}^{\prime }+m_{1}^{\prime \prime }+n}}\left[
\left( 1-\left\vert z\right\vert ^{2}\right) \sqrt{\bigtriangleup \left( \xi
,z\right) }\right] ^{m_{1}^{\prime \prime }}  \label{firstdom} \\
&&\times \left\{ \left[ \left( 1-\left\vert z\right\vert ^{2}\right) \sqrt{%
\bigtriangleup \left( \xi ,z\right) }\right] ^{m_{1}^{\prime
}}+\bigtriangleup \left( \xi ,z\right) ^{m_{1}^{\prime }}\right\} \left\vert 
\frac{1-\left\vert \xi \right\vert ^{2}}{1-\xi \overline{z}}\right\vert
^{s_{1}-n}  \notag \\
&&\times \int_{\mathbb{B}_{n}}\frac{\left( 1-\left\vert \xi \right\vert
^{2}\right) ^{m_{2}^{\prime \prime }}\left\vert 1-w\overline{\xi }%
\right\vert ^{n-2}\left( 1-|w|^{2}\right) }{\bigtriangleup \left( w,\xi
\right) ^{m_{2}^{\prime }+m_{2}^{\prime \prime }+n}}\left( \frac{\sqrt{%
\bigtriangleup \left( \xi ,z\right) }}{1-\left\vert \xi \right\vert ^{2}}%
\right) ^{m_{2}^{\prime }}  \notag \\
&&\times \left[ \left( 1-\left\vert \xi \right\vert ^{2}\right) \sqrt{%
\bigtriangleup \left( w,\xi \right) }\right] ^{m_{2}^{\prime \prime
}}\left\{ \left[ \left( 1-\left\vert \xi \right\vert ^{2}\right) \sqrt{%
\bigtriangleup \left( w,\xi \right) }\right] ^{m_{2}^{\prime
}}+\bigtriangleup \left( w,\xi \right) ^{m_{2}^{\prime }}\right\}  \notag \\
&&\times \left\vert \frac{1-\left\vert w\right\vert ^{2}}{1-w\overline{\xi }}%
\right\vert ^{s_{2}-n}\left( \frac{\sqrt{\bigtriangleup \left( w,\xi \right) 
}}{1-\left\vert w\right\vert ^{2}}\right) ^{m_{3}^{\prime }}\left( \frac{%
\sqrt{\bigtriangleup \left( w,\xi \right) }}{1-\left\vert w\right\vert ^{2}}%
\right) ^{2}  \notag \\
&&\times \left\vert \left( 1-\left\vert w\right\vert ^{2}\right)
^{m_{3}^{\prime \prime }}R^{m_{3}^{\prime \prime }}\overline{%
D^{m_{3}^{\prime }}}\left( \widehat{\Omega _{2}^{3}}h\right) \left( w\right)
\right\vert dV\left( w\right) dV\left( \xi \right) ,  \notag
\end{eqnarray}%
and the result of this application on the second integral is dominated by%
\begin{eqnarray}
&&\int_{\mathbb{B}_{n}}\frac{\left( 1-\left\vert z\right\vert ^{2}\right)
^{\sigma }\left\vert 1-\xi \overline{z}\right\vert ^{n-1}}{\bigtriangleup
\left( \xi ,z\right) ^{m_{1}^{\prime }+m_{1}^{\prime \prime }+2}}\left[
\left( 1-\left\vert z\right\vert ^{2}\right) \sqrt{\bigtriangleup \left( \xi
,z\right) }\right] ^{m_{1}^{\prime \prime }}  \label{seconddom} \\
&&\times \left\{ \left[ \left( 1-\left\vert z\right\vert ^{2}\right) \sqrt{%
\bigtriangleup \left( \xi ,z\right) }\right] ^{m_{1}^{\prime
}}+\bigtriangleup \left( \xi ,z\right) ^{m_{1}^{\prime }}\right\} \left\vert 
\frac{1-\left\vert \xi \right\vert ^{2}}{1-\xi \overline{z}}\right\vert
^{s_{1}-n}  \notag \\
&&\times \int_{\mathbb{B}_{n}}\frac{\left( 1-\left\vert \xi \right\vert
^{2}\right) ^{m_{2}^{\prime \prime }}\left\vert 1-w\overline{\xi }%
\right\vert ^{n-2}\left( 1-|w|^{2}\right) }{\bigtriangleup \left( w,\xi
\right) ^{m_{2}^{\prime }+m_{2}^{\prime \prime }+n}}\left( \frac{\sqrt{%
\bigtriangleup \left( \xi ,z\right) }}{1-\left\vert \xi \right\vert ^{2}}%
\right) ^{m_{2}^{\prime }}  \notag \\
&&\times \left[ \left( 1-\left\vert \xi \right\vert ^{2}\right) \sqrt{%
\bigtriangleup \left( w,\xi \right) }\right] ^{m_{2}^{\prime \prime
}}\left\{ \left[ \left( 1-\left\vert \xi \right\vert ^{2}\right) \sqrt{%
\bigtriangleup \left( w,\xi \right) }\right] ^{m}+\bigtriangleup \left(
w,\xi \right) ^{m_{2}^{\prime }}\right\}  \notag \\
&&\times \left\vert \frac{1-\left\vert w\right\vert ^{2}}{1-w\overline{\xi }}%
\right\vert ^{s_{2}-n}\left( \frac{\sqrt{\bigtriangleup \left( w,\xi \right) 
}}{1-\left\vert w\right\vert ^{2}}\right) ^{m_{3}^{\prime }}\left( \frac{%
\sqrt{\bigtriangleup \left( w,z\right) }}{1-\left\vert w\right\vert ^{2}}%
\right) \left( \frac{\sqrt{\bigtriangleup \left( w,\xi \right) }}{%
1-\left\vert w\right\vert ^{2}}\right)  \notag \\
&&\times \left\vert \left( 1-\left\vert w\right\vert ^{2}\right)
^{m_{3}^{\prime \prime }}R^{m_{3}^{\prime \prime }}\overline{%
D^{m_{3}^{\prime }}}\left( \widehat{\Omega _{2}^{3}}h\right) \left( w\right)
\right\vert dV\left( w\right) dV\left( \xi \right) ,  \notag
\end{eqnarray}%
The only difference between these two iterated integrals is that one of the
factors $\frac{\sqrt{\bigtriangleup \left( w,\xi \right) }}{1-\left\vert
w\right\vert ^{2}}$ that occur in the first is replaced by the factor $\frac{%
\sqrt{\bigtriangleup \left( w,z\right) }}{1-\left\vert w\right\vert ^{2}}$
in the second. Note that the ignored wedge products have now been reinstated
in $\widehat{\Omega _{2}^{3}}$.

Now for the iterated integral in (\ref{firstdom}), we can separate it into
the composition of two operators of the form treated previously. One factor
is the operator%
\begin{eqnarray}
&&\int_{\mathbb{B}_{n}}\frac{\left( 1-\left\vert z\right\vert ^{2}\right)
^{\sigma }\left\vert 1-\xi \overline{z}\right\vert ^{n-1}}{\bigtriangleup
\left( \xi ,z\right) ^{m_{1}^{\prime }+m_{1}^{\prime \prime }+n}}\left[
\left( 1-\left\vert z\right\vert ^{2}\right) \sqrt{\bigtriangleup \left( \xi
,z\right) }\right] ^{m_{1}^{\prime \prime }}  \label{firstfact} \\
&&\times \left\{ \left[ \left( 1-\left\vert z\right\vert ^{2}\right) \sqrt{%
\bigtriangleup \left( \xi ,z\right) }\right] ^{m_{1}^{\prime
}}+\bigtriangleup \left( \xi ,z\right) ^{m_{1}^{\prime }}\right\}  \notag \\
&&\times \left( \frac{\sqrt{\bigtriangleup \left( \xi ,z\right) }}{%
1-\left\vert \xi \right\vert ^{2}}\right) ^{m_{2}^{\prime }}\left\vert \frac{%
1-\left\vert \xi \right\vert ^{2}}{1-\xi \overline{z}}\right\vert
^{s_{1}-n}\left( 1-\left\vert \xi \right\vert ^{2}\right) ^{-\sigma }F\left(
\xi \right) dV\left( \xi \right) ,  \notag
\end{eqnarray}%
and the other factor is the operator%
\begin{eqnarray}
&&F\left( \xi \right) =\int_{\mathbb{B}_{n}}\frac{\left( 1-\left\vert \xi
\right\vert ^{2}\right) ^{\sigma }\left\vert 1-w\overline{\xi }\right\vert
^{n-2}\left( 1-|w|^{2}\right) }{\bigtriangleup \left( w,\xi \right)
^{m_{2}^{\prime }+m_{2}^{\prime \prime }+n}}\left[ \left( 1-\left\vert \xi
\right\vert ^{2}\right) \sqrt{\bigtriangleup \left( w,\xi \right) }\right]
^{m_{2}^{\prime \prime }}  \label{secondfact} \\
&&\times \left\{ \left[ \left( 1-\left\vert \xi \right\vert ^{2}\right) 
\sqrt{\bigtriangleup \left( w,\xi \right) }\right] ^{m_{2}^{\prime
}}+\bigtriangleup \left( w,\xi \right) ^{m_{2}^{\prime }}\right\} \left\vert 
\frac{1-\left\vert w\right\vert ^{2}}{1-w\overline{\xi }}\right\vert
^{s_{2}-n}  \notag \\
&&\times \left( \frac{\sqrt{\bigtriangleup \left( w,\xi \right) }}{%
1-\left\vert w\right\vert ^{2}}\right) ^{m_{3}^{\prime }+2}\left(
1-\left\vert w\right\vert ^{2}\right) ^{-\sigma }f\left( w\right) dV\left(
w\right) ,  \notag
\end{eqnarray}%
where $f\left( w\right) =\left( 1-\left\vert w\right\vert ^{2}\right)
^{\sigma }\left\vert \left( 1-\left\vert w\right\vert ^{2}\right)
^{m_{3}^{\prime \prime }}R^{m_{3}^{\prime \prime }}\overline{%
D^{m_{3}^{\prime }}}\left( \widehat{\Omega _{2}^{3}}h\right) \left( w\right)
\right\vert $. We now show how Lemma \ref{Zlemma} applies to obtain the
appropriate boundedness.

We will in fact compare the corresponding kernels to that in (\ref{altL}).
When we consider the summand $\bigtriangleup \left( \xi ,z\right)
^{m_{1}^{\prime }}$ in the middle line of (\ref{firstfact}), the first
operator has kernel%
\begin{eqnarray}
&&\frac{\left( 1-\left\vert z\right\vert ^{2}\right) ^{\sigma +m_{1}^{\prime
\prime }}\left( 1-\left\vert \xi \right\vert ^{2}\right)
^{s_{1}-n-m_{2}^{\prime }-\sigma }}{\left\vert 1-\xi \overline{z}\right\vert
^{s_{1}-2n+1}\bigtriangleup \left( \xi ,z\right) ^{m_{1}^{\prime
}+m_{1}^{\prime \prime }+n-\frac{m_{1}^{\prime \prime }+2m_{1}^{\prime
}+m_{2}^{\prime }}{2}}}  \label{firstker} \\
&=&\frac{\left( 1-\left\vert z\right\vert ^{2}\right) ^{\sigma
+m_{1}^{\prime \prime }}\left( 1-\left\vert \xi \right\vert ^{2}\right)
^{s_{1}-3n-m_{1}^{\prime \prime }-\sigma }}{\left\vert 1-\xi \overline{z}%
\right\vert ^{s_{1}-2n+1}},  \notag
\end{eqnarray}%
if we choose\thinspace $m_{2}^{\prime }=m_{1}^{\prime \prime }+2n$ so that
the factor $\bigtriangleup \left( \xi ,z\right) $ disappears. This is
exactly the same as the kernel of the operator in (\ref{altL}) in the
previous alternative argument but with $m_{1}^{\prime \prime }$ in place of $%
m_{1}^{\prime }$ there. When we consider instead the summand $\left[ \left(
1-\left\vert z\right\vert ^{2}\right) \sqrt{\bigtriangleup \left( \xi
,z\right) }\right] ^{m_{1}^{\prime }}$ in the middle line of (\ref{firstfact}%
), we obtain the kernel in (\ref{firstker}) but with $m_{1}^{\prime \prime
}+m_{1}^{\prime }$ in place of $m_{1}^{\prime \prime }$.

When we consider the summand $\bigtriangleup \left( w,\xi \right)
^{m_{2}^{\prime }}$ in the middle line of (\ref{secondfact}), the second
operator has kernel%
\begin{eqnarray}
&&\frac{\left( 1-\left\vert \xi \right\vert ^{2}\right) ^{m_{2}^{\prime
\prime }+\sigma }\left( 1-|w|^{2}\right) ^{1+s_{2}-n-m_{3}^{\prime
}-2-\sigma }}{\left\vert 1-w\overline{\xi }\right\vert
^{s_{2}-2n+2}\bigtriangleup \left( w,\xi \right) ^{m_{2}^{\prime
}+m_{2}^{\prime \prime }+n-\frac{m_{2}^{\prime \prime }+2m_{2}^{\prime
}+m_{3}^{\prime }+2}{2}}}  \label{secondker} \\
&=&\frac{\left( 1-\left\vert \xi \right\vert ^{2}\right) ^{m_{2}^{\prime
\prime }+\sigma }\left( 1-|w|^{2}\right) ^{s_{2}-3n+1-m_{2}^{\prime \prime
}-\sigma }}{\left\vert 1-w\overline{\xi }\right\vert ^{s_{2}-2n+2}}.  \notag
\end{eqnarray}%
if we choose $m_{3}^{\prime }=m_{2}^{\prime \prime }+2n-2$, and this is also
bounded on $L^{p}\left( d\lambda _{n}\right) $ for $m_{2}^{\prime \prime }$
and $s_{2}$ sufficiently large.

\bigskip

\textbf{Note}: It is here in choosing $m_{2}^{\prime \prime }$ large that we
are using the full force of Corollary \ref{IBP2iter} to perform integration
by parts in the \emph{radial} derivative $m_{2}^{\prime \prime }$ times in
the first iterated integral.

When we consider instead the summand $\left[ \left( 1-\left\vert
z\right\vert ^{2}\right) \sqrt{\bigtriangleup \left( \xi ,z\right) }\right]
^{m_{2}^{\prime }}$ in the middle line of (\ref{secondfact}), we obtain the
kernel in (\ref{secondker}) but with $m_{2}^{\prime \prime }+m_{2}^{\prime }$
in place of $m_{2}^{\prime \prime }$.

\bigskip

To handle the iterated integral in (\ref{seconddom}) we must first deal with
the \emph{rogue} factor $\sqrt{\bigtriangleup \left( w,z\right) }$ whose
variable pair $\left( w,z\right) $ doesn't match that of either of the
denominators $\bigtriangleup \left( \xi ,z\right) $ or $\bigtriangleup
\left( w,\xi \right) $. For this we use the fact that%
\begin{equation*}
\sqrt{\bigtriangleup \left( w,z\right) }=\left\vert 1-w\overline{z}%
\right\vert \left\vert \varphi _{z}\left( w\right) \right\vert =\delta
\left( w,z\right) ^{2}\rho \left( w,z\right) ,
\end{equation*}%
where $\rho \left( w,z\right) =\left\vert \varphi _{z}\left( w\right)
\right\vert $ is the invariant pseudohyperbolic metric on the ball
(Corollary 1.22 in \cite{Zhu}) and where $\delta \left( w,z\right)
=\left\vert 1-w\overline{z}\right\vert ^{\frac{1}{2}}$ satisfies the
triangle inequality on the ball (Proposition 5.1.2 in \cite{Rud}). Thus we
have%
\begin{eqnarray*}
\rho \left( w,z\right) &\leq &\rho \left( \xi ,z\right) +\rho \left( w,\xi
\right) , \\
\delta \left( w,z\right) &\leq &\delta \left( \xi ,z\right) +\delta \left(
w,\xi \right) ,
\end{eqnarray*}%
and so also%
\begin{eqnarray*}
\sqrt{\bigtriangleup \left( w,z\right) } &\leq &2\left[ \delta \left( \xi
,z\right) ^{2}+\delta \left( w,\xi \right) ^{2}\right] \left( \left\vert
\varphi _{z}\left( \xi \right) \right\vert +\left\vert \varphi _{\xi }\left(
w\right) \right\vert \right) \\
&=&2\left( 1+\frac{\left\vert 1-w\overline{\xi }\right\vert }{\left\vert
1-\xi \overline{z}\right\vert }\right) \sqrt{\bigtriangleup \left( \xi
,z\right) }+2\left( 1+\frac{\left\vert 1-\xi \overline{z}\right\vert }{%
\left\vert 1-w\overline{\xi }\right\vert }\right) \sqrt{\bigtriangleup
\left( w,\xi \right) }.
\end{eqnarray*}%
Thus we can write%
\begin{eqnarray}
&&\frac{\sqrt{\bigtriangleup \left( w,z\right) }}{1-\left\vert w\right\vert
^{2}}  \label{canwrite} \\
&\lesssim &\frac{1-\left\vert \xi \right\vert ^{2}}{1-\left\vert
w\right\vert ^{2}}\frac{\sqrt{\bigtriangleup \left( \xi ,z\right) }}{%
1-\left\vert \xi \right\vert ^{2}}+\frac{\left\vert 1-w\overline{\xi }%
\right\vert }{1-\left\vert w\right\vert ^{2}}\frac{1-\left\vert \xi
\right\vert ^{2}}{\left\vert 1-\xi \overline{z}\right\vert }\frac{\sqrt{%
\bigtriangleup \left( \xi ,z\right) }}{1-\left\vert \xi \right\vert ^{2}} 
\notag \\
&&+\frac{\sqrt{\bigtriangleup \left( w,\xi \right) }}{1-\left\vert
w\right\vert ^{2}}+\frac{\left\vert 1-\xi \overline{z}\right\vert }{%
1-\left\vert \xi \right\vert ^{2}}\frac{1-\left\vert \xi \right\vert ^{2}}{%
\left\vert 1-w\overline{\xi }\right\vert }\frac{\sqrt{\bigtriangleup \left(
w,\xi \right) }}{1-\left\vert w\right\vert ^{2}}.  \notag
\end{eqnarray}%
All of the terms on the right hand side of (\ref{canwrite}) are of an
appropriate form to distribute throughout the iterated integral, and again
Lemma \ref{Zlemma} applies to obtain the appropriate boundedness.

For example, the final two terms on the right side of (\ref{canwrite}) that
involve $\frac{\sqrt{\bigtriangleup \left( w,\xi \right) }}{1-\left\vert
w\right\vert ^{2}}$ are handled in the same way as the operator in (\ref%
{firstdom}) by taking $m_{3}^{\prime }=m_{2}^{\prime \prime }+2n-2$ and $%
m_{2}^{\prime }=m_{1}^{\prime \prime }+2n$, and taking $s_{1}$ and $s_{2}$
large as required by the extra factors $\frac{\left\vert 1-\xi \overline{z}%
\right\vert }{1-\left\vert \xi \right\vert ^{2}}\frac{1-\left\vert \xi
\right\vert ^{2}}{\left\vert 1-w\overline{\xi }\right\vert }$. With these
choices the first two terms on the right side of (\ref{canwrite}) that
involve $\frac{\sqrt{\bigtriangleup \left( \xi ,z\right) }}{1-\left\vert \xi
\right\vert ^{2}}$ are then handled using Lemma \ref{Zlemma} with $c=\pm 1$
as follows.

If we substitute the first term $\frac{1-\left\vert \xi \right\vert ^{2}}{%
1-\left\vert w\right\vert ^{2}}\frac{\sqrt{\bigtriangleup \left( \xi
,z\right) }}{1-\left\vert \xi \right\vert ^{2}}$ on the right in (\ref%
{canwrite}) for the factor $\frac{\sqrt{\bigtriangleup \left( w,z\right) }}{%
1-\left\vert w\right\vert ^{2}}$ in (\ref{seconddom}) we get a composition
of two operators as in (\ref{firstfact}) and (\ref{secondfact}) but with the
kernel in (\ref{firstfact}) multiplied by $\frac{\sqrt{\bigtriangleup \left(
\xi ,z\right) }}{1-\left\vert \xi \right\vert ^{2}}$ and the kernel in (\ref%
{secondfact}) multiplied by $\frac{1-\left\vert \xi \right\vert ^{2}}{%
1-\left\vert w\right\vert ^{2}}$ and divided by $\frac{\sqrt{\bigtriangleup
\left( w,\xi \right) }}{1-\left\vert w\right\vert ^{2}}$. If we consider the
summand $\bigtriangleup \left( \xi ,z\right) ^{m_{1}^{\prime }}$ in the
middle line of (\ref{firstfact}), and with the choice $m_{2}^{\prime
}=m_{1}^{\prime \prime }+2n$ already made, the first operator then has kernel%
\begin{eqnarray*}
&&\frac{\sqrt{\bigtriangleup \left( \xi ,z\right) }}{1-\left\vert \xi
\right\vert ^{2}}\times \frac{\left( 1-\left\vert z\right\vert ^{2}\right)
^{\sigma +m_{1}^{\prime \prime }}\left( 1-\left\vert \xi \right\vert
^{2}\right) ^{s_{1}-3n-m_{1}^{\prime \prime }-\sigma }}{\left\vert 1-\xi 
\overline{z}\right\vert ^{s_{1}-2n+1}} \\
&=&\frac{\left( 1-\left\vert z\right\vert ^{2}\right) ^{m_{1}^{\prime \prime
}+\sigma }\left( 1-\left\vert \xi \right\vert ^{2}\right)
^{s_{1}-m_{1}^{\prime \prime }-3n-1-\sigma }\sqrt{\bigtriangleup \left( \xi
,z\right) }}{\left\vert 1-\xi \overline{z}\right\vert ^{s_{1}-2n+1}},
\end{eqnarray*}%
and hence is of the form \thinspace $T_{a,b,c}$ with 
\begin{eqnarray*}
a &=&m_{1}^{\prime \prime }+\sigma , \\
b &=&s_{1}-3n-1-m_{1}^{\prime \prime }-\sigma , \\
c &=&1,
\end{eqnarray*}%
since $a+b+c+n+1=s_{1}-n-1$. Now we apply Lemma \ref{Zlemma} to conclude
that this operator is bounded on $L^{p}\left( \lambda _{n}\right) $ if and
only if%
\begin{equation*}
-p\left( m_{1}^{\prime \prime }+\sigma \right) <-n<p\left(
s_{1}-3n-m_{1}^{\prime \prime }-\sigma \right) ,
\end{equation*}%
i.e. $m_{1}^{\prime \prime }+\sigma >\frac{n}{p}$ and $s_{1}>m_{1}^{\prime
\prime }+\sigma +3n-\frac{n}{p}$.

If we consider the summand $\bigtriangleup \left( w,\xi \right)
^{m_{2}^{\prime }}$ in the middle line of (\ref{secondfact}), and with the
choice $m_{3}^{\prime }=m_{2}^{\prime \prime }+2n-2$ already made, the
second operator has kernel%
\begin{eqnarray*}
&&\frac{1-\left\vert \xi \right\vert ^{2}}{1-\left\vert w\right\vert ^{2}}%
\times \left( \frac{\sqrt{\bigtriangleup \left( w,\xi \right) }}{%
1-\left\vert w\right\vert ^{2}}\right) ^{-1}\times \frac{\left( 1-\left\vert
\xi \right\vert ^{2}\right) ^{m_{2}^{\prime \prime }+\sigma }\left(
1-|w|^{2}\right) ^{s_{2}-3n+1-m_{2}^{\prime \prime }-\sigma }}{\left\vert 1-w%
\overline{\xi }\right\vert ^{s_{2}-2n+2}} \\
&=&\frac{\left( 1-\left\vert \xi \right\vert ^{2}\right) ^{m_{2}^{\prime
\prime }+\sigma +1}\left( 1-|w|^{2}\right) ^{s_{2}-3n+1-m_{2}^{\prime \prime
}-\sigma }\sqrt{\bigtriangleup \left( w,\xi \right) }^{-1}}{\left\vert 1-w%
\overline{\xi }\right\vert ^{s_{2}-2n+2}},
\end{eqnarray*}%
and hence is of the form \thinspace $T_{a,b,c}$ with 
\begin{eqnarray*}
a &=&m_{2}^{\prime \prime }+\sigma +1, \\
b &=&s_{2}-3n+1-m_{2}^{\prime \prime }-\sigma , \\
c &=&-1.
\end{eqnarray*}%
This operator is bounded on $L^{p}\left( \lambda _{n}\right) $ if and only if%
\begin{equation*}
-p\left( m_{2}^{\prime \prime }+\sigma +1\right) <-n<p\left(
s_{2}-3n+2-m_{2}^{\prime \prime }-\sigma \right) ,
\end{equation*}%
i.e. $m_{2}^{\prime \prime }+\sigma >\frac{n}{p}-1$ and $s_{2}>m_{2}^{\prime
\prime }+\sigma +3n-2-\frac{n}{p}$.

If we now substitute the second term $\frac{\left\vert 1-w\overline{\xi }%
\right\vert }{1-\left\vert w\right\vert ^{2}}\frac{1-\left\vert \xi
\right\vert ^{2}}{\left\vert 1-\xi \overline{z}\right\vert }\frac{\sqrt{%
\bigtriangleup \left( \xi ,z\right) }}{1-\left\vert \xi \right\vert ^{2}}$
on the right in (\ref{canwrite}) for the factor $\frac{\sqrt{\bigtriangleup
\left( w,z\right) }}{1-\left\vert w\right\vert ^{2}}$ in (\ref{seconddom})
we similarly get a composition of two operators that are each bounded on $%
L^{p}\left( \lambda _{n}\right) $ for $m_{i}$ and $s_{i}$ chosen large
enough.

\subsubsection{Boundary terms for $\mathcal{F}^{2}$\label{Boundary terms for
F2}}

Now we must address in $\mathcal{F}^{2}$ the boundary terms that arise in
the integration by parts formula (\ref{mformula'}). Suppose the first
operator $\mathcal{C}_{n,s_{1}}^{0,0}$ is replaced by a boundary term, but
not the second. We proceed by applying Corollary \ref{IBP2iter} to the
boundary term. Since the differential operator $\left( 1-\left\vert
z\right\vert ^{2}\right) ^{m_{1}+\sigma }R^{m_{1}}$ hits only the kernel of
the boundary term, we can apply Remark \ref{Zlemma'} to the first iterated
integral and Lemma \ref{Zlemma} to the second iterated integral in the
manner indicated in the above arguments. If the second operator $\mathcal{C}%
_{n,s_{2}}^{0,1}$ is replaced by a boundary term, then as mentioned in
Remark \ref{conjugatediff}, the operators $\overline{D}^{m_{2}}$ hit only
the factors $\overline{\mathcal{D}}^{m_{3}}$, and this produces \emph{rogue}
terms that are handled as above. If the first operator $\mathcal{C}%
_{n,s_{1}}^{0,0}$ was also replaced by a boundary term, then in addition we
would have radial derivatives $R^{m}$ hitting the second boundary term.
Since radial derivatives are holomorphic, they hit only the holomorphic
kernel and not the antiholomorphic factors in $\overline{\mathcal{D}}%
^{m_{3}} $, and so these terms can also be handled as above.

\subsection{The estimates for general $\mathcal{F}^{\protect\mu }$\label{The
estimates for general}}

In view of inequality (\ref{Omegabound}), it suffices to establish the
following inequality:%
\begin{eqnarray}
&&\left\Vert \mathcal{F}^{\mu }\right\Vert _{B_{p}^{\sigma }\left( \mathbb{B}%
_{n}\right) }^{p}  \label{want1} \\
&=&\int_{\mathbb{B}_{n}}\left\vert \left( 1-\left\vert z\right\vert
^{2}\right) ^{m_{1}+\sigma }R^{m_{1}}\Lambda _{g}\mathcal{C}%
_{n,s_{1}}^{0,0}...\Lambda _{g}\mathcal{C}_{n,s_{\mu }}^{0,\mu -1}\Omega
_{\mu }^{\mu +1}h\right\vert ^{p}d\lambda _{n}\left( z\right)  \notag \\
&\leq &C_{\sigma ,n,p,\delta }\int_{\mathbb{B}_{n}}\left\vert \left(
1-\left\vert z\right\vert ^{2}\right) ^{\sigma }\mathcal{X}^{m_{\mu }}\left( 
\widehat{\Omega _{\mu }^{\mu +1}}h\right) \left( z\right) \right\vert
^{p}d\lambda _{n}\left( z\right) .  \notag
\end{eqnarray}%
Recall that the absolute value $\left\vert F\right\vert $ of an element $F$
in the exterior algebra is the square root of the sum of the squares of the
coefficients of $F$ in the standard basis.

The case $\mu >2$ involves no new ideas, and is merely complicated by
straightforward algebra. The reason is that the solution operator $\Lambda
_{g}\mathcal{C}_{n,s_{1}}^{0,0}...\Lambda _{g}\mathcal{C}_{n,s_{\mu
}}^{0,\mu -1}$ acts \emph{separately} in each entry of the form $\Omega
_{\mu }^{\mu +1}h$, an element of the exterior algebra of $\mathbb{C}%
^{\infty }\otimes \mathbb{C}^{n}$ which we view as an alternating $\ell ^{2}$%
-tensor of $\left( 0,\mu \right) $ forms in $\mathbb{C}^{n}$. These
operators decompose as a sum of simpler operators with the basic property
that their kernels are \emph{identical}, except that the rogue factors in
each kernel differ according to the entry. Nevertheless, there are always
exactly $\mu $ distinct rogue factors in each kernel and after splitting,
the $\mu $\ rogue factors can be associated in one-to-one fashion with each
of the $\frac{\partial }{\partial \overline{w_{j}}}$ derivatives in the
corresponding entry of%
\begin{equation*}
\Omega _{\mu }^{\mu +1}h=-\left( \mu +1\right) \left( \sum_{k_{0}=1}^{\infty
}\frac{\overline{g_{k_{0}}}}{\left\vert g\right\vert ^{2}}e_{k_{0}}\right)
\wedge \dbigwedge\limits_{i=1}^{\mu }\left( \sum_{k_{i}=1}^{\infty }\frac{%
\overline{\partial g_{k_{i}}}}{\left\vert g\right\vert ^{2}}e_{k_{i}}\right)
h.
\end{equation*}%
After applying the crucial inequalities, this effectively results in
replacing each derivative $\frac{\partial }{\partial \overline{w_{j}}}$ by
the derivative $\overline{D_{j}}$, and consequently we can write the
resulting form as $\widehat{\Omega _{\mu }^{\mu +1}}h$.

\bigskip

This completes our proof of Theorem \ref{baby}.

\section{Appendix}

Here in the appendix we collect proofs of formulas and modifications of
arguments already in the literature that would otherwise interrupt the main
flow of the paper.

\subsection{Charpentier's solution kernels\label{Charpentier's solution
kernels}}

Here we prove Theorem \ref{explicit}. In the computation of the Cauchy
kernel $\mathcal{C}_{n}(w,z)$, we need to compute the full exterior
derivative of the section $s(w,z)$. By definition one has, 
\begin{eqnarray*}
s_{i}(w,z) &=&\overline{w_{i}}(1-w\overline{z})-\overline{z_{i}}%
(1-\left\vert w\right\vert ^{2}), \\
ds_{i}(w,z) &\equiv &(\partial _{w}+\overline{\partial }_{w}+\partial _{z}+%
\overline{\partial }_{z})s_{i}(w,z)
\end{eqnarray*}%
Straightforward computations show that%
\begin{eqnarray}
\partial _{w}s_{i}\left( w,z\right) &=&\sum_{j=1}^{n}\left( \overline{z}_{i}%
\overline{w}_{j}-\overline{w}_{i}\overline{z}_{j}\right) dw_{j}  \label{ds}
\\
\overline{\partial }_{w}s_{i}\left( w,z\right) &=&\left( 1-w\overline{z}%
\right) d\overline{w}_{i}+\sum_{j=1}^{n}w_{j}\overline{z}_{i}d\overline{w}%
_{j}  \notag \\
\overline{\partial }_{z}s_{i}\left( w,z\right) &=&-\sum_{j=1}^{n}\overline{w}%
_{i}w_{j}d\overline{z}_{j}-\left( 1-|w|^{2}\right) d\overline{z}_{i}  \notag
\\
\partial _{z}s_{i}\left( w,z\right) &=&0,  \notag
\end{eqnarray}%
as well as 
\begin{eqnarray*}
\overline{\partial }_{w}s_{k} &=&(1-w\overline{z})d\overline{w}_{k}+%
\overline{z}_{k}\overline{\partial }_{w}|w|^{2} \\
\overline{\partial }_{z}s_{k} &=&-(1-|w|^{2})d\overline{z}_{k}-\overline{w}%
_{k}\overline{\partial }_{z}(w\overline{z}).
\end{eqnarray*}

We also have the following representations of $s_{k}$, again following by
simple computation. Recall from Notation \ref{perm}\ that $\left\{
1,2,...,n\right\} =\left\{ i_{\nu }\right\} \cup J_{\nu }\cup L_{\nu }$
where $J_{\nu }$ and $L_{\nu }$ are increasing multi-indices of lengths $%
n-q-1$ and $q$. We will use the following with $k=i_{\nu }$.%
\begin{eqnarray*}
s_{k} &=&(\overline{w}_{k}-\overline{z}_{k})+\sum_{l\neq k}w_{l}(\overline{w}%
_{l}\overline{z}_{k}-\overline{w}_{k}\overline{z}_{l}) \\
&=&(\overline{w}_{k}-\overline{z}_{k})+\sum_{j\in J_{\nu }}w_{j}(\overline{w}%
_{j}\overline{z}_{k}-\overline{w}_{k}\overline{z}_{j})+\sum_{l\in L_{\nu
}}w_{l}(\overline{w}_{l}\overline{z}_{k}-\overline{w}_{k}\overline{z}_{l}) \\
&=&(\overline{w}_{k}-\overline{z}_{k})+\overline{z}_{k}\sum_{j\in J_{\nu
}}|w_{j}|^{2}-\overline{w}_{k}\sum_{j\in J_{\nu }}w_{j}\overline{z}_{j}+%
\overline{z}_{k}\sum_{l\in L_{\nu }}|w_{l}|^{2}-\overline{w}_{k}\sum_{l\in
L_{\nu }}w_{l}\overline{z}_{l}.
\end{eqnarray*}

\begin{remark}
Since $A\wedge A=0$ for any form, we have in particular that $\overline{%
\partial _{w}}\left\vert w\right\vert ^{2}\wedge \overline{\partial _{w}}%
\left\vert w\right\vert ^{2}=0$ and $\overline{\partial _{z}}\left( w%
\overline{z}\right) \wedge \overline{\partial _{z}}\left( w\overline{z}%
\right) =0$.
\end{remark}

Using this remark we next compute $\bigwedge_{j\in J_{\nu }}\overline{%
\partial }_{w}s_{j}$. We identify $J_{\nu }$ as $j_{1}<j_{2}<\cdots
<j_{n-q-1}$ and define a map $\imath (j_{r})=r$, namely $\imath $ says where 
$j_{r}$ occurs in the multi-index. We will frequently abuse notation and
simply write $\imath (j)$. Because $\overline{\partial }_{w}|w|^{2}\wedge 
\overline{\partial }_{w}|w|^{2}=0$ it is easy to conclude that we can not
have any term in $\overline{\partial }_{w}|w|^{2}$ of degree greater than
one when expanding the wedge product of the $\overline{\partial _{w}}s_{j}$. 
\begin{eqnarray*}
\bigwedge_{j\in J_{\nu }}\overline{\partial }_{w}s_{j} &=&\bigwedge_{j\in
J_{\nu }}\left\{ (1-w\overline{z})d\overline{w}_{j}+\overline{z}_{j}%
\overline{\partial }_{w}|w|^{2}\right\} \\
&=&(1-w\overline{z})^{n-q-1}\bigwedge_{j\in J_{\nu }}d\overline{w}_{j}+(1-w%
\overline{z})^{n-q-2}\sum_{j\in J_{\nu }}(-1)^{\imath (j)-1}\overline{z}_{j}%
\overline{\partial }_{w}|w|^{2}\wedge \bigwedge_{j^{\prime }\in J_{\nu
}\setminus \{j\}}d\overline{w}_{j^{\prime }} \\
&=&(1-w\overline{z})^{n-q-2} \\
&&\left( \left( 1-w\overline{z}+\sum_{j\in J_{\nu }}w_{j}\overline{z}%
_{j}\right) \bigwedge_{j\in J_{\nu }}d\overline{w}_{j}+\sum_{j\in J_{\nu
}}(-1)^{\imath (j)-1}\overline{z}_{j}\sum_{k\in L_{\nu }\cup \{i_{\nu
}\}}w_{k}d\overline{w}_{k}\bigwedge_{j^{\prime }\in J_{\nu }\setminus \{j\}}d%
\overline{w}_{j^{\prime }}\right) .
\end{eqnarray*}

The last line follows by direct computation using 
\begin{equation*}
\overline{\partial }_{w}\left\vert w\right\vert ^{2}=\sum_{j\in J_{\nu
}}w_{j}d\overline{w}_{j}+\sum_{k\in L_{\nu }\cup \{i_{\nu }\}}w_{k}d%
\overline{w}_{k}.
\end{equation*}%
A similar computation yields that

\begin{eqnarray*}
&&\bigwedge_{l\in L_{\nu }}\overline{\partial }_{z}s_{l} \\
&=&(-1)^{q}\bigwedge_{l\in L_{\nu }}\left\{ (1-|w|^{2})d\overline{z}_{l}+%
\overline{w}_{l}\overline{\partial }_{z}(w\overline{z})\right\} \\
&=&(-1)^{q}\left( (1-|w|^{2})^{q}\bigwedge_{l\in L_{\nu }}d\overline{z}%
_{l}+(1-|w|^{2})^{q-1}\sum_{l\in L_{\nu }}(-1)^{\imath (l)-1}\overline{w}_{l}%
\overline{\partial }_{z}(w\overline{z})\wedge \bigwedge_{l^{\prime }\in
L_{\nu }\setminus \{l\}}d\overline{z}_{l^{\prime }}\right) \\
&=&(-1)^{q}(1-|w|^{2})^{q-1} \\
&&\left( \left( 1-|w|^{2}+\sum_{l\in L_{\nu }}|w_{l}|^{2}\right)
\bigwedge_{l\in L_{\nu }}d\overline{z}_{l}+\sum_{l\in L_{\nu }}(-1)^{\imath
(l)-1}\overline{w}_{l}\sum_{k\in J_{\nu }\cup \{i_{\nu }\}}w_{k}d\overline{z}%
_{k}\bigwedge_{l^{\prime }\in L_{\nu }\setminus \{l\}}d\overline{z}%
_{l^{\prime }}\right) .
\end{eqnarray*}

An important remark at this point is that the multi-index $J_{\nu }$ or $%
L_{\nu }$ can only appear in the first term of the last line above. The
terms after the plus sign have multi-indices that are related to $J_{\nu }$
and $L_{\nu }$, but differ by one element. This fact will play a role later.

Combining things, we see that 
\begin{equation*}
\bigwedge_{j\in J_{\nu }}\overline{\partial }_{w}s_{j}\bigwedge_{l\in L_{\nu
}}\overline{\partial }_{z}s_{j}=(-1)^{q}(1-w\overline{z}%
)^{n-q-2}(1-|w|^{2})^{q-1}\left( I_{\nu }+II_{\nu }+III_{\nu }+IV_{\nu
}\right) ,
\end{equation*}%
where 
\begin{equation*}
I_{\nu }=\left( 1-w\overline{z}+\sum_{j\in J_{\nu }}w_{j}\overline{z}%
_{j}\right) \left( 1-|w|^{2}+\sum_{l\in L_{\nu }}|w_{l}|^{2}\right)
\bigwedge_{j\in J_{\nu }}d\overline{w}_{j}\bigwedge_{l\in L_{\nu }}d%
\overline{z}_{l},
\end{equation*}%
\begin{equation*}
II_{\nu }=\left( 1-w\overline{z}+\sum_{j\in J_{\nu }}w_{j}\overline{z}%
_{j}\right) \bigwedge_{j\in J_{\nu }}d\overline{w}_{j}\left( \sum_{l\in
L_{\nu }}(-1)^{\imath (l)-1}\overline{w}_{l}\sum_{k\in J_{\nu }\cup \{i_{\nu
}\}}w_{k}d\overline{z}_{k}\bigwedge_{l^{\prime }\in L_{\nu }\setminus \{l\}}d%
\overline{z}_{l^{\prime }}\right) ,
\end{equation*}%
\begin{equation*}
III_{\nu }=\left( \sum_{j\in J_{\nu }}(-1)^{\imath (j)-1}\overline{z}%
_{j}\sum_{k\in L_{\nu }\cup \{i_{\nu }\}}w_{k}d\overline{w}%
_{k}\bigwedge_{j^{\prime }\in J_{\nu }\setminus \{j\}}d\overline{w}%
_{j^{\prime }}\right) \left( 1-|w|^{2}+\sum_{l\in L_{\nu
}}|w_{l}|^{2}\right) \bigwedge_{l\in L_{\nu }}d\overline{z}_{l},
\end{equation*}%
\begin{eqnarray*}
IV_{\nu } &=&\left( \sum_{j\in J_{\nu }}(-1)^{\imath (j)-1}\overline{z}%
_{j}\sum_{k\in L_{\nu }\cup \{i_{\nu }\}}w_{k}d\overline{w}%
_{k}\bigwedge_{j^{\prime }\in J_{\nu }\setminus \{j\}}d\overline{w}%
_{j^{\prime }}\right) \\
&&\times \left( \sum_{l\in L_{\nu }}(-1)^{\imath (l)-1}\overline{w}%
_{l}\sum_{k\in J_{\nu }\cup \{i_{\nu }\}}w_{k}d\overline{z}%
_{k}\bigwedge_{l^{\prime }\in L_{\nu }\setminus \{l\}}d\overline{z}%
_{l^{\prime }}\right) .
\end{eqnarray*}

We next introduce a little more notation to aid in the computation of the
kernel $\mathcal{C}_{n}^{0,q}(w,z)$. For $1\leq k\leq n$ we let $%
P_{n}^{q}(k)=\{\nu \in P_{n}^{q}:\nu (1)=i_{\nu }=k\}$. This divides the set 
$P_{n}^{q}$ into $n$ classes with $\frac{(n-1)!}{(n-q-1)!q!}$ elements. At
this point, with the notation introduced in Notation \ref{perm} and
computations performed above, we have reduced the calculation of $\mathcal{C}%
_{n}^{0,q}(w,z)$ to 
\begin{eqnarray*}
\mathcal{C}_{n}^{0,q}(w,z) &=&\frac{1}{\bigtriangleup (w,z)^{n}}\sum_{\nu
\in P_{n}^{q}}\epsilon _{\nu }s_{i_{\nu }}\bigwedge_{j\in J_{\nu }}\overline{%
\partial }_{w}s_{j}\bigwedge_{l\in L_{\nu }}\overline{\partial }%
_{z}s_{l}\wedge \omega (w) \\
&=&\frac{(-1)^{q}(1-w\overline{z})^{n-q-2}(1-|w|^{2})^{q-1}}{\bigtriangleup
(w,z)^{n}}\sum_{k=1}^{n}s_{k}\sum_{\nu \in P_{n}^{q}(k)}\epsilon _{\nu
}(I_{\nu }+II_{\nu }+III_{\nu }+IV_{\nu }) \\
&=&\frac{(-1)^{q}(1-w\overline{z})^{n-q-2}(1-|w|^{2})^{q-1}}{\bigtriangleup
(w,z)^{n}}\sum_{k=1}^{n}s_{k}(I(k)+II(k)+III(k)+IV(k)) \\
&=&\frac{(-1)^{q}(1-w\overline{z})^{n-q-2}(1-|w|^{2})^{q-1}}{\bigtriangleup
(w,z)^{n}}\sum_{k=1}^{n}s_{k}C(k).
\end{eqnarray*}%
Here we have defined $C(k)\equiv I(k)+II(k)+III(k)+IV(k)$, and 
\begin{eqnarray*}
I(k)\equiv \sum_{\nu \in P_{n}^{q}(k)}\epsilon _{\nu }I_{\nu } &\quad
&II(k)\equiv \sum_{\nu \in P_{n}^{q}(k)}\epsilon _{\nu }II_{\nu } \\
III(k)\equiv \sum_{\nu \in P_{n}^{q}(k)}\epsilon _{\nu }III_{\nu } &\quad
&IV(k)\equiv \sum_{\nu \in P_{n}^{q}(k)}\epsilon _{\nu }IV_{\nu }.
\end{eqnarray*}

For a fixed $\tau \in P_{n}^{q}$ we will compute the coefficient of $%
\bigwedge_{j\in J_{\tau }}d\overline{w}_{j}\bigwedge_{l\in L_{\tau }}d%
\overline{z}_{l}$. We will ignore the functional coefficient in front of the
sum since it only needs to be taken into consideration at the final stage.
We will show that for this fixed $\tau $ the sum on $k$ of $s_{k}$ times $%
I(k)$, $II(k)$, $III(k)$ and $IV(k)$ can be replaced by $\epsilon _{\tau
}(1-w\overline{z})(1-|w|^{2})(\overline{w}_{i_{\tau }}-\overline{z}_{i_{\tau
}})\bigwedge_{j\in J_{\tau }}d\overline{w}_{j}\bigwedge_{l\in L_{\tau }}d%
\overline{z}_{l}$. There will also be other terms that appear in this
expression that arise from multi-indices $J$ and $I$ that are not disjoint.
Using the computations below it can be seen that these terms actually vanish
and hence provide no contribution for $\mathcal{C}_{n}^{0,q}(w,z)$. Since $%
\tau $ is an arbitrary element of $P_{n}^{q}$ this will then complete the
computation of the kernel.

Note that when $k=i_{\tau }$ then we have the following contributions. It is
easy to see that $II(i_{\tau })=III(i_{\tau })=0$. It is also easy to see
that 
\begin{eqnarray*}
I(i_{\tau }) &=&\epsilon _{\tau }\left( 1-w\overline{z}+\sum_{j\in J_{\tau
}}w_{j}\overline{z}_{j}\right) \left( 1-|w|^{2}+\sum_{l\in L_{\tau
}}|w_{l}|^{2}\right) \bigwedge_{j\in J_{\tau }}d\overline{w}%
_{j}\bigwedge_{l\in L_{\tau }}d\overline{z}_{l} \\
&=&\epsilon _{\tau }(1-w\overline{z})(1-|w|^{2})\bigwedge_{j\in J_{\tau }}d%
\overline{w}_{j}\bigwedge_{l\in L_{\tau }}d\overline{z}_{l} \\
&&+\left( (1-w\overline{z})\sum_{l\in L_{\tau
}}|w_{l}|^{2}+(1-|w|^{2})\sum_{j\in J_{\tau }}w_{j}\overline{z}%
_{j}+\sum_{l\in L_{\tau }}|w_{l}|^{2}\sum_{j\in J_{\tau }}w_{j}\overline{z}%
_{j}\right) \bigwedge_{j\in J_{\tau }}d\overline{w}_{j}\bigwedge_{l\in
L_{\tau }}d\overline{z}_{l}.
\end{eqnarray*}

We also receive a contribution from term $IV(i_{\tau })$ is this case. This
happens by interchanging an index in the multi-index $J_{\tau }$ with one in 
$L_{\tau }$. Namely, we consider the permutations $\nu :\{1,\ldots
,n\}\rightarrow \{i_{\tau },(J_{\tau }\setminus \{j\})\cup \{l\},(L_{\tau
}\setminus \{l\})\cup \{j\}\}.$ This permutation contributes the term $%
\overline{z}_{l}w_{l}\overline{w}_{j}w_{j}$. After summing over all these
possible permutations, we arrive at the simplified formula, 
\begin{equation*}
IV(i_{\tau })=-\epsilon _{\tau }\left( \sum_{j\in J_{\tau
}}|w_{j}|^{2}\right) \left( \sum_{l\in L_{\tau }}w_{l}\overline{z}%
_{l}\right) \bigwedge_{j\in J_{\tau }}d\overline{w}_{j}\bigwedge_{l\in
L_{\tau }}d\overline{z}_{l}.
\end{equation*}

Collecting all these terms, when $k=i_{\tau }$ we have that the coefficient
of $\epsilon _{\tau }\bigwedge_{j\in J_{\tau }}d\overline{w}%
_{j}\bigwedge_{l\in L_{\tau }}d\overline{z}_{l}$ is: 
\begin{eqnarray*}
C(i_{\tau }) &=&(1-w\overline{z})(1-|w|^{2})+(1-w\overline{z}+\sum_{j\in
J_{\tau }}w_{j}\overline{z}_{j})\sum_{l\in L_{\tau }}|w_{l}|^{2} \\
&&+(1-|w|^{2}+\sum_{l\in L_{\tau }}|w_{l}|^{2})\sum_{j\in J_{\tau }}w_{j}%
\overline{z}_{j}-\sum_{l\in L_{\tau }}|w_{l}|^{2}\sum_{j\in J_{\tau }}w_{j}%
\overline{z}_{j}-\sum_{j\in J_{\tau }}|w_{j}|^{2}\sum_{l\in L_{\tau }}w_{l}%
\overline{z}_{l}.
\end{eqnarray*}

We next note that when $k\neq i_{\tau }$ it is still possible to have terms
which contribute to the coefficient of $\bigwedge_{j\in J_{\tau }}d\overline{%
w}_{j}\bigwedge_{l\in L_{\tau }}d\overline{z}_{l}$. To see this we further
split the conditions on $k$ into the situations where $k\in J_{\tau }$ and $%
k\in L_{\tau }$. First, observe in this situation that if $k\neq i_{\tau }$
then term $I(k)$ can never contribute. So all contributions must come from
terms $II(k)$, $III(k)$, and $IV(k)$. In these terms it is possible to
obtain the term $\bigwedge_{j\in J_{\tau }}d\overline{w}_{j}\bigwedge_{l\in
L_{\tau }}d\overline{z}_{l}$ by replacing some index in $\nu $. Namely, it
is possible to have $\nu $ and $\tau $ differ by one index from each other,
or one by replacing an index with $i_{\tau }$.

Next, observe that when $k\in L_\tau$ there exists a unique $\nu\in P^q_n(k)$
such that $J_\nu=J_\tau$. Namely, we have that $\nu:\{1,\ldots, n\}\to\{k,
J_\tau, (L_\tau\setminus\{k\})\cup i_\tau\}$. Here, we used that $i_\nu=k$.
Terms of this type will contribute to term $II(k)$ but will give no
contribution to term $III(k)$. However, they will give a contribution to
term $IV(k)$.

Similarly, when $k\in J_\tau$ there will exist a unique $\mu\in P^q_n(k)$
with $L_\mu=L_\tau$. This happens with $\mu:\{1,\ldots, n\}\to \{k,
(J_\tau\setminus\{k\})\cup i_\tau, L_\tau\}$. Here we used that $i_\mu=k$.
Again, we get a contribution to term $III(k)$ and $IV(k)$ and they give no
contribution to the term $II(k)$.

Using these observations when $k\in L_{\tau }$ we arrive at the following
for $I(k)$, $II(k)$, $III(k)$, and $IV(k)$: 
\begin{eqnarray*}
I(k) &=&0 \\
II(k) &=&-\epsilon _{\tau }\left( 1-w\overline{z}+\sum_{j\in J_{\tau }}w_{j}%
\overline{z}_{j}\right) \overline{w}_{i_{\tau }}w_{k}\bigwedge_{j\in J_{\tau
}}d\overline{w}_{j}\bigwedge_{l\in L_{\tau }}d\overline{z}_{l} \\
III(k) &=&0 \\
IV(k) &=&\epsilon _{\tau }\overline{z}_{i_{\tau }}w_{k}\left( \sum_{j\in
J_{\tau }}|w_{j}|^{2}\right) \bigwedge_{j\in J_{\tau }}d\overline{w}%
_{j}\bigwedge_{l\in L_{\tau }}d\overline{z}_{l}.
\end{eqnarray*}

Similarly, when $k\in J_{\tau }$ we arrive at the following for $I(k)$, $%
II(k)$, $III(k)$, and $IV(k)$: 
\begin{eqnarray*}
I(k) &=&0 \\
II(k) &=&0 \\
III(k) &=&-\epsilon _{\tau }\left( 1-|w|^{2}+\sum_{l\in L_{\tau
}}|w_{l}|^{2}\right) \overline{z}_{i_{\tau }}w_{k}\bigwedge_{j\in J_{\tau }}d%
\overline{w}_{j}\bigwedge_{l\in L_{\tau }}d\overline{z}_{l} \\
IV(k) &=&\epsilon _{\tau }\overline{w}_{i_{\tau }}w_{k}\left( \sum_{l\in
L_{\tau }}w_{l}\overline{z}_{l}\right) \bigwedge_{j\in J_{\tau }}d\overline{w%
}_{j}\bigwedge_{l\in L_{\tau }}d\overline{z}_{l}.
\end{eqnarray*}

Collecting these terms, we see the following for the coefficient of $%
\epsilon _{\tau }\bigwedge_{j\in J_{\tau }}d\overline{w}_{j}\bigwedge_{l\in
L_{\tau }}d\overline{z}_{l}$: 
\begin{equation*}
\begin{array}{ccc}
C(k)=-w_{k}\left( \overline{z}_{i_{\tau }}\left( 1-|w|^{2}+\sum_{l\in
L_{\tau }}|w_{l}|^{2}\right) -\overline{w}_{i_{\tau }}\left( \sum_{l\in
L_{\tau }}w_{l}\overline{z}_{l}\right) \right) &  & \forall k\in J_{\tau },
\\ 
&  &  \\ 
C(k)=-w_{k}\left( \overline{w}_{i_{\tau }}\left( 1-w\overline{z}+\sum_{j\in
J_{\tau }}w_{j}\overline{z}_{j}\right) -\overline{z}_{i_{\tau }}\left(
\sum_{j\in J_{\tau }}|w_{j}|^{2}\right) \right) &  & \forall k\in L_{\tau }.%
\end{array}%
\end{equation*}

This then implies that the \textit{total} coefficient of $\epsilon _{\tau
}\bigwedge_{j\in J_{\tau }}d\overline{w}_{j}\bigwedge_{l\in L_{\tau }}d%
\overline{z}_{l}$ is given by 
\begin{equation*}
s_{i_{\tau }}C(i_{\tau })+\sum_{k\in J_{\tau }}s_{k}C(k)+\sum_{k\in L_{\tau
}}s_{k}C(k).
\end{equation*}

At this point the remainder of the proof of the Theorem \ref{explicit}
reduces to tedious algebra. The term $s_{i_{\tau }}C(i_{\tau })$ will
contribute the term $(1-w\overline{z})(1-|w|^{2})(\overline{w}_{i_{\tau }}-%
\overline{z}_{i_{\tau }})$ and a remainder term. The remainder term will
cancel with the terms $\sum_{k\neq i_{\tau }}s_{k}C(k)$.

We first compute the term $s_{k}C(k)$ for $k\in J_{\tau }$. Note that in
this case, we have that 
\begin{eqnarray*}
C(k) &=&w_{k}\left( \overline{w}_{i_{\tau }}\left( \sum_{l\in L_{\tau }}w_{l}%
\overline{z}_{l}\right) -\overline{z}_{i_{\tau }}\left( 1-|w|^{2}+\sum_{l\in
L_{\tau }}|w_{l}|^{2}\right) \right) \\
&=&w_{k}\left( \overline{w}_{i_{\tau }}\left( \sum_{l\in L_{\tau }}w_{l}%
\overline{z}_{l}\right) -\overline{z}_{i_{\tau }}\left( 1-\sum_{l\in J_{\tau
}}|w_{l}|^{2}\right) \right) +w_{k}\overline{z}_{i_{\tau }}|w_{i_{\tau
}}|^{2}.
\end{eqnarray*}

Multiplying this by $s_{k}$ we see that 
\begin{eqnarray*}
s_{k}C(k) &=&(1-w\overline{z})\left( \overline{w}_{i_{\tau }}\left(
\sum_{l\in L_{\tau }}w_{l}\overline{z}_{l}\right) -\overline{z}_{i_{\tau
}}\left( 1-\sum_{l\in J_{\tau }}|w_{l}|^{2}\right) \right) |w_{k}|^{2} \\
&&-(1-|w|^{2})\left( \overline{w}_{i_{\tau }}\left( \sum_{l\in L_{\tau
}}w_{l}\overline{z}_{l}\right) -\overline{z}_{i_{\tau }}\left( 1-\sum_{l\in
J_{\tau }}|w_{l}|^{2}\right) \right) w_{k}\overline{z}_{k} \\
&&+(1-w\overline{z})\overline{z}_{i_{\tau }}|w_{i_{\tau
}}|^{2}|w_{k}|^{2}-(1-|w|^{2})\overline{z}_{i_{\tau }}|w_{i_{\tau
}}|^{2}w_{k}\overline{z}_{k}.
\end{eqnarray*}

Upon summing in $k\in J_{\tau }$ we find that 
\begin{eqnarray*}
\sum_{k\in J_{\tau }}s_{k}C(k) &=&(1-w\overline{z})\left( \overline{w}%
_{i_{\tau }}\left( \sum_{l\in L_{\tau }}w_{l}\overline{z}_{l}\right) -%
\overline{z}_{i_{\tau }}\left( 1-\sum_{j\in J_{\tau }}|w_{j}|^{2}\right)
\right) \sum_{k\in J_{\tau }}|w_{k}|^{2} \\
&&-(1-|w|^{2})\left( \overline{w}_{i_{\tau }}\left( \sum_{l\in L_{\tau
}}w_{l}\overline{z}_{l}\right) -\overline{z}_{i_{\tau }}\left( 1-\sum_{j\in
J_{\tau }}|w_{j}|^{2}\right) \right) \sum_{k\in J_{\tau }}w_{k}\overline{z}%
_{k} \\
&&+(1-w\overline{z})\overline{z}_{i_{\tau }}|w_{i_{\tau }}|^{2}\sum_{k\in
J_{\tau }}|w_{k}|^{2}-(1-|w|^{2})\overline{z}_{i_{\tau }}|w_{i_{\tau
}}|^{2}\sum_{k\in J_{\tau }}w_{k}\overline{z}_{k}.
\end{eqnarray*}

Performing similar computations for $k\in L_{\tau }$ we find, 
\begin{eqnarray*}
\sum_{k\in L_{\tau }}s_{k}C(k) &=&(1-w\overline{z})\left( \overline{z}%
_{i_{\tau }}\left( \sum_{k\in J_{\tau }}|w_{j}|^{2}\right) -\overline{w}%
_{i_{\tau }}\left( 1-\sum_{l\in L_{\tau }}w_{l}\overline{z}_{l}\right)
\right) \sum_{k\in L_{\tau }}|w_{k}|^{2} \\
&&-(1-|w|^{2})\left( \overline{z}_{i_{\tau }}\left( \sum_{k\in J_{\tau
}}|w_{j}|^{2}\right) -\overline{w}_{i_{\tau }}\left( 1-\sum_{l\in L_{\tau
}}w_{l}\overline{z}_{l}\right) \right) \sum_{k\in L_{\tau }}w_{k}\overline{z}%
_{k} \\
&&+(1-w\overline{z})\overline{z}_{i_{\tau }}|w_{i_{\tau }}|^{2}\sum_{k\in
L_{\tau }}|w_{k}|^{2}-(1-|w|^{2})\overline{z}_{i_{\tau }}|w_{i_{\tau
}}|^{2}\sum_{k\in L_{\tau }}w_{k}\overline{z}_{k}.
\end{eqnarray*}%
Putting this all together we find that 
\begin{eqnarray*}
&&\sum_{k\neq i_{\tau }}s_{k}C(k) \\
&=&\overline{w}_{i_{\tau }}(1-w\overline{z})\left( \left( \sum_{k\in L_{\tau
}}w_{l}\overline{z}_{l}\right) \left( \sum_{k\in J_{\tau
}}|w_{k}|^{2}\right) -\left( 1-\sum_{k\in L_{\tau }}w_{k}\overline{z}%
_{k}-w_{i_{\tau }}\overline{z}_{i_{\tau }}\right) \left( \sum_{k\in L_{\tau
}}|w_{k}|^{2}\right) \right) \\
&&+\overline{z}_{i_{\tau }}(1-|w|^{2})\left( \left( 1-\sum_{k\in J_{\tau
}}|w_{k}|^{2}-|w_{i_{\tau }}|^{2}\right) \left( \sum_{k\in J_{\tau }}w_{k}%
\overline{z}_{k}\right) -\left( \sum_{k\in J_{\tau }}|w_{j}|^{2}\right)
\left( \sum_{k\in L_{\tau }}w_{k}\overline{z}_{k}\right) \right) \\
&&-\overline{z}_{i_{\tau }}(1-w\overline{z})(1-|w|^{2})\left( \sum_{k\in
J_{\tau }}|w_{j}|^{2}\right) +\overline{w}_{i_{\tau }}(1-w\overline{z}%
)(1-|w|^{2})\left( \sum_{k\in L_{\tau }}w_{k}\overline{z}_{k}\right) .
\end{eqnarray*}

We next compute the term $s_{i_{\tau }}C(i_{\tau })$. Using the properties
of $s_{k}$ we have that $s_{i_{\tau }}C(i_{\tau })$ is%
\begin{eqnarray*}
&&\left( \overline{w}_{i_{\tau }}-\overline{z}_{i_{\tau }}\right) (1-w%
\overline{z})(1-|w|^{2}) \\
&&+\overline{z}_{i_{\tau }}(1-w\overline{z})(1-|w|^{2})\left( \sum_{k\in
J_{\tau }}|w_{k}|^{2}\right) -\overline{w}_{i_{\tau }}(1-w\overline{z}%
)(1-|w|^{2})\left( \sum_{k\in L_{\tau }}w_{k}\overline{z}_{k}\right) \\
&&+\overline{w}_{i_{\tau }}(1-w\overline{z})\left\{ (1-w\overline{z})\left(
\sum_{k\in L_{\tau }}|w_{k}|^{2}\right) +\left( \sum_{k\in L_{\tau
}}|w_{k}|^{2}\right) \left( \sum_{k\in J_{\tau }}w_{k}\overline{z}%
_{k}\right) \right. \\
&&\ \ \ \ \ \ \ \ \ \ \ \ \ \ \ \ \ \ \ \ \ \ \ \ \ \ \ \ \ \ \left. -\left(
\sum_{k\in J_{\tau }}|w_{k}|^{2}\right) \left( \sum_{k\in L_{\tau }}w_{k}%
\overline{z}_{k}\right) \right\} \\
&&+\overline{z}_{i_{\tau }}(1-|w|^{2})\left\{ -(1-|w|^{2})\left( \sum_{k\in
J_{\tau }}w_{k}\overline{z}_{k}\right) -\left( \sum_{k\in L_{\tau
}}|w_{k}|^{2}\right) \left( \sum_{k\in J_{\tau }}w_{k}\overline{z}%
_{k}\right) \right. \\
&&\ \ \ \ \ \ \ \ \ \ \ \ \ \ \ \ \ \ \ \ \ \ \ \ \ \ \ \ \ \ \left. +\left(
\sum_{k\in J_{\tau }}|w_{k}|^{2}\right) \left( \sum_{k\in L_{\tau }}w_{k}%
\overline{z}_{k}\right) \right\} .
\end{eqnarray*}

From this point on it is simple to see that the remainder of the term $%
s_{i_{\tau }}C(i_{\tau })$ cancels with $\sum_{k\neq i_{\tau }}s_{k}C(k)$.
One simply adds and subtracts a common term in parts of $\sum_{k\neq i_{\tau
}}s_{k}C(k)$. The only term that remains is $(\overline{w}_{i_{\tau }}-%
\overline{z}_{i\tau })(1-w\overline{z})(1-|w|^{2})$. Thus, we see that the
term corresponding to $\tau $ in the sum $\mathcal{C}_{n}^{0,q}(w,z)$ is 
\begin{equation*}
\epsilon _{\tau }\frac{(-1)^{q}(1-w\overline{z})^{n-q-2}(1-|w|^{2})^{q-1}}{%
\bigtriangleup (w,z)^{n}}(1-w\overline{z})(1-|w|^{2})(\overline{w}_{i_{\tau
}}-\overline{z}_{i_{\tau }})\bigwedge_{j\in J_{\tau }}d\overline{w}%
_{j}\bigwedge_{l\in L_{\tau }}d\overline{z}_{l}\wedge \omega (w).
\end{equation*}%
Since $\tau $ was arbitrary we conclude that $\mathcal{C}_{n}^{0,q}(w,z)$
equals%
\begin{equation*}
\frac{\left( 1-w\overline{z}\right) ^{n-q-1}\left( 1-\left\vert w\right\vert
^{2}\right) ^{q}}{\bigtriangleup \left( w,z\right) ^{n}}
\end{equation*}%
times 
\begin{equation*}
\sum_{\nu \in P_{n}^{q}}\epsilon _{\nu }(\overline{w}_{i_{\nu }}-\overline{z}%
_{i_{\nu }})\bigwedge_{j\in J_{\nu }}d\overline{w}_{j}\bigwedge_{l\in L_{\nu
}}d\overline{z}_{l}\wedge \omega (w),
\end{equation*}%
which completes the proof of Theorem \ref{explicit}.

\subsubsection{Explicit formulas for kernels in $n=2$ and $3$ dimensions 
\label{Explicit formulas}}

Using the above computations and simplifying algebra we obtain the formulas 
\begin{eqnarray}
&&\mathcal{C}_{2}^{0,0}(w,z)  \label{C002} \\
&=&\frac{(1-w\overline{z})}{\bigtriangleup (w,z)^{2}}\left[ (\overline{z}%
_{2}-\overline{w}_{2})d\overline{w}_{1}\wedge dw_{1}\wedge dw_{2}-(\overline{%
z}_{1}-\overline{w}_{1})d\overline{w}_{2}\wedge dw_{1}\wedge dw_{2}\right] ,
\notag
\end{eqnarray}

and 
\begin{eqnarray}
&&\mathcal{C}_{2}^{0,1}\left( w,z\right)  \label{C012} \\
&=&\frac{(1-|w|^{2})}{\bigtriangleup (w,z)^{2}}\left[ (\overline{w}_{2}-%
\overline{z}_{2})d\overline{z}_{1}\wedge dw_{1}\wedge dw_{2}-(\overline{w}%
_{1}-\overline{z}_{1})d\overline{z}_{2}\wedge dw_{1}\wedge dw_{2}\right] , 
\notag
\end{eqnarray}

and%
\begin{eqnarray}
&&\mathcal{C}_{3}^{0,q}\left( w,z\right)  \label{C0q3} \\
&=&\sum_{\sigma \in \mathcal{S}_{3}}sgn\left( \sigma \right) \frac{\left( 1-w%
\overline{z}\right) ^{2-q}\left( 1-\left\vert w\right\vert ^{2}\right)
^{q}\left( \overline{z_{\sigma \left( 1\right) }}-\overline{w_{\sigma \left(
1\right) }}\right) }{\bigtriangleup \left( w,z\right) ^{3}}d\overline{\zeta
_{\sigma \left( 2\right) }}\wedge d\overline{\zeta _{\sigma \left( 3\right) }%
}\wedge \omega _{3}\left( w\right) ,  \notag
\end{eqnarray}%
where $\mathcal{S}_{3}$ denotes the group of permutations on $\left\{
1,2,3\right\} $ and $q$ determines the number of $d\overline{z_{i}}$ in the
form $d\overline{\zeta _{\sigma \left( 2\right) }}\wedge d\overline{\zeta
_{\sigma \left( 3\right) }}$: 
\begin{equation*}
d\overline{\zeta _{\sigma \left( 2\right) }}\wedge d\overline{\zeta _{\sigma
\left( 3\right) }}=\left\{ 
\begin{array}{ccc}
d\overline{w_{\sigma \left( 2\right) }}\wedge d\overline{w_{\sigma \left(
3\right) }} & \text{ if } & q=0 \\ 
d\overline{z_{\sigma \left( 2\right) }}\wedge d\overline{w_{\sigma \left(
3\right) }} & \text{ if } & q=1 \\ 
d\overline{z_{\sigma \left( 2\right) }}\wedge d\overline{z_{\sigma \left(
3\right) }} & \text{ if } & q=2%
\end{array}%
\right. .
\end{equation*}

\subsubsection{Integrating in higher dimensions\label{Integrating in higher
dimensions}}

Here we give the proof of Lemma \ref{amelcoeff}. Let%
\begin{equation*}
B\equiv \frac{\left( 1-\left\vert z\right\vert ^{2}\right) }{\left\vert 1-w%
\overline{z}\right\vert ^{2}}\text{ and }R\equiv \sqrt{1-|w|^{2}},
\end{equation*}%
so that%
\begin{equation*}
BR^{2}=\frac{\left( 1-|w|^{2}\right) \left( 1-|z|^{2}\right) }{\left\vert 1-w%
\overline{z}\right\vert ^{2}}=1-\left\vert \varphi _{w}(z)\right\vert ^{2}.
\end{equation*}%
Then with the change of variable $\rho =Br^{2}$ we obtain%
\begin{eqnarray*}
&&\left( 1-w\overline{z}\right) ^{s-q-1}\int_{\sqrt{1-\left\vert
w\right\vert ^{2}}\mathbb{B}_{k}}\frac{(1-\left\vert w\right\vert
^{2}-\left\vert w^{\prime }\right\vert ^{2})^{q}}{\bigtriangleup (\left(
w,w^{\prime }\right) ,\left( z,0\right) )^{s}}dV(w^{\prime }) \\
&=&\frac{\left( 1-w\overline{z}\right) ^{s-q-1}}{\left\vert 1-w\overline{z}%
\right\vert ^{2s}}\int_{\sqrt{1-|w|^{2}}\mathbb{B}_{k}}\frac{\left(
1-|w|^{2}-\left\vert w^{\prime }\right\vert ^{2}\right) ^{q}}{\left( 1-\frac{%
(1-|z|^{2})}{\left\vert 1-w\overline{z}\right\vert ^{2}}\left( 1-\left\vert
w\right\vert ^{2}-\left\vert w^{\prime }\right\vert ^{2}\right) \right) ^{s}}%
dV\left( w^{\prime }\right) \\
&=&\frac{\left( 1-w\overline{z}\right) ^{s-q-1}}{\left\vert 1-w\overline{z}%
\right\vert ^{2s}}\int_{0}^{R}\frac{\left( R^{2}-r^{2}\right) ^{q}}{\left(
1-BR^{2}+Br^{2}\right) ^{s}}r^{2k-1}dr \\
&=&\frac{\left( 1-w\overline{z}\right) ^{s-q-1}}{2B^{k+q}\left\vert 1-w%
\overline{z}\right\vert ^{2s}}\int_{0}^{BR^{2}}\frac{\left( BR^{2}-\rho
\right) ^{q}}{\left( 1-BR^{2}+\rho \right) ^{s}}\rho ^{k-1}d\rho ,
\end{eqnarray*}%
which with%
\begin{equation*}
\Psi _{n,k}^{0,q}\left( t\right) =\frac{\left( 1-t\right) ^{n}}{t^{k}}%
\int_{0}^{t}\frac{\left( t-\rho \right) ^{q}}{\left( 1-t+\rho \right) ^{n+k}}%
\rho ^{k-1}d\rho ,
\end{equation*}%
we rewrite as%
\begin{eqnarray*}
&&\frac{\left( 1-w\overline{z}\right) ^{s-q-1}}{2B^{k+q}\left\vert 1-w%
\overline{z}\right\vert ^{2s}}\frac{(BR^{2})^{k}}{\left\vert \varphi
_{w}(z)\right\vert ^{2n}}\Psi _{n,k}^{0,q}\left( BR^{2}\right) \\
&=&\frac{\left( 1-w\overline{z}\right) ^{s-q-1}\left( 1-\left\vert
w\right\vert ^{2}\right) ^{k}}{2\left( 1-\left\vert z\right\vert ^{2}\right)
^{q}\left\vert 1-w\overline{z}\right\vert ^{2s}}\frac{\left\vert 1-w%
\overline{z}\right\vert ^{2q}}{\left\vert \varphi _{w}(z)\right\vert ^{2n}}%
\Psi _{n,k}^{0,q}\left( BR^{2}\right) \\
&=&\frac{\left( 1-w\overline{z}\right) ^{s-q-1}\left( 1-\left\vert
w\right\vert ^{2}\right) ^{k}}{2\left( 1-\left\vert z\right\vert ^{2}\right)
^{q}}\frac{\left\vert 1-w\overline{z}\right\vert ^{2q-2k}}{\bigtriangleup
\left( w,z\right) ^{n}}\Psi _{n,k}^{0,q}\left( BR^{2}\right) \\
&=&\frac{1}{2}\Phi _{n}^{q}\left( w,z\right) \left( \frac{1-\left\vert
w\right\vert ^{2}}{1-\overline{w}z}\right) ^{k-q}\left( \frac{1-w\overline{z}%
}{1-\left\vert z\right\vert ^{2}}\right) ^{q}\Psi _{n,k}^{0,q}\left(
BR^{2}\right) .
\end{eqnarray*}%
since $\Phi _{n}^{q}\left( w,z\right) =\frac{\left( 1-w\overline{z}\right)
^{n-1-q}\left( 1-\left\vert w\right\vert ^{2}\right) ^{q}}{\bigtriangleup
\left( w,z\right) ^{n}}$.

At this point we claim that%
\begin{equation}
\Psi _{n,k}^{0,q}\left( t\right) =\frac{\left( 1-t\right) ^{n}}{t^{k}}%
\int_{0}^{t}\frac{\left( t-r\right) ^{q}}{\left( 1-t+r\right) ^{n+k}}%
r^{k-1}dr  \label{defpsi}
\end{equation}%
is a polynomial in%
\begin{equation*}
t=BR^{2}=1-\left\vert \varphi _{w}(z)\right\vert ^{2}
\end{equation*}%
of degree $n-1$ that vanishes to order $q$ at $t=0$, so that 
\begin{equation*}
\Psi _{n,k}^{0,q}\left( t\right) =\sum_{j=q}^{n-1}c_{j,n,s}\left( \frac{%
\left( 1-|w|^{2}\right) \left( 1-|z|^{2}\right) }{\left\vert 1-w\overline{z}%
\right\vert ^{2}}\right) ^{j},
\end{equation*}%
With this claim established, the proof of Lemma \ref{amelcoeff} is complete.

To see that $\Psi _{n,k}^{0,q}$ vanishes of order $q$ at $t=0$ is easy since
for $t$ small (\ref{defpsi}) yields%
\begin{equation*}
\left\vert \Psi _{n,k}^{0,q}\left( t\right) \right\vert \leq
Ct^{-k}\int_{0}^{t}\frac{t^{q}}{C}r^{k-1}dr\leq Ct^{q}.
\end{equation*}%
To see that$\Psi _{n,k}^{0,q}$ is a polynomial of degree $n-1$ we prove two
recursion formulas valid for $0\leq t<1$ (we let $t\rightarrow 1$ at the end
of the argument):%
\begin{eqnarray}
\Psi _{n,k}^{0,q}\left( t\right) -\Psi _{n,k}^{0,q+1}\left( t\right)
&=&\left( 1-t\right) \Psi _{n-1,k}^{0,q}\left( t\right) ,  \label{tworecur}
\\
\Psi _{n,k}^{0,0}\left( t\right) &=&\frac{1}{k}\left( 1-t\right) ^{n}+\frac{%
n+k}{k}t\Psi _{n,k+1}^{0,0}\left( t\right) .  \notag
\end{eqnarray}%
The first formula follows from%
\begin{equation*}
\left( t-r\right) ^{q}-\left( t-r\right) ^{q+1}=\left( t-r\right) ^{q}\left(
1-t+r\right) ,
\end{equation*}%
while the second is an integration by parts:%
\begin{eqnarray*}
\int_{0}^{t}\frac{r^{k-1}}{\left( 1-t+r\right) ^{n+k}}dr &=&\frac{1}{k}\frac{%
r^{k}}{\left( 1-t+r\right) ^{n+k}}\mid _{0}^{t} \\
&&+\frac{n+k}{k}\int_{0}^{t}\frac{r^{k}}{\left( 1-t+r\right) ^{n+k+1}}dr \\
&=&\frac{1}{k}t^{k}+\frac{n+k}{k}\int_{0}^{t}\frac{r^{k}}{\left(
1-t+r\right) ^{n+k+1}}dr.
\end{eqnarray*}%
If we multiply this equality through by $\frac{(1-t)^{n}}{t^{k}}$ we obtain
the second formula in (\ref{tworecur}).

The recursion formulas in (\ref{tworecur}) reduce matters to proving that $%
\Psi _{n,1}^{0,0}$ is a polynomial of degree $n-1$. Indeed, once we know
that $\Psi _{n,1}^{0,0}$ is a polynomial of degree $n-1$, then the second
formula in (\ref{tworecur}) and induction on $k$ shows that $\Psi
_{n,k}^{0,0}$ is as well. Then the first formula and induction on $q$ then
shows that $\Psi _{n,k}^{0,q}$ is also. To see that $\Psi _{n,1}^{0,0}$ is a
polynomial of degree $n-1$ we compute%
\begin{eqnarray*}
\Psi _{n,1}^{0,0}\left( t\right) &=&\frac{(1-t)^{n}}{t}\int_{0}^{t}\frac{1}{%
\left( 1-t+r\right) ^{n+1}}dr \\
&=&\frac{(1-t)^{n}}{t}\left\{ -\frac{1}{n\left( 1-t+r\right) ^{n}}\right\}
\mid _{0}^{t} \\
&=&\frac{1-(1-t)^{n}}{nt},
\end{eqnarray*}%
which is a polynomial of degree $n-1$. This finishes the proof of the claim,
and hence that of Lemma \ref{amelcoeff} as well.

\subsection{Integration by parts formulas in the ball\label{Integration by
parts formulas}}

We begin by proving the generalized analogue of the integration by parts
formula of Ortega and Fabrega \cite{OrFa} as given in Lemma \ref{IBP1'}. For
this we will use the following identities.

\begin{lemma}
\label{collect}For $\ell \in \mathbb{Z}$, we have 
\begin{eqnarray}
\overline{\mathcal{Z}}\left\{ \bigtriangleup \left( w,z\right) ^{\ell
}\right\} &=&\ell \bigtriangleup \left( w,z\right) ^{\ell },  \label{Zbar1}
\\
\overline{\mathcal{Z}}\left\{ \left( 1-w\overline{z}\right) ^{\ell }\right\}
&=&0,  \notag \\
\overline{\mathcal{Z}}\left\{ \left( 1-\left\vert w\right\vert ^{2}\right)
^{\ell }\right\} &=&\ell \left( 1-\left\vert w\right\vert ^{2}\right) ^{\ell
}-\ell \left( 1-\left\vert w\right\vert ^{2}\right) ^{\ell -1}\left( 1-%
\overline{z}w\right) .  \notag
\end{eqnarray}
\end{lemma}

\textbf{Proof}: (of Lemma \ref{collect}) The computation 
\begin{eqnarray*}
\frac{\partial \bigtriangleup }{\partial \overline{w_{j}}} &=&\frac{\partial 
}{\partial \overline{w_{j}}}\left\{ \left\vert 1-w\overline{z}\right\vert
^{2}-\left( 1-\left\vert w\right\vert ^{2}\right) \left( 1-\left\vert
z\right\vert ^{2}\right) \right\} \\
&=&\left( w\overline{z}-1\right) z_{j}+\left( 1-\left\vert z\right\vert
^{2}\right) w_{j},
\end{eqnarray*}%
shows that $\overline{\mathcal{Z}}\bigtriangleup =\bigtriangleup $:%
\begin{eqnarray*}
\overline{\mathcal{Z}}\bigtriangleup \left( w,z\right) &=&\left(
\sum_{j=1}^{n}\left( \overline{w_{j}}-\overline{z_{j}}\right) \frac{\partial 
}{\partial \overline{w_{j}}}\right) \left\{ \left\vert 1-w\overline{z}%
\right\vert ^{2}-\left( 1-\left\vert w\right\vert ^{2}\right) \left(
1-\left\vert z\right\vert ^{2}\right) \right\} \\
&=&\sum_{j=1}^{n}\left( \overline{w_{j}}-\overline{z_{j}}\right) \left\{
\left( w\overline{z}-1\right) z_{j}+\left( 1-\left\vert z\right\vert
^{2}\right) w_{j}\right\} \\
&=&\left( \overline{w}z-\left\vert z\right\vert ^{2}\right) \left( w%
\overline{z}-1\right) +\left( 1-\left\vert z\right\vert ^{2}\right) \left(
\left\vert w\right\vert ^{2}-\overline{z}w\right) \\
&=&-\overline{w\overline{z}}+\left\vert z\right\vert ^{2}+\left\vert 
\overline{w}z\right\vert ^{2}-\left\vert z\right\vert ^{2}w\overline{z}%
+\left\vert w\right\vert ^{2}-w\overline{z}-\left\vert z\right\vert
^{2}\left\vert w\right\vert ^{2}+\left\vert z\right\vert ^{2}w\overline{z} \\
&=&-2\func{Re}w\overline{z}+\left\vert z\right\vert ^{2}+\left\vert w%
\overline{z}\right\vert ^{2}+\left\vert w\right\vert ^{2}-\left\vert
z\right\vert ^{2}\left\vert w\right\vert ^{2} \\
&=&\left\vert w-z\right\vert ^{2}+\left\vert w\overline{z}\right\vert
^{2}-\left\vert z\right\vert ^{2}\left\vert w\right\vert ^{2}=\bigtriangleup
\left( w,z\right)
\end{eqnarray*}%
by the second line in (\ref{Dis}) above. Iteration then gives the first line
in (\ref{Zbar1}). The second line is trivial since $1-w\overline{z}$ is
holomorphic in $w$. The third line follows by iterating%
\begin{equation*}
\overline{\mathcal{Z}}\left( 1-\left\vert w\right\vert ^{2}\right) =%
\overline{z}w-\left\vert w\right\vert ^{2}=\left( 1-\left\vert w\right\vert
^{2}\right) -\left( 1-\overline{z}w\right) .
\end{equation*}

\textbf{Proof of Lemma \ref{IBP1'}}: We use the general formula (\ref%
{gensolutionker'}) for the solution kernels $\mathcal{C}_{n}^{0,q}$ to prove
(\ref{mformula'}) by induction on $m$. For $m=0$ we obtain%
\begin{equation}
\mathcal{C}_{n}^{0,q}\eta \left( z\right) =c_{0}\int_{\mathbb{B}_{n}}\Phi
_{n}^{q}\left( w,z\right) \left\{ \sum_{\left\vert J\right\vert =q}\overline{%
\mathcal{D}^{0}}\left( \eta \lrcorner d\overline{w}^{J}\right) d\overline{z}%
^{J}\right\} dV\left( w\right) \equiv c_{0}\Phi _{n}^{q}\left( \overline{%
\mathcal{D}^{0}}\eta \right) \left( z\right) ,  \label{m=0'}
\end{equation}%
from (\ref{key}) and the following calculation using (\ref{gensolutionker}): 
\begin{eqnarray*}
&&\mathcal{C}_{n}^{0,q}\eta \left( z\right) \\
&\equiv &\int_{\mathbb{B}_{n}}\mathcal{C}_{n}^{0,q}\left( w,z\right) \wedge
\eta \left( w\right) \\
&=&\int_{\mathbb{B}_{n}}\sum_{\left\vert J\right\vert =q}\Phi _{n}^{q}\left(
w,z\right) \sum_{k\notin J}\left( -1\right) ^{\mu \left( k,J\right) }\left( 
\overline{z_{k}}-\overline{\eta _{k}}\right) d\overline{z}^{J}\wedge d%
\overline{w}^{\left( J\cup \left\{ k\right\} \right) ^{c}}\wedge \omega
_{n}\left( w\right) \wedge \left( \sum_{\left\vert I\right\vert =q+1}\eta
_{I}d\overline{w_{I}}\right) \\
&=&\left\{ \int_{\mathbb{B}_{n}}\Phi _{n}^{q}\left( w,z\right) \left[
\sum_{\left\vert J\right\vert =q}\sum_{k\notin J}\left( -1\right) ^{\mu
\left( k,J\right) }\left( \overline{z_{k}}-\overline{w_{k}}\right) \eta
_{J\cup \left\{ k\right\} }d\overline{z}^{J}\right] dV\left( w\right)
\right\} .
\end{eqnarray*}

Now we consider the case $m=1$. First we note that for each $J$ with $%
\left\vert J\right\vert =q$,%
\begin{equation}
\overline{\mathcal{Z}}\overline{\mathcal{D}^{0}}\left( \eta \lrcorner d%
\overline{w}^{J}\right) -\overline{\mathcal{D}^{0}}\left( \eta \lrcorner d%
\overline{w}^{J}\right) =\overline{\mathcal{D}^{1}}\left( \eta \lrcorner d%
\overline{w}^{J}\right) .  \label{indint}
\end{equation}%
Indeed, we compute%
\begin{eqnarray*}
\overline{\mathcal{Z}}\overline{\mathcal{D}^{0}}\left( \eta \lrcorner d%
\overline{w}^{J}\right) &=&\left( \sum_{j=1}^{n}\left( \overline{w_{j}}-%
\overline{z_{j}}\right) \frac{\partial }{\partial \overline{w_{j}}}\right)
\left( \sum_{k\notin J}\left( \overline{w_{k}}-\overline{z_{k}}\right)
\sum_{I\setminus J=\left\{ k\right\} }\left( -1\right) ^{\mu \left(
k,J\right) }\eta _{I}\right) \\
&=&\sum_{j=1}^{n}\sum_{k\notin J}\sum_{I\setminus J=\left\{ k\right\}
}\left( -1\right) ^{\mu \left( k,J\right) }\left( \overline{w_{j}}-\overline{%
z_{j}}\right) \left( \overline{w_{k}}-\overline{z_{k}}\right) \frac{\partial 
}{\partial \overline{w_{j}}}\eta _{I} \\
&&+\sum_{k\notin J}\left( \overline{w_{k}}-\overline{z_{k}}\right)
\sum_{I\setminus J=\left\{ k\right\} }\left( -1\right) ^{\mu \left(
k,J\right) }\eta _{I},
\end{eqnarray*}%
so that%
\begin{eqnarray*}
&&\overline{\mathcal{Z}}\overline{\mathcal{D}^{0}}\left( \eta \lrcorner d%
\overline{w}^{J}\right) -\overline{\mathcal{D}^{0}}\left( \eta \lrcorner d%
\overline{w}^{J}\right) \\
&=&\sum_{j=1}^{n}\sum_{k\notin J}\sum_{I\setminus J=\left\{ k\right\}
}\left( -1\right) ^{\mu \left( k,J\right) }\left( \overline{w_{j}}-\overline{%
z_{j}}\right) \left( \overline{w_{k}}-\overline{z_{k}}\right) \frac{\partial 
}{\partial \overline{w_{j}}}\eta _{I}=\overline{\mathcal{D}^{1}}\left( \eta
\lrcorner d\overline{w}^{J}\right) .
\end{eqnarray*}

For $\left\vert J\right\vert =q$ and $0\leq \ell \leq q$ define%
\begin{equation*}
\mathcal{I}_{J}^{\ell }\equiv \sum_{j=1}^{n}\int_{\mathbb{B}_{n}}\frac{%
\partial }{\partial \overline{w_{j}}}\left\{ \frac{\left( 1-w\overline{z}%
\right) ^{n-1-\ell }\left( 1-\left\vert w\right\vert ^{2}\right) ^{\ell }}{%
\bigtriangleup \left( w,z\right) ^{n}}\left( \overline{w_{j}}-\overline{z_{j}%
}\right) \overline{\mathcal{D}^{0}}\left( \eta \lrcorner d\overline{w}%
^{J}\right) \right\} \omega \left( \overline{w}\right) \wedge \omega \left(
w\right) .
\end{equation*}%
By (3) and (4) of Proposition 16.4.4 in \cite{Rud} we have%
\begin{equation*}
\sum_{j=1}^{n}\left( -1\right) ^{j-1}\left( \overline{w_{j}}-\overline{z_{j}}%
\right) \dbigwedge\limits_{k\neq j}d\overline{w_{k}}\wedge \omega \left(
w\right) \mid _{\partial \mathbb{B}_{n}}=c\left( 1-\overline{z}w\right)
d\sigma \left( w\right) ,
\end{equation*}%
and Stokes' theorem then yields%
\begin{equation*}
\mathcal{I}_{J}^{\ell }=c\int_{\partial \mathbb{B}_{n}}\frac{\left( 1-w%
\overline{z}\right) ^{n-\ell }\left( 1-\left\vert w\right\vert ^{2}\right)
^{\ell }}{\bigtriangleup \left( w,z\right) ^{n}}\overline{\mathcal{D}^{0}}%
\left( \eta \lrcorner d\overline{w}^{J}\right) d\sigma \left( w\right) =0,
\end{equation*}%
since $\ell \geq 1$ and $1-\left\vert w\right\vert ^{2}$ vanishes on $%
\partial \mathbb{B}_{n}$. Moreover, from Lemma \ref{collect} we obtain%
\begin{eqnarray*}
\mathcal{I}_{J}^{\ell } &=&n\int_{\mathbb{B}_{n}}\frac{\left( 1-w\overline{z}%
\right) ^{n-1-\ell }\left( 1-\left\vert w\right\vert ^{2}\right) ^{\ell }}{%
\bigtriangleup \left( z,w\right) ^{n}}\overline{\mathcal{D}^{0}}\left( \eta
\lrcorner d\overline{w}^{J}\right) dV\left( w\right) \\
&&+\int_{\mathbb{B}_{n}}\overline{\mathcal{Z}}\left\{ \frac{\left( 1-w%
\overline{z}\right) ^{n-1-\ell }\left( 1-\left\vert w\right\vert ^{2}\right)
^{\ell }}{\bigtriangleup \left( z,w\right) ^{n}}\overline{\mathcal{D}^{0}}%
\left( \eta \lrcorner d\overline{w}^{J}\right) \right\} dV\left( w\right) \\
&=&\int_{\mathbb{B}_{n}}\frac{\left( 1-w\overline{z}\right) ^{n-1-\ell
}\left( 1-\left\vert w\right\vert ^{2}\right) ^{\ell }}{\bigtriangleup
\left( z,w\right) ^{n}}\overline{\mathcal{Z}}\overline{\mathcal{D}^{0}}%
\left( \eta \lrcorner d\overline{w}^{J}\right) dV\left( w\right) \\
&&+\ell \int_{\mathbb{B}_{n}}\frac{\left( 1-w\overline{z}\right) ^{n-1-\ell
}\left( 1-\left\vert w\right\vert ^{2}\right) ^{\ell }}{\bigtriangleup
\left( z,w\right) ^{n}}\overline{\mathcal{D}^{0}}\left( \eta \lrcorner d%
\overline{w}^{J}\right) dV\left( w\right) \\
&&-\ell \int_{\mathbb{B}_{n}}\frac{\left( 1-w\overline{z}\right) ^{n-\ell
}\left( 1-\left\vert w\right\vert ^{2}\right) ^{\ell -1}}{\bigtriangleup
\left( z,w\right) ^{n}}\overline{\mathcal{D}^{0}}\left( \eta \lrcorner d%
\overline{w}^{J}\right) dV\left( w\right) .
\end{eqnarray*}%
Combining this with (\ref{indint})\ and (\ref{m=0'}) yields%
\begin{eqnarray*}
\Phi _{n}^{\ell }\left( \overline{\mathcal{D}^{0}}\eta \right) \left(
z\right) &=&\sum_{J}\int_{\mathbb{B}_{n}}\Phi _{n}^{\ell }\left( w,z\right) 
\overline{\mathcal{D}^{0}}\left( \eta \lrcorner d\overline{w}^{J}\right)
dV\left( w\right) d\overline{z}^{J} \\
&=&\sum_{J}\int_{\mathbb{B}_{n}}\Phi _{n}^{\ell }\left( w,z\right) \overline{%
\mathcal{Z}}\overline{\mathcal{D}^{0}}\left( \eta \lrcorner d\overline{w}%
^{J}\right) dV\left( w\right) d\overline{z}^{J} \\
&&-\sum_{J}\int_{\mathbb{B}_{n}}\Phi _{n}^{\ell }\left( w,z\right) \overline{%
\mathcal{D}^{1}}\left( \eta \lrcorner d\overline{w}^{J}\right) dV\left(
w\right) d\overline{z}^{J} \\
&=&-\sum_{J}\int_{\mathbb{B}_{n}}\Phi _{n}^{\ell }\left( w,z\right) 
\overline{\mathcal{D}^{1}}\left( \eta \lrcorner d\overline{w}^{J}\right)
dV\left( w\right) d\overline{z}^{J} \\
&&-\ell \sum_{J}\int_{\mathbb{B}_{n}}\Phi _{n}^{\ell }\left( w,z\right) 
\overline{\mathcal{D}^{0}}\left( \eta \lrcorner d\overline{w}^{J}\right)
dV\left( w\right) d\overline{z}^{J} \\
&&+\ell \sum_{J}\int_{\mathbb{B}_{n}}\Phi _{n}^{\ell -1}\left( w,z\right) 
\overline{\mathcal{D}^{0}}\left( \eta \lrcorner d\overline{w}^{J}\right)
dV\left( w\right) d\overline{z}^{J} \\
&=&-\Phi _{n}^{\ell }\left( \overline{\mathcal{D}^{1}}\eta \right) \left(
z\right) -\ell \Phi _{n}^{\ell }\left( \overline{\mathcal{D}^{0}}\eta
\right) \left( z\right) +\ell \Phi _{n}^{\ell -1}\left( \overline{\mathcal{D}%
^{0}}\eta \right) \left( z\right) .
\end{eqnarray*}

Thus we have%
\begin{equation}
\Phi _{n}^{\ell }\left( \overline{\mathcal{D}^{0}}\eta \right) \left(
z\right) =-\frac{1}{\ell +1}\Phi _{n}^{\ell }\left( \overline{\mathcal{D}^{1}%
}\eta \right) \left( z\right) +\frac{\ell }{\ell +1}\Phi _{n}^{\ell
-1}\left( \overline{\mathcal{D}^{0}}\eta \right) \left( z\right) .
\label{once}
\end{equation}%
From (\ref{m=0'}) and then iterating (\ref{once}) we obtain%
\begin{eqnarray}
&&\mathcal{C}_{n}^{\left( 0,q\right) }\eta \left( z\right)  \label{onceiter}
\\
&=&\Phi _{n}^{q}\left( \overline{\mathcal{D}^{0}}\eta \right) \left(
z\right) =-\frac{1}{q+1}\Phi _{n}^{q}\left( \overline{\mathcal{D}^{1}}\eta
\right) \left( z\right) +\frac{q}{q+1}\Phi _{n}^{q-1}\left( \overline{%
\mathcal{D}^{0}}\eta \right) \left( z\right)  \notag \\
&=&-\frac{1}{q+1}\Phi _{n}^{q}\left( \overline{\mathcal{D}^{1}}\eta \right)
\left( z\right) +\frac{q}{q+1}\left\{ -\frac{1}{q}\Phi _{n}^{q-1}\left( 
\overline{\mathcal{D}^{1}}\eta \right) \left( z\right) +\frac{q-1}{q}\Phi
_{n}^{q-2}\left( \overline{\mathcal{D}^{0}}\eta \right) \left( z\right)
\right\}  \notag \\
&=&-\frac{1}{q+1}\sum_{\ell =1}^{q}\Phi _{n}^{\ell }\left( \overline{%
\mathcal{D}^{1}}\eta \right) \left( z\right) +boundary\ term.  \notag
\end{eqnarray}%
Thus we have obtained the second sum in (\ref{mformula'}) with $c_{\ell }=-%
\frac{1}{q+1}$ for $1\leq \ell \leq q$ in the case $m=1$.

We have included $boundary\ term$ in (\ref{onceiter}) since when we use
Stokes' theorem on $\Phi _{n}^{0}\left( \overline{\mathcal{D}^{0}}\eta
\right) $ the boundary integral no longer vanishes. In fact when $\ell =0$
the boundary term in Stokes' theorem is%
\begin{eqnarray*}
\mathcal{I}_{J}^{0} &=&c\int_{\partial \mathbb{B}_{n}}\frac{\left( 1-\zeta 
\overline{z}\right) ^{n}}{\bigtriangleup \left( \zeta ,z\right) ^{n}}%
\overline{\mathcal{D}^{0}}\left( \eta \lrcorner d\overline{w}^{J}\right)
d\sigma \left( \zeta \right) \\
&=&c\int_{\partial \mathbb{B}_{n}}\frac{1}{\left( 1-\overline{\zeta }%
z\right) ^{n}}\overline{\mathcal{D}^{0}}\left( \eta \lrcorner d\overline{w}%
^{J}\right) d\sigma \left( \zeta \right) ,
\end{eqnarray*}%
since from (\ref{manyfaces}) we have%
\begin{equation*}
\frac{\left( 1-w\overline{z}\right) ^{n}}{\bigtriangleup \left( z,w\right)
^{n}}=\frac{\left( 1-w\overline{z}\right) ^{n}}{\left\vert 1-w\overline{z}%
\right\vert ^{2n}\left\vert \varphi _{z}\left( w\right) \right\vert ^{2n}}=%
\frac{1}{\left( 1-\overline{w}z\right) ^{n}},\ \ \ \ \ w\in \partial \mathbb{%
B}_{n}.
\end{equation*}%
Thus the boundary term in (\ref{onceiter}) is 
\begin{equation*}
c\sum_{J}\int_{\partial \mathbb{B}_{n}}\frac{1}{\left( 1-\overline{\zeta }%
z\right) ^{n}}\overline{\mathcal{D}^{0}}\left( \eta \lrcorner d\overline{w}%
^{J}\right) d\sigma \left( \zeta \right) d\overline{z}^{J}=c\mathcal{S}%
_{n}\left( \overline{\mathcal{D}^{0}}\eta \right) \left( z\right) .
\end{equation*}%
This completes the proof of (\ref{mformula'}) in the case $m=1$. Now we
proceed by induction on $m$ to complete the proof of Lemma \ref{IBP1'}.

\bigskip

Finally here is the simple proof of the integration by parts formula for the
radial derivative in Lemma \ref{IBP2}.

\textbf{Proof\ of Lemma \ref{IBP2}}: Since $\left( 1-\left\vert w\right\vert
^{2}\right) ^{b+1}$ vanishes on the boundary for $b>-1$, and since%
\begin{equation*}
R\left( 1-\left\vert w\right\vert ^{2}\right) ^{b+1}=\sum_{j=1}^{n}w_{j}%
\frac{\partial }{\partial w_{j}}\left( 1-\left\vert w\right\vert ^{2}\right)
^{b+1}=-\left( b+1\right) \left( 1-\left\vert w\right\vert ^{2}\right)
^{b}\left\vert w\right\vert ^{2},
\end{equation*}%
the divergence theorem yields 
\begin{eqnarray*}
0 &=&\int_{\partial \mathbb{B}_{n}}\left( 1-\left\vert w\right\vert
^{2}\right) ^{b+1}\Psi \left( w\right) w\cdot \mathbf{n}d\sigma \left(
w\right) \\
&=&\int_{\mathbb{B}_{n}}\sum_{j=1}^{n}\frac{\partial }{\partial w_{j}}%
\left\{ w_{j}\left( 1-\left\vert w\right\vert ^{2}\right) ^{b+1}\Psi \left(
w\right) \right\} dV\left( w\right) \\
&=&n\int_{\mathbb{B}_{n}}\left( 1-\left\vert w\right\vert ^{2}\right)
^{b+1}\Psi \left( w\right) dV\left( w\right) \\
&&+\left( b+1\right) \int_{\mathbb{B}_{n}}\left( 1-\left\vert w\right\vert
^{2}\right) ^{b}\left( -\left\vert w\right\vert ^{2}\right) \Psi \left(
w\right) dV\left( w\right) \\
&&+\int_{\mathbb{B}_{n}}\left( 1-\left\vert w\right\vert ^{2}\right)
^{b+1}R\Psi \left( w\right) dV\left( w\right) ,
\end{eqnarray*}%
which after rearranging becomes%
\begin{eqnarray*}
&&\left( n+b+1\right) \int_{\mathbb{B}_{n}}\left( 1-\left\vert w\right\vert
^{2}\right) ^{b+1}\Psi \left( w\right) dV\left( w\right) \\
&&+\int_{\mathbb{B}_{n}}\left( 1-\left\vert w\right\vert ^{2}\right)
^{b+1}R\Psi \left( w\right) dV\left( w\right) . \\
&=&\left( b+1\right) \int_{\mathbb{B}_{n}}\left( 1-\left\vert w\right\vert
^{2}\right) ^{b}\Psi \left( w\right) dV\left( w\right) .
\end{eqnarray*}

\subsection{Equivalent seminorms on Besov-Sobolev spaces\label{Equivalent
seminorms}}

It is a routine matter to take known scalar-valued proofs of the results in
this section and replace the scalars with vectors in $\ell ^{2}$ to obtain
proofs for the $\ell ^{2}$-valued versions. We begin illustrating this
process by proving the equivalence of norms in Proposition \ref{Bequiv}.

\textbf{Proof of Proposition \ref{Bequiv}}: First we note the equivalence of
the following two conditions (the case $\sigma =0$ is Theorem 6.1 of \cite%
{Zhu}):

\begin{enumerate}
\item The functions%
\begin{equation*}
\left( 1-\left\vert z\right\vert ^{2}\right) ^{\left\vert k\right\vert
+\sigma }\frac{\partial ^{\left\vert k\right\vert }}{\partial z^{k}}f\left(
z\right) ,\ \ \ \ \ \left\vert k\right\vert =N
\end{equation*}%
are in $L^{p}\left( d\lambda _{n};\ell ^{2}\right) $ for some $N>\frac{n}{p}%
-\sigma $,

\item The functions%
\begin{equation*}
\left( 1-\left\vert z\right\vert ^{2}\right) ^{\left\vert k\right\vert
+\sigma }\frac{\partial ^{\left\vert k\right\vert }}{\partial z^{k}}f\left(
z\right) ,\ \ \ \ \ \left\vert k\right\vert =N
\end{equation*}%
are in $L^{p}\left( d\lambda _{n};\ell ^{2}\right) $ for every $N>\frac{n}{p}%
-\sigma $.
\end{enumerate}

Indeed, $L^{p}\left( d\lambda _{n};\ell ^{2}\right) =L^{p}\left( \nu
_{-n-1};\ell ^{2}\right) $ and $\left( 1-\left\vert z\right\vert ^{2}\right)
^{\left\vert k\right\vert +\sigma }\frac{\partial ^{\left\vert k\right\vert }%
}{\partial z^{k}}f\left( z\right) \in L^{p}\left( \nu _{-n-1};\ell
^{2}\right) $ if and only if $\frac{\partial ^{\left\vert k\right\vert }}{%
\partial z^{k}}f\left( z\right) \in L^{p}\left( \nu _{p\left( \left\vert
k\right\vert +\sigma \right) -n-1};\ell ^{2}\right) $. Provided $p\left(
\left\vert k\right\vert +\sigma \right) -n-1>-1$, Theorem 2.17 of \cite{Zhu}
shows that $\left( 1-\left\vert z\right\vert ^{2}\right) ^{\ell }\frac{%
\partial ^{\left\vert \ell \right\vert }}{\partial z^{\ell }}\left( \frac{%
\partial ^{\left\vert k\right\vert }}{\partial z^{k}}f\right) \left(
z\right) \in L^{p}\left( \nu _{p\left( \left\vert k\right\vert +\sigma
\right) -n-1};\ell ^{2}\right) $, which shows that (2) follows from (1).

From the equivalence of (1) and (2) we obtain the equivalence of the first
two conditions in Proposition \ref{Bequiv}. The equivalence with the next
two conditions follows from the corresponding generalization to $\sigma >0$
of Theorem 6.4 in \cite{Zhu}, which in turn is achieved by arguing as in the
previous paragraph.

\bigskip

Next we prove Lemma \ref{radinv} by adapting the proof of Lemma 6.4 in \cite%
{ArRoSa}.

\textbf{Proof of Lemma \ref{radinv}}: We have 
\begin{equation}
\left\vert D_{a}f\left( z\right) \right\vert =\left\vert f^{\prime }\left(
z\right) \left\{ \left( 1-\left\vert a\right\vert ^{2}\right) P_{a}+\left(
1-\left\vert a\right\vert ^{2}\right) ^{\frac{1}{2}}Q_{a}\right\}
\right\vert \geq \left\vert \left( 1-\left\vert a\right\vert ^{2}\right)
f^{\prime }\left( z\right) \right\vert ,  \label{fin}
\end{equation}%
and iterating with $f$ replaced by (the components of) $D_{a}f$ in (\ref{fin}%
), we obtain 
\begin{equation*}
\left\vert D_{a}^{2}f\left( z\right) \right\vert \geq \left\vert \left(
1-\left\vert a\right\vert ^{2}\right) \left( D_{a}f\right) ^{\prime }\left(
z\right) \right\vert .
\end{equation*}%
Applying (\ref{fin}) once more with $f$ replaced by (the components of) $%
f^{\prime }$, we get 
\begin{equation*}
\left\vert \left( 1-\left\vert a\right\vert ^{2}\right) \left( D_{a}f\right)
^{\prime }\left( z\right) \right\vert =\left\vert \left( 1-\left\vert
a\right\vert ^{2}\right) D_{a}\left( f^{\prime }\right) \left( z\right)
\right\vert \geq \left\vert \left( 1-\left\vert a\right\vert ^{2}\right)
^{2}f^{\prime \prime }\left( z\right) \right\vert ,
\end{equation*}%
which when combined with the previous inequality yields 
\begin{equation*}
\left\vert D_{a}^{2}f\left( z\right) \right\vert \geq \left\vert \left(
1-\left\vert a\right\vert ^{2}\right) ^{2}f^{\prime \prime }\left( z\right)
\right\vert .
\end{equation*}%
Continuing by induction we have 
\begin{equation}
\left\vert D_{a}^{m}f\left( z\right) \right\vert \geq \left\vert \left(
1-\left\vert a\right\vert ^{2}\right) ^{m}f^{\left( m\right) }\left(
z\right) \right\vert ,\;\;\;\;\;m\geq 1.  \label{allm}
\end{equation}%
Proposition \ref{Bequiv} and (\ref{allm}) now show that 
\begin{eqnarray*}
&&\left( \int_{\mathbb{B}_{n}}\left\vert \left( 1-\left\vert z\right\vert
^{2}\right) ^{m+\sigma }R^{0,m}f\left( z\right) \right\vert ^{p}d\lambda
_{n}\left( z\right) \right) ^{\frac{1}{p}} \\
&\leq &C\left( \int_{\mathbb{B}_{n}}\left\vert \left( 1-\left\vert
z\right\vert ^{2}\right) ^{m+\sigma }f^{\left( m\right) }\left( z\right)
\right\vert ^{p}d\lambda _{n}\left( z\right) \right) ^{\frac{1}{p}%
}+\sum_{j=0}^{m-1}\left\vert \nabla ^{j}f\left( 0\right) \right\vert \\
&\leq &C\left( \sum_{\alpha \in \mathcal{T}_{n}}\int_{B_{\beta }\left(
c_{\alpha },C_{2}\right) }\left\vert \left( 1-\left\vert z\right\vert
^{2}\right) ^{m+\sigma }f^{\left( m\right) }\left( z\right) \right\vert
^{p}d\lambda _{n}\left( z\right) \right) ^{\frac{1}{p}}+\sum_{j=0}^{m-1}%
\left\vert \nabla ^{j}f\left( 0\right) \right\vert \\
&\leq &C\left( \sum_{\alpha \in \mathcal{T}_{n}}\int_{B_{\beta }\left(
c_{\alpha },C_{2}\right) }\left\vert \left( 1-\left\vert c_{\alpha
}\right\vert ^{2}\right) ^{m+\sigma }f^{\left( m\right) }\left( z\right)
\right\vert ^{p}d\lambda _{n}\left( z\right) \right) ^{\frac{1}{p}%
}+\sum_{j=0}^{m-1}\left\vert \nabla ^{j}f\left( 0\right) \right\vert \\
&\leq &C\left( \sum_{\alpha \in \mathcal{T}_{n}}\int_{B_{\beta }\left(
c_{\alpha },C_{2}\right) }\left\vert \left( 1-\left\vert z\right\vert
^{2}\right) ^{\sigma }D_{c_{\alpha }}^{m}f\left( z\right) \right\vert
^{p}d\lambda _{n}\left( z\right) \right) ^{\frac{1}{p}}+\sum_{j=0}^{m-1}%
\left\vert \nabla ^{j}f\left( 0\right) \right\vert \\
&=&C\left\Vert f\right\Vert _{B_{p,m}^{\sigma }\left( \mathbb{B}_{n}\right)
}^{\ast }+\sum_{j=0}^{m-1}\left\vert \nabla ^{j}f\left( 0\right) \right\vert
.
\end{eqnarray*}

For the opposite inequality, just as in \cite{ArRoSa}, we employ some of the
ideas in the proofs of Theorem 6.11 and Lemma 3.3 in \cite{Zhu}, where the
case $\sigma =0$ and $m=1>\frac{2n}{p}$ is proved. Suppose $f\in H\left( 
\mathbb{B}_{n}\right) $ and that the right side of (\ref{norms}) is finite.
By Proposition \ref{Bequiv} and the proof of Theorem 6.7 of \cite{Zhu} we
have 
\begin{equation}
f\left( z\right) =c\int_{\mathbb{B}_{n}}\frac{g\left( w\right) }{\left( 1-%
\overline{w}z\right) ^{n+1+\sigma }}dV\left( w\right) ,\;\;\;\;\;z\in 
\mathbb{B}_{n},  \label{repproj}
\end{equation}%
for some $g\in L^{p}\left( \lambda _{n}\right) $ where 
\begin{equation}
\left\Vert g\right\Vert _{L^{p}\left( \lambda _{n}\right) }\approx
\sum_{j=0}^{m-1}\left\vert \nabla ^{j}f\left( 0\right) \right\vert +\left(
\int_{\mathbb{B}_{n}}\left\vert \left( 1-\left\vert z\right\vert ^{2}\right)
^{m+\sigma }R^{\sigma ,m}f\left( z\right) \right\vert ^{p}d\lambda
_{n}\left( z\right) \right) ^{\frac{1}{p}}.  \label{gnorm}
\end{equation}%
Indeed, Proposition \ref{Bequiv} shows that 
\begin{eqnarray*}
f &\in &B_{p}^{\sigma }\left( \mathbb{B}_{n}\right) \Leftrightarrow \left(
1-\left\vert z\right\vert ^{2}\right) ^{m+\sigma }R^{\sigma ,m}f\left(
z\right) \in L^{p}\left( \lambda _{n}\right) \\
&\Leftrightarrow &R^{\sigma ,m}f\left( z\right) \in L^{p}\left( \nu
_{p\left( m+\sigma \right) -n-1}\right) \cap H\left( \mathbb{B}_{n}\right) ,
\end{eqnarray*}%
where as in \cite{Zhu} we write $d\nu _{\alpha }\left( z\right) =\left(
1-\left\vert z\right\vert ^{2}\right) ^{\alpha }dV\left( z\right) $. Now
Lemma \ref{Zlemma} above (see also Proposition 2.11 in \cite{Zhu}) shows that%
\begin{equation*}
T_{0,\beta ,0}L^{p}\left( \nu _{\gamma }\right) =L^{p}\left( \nu _{\gamma
}\right) \cap H\left( \mathbb{B}_{n}\right)
\end{equation*}%
if and only if $p\left( \beta +1\right) >\gamma +1$. Choosing $\beta
=m+\sigma $ and $\gamma =p\left( m+\sigma \right) -n-1$ we see that $p\left(
\beta +1\right) >\gamma +1$ and so $f\in B_{p}^{\sigma }\left( \mathbb{B}%
_{n}\right) $ if and only if%
\begin{equation*}
R^{\sigma ,m}f\left( z\right) =c\int_{\mathbb{B}_{n}}\frac{\left(
1-\left\vert w\right\vert ^{2}\right) ^{m+\sigma }h\left( w\right) }{\left(
1-\overline{w}z\right) ^{n+1+m+\sigma }}dV\left( w\right)
\end{equation*}%
for some $h\in L^{p}\left( \nu _{p\left( m+\sigma \right) -n-1}\right) $. If
we set $g\left( w\right) =\left( 1-\left\vert w\right\vert ^{2}\right)
^{m+\sigma }h\left( w\right) $ we obtain%
\begin{equation}
R^{\sigma ,m}f\left( z\right) =c\int_{\mathbb{B}_{n}}\frac{g\left( w\right) 
}{\left( 1-\overline{w}z\right) ^{n+1+m+\sigma }}dV\left( w\right)
\label{Rf}
\end{equation}%
with $g\in L^{p}\left( \lambda _{n}\right) $. Now apply the inverse operator 
$R_{\sigma ,m}=\left( R^{\sigma ,m}\right) ^{-1}$ to both sides of (\ref{Rf}%
) and use (\ref{Zhuidentity}),%
\begin{equation*}
R_{\sigma ,m}\left( \frac{1}{\left( 1-\overline{w}z\right) ^{n+1+m+\sigma }}%
\right) =\frac{1}{\left( 1-\overline{w}z\right) ^{n+1+\sigma }},
\end{equation*}%
to obtain (\ref{repproj}) and (\ref{gnorm}).

Fix $\alpha \in \mathcal{T}_{n}$ and let $a=c_{\alpha }\in \mathbb{B}_{n}$.
We claim that 
\begin{equation}
\left\vert D_{a}^{m}f\left( z\right) \right\vert \leq C_{m}\left(
1-\left\vert a\right\vert ^{2}\right) ^{\frac{m}{2}}\int_{\mathbb{B}_{n}}%
\frac{\left\vert g\left( w\right) \right\vert }{\left\vert 1-\overline{w}%
z\right\vert ^{n+1+\frac{m}{2}+\sigma }}dV\left( w\right) ,\;\;\;\;\;m\geq
1,z\in B_{\beta }\left( a,C\right) .  \label{mclaim}
\end{equation}%
To see this we compute $D_{a}^{m}f\left( z\right) $ for $z\in B_{\beta
}\left( a,C\right) $, beginning with the case $m=1$. Since 
\begin{eqnarray*}
D_{a}\left( \overline{w}z\right) &=&\left( \overline{w}z\right) ^{\prime
}\varphi _{a}^{\prime }\left( 0\right) =-\overline{w}^{t}\left\{ \left(
1-\left\vert a\right\vert ^{2}\right) P_{a}+\left( 1-\left\vert a\right\vert
^{2}\right) ^{\frac{1}{2}}Q_{a}\right\} \\
&=&-\overline{\left\{ \left( 1-\left\vert a\right\vert ^{2}\right)
P_{a}w+\left( 1-\left\vert a\right\vert ^{2}\right) ^{\frac{1}{2}%
}Q_{a}w\right\} }^{t},
\end{eqnarray*}%
we have 
\begin{eqnarray}
&&D_{a}f\left( z\right)  \label{mequ} \\
&=&c_{n}\int_{\mathbb{B}_{n}}D_{a}\left( 1-\overline{w}z\right) ^{-\left(
n+1+\sigma \right) }g\left( w\right) dV\left( w\right)  \notag \\
&=&c_{n}\int_{\mathbb{B}_{n}}\left( 1-\overline{w}z\right) ^{-\left(
n+2+\sigma \right) }D_{a}\left( \overline{w}z\right) g\left( w\right)
dV\left( w\right)  \notag \\
&=&c_{n}\int_{\mathbb{B}_{n}}\left( 1-\overline{w}z\right) ^{-\left(
n+2+\sigma \right) }\overline{\left\{ \left( 1-\left\vert a\right\vert
^{2}\right) P_{a}w+\left( 1-\left\vert a\right\vert ^{2}\right) ^{\frac{1}{2}%
}Q_{a}w\right\} }^{t}g\left( w\right) dV\left( w\right) .  \notag
\end{eqnarray}%
Taking absolute values inside, we obtain 
\begin{equation}
\left\vert D_{a}f\left( z\right) \right\vert \leq C\left( 1-\left\vert
a\right\vert ^{2}\right) ^{\frac{1}{2}}\int_{\mathbb{B}_{n}}\frac{\left(
1-\left\vert a\right\vert ^{2}\right) ^{\frac{1}{2}}\left\vert
P_{a}w\right\vert +\left\vert Q_{a}w\right\vert }{\left\vert 1-\overline{w}%
z\right\vert ^{n+2+\sigma }}\left\vert g\left( w\right) \right\vert dV\left(
w\right) .  \label{daf1}
\end{equation}%
From the following elementary inequalities 
\begin{eqnarray}
\left\vert Q_{a}w\right\vert ^{2} &=&\left\vert Q_{a}\left( w-a\right)
\right\vert ^{2}\leq \left\vert w-a\right\vert ^{2},  \label{elineq} \\
&=&\left\vert w\right\vert ^{2}+\left\vert a\right\vert ^{2}-2\func{Re}%
\left( w\overline{a}\right)  \notag \\
&\leq &2\func{Re}\left( 1-w\overline{a}\right) \leq 2\left\vert 1-w\overline{%
a}\right\vert ,  \notag
\end{eqnarray}%
we obtain that $\left\vert Q_{a}w\right\vert \leq C\left\vert 1-\overline{w}%
a\right\vert ^{\frac{1}{2}}$. Now 
\begin{equation*}
\left\vert 1-w\overline{a}\right\vert \approx \left\vert 1-w\overline{z}%
\right\vert \geq \frac{1}{2}\left( 1-\left\vert z\right\vert ^{2}\right)
\approx \left( 1-\left\vert a\right\vert ^{2}\right) ,\ \ \ \ \ z\in
B_{\beta }\left( a,C\right)
\end{equation*}%
shows that%
\begin{equation*}
\left( 1-\left\vert a\right\vert ^{2}\right) ^{\frac{1}{2}}+\left\vert 1-%
\overline{w}a\right\vert ^{\frac{1}{2}}\leq C\left\vert 1-\overline{w}%
z\right\vert ^{\frac{1}{2}},\ \ \ \ \ z\in B_{\beta }\left( a,C\right) ,
\end{equation*}%
and so we see that 
\begin{equation*}
\frac{\left( 1-\left\vert a\right\vert ^{2}\right) ^{\frac{1}{2}}\left\vert
P_{a}w\right\vert +\left\vert Q_{a}w\right\vert }{\left\vert 1-\overline{w}%
z\right\vert ^{n+2}}\leq \frac{C}{\left\vert 1-\overline{w}z\right\vert ^{n+%
\frac{3}{2}}},\;\;\;\;\;z\in B_{\beta }\left( a,C\right) .
\end{equation*}%
Plugging this estimate into (\ref{daf1}) yields 
\begin{equation*}
\left\vert D_{a}f\left( z\right) \right\vert \leq C\left( 1-\left\vert
a\right\vert ^{2}\right) ^{\frac{1}{2}}\int_{\mathbb{B}_{n}}\frac{\left\vert
g\left( w\right) \right\vert }{\left\vert 1-\overline{w}z\right\vert ^{n+%
\frac{3}{2}+\sigma }}dV\left( w\right) ,
\end{equation*}%
which is the case $m=1$ of (\ref{mclaim}).

To obtain the case $m=2$ of (\ref{mclaim}), we differentiate (\ref{mequ})
again to get 
\begin{equation*}
D_{a}^{2}f\left( z\right) =c\int_{\mathbb{B}_{n}}\left( 1-\overline{w}%
z\right) ^{-\left( n+3+\sigma \right) }W\overline{W}^{t}g\left( w\right)
dV\left( w\right) .
\end{equation*}%
where we have written $W=\left\{ \left( 1-\left\vert a\right\vert
^{2}\right) P_{a}w+\left( 1-\left\vert a\right\vert ^{2}\right) ^{\frac{1}{2}%
}Q_{a}w\right\} $ for convenience. Again taking absolute values inside, we
obtain 
\begin{equation*}
\left\vert D_{a}^{2}f\left( z\right) \right\vert \leq C\left( 1-\left\vert
a\right\vert ^{2}\right) \int_{\mathbb{B}_{n}}\frac{\left( \left(
1-\left\vert a\right\vert ^{2}\right) ^{\frac{1}{2}}\left\vert
P_{a}w\right\vert +\left\vert Q_{a}w\right\vert \right) ^{2}}{\left\vert 1-%
\overline{w}z\right\vert ^{n+3+\sigma }}\left\vert g\left( w\right)
\right\vert dV\left( w\right) .
\end{equation*}%
Once again, using $\left\vert Q_{a}w\right\vert \leq C\left\vert 1-\overline{%
w}a\right\vert ^{\frac{1}{2}}$ and $\left( 1-\left\vert a\right\vert
^{2}\right) ^{\frac{1}{2}}+\left\vert 1-\overline{w}a\right\vert ^{\frac{1}{2%
}}\leq C\left\vert 1-\overline{w}z\right\vert ^{\frac{1}{2}}$ for $z\in
B_{\beta }\left( a,C\right) $, we see that 
\begin{equation*}
\frac{\left( \left( 1-\left\vert a\right\vert ^{2}\right) ^{\frac{1}{2}%
}\left\vert P_{a}w\right\vert +\left\vert Q_{a}w\right\vert \right) ^{2}}{%
\left\vert 1-\overline{w}z\right\vert ^{n+3+\sigma }}\leq \frac{C}{%
\left\vert 1-\overline{w}z\right\vert ^{n+2+\sigma }},\;\;\;\;\;z\in
B_{\beta }\left( a,C\right) ,
\end{equation*}%
which yields the case $m=2$ of (\ref{mclaim}). The general case of (\ref%
{mclaim}) follows by induction on $m$.

The inequality (\ref{mclaim}) shows that $\left( 1-\left\vert z\right\vert
^{2}\right) ^{\sigma }\left\vert D_{c_{\alpha }}^{m}f\left( z\right)
\right\vert \leq C_{m}S\left\vert g\right\vert \left( z\right) $ for $z\in
B_{\beta }\left( c_{\alpha },C\right) $, where%
\begin{equation*}
Sg\left( z\right) =\int_{\mathbb{B}_{n}}\frac{\left( 1-\left\vert
z\right\vert ^{2}\right) ^{\frac{m}{2}+\sigma }}{\left\vert 1-\overline{w}%
z\right\vert ^{n+1+\frac{m}{2}+\sigma }}g\left( w\right) dV\left( w\right) .
\end{equation*}%
We will now use the symbol $a$ differently than before. The operator $S$ is
the operator $T_{a,b,c}$ in Lemma \ref{Zlemma} above (see also Theorem 2.10
of \cite{Zhu}) with parameters $a=\frac{m}{2}+\sigma $ and $b=c=0$. Now with 
$t=-n-1$, our assumption that $m>2\left( \frac{n}{p}-\sigma \right) $ yields 
$-p\left( \frac{m}{2}+\sigma \right) <-n<p\left( 0+1\right) $, i.e. 
\begin{equation*}
-pa<t+1<p\left( b+1\right) .
\end{equation*}%
Thus the bounded overlap property of the balls $B_{\beta }\left( c_{\alpha
},C_{2}\right) $ together with Lemma \ref{Zlemma} above yields 
\begin{eqnarray*}
\left\Vert f\right\Vert _{B_{p,m}^{\sigma }\left( \mathbb{B}_{n}\right)
}^{\ast } &=&\left( \sum_{\alpha \in \mathcal{T}_{n}}\int_{B_{\beta }\left(
c_{\alpha },C_{2}\right) }\left\vert \left( 1-\left\vert z\right\vert
^{2}\right) ^{\sigma }D_{c_{\alpha }}^{m}f\left( z\right) \right\vert
^{p}d\lambda _{n}\left( z\right) \right) ^{\frac{1}{p}} \\
&\leq &C_{m}\left( \int_{\mathbb{B}_{n}}\left\vert Sg\left( z\right)
\right\vert ^{p}d\lambda _{n}\left( z\right) \right) ^{\frac{1}{p}} \\
&\leq &C_{m}^{\prime }\left( \int_{\mathbb{B}_{n}}\left\vert g\left(
z\right) \right\vert ^{p}d\lambda _{n}\left( z\right) \right) ^{\frac{1}{p}}
\\
&\leq &C_{m}^{\prime \prime }\left( \int_{\mathbb{B}_{n}}\left\vert \left(
1-\left\vert z\right\vert ^{2}\right) ^{m+\sigma }R^{\sigma ,m}f\left(
z\right) \right\vert ^{p}d\lambda _{n}\left( z\right) \right) ^{\frac{1}{p}}
\end{eqnarray*}%
by (\ref{gnorm}). This completes the proof of Lemma \ref{radinv}.

\subsubsection{Multilinear inequalities\label{Multilinear inequalities}}

Proposition \ref{multilinear} is proved by adapting the proof of Theorem 3.5
in Ortega and Fabrega \cite{OrFa} to $\ell ^{2}$-valued functions. This
argument uses the complex interpolation theorem of Beatrous \cite{Bea} and
Ligocka \cite{Lig}, which extends to Hilbert space valued functions with the
same proof.\ In order to apply this extension we will need the following
operator norm inequality.

If $\varphi \in M_{B_{p}^{\sigma }\left( \mathbb{B}_{n}\right) \rightarrow
B_{p}^{\sigma }\left( \mathbb{B}_{n};\ell ^{2}\right) }$ and $%
f=\sum_{\left\vert I\right\vert =\kappa }f_{I}e_{I}\in B_{p}^{\sigma }\left( 
\mathbb{B}_{n};\otimes ^{\kappa -1}\ell ^{2}\right) $, we define%
\begin{equation*}
\mathbb{M}_{\varphi }f=\varphi \otimes f=\varphi \otimes \left(
\sum_{\left\vert I\right\vert =\kappa -1}f_{I}e_{I}\right) =\sum_{\left\vert
I\right\vert =\kappa -1}\left( \varphi f_{I}\right) \otimes e_{I},
\end{equation*}%
where $I=\left( i_{1},...,i_{\kappa -1}\right) \in \mathbb{N}^{\kappa -1}$
and $e_{I}=e_{i_{1}}\otimes ...\otimes e_{i_{\kappa -1}}$.

\begin{lemma}
\label{tensornorms}Suppose that $\sigma \geq 0$, $1<p<\infty $ and $\kappa
\geq 1$. Then there is a constant $C_{n,\sigma ,p,\kappa }$ such that 
\begin{equation}
\left\Vert \mathbb{M}_{g}\right\Vert _{B_{p}^{\sigma }\left( \mathbb{B}%
_{n};\otimes ^{\kappa -1}\ell ^{2}\right) \rightarrow B_{p}^{\sigma }\left( 
\mathbb{B}_{n};\otimes ^{\kappa }\ell ^{2}\right) }\leq C_{n,\sigma
,p,\kappa }\left\Vert \mathbb{M}_{g}\right\Vert _{B_{p}^{\sigma }\left( 
\mathbb{B}_{n};\ell ^{2}\right) \rightarrow B_{p}^{\sigma }\left( \mathbb{B}%
_{n};\ell ^{2}\right) }.  \label{ktensor}
\end{equation}%
In the case $p=2$ we have equality: 
\begin{equation}
\left\Vert \mathbb{M}_{\varphi }\right\Vert _{B_{2}^{\sigma }\left( \mathbb{B%
}_{n};\otimes ^{\kappa -1}\ell ^{2}\right) \rightarrow B_{2}^{\sigma }\left( 
\mathbb{B}_{n};\otimes ^{\kappa }\ell ^{2}\right) }=\left\Vert \mathbb{M}%
_{\varphi }\right\Vert _{B_{2}^{\sigma }\left( \mathbb{B}_{n}\right)
\rightarrow B_{2}^{\sigma }\left( \mathbb{B}_{n};\ell ^{2}\right) }.
\label{tensor equality}
\end{equation}
\end{lemma}

It turns out that in order to prove (\ref{ktensor}) for $p\neq 2$ we will
need the case $M=1$ of Proposition \ref{multilinear}. Fortunately, the case $%
M=1$ does not require inequality (\ref{ktensor}), thus avoiding circularity.

\bigskip

\textbf{Proof of Proposition \ref{multilinear} and Lemma \ref{tensornorms}}:
We begin with the proof of the case $M=1$ of Proposition \ref{multilinear}.
We will show that for $m=\ell +k$, 
\begin{equation}
\int_{\mathbb{B}_{n}}\left\vert \left( 1-\left\vert z\right\vert ^{2}\right)
^{\sigma }\left( \mathcal{Y}^{\ell }g\right) \left( \mathcal{Y}^{k}h\right)
\right\vert ^{p}d\lambda _{n}\left( z\right) \leq C_{n,\sigma ,p}\left\Vert 
\mathbb{M}_{g}\right\Vert _{B_{p}^{\sigma }\left( \mathbb{B}_{n}\right)
\rightarrow B_{p}^{\sigma }\left( \mathbb{B}_{n};\ell ^{2}\right)
}^{p}\left\Vert h\right\Vert _{B_{p}^{\sigma }\left( \mathbb{B}_{n}\right)
}^{p}.  \label{N=1}
\end{equation}%
Following the proof of Theorem 3.1 in \cite{OrFa} we first convert the
Leibniz formula%
\begin{equation*}
\left( \mathcal{Y}^{\ell }g\right) \left( \mathcal{Y}^{k}h\right) =\mathcal{Y%
}^{\ell }\left( g\mathcal{Y}^{k}h\right) -\sum_{\alpha =0}^{\ell -1}\left( 
\begin{array}{c}
\ell \\ 
\alpha%
\end{array}%
\right) \left( \mathcal{Y}^{\alpha }g\right) \left( \mathcal{Y}^{k+\ell
-\alpha }h\right)
\end{equation*}%
to "divergence form"%
\begin{equation*}
\left( \mathcal{Y}^{\ell }g\right) \left( \mathcal{Y}^{k}h\right)
=\sum_{\alpha =0}^{\ell }\left( -1\right) ^{\alpha }\left( 
\begin{array}{c}
\ell \\ 
\ell -\alpha%
\end{array}%
\right) \mathcal{Y}^{\ell -\alpha }\left( g\mathcal{Y}^{k+\alpha }h\right) .
\end{equation*}%
This is easily established by induction on $\ell $ with $k$ held fixed and
can be stated as%
\begin{equation}
\left( \mathcal{Y}^{\ell }g\right) \left( \mathcal{Y}^{k}h\right)
=\sum_{\alpha =0}^{\ell }c_{\alpha }^{\ell }\mathcal{Y}^{\alpha }\left( g%
\mathcal{Y}^{k+\ell -\alpha }h\right) .  \label{revLieb}
\end{equation}

Next we note that for $s>\frac{n}{p}$, $B_{p}^{s}\left( \mathbb{B}_{n};\ell
^{2}\right) $ is a Bergman space, hence $M_{B_{p}^{s}\left( \mathbb{B}%
_{n}\right) \rightarrow B_{p}^{s}\left( \mathbb{B}_{n};\ell ^{2}\right)
}=H^{\infty }\left( \mathbb{B}_{n};\ell ^{2}\right) $. Thus using (\ref%
{multarebounded}) we have for $s>\frac{n}{p}$, 
\begin{equation*}
g\in M_{B_{p}^{\sigma }\left( \mathbb{B}_{n}\right) \rightarrow
B_{p}^{\sigma }\left( \mathbb{B}_{n};\ell ^{2}\right) }\cap H^{\infty
}\left( \mathbb{B}_{n};\ell ^{2}\right) =M_{B_{p}^{\sigma }\left( \mathbb{B}%
_{n}\right) \rightarrow B_{p}^{\sigma }\left( \mathbb{B}_{n};\ell
^{2}\right) }\cap M_{B_{p}^{s}\left( \mathbb{B}_{n}\right) \rightarrow
B_{p}^{s}\left( \mathbb{B}_{n};\ell ^{2}\right) }.
\end{equation*}%
Then, still following the argument in \cite{OrFa}, we use the complex
interpolation theorem of Beatrous \cite{Bea} and Ligocka \cite{Lig} (they
prove only the scalar-valued version but the Hilbert space valued version
has the same proof),%
\begin{eqnarray*}
\left( B_{p}^{\sigma }\left( \mathbb{B}_{n}\right) ,B_{p}^{\frac{n}{p}%
+\varepsilon }\left( \mathbb{B}_{n}\right) \right) _{\theta }
&=&B_{p}^{\left( 1-\theta \right) \sigma +\theta \left( \frac{n}{p}%
+\varepsilon \right) }\left( \mathbb{B}_{n}\right) ,\ \ \ \ \ 0\leq \theta
\leq 1, \\
\left( B_{p}^{\sigma }\left( \mathbb{B}_{n};\ell ^{2}\right) ,B_{p}^{\frac{n%
}{p}+\varepsilon }\left( \mathbb{B}_{n};\ell ^{2}\right) \right) _{\theta }
&=&B_{p}^{\left( 1-\theta \right) \sigma +\theta \left( \frac{n}{p}%
+\varepsilon \right) }\left( \mathbb{B}_{n};\ell ^{2}\right) ,\ \ \ \ \
0\leq \theta \leq 1,
\end{eqnarray*}%
to conclude that $g\in M_{B_{p}^{s}\left( \mathbb{B}_{n}\right) \rightarrow
B_{p}^{s}\left( \mathbb{B}_{n};\ell ^{2}\right) }$ for all $s\geq \sigma $,
and with multiplier norm $\left\Vert \mathbb{M}_{g}\right\Vert
_{B_{p}^{s}\left( \mathbb{B}_{n}\right) \rightarrow B_{p}^{s}\left( \mathbb{B%
}_{n};\ell ^{2}\right) }$ bounded by $\left\Vert \mathbb{M}_{g}\right\Vert
_{B_{p}^{\sigma }\left( \mathbb{B}_{n}\right) \rightarrow B_{p}^{\sigma
}\left( \mathbb{B}_{n};\ell ^{2}\right) }$. Recall now that%
\begin{equation*}
\left\Vert h\right\Vert _{B_{p}^{\sigma }\left( \mathbb{B}_{n}\right)
}^{p}=\int_{\mathbb{B}_{n}}\left\vert \left( 1-\left\vert z\right\vert
^{2}\right) ^{\sigma }\mathcal{Y}^{m}h\left( z\right) \right\vert
^{p}d\lambda _{n}\left( z\right) ,
\end{equation*}%
and similarly for $\left\Vert f\right\Vert _{B_{p}^{\sigma }\left( \mathbb{B}%
_{n};\ell ^{2}\right) }^{p}$, provided $m$ satisfies%
\begin{equation}
\left( \sigma +\frac{m}{2}\right) p>n,  \label{m}
\end{equation}%
where $\frac{m}{2}$ appears in the inequality since the derivatives $D$ that
can appear in $\mathcal{Y}^{m}$ only contribute $\left( 1-\left\vert
z\right\vert ^{2}\right) ^{\frac{1}{2}}$\ to the power of $1-\left\vert
z\right\vert ^{2}$ in the integral (see Section \ref{An almost invariant}).

\begin{remark}
\label{D}At this point we recall the convention established in Definitions %
\ref{Dpowers} and \ref{calX} that the factors of $1-\left\vert z\right\vert
^{2}$ that are embedded in the notation for the derivative $\mathcal{Y}%
^{\alpha }$ are treated as constants relative to the actual
differentiations. In the calculations below, we will adopt \textbf{the same
convention} for the factors $\left( 1-\left\vert z\right\vert ^{2}\right)
^{s}$ that we introduce into the integrals. Alternatively, the reader may
wish to write out all the derivatives explicitly with the appropriate power
of $1-\left\vert z\right\vert ^{2}$ set aside as is done in \cite{OrFa}.
\end{remark}

So we have, keeping in mind Remark \ref{D},%
\begin{eqnarray*}
&&\int_{\mathbb{B}_{n}}\left\vert \left( 1-\left\vert z\right\vert
^{2}\right) ^{\sigma }\mathcal{Y}^{\alpha }\left( g\left( z\right) \mathcal{Y%
}^{k+\ell -\alpha }h\left( z\right) \right) \right\vert ^{p}d\lambda _{n} \\
&&\ \ \ \ \ \ \ \ \ \ =\int_{\mathbb{B}_{n}}\left\vert \left( 1-\left\vert
z\right\vert ^{2}\right) ^{s}\mathcal{Y}^{\alpha }\left\{ g\left( z\right)
\left( 1-\left\vert z\right\vert ^{2}\right) ^{\sigma -s}\mathcal{Y}^{k+\ell
-\alpha }h\left( z\right) \right\} \right\vert ^{p}d\lambda _{n} \\
&&\ \ \ \ \ \ \ \ \ \ =\left\Vert g\left( z\right) \left( 1-\left\vert
z\right\vert ^{2}\right) ^{\sigma -s}\mathcal{Y}^{k+\ell -\alpha
}h\right\Vert _{B_{p,\alpha }^{s}\left( \mathbb{B}_{n};\ell ^{2}\right)
}^{p}.
\end{eqnarray*}%
Here the function%
\begin{equation*}
H\left( z\right) =\left( 1-\left\vert z\right\vert ^{2}\right) ^{\sigma -s}%
\mathcal{Y}^{k+\ell -\alpha }h\left( z\right)
\end{equation*}%
is \emph{not} holomorphic, but we have defined the norm $\left\Vert \cdot
\right\Vert _{B_{p,\alpha }^{s}\left( \mathbb{B}_{n};\ell ^{2}\right) }$ on
smooth functions anyway. Now we would like to apply a multiplier property of 
$g$, and for this we must be acting on a Besov-Sobolev space of $\emph{%
holomorphic}$ functions, since that is what we get from the complex
interpolation of Beatrous and Ligocka. Precisely, we get that $\mathbb{M}%
_{g} $ is a bounded operator from $B_{p}^{s}\left( \mathbb{B}_{n}\right) $
to $B_{p}^{s}\left( \mathbb{B}_{n};\ell ^{2}\right) $ for \emph{all} $s\geq
\sigma $.

Now we express $\mathcal{Y}^{k+\ell -\alpha }h\left( z\right) $ as a sum of
terms that are products of a power of $1-\left\vert z\right\vert ^{2}$ and a
derivative $R^{i}L^{j}h\left( z\right) $ where $i+j=k+\ell -\alpha $ and $R$
is the radial derivative and $L$ denotes a complex tangential derivative $%
\frac{\partial }{\partial z_{j}}-\overline{z_{j}}R$ as in \cite{OrFa}.
However, the operators $R^{i}L^{j}$ have different weights in the sense that
the power of $1-\left\vert z\right\vert ^{2}$ that is associated with $%
R^{i}L^{j}$ is $\left( 1-\left\vert z\right\vert ^{2}\right) ^{i+\frac{j}{2}%
} $, i.e.%
\begin{equation*}
\mathcal{Y}^{k+\ell -\alpha }h\left( z\right) =\sum \left( 1-\left\vert
z\right\vert ^{2}\right) ^{i+\frac{j}{2}}R^{i}L^{j}h\left( z\right) .
\end{equation*}%
It turns out that to handle the term $\left( 1-\left\vert z\right\vert
^{2}\right) ^{i+\frac{j}{2}}R^{i}L^{j}h\left( z\right) $ we will use that $g$
is a multiplier on $B_{p}^{s}\left( \mathbb{B}_{n}\right) $ with%
\begin{equation*}
s=\sigma +i+\frac{j}{2},
\end{equation*}%
an exponent that depends on $i+\frac{j}{2}$ and \emph{not} on $i+j=k+\ell
-\alpha $.

Indeed, we have using our "convention" that%
\begin{eqnarray*}
&&\left\Vert g\left( z\right) \left( 1-\left\vert z\right\vert ^{2}\right)
^{\sigma -s}\left( 1-\left\vert z\right\vert ^{2}\right) ^{i+\frac{j}{2}%
}R^{i}L^{j}h\left( z\right) \right\Vert _{B_{p,\alpha }^{s}\left( \mathbb{B}%
_{n};\ell ^{2}\right) }^{p} \\
&=&\int_{\mathbb{B}_{n}}\left\vert \left( 1-\left\vert z\right\vert
^{2}\right) ^{s}\mathcal{Y}^{\alpha }\left\{ g\left( z\right) \left(
1-\left\vert z\right\vert ^{2}\right) ^{\sigma -s}\left( 1-\left\vert
z\right\vert ^{2}\right) ^{i+\frac{j}{2}}R^{i}L^{j}h\left( z\right) \right\}
\right\vert ^{p}d\lambda _{n} \\
&=&\int_{\mathbb{B}_{n}}\left\vert \left( 1-\left\vert z\right\vert
^{2}\right) ^{\sigma +i+\frac{j}{2}}\mathcal{Y}^{\alpha }\left\{ g\left(
z\right) R^{i}L^{j}h\left( z\right) \right\} \right\vert ^{p}d\lambda _{n} \\
&=&\left\Vert g\left( z\right) R^{i}L^{j}h\left( z\right) \right\Vert
_{B_{p,\alpha }^{s}\left( \mathbb{B}_{n};\ell ^{2}\right) }^{p}.
\end{eqnarray*}%
Now the function $g\left( z\right) R^{i}L^{j}h\left( z\right) $ is
holomorphic and $s=\sigma +i+\frac{j}{2}\geq \sigma $ so that we can use
that $g$ is a multiplier on $B_{p}^{s}\left( \mathbb{B}_{n}\right)
=B_{p,\alpha }^{s}\left( \mathbb{B}_{n}\right) $ (this latter equality holds
because $\left( s+\frac{\alpha }{2}\right) p>n$ by (\ref{m})). The result is
that%
\begin{eqnarray*}
&&\left\Vert g\left( z\right) R^{i}L^{j}h\left( z\right) \right\Vert
_{B_{p}^{s}\left( \mathbb{B}_{n};\ell ^{2}\right) }^{p} \\
&\leq &\left\Vert \mathbb{M}_{g}\right\Vert _{B_{p}^{s}\left( \mathbb{B}%
_{n}\right) \rightarrow B_{p}^{s}\left( \mathbb{B}_{n};\ell ^{2}\right)
}^{p}\left\Vert R^{i}L^{j}h\left( z\right) \right\Vert _{B_{p,\alpha
}^{s}\left( \mathbb{B}_{n}\right) }^{p} \\
&\leq &\left\Vert \mathbb{M}_{g}\right\Vert _{B_{p}^{s}\left( \mathbb{B}%
_{n}\right) \rightarrow B_{p}^{s}\left( \mathbb{B}_{n};\ell ^{2}\right)
}^{p}\int_{\mathbb{B}_{n}}\left\vert \left( 1-\left\vert z\right\vert
^{2}\right) ^{\sigma +i+\frac{j}{2}}\mathcal{Y}^{\alpha }R^{i}L^{j}h\left(
z\right) \right\vert ^{p}d\lambda _{n} \\
&=&\left\Vert \mathbb{M}_{g}\right\Vert _{B_{p}^{s}\left( \mathbb{B}%
_{n}\right) \rightarrow B_{p}^{s}\left( \mathbb{B}_{n};\ell ^{2}\right)
}^{p}\int_{\mathbb{B}_{n}}\left\vert \left( 1-\left\vert z\right\vert
^{2}\right) ^{\sigma }\mathcal{Y}^{\alpha }\left[ \left( 1-\left\vert
z\right\vert ^{2}\right) R\right] ^{i}\left[ \sqrt{1-\left\vert z\right\vert
^{2}}L\right] ^{j}h\left( z\right) \right\vert ^{p}d\lambda _{n} \\
&\leq &\left\Vert \mathbb{M}_{g}\right\Vert _{B_{p}^{s}\left( \mathbb{B}%
_{n}\right) \rightarrow B_{p}^{s}\left( \mathbb{B}_{n};\ell ^{2}\right)
}^{p}\int_{\mathbb{B}_{n}}\left\vert \left( 1-\left\vert z\right\vert
^{2}\right) ^{\sigma }\mathcal{Y}^{\alpha +i+j}h\left( z\right) \right\vert
^{p}d\lambda _{n} \\
&=&\left\Vert \mathbb{M}_{g}\right\Vert _{B_{p}^{s}\left( \mathbb{B}%
_{n}\right) \rightarrow B_{p}^{s}\left( \mathbb{B}_{n};\ell ^{2}\right)
}^{p}\int_{\mathbb{B}_{n}}\left\vert \left( 1-\left\vert z\right\vert
^{2}\right) ^{\sigma }\mathcal{Y}^{m}h\left( z\right) \right\vert
^{p}d\lambda _{n} \\
&\leq &\left\Vert \mathbb{M}_{g}\right\Vert _{B_{p}^{\sigma }\left( \mathbb{B%
}_{n}\right) \rightarrow B_{p}^{\sigma }\left( \mathbb{B}_{n};\ell
^{2}\right) }^{p}\left\Vert h\right\Vert _{B_{p}^{\sigma }\left( \mathbb{B}%
_{n}\right) }^{p},
\end{eqnarray*}%
and the case $M=1$ of Proposition \ref{multilinear} is proved.

\bigskip

Now we turn to the proof of the operator norm inequality (\ref{ktensor}) in
Lemma \ref{tensornorms}. The case $p=2$ is particularly easy:%
\begin{eqnarray*}
\left\Vert \mathbb{M}_{\varphi }f\right\Vert _{B_{2}^{\sigma }\left( \mathbb{%
B}_{n};\otimes ^{\kappa }\ell ^{2}\right) }^{2} &=&\int_{\mathbb{B}%
_{n}}\left( 1-\left\vert z\right\vert ^{2}\right) ^{2\sigma
}\sum_{\left\vert I\right\vert =\kappa -1}\left\vert \mathcal{Y}^{m}\left(
\varphi f_{I}\right) \right\vert ^{2}d\lambda _{n} \\
&=&\sum_{\left\vert I\right\vert =\kappa -1}\left\Vert \mathbb{M}_{\varphi
}f_{I}\right\Vert _{B_{2}^{\sigma }\left( \mathbb{B}_{n};\ell ^{2}\right)
}^{2} \\
&\leq &\left\Vert \mathbb{M}_{\varphi }\right\Vert _{B_{2}^{\sigma }\left( 
\mathbb{B}_{n}\right) \rightarrow B_{2}^{\sigma }\left( \mathbb{B}_{n};\ell
^{2}\right) }^{2}\sum_{\left\vert I\right\vert =\kappa -1}\left\Vert
f_{I}\right\Vert _{B_{2}^{\sigma }\left( \mathbb{B}_{n}\right) }^{2} \\
&=&\left\Vert \mathbb{M}_{\varphi }\right\Vert _{B_{2}^{\sigma }\left( 
\mathbb{B}_{n}\right) \rightarrow B_{2}^{\sigma }\left( \mathbb{B}_{n};\ell
^{2}\right) }^{2}\int_{\mathbb{B}_{n}}\left( 1-\left\vert z\right\vert
^{2}\right) ^{2\sigma }\sum_{\left\vert I\right\vert =\kappa -1}\left\vert 
\mathcal{Y}^{m}f_{I}\right\vert ^{2}d\lambda _{n} \\
&=&\left\Vert \mathbb{M}_{\varphi }\right\Vert _{B_{2}^{\sigma }\left( 
\mathbb{B}_{n}\right) \rightarrow B_{2}^{\sigma }\left( \mathbb{B}_{n};\ell
^{2}\right) }^{2}\left\Vert f\right\Vert _{B_{2}^{\sigma }\left( \mathbb{B}%
_{n};\otimes ^{\kappa -1}\ell ^{2}\right) }^{2},
\end{eqnarray*}%
and from this we easily obtain (\ref{tensor equality}).

For $p\neq 2$ it suffices to show that%
\begin{equation}
\left\Vert \mathbb{M}_{\varphi }\right\Vert _{B_{p}^{\sigma }\left( \mathbb{B%
}_{n};\mathbb{C}^{\nu }\right) \rightarrow B_{p}^{\sigma }\left( \mathbb{B}%
_{n};\mathbb{C}^{\mu }\otimes \mathbb{C}^{\nu }\right) }\leq C_{n,\sigma
,p}\left\Vert \mathbb{M}_{\varphi }\right\Vert _{B_{p}^{\sigma }\left( 
\mathbb{B}_{n}\right) \rightarrow B_{p}^{\sigma }\left( \mathbb{B}_{n};%
\mathbb{C}^{\mu }\right) }  \label{itsuff}
\end{equation}%
for all $\mu ,\nu \geq 1$ where the constant $C_{n,\sigma ,p}$ is
independent of $\mu ,\nu $. Indeed, both $\ell ^{2}$ and $\otimes ^{\kappa
-1}\ell ^{2}$ are separable Hilbert spaces and so can be appropriately
approximated by $\mathbb{C}^{\mu }$ and $\mathbb{C}^{\nu }$ respectively.
For each $z\in \mathbb{B}_{n}$ we will view $\varphi \left( z\right) \in 
\mathbb{C}^{\mu }$ as a column vector and $f\left( z\right) \in \mathbb{C}%
^{\nu }$ as a row vector so that $\left( \mathbb{M}_{\varphi }f\right)
\left( z\right) $ is the rank one $\mu \times \nu $ matrix%
\begin{equation*}
\left( \mathbb{M}_{\varphi }f\right) \left( z\right) =\left[ 
\begin{array}{ccc}
\left( \varphi _{1}f_{1}\right) \left( z\right) & \cdots & \left( \varphi
_{1}f_{\nu }\right) \left( z\right) \\ 
\vdots & \ddots & \vdots \\ 
\left( \varphi _{\mu }f_{1}\right) \left( z\right) & \cdots & \left( \varphi
_{\mu }f_{\nu }\right) \left( z\right)%
\end{array}%
\right] =\varphi \left( z\right) \odot f\left( z\right) ,
\end{equation*}%
where we have inserted the symbol $\odot $ simply to remind the reader that
this is \emph{not} the dot product $\varphi \left( z\right) \cdot f\left(
z\right) =f\left( z\right) \varphi \left( z\right) $ of the vectors $\varphi
\left( z\right) $ and $f\left( z\right) $.

Now we consider a single component $X^{m}$ of the vector differential
operator $\mathcal{Y}^{m}$ for some $m>2\left( \frac{n}{p}-\sigma \right) $,
which can be chosen independent of $\mu $ and $\nu $. The main point in the
proof of the lemma is that the matrix $X^{m}\left( \mathbb{M}_{\varphi
}f\right) \left( z\right) $ has rank at most $m+1$ independent of $\mu $ and 
$\nu $. Indeed, the Leibniz formula yields%
\begin{equation*}
X^{m}\left( \mathbb{M}_{\varphi }f\right) \left( z\right) =X^{m}\left(
\varphi \left( z\right) \odot f\left( z\right) \right) =\sum_{\ell
=0}^{m}c_{\ell ,m}X^{m-\ell }\varphi \left( z\right) \odot X^{\ell }f\left(
z\right) ,
\end{equation*}%
where each matrix $X^{m-\ell }\varphi \left( z\right) \odot X^{\ell }f\left(
z\right) $ is rank one, and where the Hilbert Schmidt norm is multiplicative:%
\begin{equation*}
\left\vert X^{m-\ell }\varphi \left( z\right) \odot X^{\ell }f\left(
z\right) \right\vert =\left\vert X^{m-\ell }\varphi \left( z\right)
\right\vert \left\vert X^{\ell }f\left( z\right) \right\vert .
\end{equation*}

Momentarily fix $0\leq \ell \leq m$ and define%
\begin{eqnarray*}
T^{\ell }h\left( z\right) &=&X^{m-\ell }\varphi \left( z\right) h\left(
z\right) ,\ \ \ \ \ h\left( z\right) \in \mathbb{C}, \\
T^{\ell }g\left( z\right) &=&X^{m-\ell }\varphi \left( z\right) \odot
g\left( z\right) ,\ \ \ \ \ g\left( z\right) \in \mathbb{C}^{\nu }.
\end{eqnarray*}%
For $x\in \mathbb{\partial B}_{\mu }$, which we view as a row vector,\ define%
\begin{equation*}
T_{x}^{\ell }g\left( z\right) =xT^{\ell }g\left( z\right) =x\left( X^{m-\ell
}\varphi \right) \left( z\right) \odot g\left( z\right) .
\end{equation*}%
Now choose $x\left( z\right) \in \mathbb{\partial B}_{\mu }$\ such that $%
x\left( z\right) \left( X^{m-\ell }\varphi \right) \left( z\right)
=\left\vert X^{m-\ell }\varphi \left( z\right) \right\vert $ so that%
\begin{equation*}
T_{x\left( z\right) }^{\ell }g\left( z\right) =x\left( z\right) \left(
X^{m-\ell }\varphi \right) \left( z\right) \odot g\left( z\right)
=\left\vert X^{m-\ell }\varphi \left( z\right) \right\vert g\left( z\right) ,
\end{equation*}%
and hence%
\begin{equation*}
\left\vert T_{x\left( z\right) }^{\ell }\left( X^{\ell }f\right) \left(
z\right) \right\vert =\left\vert X^{m-\ell }\varphi \left( z\right)
\right\vert \left\vert X^{\ell }f\left( z\right) \right\vert =\left\vert
X^{m-\ell }\varphi \left( z\right) \odot X^{\ell }f\left( z\right)
\right\vert =\left\vert T^{\ell }\left( X^{\ell }f\right) \left( z\right)
\right\vert .
\end{equation*}

Now we follow the well known argument on page 451 of \cite{Ste}. For $y\in
\partial \mathbb{B}_{\nu }$, which we view as a column vector, and $g\left(
z\right) \in \mathbb{C}^{\nu }$ define the scalars%
\begin{eqnarray*}
g_{y}\left( z\right) &=&g\left( z\right) y, \\
\left( T_{x\left( z\right) }^{\ell }g\right) _{y}\left( z\right)
&=&T_{x\left( z\right) }^{\ell }g\left( z\right) y=x\left( z\right) \left(
X^{m-\ell }\varphi \right) \left( z\right) \odot g\left( z\right) y,
\end{eqnarray*}%
and note that%
\begin{equation*}
T_{x\left( z\right) }^{\ell }\left( X^{\ell }f\right) \left( z\right)
y=x\left( z\right) \left( X^{m-\ell }\varphi \right) \left( z\right) \odot
\left( X^{\ell }f\right) \left( z\right) y=T_{x\left( z\right) }^{\ell
}\left( X^{\ell }f\right) _{y}\left( z\right) .
\end{equation*}%
Thus we have with $d\sigma _{\nu }$ surface measure on $\mathbb{\partial B}%
_{\nu }$, 
\begin{equation*}
\int_{\mathbb{\partial B}_{\nu }}\left\vert T_{x\left( z\right) }^{\ell
}\left( X^{\ell }f\right) \left( z\right) y\right\vert ^{p}d\sigma _{\nu
}\left( y\right) =\left\vert T_{x\left( z\right) }^{\ell }\left( X^{\ell
}f\right) \left( z\right) \right\vert ^{p}\int_{\mathbb{\partial B}_{\nu
}}\left\vert \frac{T_{x\left( z\right) }^{\ell }\left( X^{\ell }f\right)
\left( z\right) }{\left\vert T_{x\left( z\right) }^{\ell }\left( X^{\ell
}f\right) \left( z\right) \right\vert }\cdot y\right\vert ^{p}d\sigma _{\nu
}\left( y\right) ,
\end{equation*}%
as well as%
\begin{equation*}
\int_{\mathbb{\partial B}_{\nu }}\left\vert \left( X^{\ell }f\right)
_{y}\left( z\right) \right\vert ^{p}d\sigma _{\nu }\left( y\right)
=\left\vert X^{\ell }f\left( z\right) \right\vert ^{p}\int_{\mathbb{\partial
B}_{\nu }}\left\vert \frac{X^{\ell }f\left( z\right) }{\left\vert X^{\ell
}f\left( z\right) \right\vert }\cdot y\right\vert ^{p}d\sigma _{\nu }\left(
y\right) .
\end{equation*}

The crucial observation now is that 
\begin{equation*}
\int_{\mathbb{\partial B}_{\nu }}\left\vert \frac{T_{x\left( z\right)
}^{\ell }\left( X^{\ell }f\right) \left( z\right) }{\left\vert T_{x\left(
z\right) }^{\ell }\left( X^{\ell }f\right) \left( z\right) \right\vert }%
\cdot y\right\vert ^{p}d\sigma _{\nu }\left( y\right) =\int_{\mathbb{%
\partial B}_{\nu }}\left\vert \frac{X^{\ell }f\left( z\right) }{\left\vert
X^{\ell }f\left( z\right) \right\vert }\cdot y\right\vert ^{p}d\sigma _{\nu
}\left( y\right) =\gamma _{p,\nu }
\end{equation*}%
is \emph{independent} of the row vector in $\mathbb{\partial B}_{\nu }$ that
is dotted with $y$. Thus we have%
\begin{eqnarray*}
\left\vert T^{\ell }\left( X^{\ell }f\right) \left( z\right) \right\vert
^{p} &=&\left\vert T_{x\left( z\right) }^{\ell }\left( X^{\ell }f\right)
\left( z\right) \right\vert ^{p}=\frac{1}{\gamma _{p,\nu }}\int_{\mathbb{%
\partial B}_{\nu }}\left\vert T_{x\left( z\right) }^{\ell }\left( X^{\ell
}f\right) \left( z\right) y\right\vert ^{p}d\sigma _{\nu }\left( y\right) ,
\\
\left\vert X^{\ell }f\left( z\right) \right\vert ^{p} &=&\frac{1}{\gamma
_{p,\nu }}\int_{\mathbb{\partial B}_{\nu }}\left\vert \left( X^{\ell
}f\right) _{y}\left( z\right) \right\vert ^{p}d\sigma _{\nu }\left( y\right)
.
\end{eqnarray*}%
So with $d\omega _{p\sigma }\left( z\right) =\left( 1-\left\vert
z\right\vert ^{2}\right) ^{p\sigma }d\lambda _{n}\left( z\right) $, we
conclude that%
\begin{eqnarray*}
&&\int_{\mathbb{B}_{n}}\left\vert X^{m}\left( \mathbb{M}_{\varphi }f\right)
\right\vert ^{p}d\omega _{p\sigma }\left( z\right) \\
&\leq &C_{n,\sigma ,p,m}\sum_{\ell =0}^{m}\int_{\mathbb{B}_{n}}\left\vert
T^{\ell }\left( X^{\ell }f\right) \left( z\right) \right\vert ^{p}d\omega
_{p\sigma }\left( z\right) \\
&=&C_{n,\sigma ,p,m}\sum_{\ell =0}^{m}\frac{1}{\gamma _{p,\nu }}\int_{%
\mathbb{\partial B}_{\nu }}\int_{\mathbb{B}_{n}}\left\vert x\left( z\right)
\left( X^{m-\ell }\varphi \right) \left( z\right) \left( X^{\ell
}f_{y}\right) \left( z\right) \right\vert ^{p}d\omega _{p\sigma }\left(
z\right) d\sigma _{\nu }\left( y\right) \\
&\leq &C_{n,\sigma ,p,m}\sum_{\ell =0}^{m}\frac{1}{\gamma _{p,\nu }}\int_{%
\mathbb{\partial B}_{\nu }}\int_{\mathbb{B}_{n}}\left\vert \left( X^{m-\ell
}\varphi \right) \left( z\right) \left( X^{\ell }f\right) _{y}\left(
z\right) \right\vert ^{p}d\omega _{p\sigma }\left( z\right) d\sigma _{\nu
}\left( y\right) \\
&\leq &C_{n,\sigma ,p,m}\sum_{\ell =0}^{m}\frac{1}{\gamma _{p,\nu }}\int_{%
\mathbb{\partial B}_{\nu }}\left\Vert \mathbb{M}_{\varphi }\right\Vert
_{B_{p}^{\sigma }\left( \mathbb{B}_{n}\right) \rightarrow B_{p}^{\sigma
}\left( \mathbb{B}_{n};\mathbb{C}^{\mu }\right) }^{p}\int_{\mathbb{B}%
_{n}}\left\vert \left( \mathcal{X}^{m}f\right) _{y}\left( z\right)
\right\vert ^{p}d\omega _{p\sigma }\left( z\right) d\sigma _{\nu }\left(
y\right)
\end{eqnarray*}%
by the case $M=1$ of Proposition \ref{multilinear}, where $\ell ^{2}$ there
is replaced by $\mathbb{C}^{\nu }$, $g$ by $\varphi $ and $h$ by $f_{y}$.
Now we use the equality%
\begin{equation*}
\int_{\mathbb{\partial B}_{\nu }}\left\vert \left( \mathcal{X}^{m}f\right)
_{y}\left( z\right) \right\vert ^{p}d\sigma _{\nu }\left( y\right) =\gamma
_{p,\nu }\left\vert \mathcal{X}^{m}f\left( z\right) \right\vert ^{p}
\end{equation*}%
to obtain%
\begin{eqnarray*}
\int_{\mathbb{B}_{n}}\left\vert X^{m}\left( \mathbb{M}_{\varphi }f\right)
\right\vert ^{p}d\omega _{p\sigma }\left( z\right) &\leq &C_{n,\sigma
,p,m}\left\Vert \mathbb{M}_{\varphi }\right\Vert _{B_{p}^{\sigma }\left( 
\mathbb{B}_{n}\right) \rightarrow B_{p}^{\sigma }\left( \mathbb{B}_{n};%
\mathbb{C}^{\mu }\right) }^{p}\int_{\mathbb{B}_{n}}\left\vert \mathcal{X}%
^{m}f\left( z\right) \right\vert ^{p}d\omega _{p\sigma }\left( z\right) \\
&\leq &C_{n,\sigma ,p,m}\left\Vert \mathbb{M}_{\varphi }\right\Vert
_{B_{p}^{\sigma }\left( \mathbb{B}_{n}\right) \rightarrow B_{p}^{\sigma
}\left( \mathbb{B}_{n};\mathbb{C}^{\mu }\right) }^{p}\left\Vert f\right\Vert
_{B_{p}^{\sigma }\left( \mathbb{B}_{n};\mathbb{C}^{\mu }\right) }^{p}.
\end{eqnarray*}%
Since $m$ depends only on $n$, $\sigma $ and $p$, this completes the proof
of (\ref{itsuff}), and hence that of Lemma \ref{tensornorms}

\bigskip

Finally we return to complete the proof of Proposition \ref{multilinear}. We
have already proved the case $M=1$. Now we sketch a proof of the case $M=2$
using the multiplier norm inequality (\ref{ktensor}) with $\kappa =2$. By
multiplicativity of $\left\vert \cdot \right\vert $ on tensors, it suffices
to show that for $m=\ell _{1}+\ell _{2}+k$, 
\begin{eqnarray}
&&\int_{\mathbb{B}_{n}}\left\vert \left( 1-\left\vert z\right\vert
^{2}\right) ^{\sigma }\left( \mathcal{Y}^{\ell _{1}}g\right) \otimes \left( 
\mathcal{Y}^{\ell _{2}}g\right) \left( \mathcal{Y}^{k}h\right) \right\vert
^{p}d\lambda _{n}\left( z\right)  \label{N=2} \\
&&\ \ \ \ \ \ \ \ \ \ \leq C_{n,\sigma ,p}\left\Vert \mathbb{M}%
_{g}\right\Vert _{B_{p}^{\sigma }\left( \mathbb{B}_{n}\right) \rightarrow
B_{p}^{\sigma }\left( \mathbb{B}_{n};\ell ^{2}\right) }^{2p}\left\Vert
h\right\Vert _{B_{p}^{\sigma }\left( \mathbb{B}_{n}\right) }^{p}.  \notag
\end{eqnarray}%
This time we write using the divergence form of Leibniz' formula on tensor
products (c.f. (\ref{revLieb})),%
\begin{eqnarray*}
\left( \mathcal{Y}^{\ell _{1}}g\right) \otimes \left( \mathcal{Y}^{\ell
_{2}}g\right) \left( \mathcal{Y}^{k}h\right) &=&\left( \mathcal{Y}^{\ell
_{1}}g\right) \otimes \left\{ \sum_{\alpha =0}^{\ell _{2}}c_{\alpha }^{\ell
_{2}}\mathcal{Y}^{\alpha }\left( g\mathcal{Y}^{k+\ell _{2}-\alpha }h\right)
\right\} \\
&=&\sum_{\alpha =0}^{\ell _{2}}c_{\alpha }^{\ell _{2}}\left( \mathcal{Y}%
^{\ell _{1}}g\right) \otimes \left[ \mathcal{Y}^{\alpha }\left( g\mathcal{Y}%
^{k+\ell _{2}-\alpha }h\right) \right] \\
&=&\sum_{\alpha =0}^{\ell _{2}}c_{\alpha }^{\ell _{2}}\left\{ \sum_{\beta
=0}^{\ell _{1}}c_{\beta }^{\ell _{1}}\mathcal{Y}^{\beta }\left( g\otimes 
\mathcal{Y}^{\alpha +\ell _{1}-\beta }\left( g\mathcal{Y}^{k+\ell
_{2}-\alpha }h\right) \right) \right\} .
\end{eqnarray*}

We first use the Hilbert space valued interpolation theorem together with
the case $\kappa =2$ of Lemma \ref{tensornorms} to show that $g\in
M_{B_{p}^{s_{1}}\left( \mathbb{B}_{n};\ell ^{2}\right) \rightarrow
B_{p}^{s_{1}}\left( \mathbb{B}_{n};\ell ^{2}\otimes \ell ^{2}\right) }$ and $%
g\in M_{B_{p}^{s_{2}}\left( \mathbb{B}_{n}\right) \rightarrow
B_{p}^{s_{2}}\left( \mathbb{B}_{n};\ell ^{2}\right) }$ for appropriate
values of $s_{1}$ and $s_{2}$. Assuming for convenience that $\mathcal{Y}%
=\left( 1-\left\vert z\right\vert ^{2}\right) R$, and keeping in mind Remark %
\ref{D}, we obtain%
\begin{eqnarray*}
&&\left\Vert g\left( z\right) \otimes \left( 1-\left\vert z\right\vert
^{2}\right) ^{\sigma -s_{1}}\mathcal{Y}^{\alpha +\ell _{1}-\beta }\left( g%
\mathcal{Y}^{k+\ell _{2}-\alpha }h\right) \right\Vert _{B_{p}^{s_{1}}\left( 
\mathbb{B}_{n};\ell ^{2}\otimes \ell ^{2}\right) }^{p} \\
&\leq &\left\Vert \mathbb{M}_{g}\right\Vert _{B_{p}^{s_{1}}\left( \mathbb{B}%
_{n};\ell ^{2}\right) \rightarrow B_{p}^{s_{1}}\left( \mathbb{B}_{n};\ell
^{2}\otimes \ell ^{2}\right) }^{p}\left\Vert \left( 1-\left\vert
z\right\vert ^{2}\right) ^{\sigma -s_{1}}\mathcal{Y}^{\alpha +\ell
_{1}-\beta }\left( g\mathcal{Y}^{k+\ell _{2}-\alpha }h\right) \right\Vert
_{B_{p}^{s_{1}}\left( \mathbb{B}_{n};\ell ^{2}\right) }^{p} \\
&=&\left\Vert \mathbb{M}_{g}\right\Vert _{B_{p}^{s_{1}}\left( \mathbb{B}%
_{n};\ell ^{2}\right) \rightarrow B_{p}^{s_{1}}\left( \mathbb{B}_{n};\ell
^{2}\otimes \ell ^{2}\right) }^{p}\int_{\mathbb{B}_{n}}\left\vert \left(
1-\left\vert z\right\vert ^{2}\right) ^{s_{1}}\mathcal{Y}^{\beta }\left(
1-\left\vert z\right\vert ^{2}\right) ^{\sigma -s_{1}}\mathcal{Y}^{\alpha
+\ell _{1}-\beta }\left( g\mathcal{Y}^{k+\ell _{2}-\alpha }h\right)
\right\vert ^{p}d\lambda _{n},
\end{eqnarray*}%
which by (\ref{ktensor}) is at most%
\begin{eqnarray*}
&&C_{n,\sigma ,p}\left\Vert \mathbb{M}_{g}\right\Vert _{B_{p}^{\sigma
}\left( \mathbb{B}_{n}\right) \rightarrow B_{p}^{\sigma }\left( \mathbb{B}%
_{n};\ell ^{2}\right) }^{p}\int_{\mathbb{B}_{n}}\left\vert \left(
1-\left\vert z\right\vert ^{2}\right) ^{s_{2}}\mathcal{Y}^{\alpha +\ell
_{1}}\left( g\left( 1-\left\vert z\right\vert ^{2}\right) ^{\sigma -s_{2}}%
\mathcal{Y}^{k+\ell _{2}-\alpha }h\right) \right\vert ^{p}d\lambda _{n} \\
&=&C_{n,\sigma ,p}\left\Vert \mathbb{M}_{g}\right\Vert _{B_{p}^{\sigma
}\left( \mathbb{B}_{n}\right) \rightarrow B_{p}^{\sigma }\left( \mathbb{B}%
_{n};\ell ^{2}\right) }^{p}\left\Vert g\left( 1-\left\vert z\right\vert
^{2}\right) ^{\sigma -s_{2}}\mathcal{Y}^{k+\ell _{2}-\alpha }h\right\Vert
_{B_{p}^{s_{2}}\left( \mathbb{B}_{n};\ell ^{2}\right) }^{p} \\
&\leq &C_{n,\sigma ,p}\left\Vert \mathbb{M}_{g}\right\Vert _{B_{p}^{\sigma
}\left( \mathbb{B}_{n}\right) \rightarrow B_{p}^{\sigma }\left( \mathbb{B}%
_{n};\ell ^{2}\right) }^{p}\left\Vert \mathbb{M}_{g}\right\Vert
_{B_{p}^{s_{2}}\left( \mathbb{B}_{n}\right) \rightarrow B_{p}^{s_{2}}\left( 
\mathbb{B}_{n};\ell ^{2}\right) }^{p}\left\Vert \left( 1-\left\vert
z\right\vert ^{2}\right) ^{\sigma -s_{2}}\mathcal{Y}^{k+\ell _{2}-\alpha
}h\right\Vert _{B_{p}^{s_{2}}\left( \mathbb{B}_{n}\right) }^{p} \\
&\leq &C_{n,\sigma ,p}\left\Vert \mathbb{M}_{g}\right\Vert _{B_{p}^{\sigma
}\left( \mathbb{B}_{n}\right) \rightarrow B_{p}^{\sigma }\left( \mathbb{B}%
_{n};\ell ^{2}\right) }^{2p}\left\Vert h\right\Vert _{B_{p}^{\sigma }\left( 
\mathbb{B}_{n}\right) }^{p}.
\end{eqnarray*}%
Summing up over $\alpha $ and $\beta $ gives (\ref{N=2}).

Repeating this procedure for $M\geq 3$ and using (\ref{ktensor}) with $%
\kappa =M$ finishes the proof of Proposition \ref{multilinear}.

\subsection{Schur's test\label{Schur's Test}}

We prove Lemma \ref{Zlemma} using Schur's Test as given in Theorem 2.9 on
page 51 of \cite{Zhu}.

\begin{lemma}
\label{Schurtest}Let $\left( X,\mu \right) $ be a measure space and $H\left(
x,y\right) $ be a nonnegative kernel. Let $1<p<\infty $ and $\frac{1}{p}+%
\frac{1}{q}=1$. Define%
\begin{eqnarray*}
Tf\left( x\right) &=&\int_{X}H\left( x,y\right) f\left( y\right) d\mu \left(
y\right) , \\
T^{\ast }g\left( y\right) &=&\int_{X}H\left( x,y\right) g\left( x\right)
d\mu \left( x\right) .
\end{eqnarray*}%
If there is a positive function $h$ on $X$ and a positive constant $A$ such
that%
\begin{eqnarray*}
Th^{q}\left( x\right) &=&\int_{X}H\left( x,y\right) h\left( y\right)
^{q}d\mu \left( y\right) \leq Ah\left( x\right) ^{q},\ \ \ \ \ \mu -a.e.x\in
X, \\
T^{\ast }h^{p}\left( y\right) &=&\int_{X}H\left( x,y\right) h\left( x\right)
^{p}d\mu \left( x\right) \leq Ah\left( y\right) ^{p},\ \ \ \ \ \mu -a.e.y\in
X,
\end{eqnarray*}%
then $T$ is bounded on $L^{p}\left( \mu \right) $ with $\left\Vert
T\right\Vert \leq A$.
\end{lemma}

Now we turn to the proof of Lemma \ref{Zlemma}. The case $c=0$ of Lemma \ref%
{Zlemma} is Lemma 2.10 in \cite{Zhu}. To minimize the clutter of indices, we
first consider the proof for the case $c\neq 0$ when $p=2$ and $t=-n-1$.
Recall that 
\begin{eqnarray*}
\sqrt{\bigtriangleup \left( w,z\right) } &=&\left\vert 1-w\overline{z}%
\right\vert \left\vert \varphi _{z}\left( w\right) \right\vert , \\
\psi _{\varepsilon }\left( \zeta \right) &=&\left( 1-\left\vert \zeta
\right\vert ^{2}\right) ^{\varepsilon },
\end{eqnarray*}%
and%
\begin{equation*}
T_{a,b,c}f\left( z\right) =\int_{\mathbb{B}_{n}}\frac{\left( 1-\left\vert
z\right\vert ^{2}\right) ^{a}\left( 1-\left\vert w\right\vert ^{2}\right)
^{b+n+1}\left( \sqrt{\bigtriangleup \left( w,z\right) }\right) ^{c}}{%
\left\vert 1-w\overline{z}\right\vert ^{n+1+a+b+c}}f\left( w\right) d\lambda
_{n}\left( w\right) .
\end{equation*}%
We will compute conditions on $a$, $b$, $c$ and $\varepsilon $ such that we
have%
\begin{equation}
T_{a,b,c}\psi _{\varepsilon }\left( z\right) \leq C\psi _{\varepsilon
}\left( z\right) \text{ and }T_{a,b,c}^{\ast }\psi _{\varepsilon }\left(
w\right) \leq C\psi _{\varepsilon }\left( w\right) ,\ \ \ \ \ z,w\in \mathbb{%
B}_{n},  \label{Schurfunctions}
\end{equation}%
where $T_{a,b,c}^{\ast }$ denotes the dual relative to $L^{2}\left( \lambda
_{n}\right) $. For this we take $\varepsilon \in \mathbb{R}$ and compute%
\begin{equation*}
T_{a,b,c}\psi _{\varepsilon }\left( z\right) =\int_{\mathbb{B}_{n}}\frac{%
\left( 1-\left\vert z\right\vert ^{2}\right) ^{a}\left( 1-\left\vert
w\right\vert ^{2}\right) ^{n+1+b+\varepsilon }\left\vert \varphi _{z}\left(
w\right) \right\vert ^{c}}{\left\vert 1-w\overline{z}\right\vert ^{n+1+a+b}}%
d\lambda _{n}\left( w\right) .
\end{equation*}

Note that the integral is finite if and only if $\varepsilon >-b-1$. Now
make the change of variable $w=\varphi _{z}\left( \zeta \right) $ and use
that $\lambda _{n}$ is invariant to obtain%
\begin{eqnarray*}
T_{a,b,c}\psi _{\varepsilon }\left( z\right) &=&\int_{\mathbb{B}_{n}}\frac{%
\left( 1-\left\vert z\right\vert ^{2}\right) ^{a}\left( 1-\left\vert
w\right\vert ^{2}\right) ^{n+1+b+\varepsilon }\left\vert \varphi _{z}\left(
w\right) \right\vert ^{c}}{\left\vert 1-w\overline{z}\right\vert ^{n+1+a+b}}%
d\lambda _{n}\left( w\right) \\
&=&\int_{\mathbb{B}_{n}}F\left( w\right) d\lambda _{n}\left( w\right) =\int_{%
\mathbb{B}_{n}}F\left( \varphi _{z}\left( \zeta \right) \right) d\lambda
_{n}\left( \zeta \right) \\
&=&\int_{\mathbb{B}_{n}}\frac{\left( 1-|z|^{2}\right) ^{a}\left(
1-\left\vert \varphi _{z}\left( \zeta \right) \right\vert ^{2}\right)
^{n+1+b+\varepsilon }\left\vert \zeta \right\vert ^{c}}{\left\vert 1-%
\overline{\varphi _{z}\left( \zeta \right) }z\right\vert ^{n+1+a+b}\left(
1-|\zeta |^{2}\right) ^{n+1}}dV\left( \zeta \right) .
\end{eqnarray*}%
From the identity (Theorem 2.2.2 in \cite{Rud}),%
\begin{equation*}
1-\left\langle \varphi _{a}\left( \beta \right) ,\varphi _{a}\left( \gamma
\right) \right\rangle =\frac{\left( 1-\left\langle a,a\right\rangle \right)
\left( 1-\left\langle \beta ,\gamma \right\rangle \right) }{\left(
1-\left\langle \beta ,a\right\rangle \right) \left( 1-\left\langle a,\gamma
\right\rangle \right) },
\end{equation*}%
we obtain the identities%
\begin{eqnarray*}
1-\varphi _{z}\left( \zeta \right) \overline{z} &=&1-\left\langle \varphi
_{z}\left( \zeta \right) ,\varphi _{z}\left( 0\right) \right\rangle =\frac{%
1-|z|^{2}}{1-\zeta \overline{z}}, \\
1-\left\vert \varphi _{z}\left( \zeta \right) \right\vert ^{2}
&=&1-\left\langle \varphi _{z}\left( \zeta \right) ,\varphi _{z}\left( \zeta
\right) \right\rangle =\frac{\left( 1-\left\vert z\right\vert ^{2}\right)
\left( 1-\left\vert \zeta \right\vert ^{2}\right) }{\left\vert 1-\zeta 
\overline{z}\right\vert ^{2}}.
\end{eqnarray*}%
Plugging these identities into the formula for $T_{a,b,c}\psi _{\varepsilon
}\left( z\right) $ we obtain%
\begin{eqnarray}
T_{a,b,c}\psi _{\varepsilon }\left( z\right) &=&\int_{\mathbb{B}_{n}}\frac{%
\left( 1-|z|^{2}\right) ^{a}\left( \frac{\left( 1-\left\vert z\right\vert
^{2}\right) \left( 1-\left\vert \zeta \right\vert ^{2}\right) }{\left\vert
1-\zeta \overline{z}\right\vert ^{2}}\right) ^{n+1+b+\varepsilon }\left\vert
\zeta \right\vert ^{c}}{\left\vert \frac{1-|z|^{2}}{1-\zeta \overline{z}}%
\right\vert ^{n+1+a+b}\left( 1-|\zeta |^{2}\right) ^{n+1}}dV\left( \zeta
\right)  \label{Tpsi} \\
&=&\psi _{\varepsilon }\left( z\right) \int_{\mathbb{B}_{n}}\frac{\left(
1-\left\vert \zeta \right\vert ^{2}\right) ^{b+\varepsilon }\left\vert \zeta
\right\vert ^{c}}{\left\vert 1-\zeta \overline{z}\right\vert
^{n+1+b-a+2\varepsilon }}dV\left( \zeta \right) .  \notag
\end{eqnarray}

Now from Theorem 1.12 in \cite{Zhu} we obtain that%
\begin{equation*}
\sup_{z\in \mathbb{B}_{n}}\int_{\mathbb{B}_{n}}\frac{\left( 1-\left\vert
\zeta \right\vert ^{2}\right) ^{\alpha }}{\left\vert 1-\zeta \overline{z}%
\right\vert ^{\beta }}dV\left( \zeta \right) <\infty
\end{equation*}%
if and only if $\beta -\alpha <n+1$. Provided $c>-2n$ it is now easy to see
that we also have%
\begin{equation*}
\sup_{z\in \mathbb{B}_{n}}\int_{\mathbb{B}_{n}}\frac{\left( 1-\left\vert
\zeta \right\vert ^{2}\right) ^{\alpha }\left\vert \zeta \right\vert ^{c}}{%
\left\vert 1-\zeta \overline{z}\right\vert ^{\beta }}dV\left( \zeta \right)
<\infty
\end{equation*}%
if and only if $\beta -\alpha <n+1$. It now follows from the above that 
\begin{equation*}
T_{a,b,c}\psi _{\varepsilon }\left( z\right) \leq C\psi _{\varepsilon
}\left( z\right) ,\ \ \ \ \ z\in \mathbb{B}_{n},
\end{equation*}%
if and only if%
\begin{equation*}
-b-1<\varepsilon <a.
\end{equation*}

Now we turn to the adjoint $T_{a,b,c}^{\ast }=T_{b+n+1,a-n-1,c}$ with
respect to the space $L^{2}\left( \lambda _{n}\right) $. With the change of
variable $z=\varphi _{w}\left( \zeta \right) $ we have 
\begin{eqnarray}
T_{a,b,c}^{\ast }\psi _{\varepsilon }\left( w\right) &=&\int_{\mathbb{B}_{n}}%
\frac{\left( 1-\left\vert z\right\vert ^{2}\right) ^{a+\varepsilon }\left(
1-\left\vert w\right\vert ^{2}\right) ^{b+n+1}\left\vert \varphi _{w}\left(
z\right) \right\vert ^{c}}{\left\vert 1-w\overline{z}\right\vert ^{n+1+a+b}}%
d\lambda _{n}\left( z\right)  \label{Tpsistar} \\
&=&\int_{\mathbb{B}_{n}}G\left( z\right) d\lambda _{n}\left( z\right) =\int_{%
\mathbb{B}_{n}}G\left( \varphi _{w}\left( \zeta \right) \right) d\lambda
n\left( \zeta \right)  \notag \\
&=&\int_{\mathbb{B}_{n}}\frac{\left( 1-|\varphi _{w}\left( \zeta \right)
|^{2}\right) ^{a+\varepsilon }\left( 1-\left\vert w\right\vert ^{2}\right)
^{b+n+1}\left\vert \zeta \right\vert ^{c}}{\left\vert 1-w\overline{\varphi
_{w}\left( \zeta \right) }\right\vert ^{n+1+a+b}\left( 1-\left\vert \zeta
\right\vert ^{2}\right) ^{n+1}}dV\left( \zeta \right)  \notag \\
&=&\int_{\mathbb{B}_{n}}\frac{\left( \frac{\left( 1-\left\vert w\right\vert
^{2}\right) \left( 1-\left\vert \zeta \right\vert ^{2}\right) }{\left\vert
1-\zeta \overline{w}\right\vert ^{2}}\right) ^{a+\varepsilon }\left(
1-\left\vert w\right\vert ^{2}\right) ^{b+n+1}\left\vert \zeta \right\vert
^{c}}{\left\vert \frac{1-|w|^{2}}{1-\zeta \overline{w}}\right\vert
^{n+1+a+b}\left( 1-\left\vert \zeta \right\vert ^{2}\right) ^{n+1}}dV\left(
\zeta \right)  \notag \\
&=&\psi _{\varepsilon }\left( w\right) \int_{\mathbb{B}_{n}}\frac{\left(
1-\left\vert \zeta \right\vert ^{2}\right) ^{a+\varepsilon -n-1}\left\vert
\zeta \right\vert ^{c}}{\left\vert 1-\zeta \overline{w}\right\vert
^{a-b+2\varepsilon -n-1}}dV\left( \zeta \right) .  \notag
\end{eqnarray}%
Arguing as above and provided $c>-2n$, we obtain%
\begin{equation*}
T_{a,b,c}^{\ast }\psi _{\varepsilon }\left( w\right) \leq C\psi
_{\varepsilon }\left( w\right) ,\ \ \ \ \ w\in \mathbb{B}_{n},
\end{equation*}%
if and only if%
\begin{equation*}
-a+n<\varepsilon <b+n+1.
\end{equation*}

Altogether then there is $\varepsilon \in \mathbb{R}$ such that $h=\sqrt{%
\psi _{\varepsilon }}$ is a Schur function for $T_{a,b,c}$ on $L^{2}\left(
\lambda _{n}\right) $ in Lemma \ref{Schurtest} if and only if%
\begin{equation*}
\max \left\{ -a+n,-b-1\right\} <\min \left\{ a,b+n+1\right\} .
\end{equation*}%
This is equivalent to $-2a<-n<2\left( b+1\right) $, which is (\ref%
{indexcondition}) in the case $p=2,t=-n-1$. Thus Lemma \ref{Schurtest}
completes the proof that this case of (\ref{indexcondition}) implies the
boundedness of $T_{a,b,c}$ on $L^{2}\left( \lambda _{n}\right) $. The
converse is easy - see for example the argument for the case $c=0$ on page
52 of \cite{Zhu}.

\bigskip

We now turn to the general case. The adjoint $T_{a,b,c}^{\ast }$ relative to
the Banach space $L^{p}\left( \nu _{t}\right) $ is easily computed to be $%
T_{a,b,c}^{\ast }=T_{b-t,a+t,c}$ (see page 52 of \cite{Zhu} for the case $%
c=0 $). Then from (\ref{Tpsi}) and (\ref{Tpsistar}) we have%
\begin{eqnarray*}
T_{a,b,c}\psi _{\varepsilon }\left( z\right) &=&\psi _{\varepsilon }\left(
z\right) \int_{\mathbb{B}_{n}}\frac{\left( 1-\left\vert \zeta \right\vert
^{2}\right) ^{b+\varepsilon }\left\vert \zeta \right\vert ^{c}}{\left\vert
1-\zeta \overline{z}\right\vert ^{n+1+b-a+2\varepsilon }}dV\left( \zeta
\right) , \\
T_{a,b,c}^{\ast }\psi _{\varepsilon }\left( w\right) &=&\psi _{\varepsilon
}\left( w\right) \int_{\mathbb{B}_{n}}\frac{\left( 1-\left\vert \zeta
\right\vert ^{2}\right) ^{a+t+\varepsilon }\left\vert \zeta \right\vert ^{c}%
}{\left\vert 1-\zeta \overline{w}\right\vert ^{a-b+2\varepsilon +t}}dV\left(
\zeta \right) .
\end{eqnarray*}%
Let $\frac{1}{p}+\frac{1}{q}=1$. We apply Schur's Lemma \ref{Schurtest} with 
$h\left( \zeta \right) =\left( 1-\left\vert \zeta \right\vert ^{2}\right)
^{s}$ and%
\begin{equation}
s\in \left( -\frac{b+1}{q},\frac{a}{q}\right) \cap \left( -\frac{a+1+t}{p},%
\frac{b-t}{p}\right) .  \label{s}
\end{equation}%
Using Theorem 1.12 in \cite{Zhu} we obtain for $h$ with $s$ as in (\ref{s})
that%
\begin{equation*}
T_{a,b,c}h^{q}\leq Ch^{q}\text{ and }T_{a,b,c}^{\ast }h^{p}\leq Ch^{p}.
\end{equation*}%
Schur's Lemma \ref{Schurtest} now shows that $T_{a,b,c}$ is bounded on $%
L^{p}\left( \nu _{t}\right) $. Again, the converse follows from the argument
for the case $c=0$ on page 52 of \cite{Zhu}.

\end{document}